\pdfoutput=1
\RequirePackage{silence}
\documentclass[a4paper,svgnames]{amsart}
\usepackage{mathtools,amsfonts,amssymb,mathrsfs,stackrel,amsmath}
\usepackage[a4paper,hmargin=26mm,vmargin=27mm]{geometry}
\usepackage{tikz,tikz-cd}
\usetikzlibrary{matrix,arrows.meta,decorations.markings,shapes.multipart,decorations.pathreplacing}
\usepackage{enumitem}
\setlist[enumerate]{itemsep=0.15cm,label=\emph{\upshape(\alph*)}}
\setlist[enumerate,2]{itemsep=0.15cm,label=\emph{\upshape(\roman*)}}
\usepackage{aliascnt}
\usepackage{bbm}

\usepackage{array}
\newcolumntype{C}{>{$}c<{$}}

\usepackage{xparse}
\usepackage{cite}

\allowdisplaybreaks

\definecolor{mygray}{gray}{0.6}
\definecolor{mygraydark}{gray}{0.4}
\definecolor{mygraylight}{gray}{0.85}
\definecolor{spinach}{RGB}{46,139,87}
\definecolor{tomato}{RGB}{255,99,71}
\definecolor{orchid}{RGB}{143,40,194}
\definecolor{neon}{RGB}{77,77,255}
\definecolor{pumpkin}{RGB}{224,180,80}
\definecolor{citron}{RGB}{190,180,90}

\definecolor{lava}{RGB}{207,16,32}
\definecolor{cream}{RGB}{255,253,208}
\definecolor{verdigris}{RGB}{67,179,174}
\definecolor{Black}{RGB}{0,0,0}
\definecolor{mydarkblue}{RGB}{10,10,170}
\definecolor{darkspinach}{RGB}{20,70,20}
\definecolor{darktomato}{RGB}{155,40,30}
\definecolor{darkorchid}{RGB}{50,10,100}
\definecolor{darklava}{RGB}{150,8,16}

\let\emph\relax
\DeclareTextFontCommand{\emph}{\bfseries\em}

\newcommand{\cf}{cf.}
\newcommand{\etc}{etc.}

\newcommand{\ver}{verbatim}
\newcommand{\vive}{vice versa}
\newcommand{\muta}{mutatis mutandis}
\newcommand{\loccit}{loc. cit.}

\renewcommand{\dots}{\text{...}}

\newcommand\mapsfrom{\mathrel{\reflectbox{\ensuremath{\mapsto}}}}

\newcommand\floor[2]{\lfloor\tfrac{#1}{#2}\rfloor}

\newcommand\bsig{{\boldsymbol\sigma}}
\newcommand\balp{{\boldsymbol\alpha}}
\newcommand\bbet{{\boldsymbol\beta}}

\newcommand\blam{{\boldsymbol\lambda}}

\newcommand\brho{{\boldsymbol\rho}}
\newcommand\bmu{{\boldsymbol\mu}}
\newcommand\bnu{{\boldsymbol\nu}}
\newcommand\bom{{\boldsymbol\omega}}
\newcommand\charge{{\boldsymbol\kappa}}
\newcommand\chargetwo{{\boldsymbol\nu}}
\newcommand\chargethree{{\boldsymbol\upsilon}}

\newcommand\ba{{\boldsymbol{a}}}
\newcommand\bb{{\boldsymbol{b}}}
\newcommand\bd{\boldsymbol{d}}
\newcommand\bx{\boldsymbol{x}}
\newcommand\by{\boldsymbol{y}}
\newcommand\bz{\boldsymbol{z}}
\newcommand\bS{{\boldsymbol{S}}}
\newcommand\bT{{\boldsymbol{T}}}
\newcommand\bU{{\boldsymbol{U}}}
\newcommand\bV{{\boldsymbol{V}}}

\newcommand\bs{\boldsymbol{s}}
\newcommand\bt{\boldsymbol{t}}

\newcommand\Dcal{\mathcal{D}}
\newcommand\Pcal{\mathcal{P}}

\newcommand\str[1][s]{\mathsf{#1}}
\newcommand\Sub[1]{\bar{#1}}
\newcommand\Submap{S}
\newcommand\Ind[1]{\mathcal{#1}^{!}}
\newcommand\Res[1]{\mathcal{#1}_{*}}
\newcommand\CoInd[1]{\mathcal{#1}^{*}}

\newcommand\Indz[1]{{#1}^{!}}
\newcommand\Resz[1]{{#1}_{*}}

\newcommand\Web{\mathbb{W}}

\NewDocumentCommand\Webaa{ O{\bx,\bi} O{\by,\bj} }{\Web_{#1}^{#2}}
\NewDocumentCommand\Webab{ O{(\charge,\bx),\bi} O{(\chargetwo,\by),\bj} }{\Web_{#1}^{#2}}
\NewDocumentCommand\Webabs{ O{\bx} O{\by} O{\bi} O{\bj} }{\Web_{\Sub{#1},\Sub{#3}}^{\Sub{#2},\Sub{#4}}}

\NewDocumentCommand\WA{ O{\beta} O{\brho} }{\mathscr{W}_{#1}^{#2}}
\NewDocumentCommand\WAc{ O{\beta} O{\brho} }{\mathscr{R}_{#1}^{#2}}
\NewDocumentCommand\WAs{ O{\beta} O{\brho} }{\mathscr{W}_{\Sub{#1}}^{{#2}}}
\NewDocumentCommand\WAsc{ O{\beta} O{\brho} }{\mathscr{R}^{\Sub{#1}}^{\Sub{#2}}}
\newcommand\WABasis{\mathcal{B}_{\beta}}
\NewDocumentCommand\SAX{ O{\Sub{\beta}} O{\brho} }{\Sub{\mathscr{W}}_{#1}^{#2}(X,\Sub{X})}
\newcommand\WAsBasis{\mathcal{B}_{\Sub{\beta}}}

\NewDocumentCommand\WAg{ O{\beta} O{\brho} }{g\mathscr{W}_{#1}^{#2}}

\NewDocumentCommand\TA{ O{\beta} O{\brho} }{\mathcal{W}_{#1}^{#2}}
\NewDocumentCommand\TAc{ O{\beta} O{\brho} }{\mathcal{R}_{#1}^{#2}}

\NewDocumentCommand\WAlam{ s O{\blam} D(){A}} {\mathscr{W}_{n}^{\IfBooleanTF{#1}{\gdom}{\gedom}#2}}

\newcommand\Oneg{\1_{\Sub{\Gamma},\Gamma}}
\newcommand\sg{\Submap_{\Gamma,\Sub{\Gamma}}}

\newcommand\affine[1]{\underline{#1}}
\newcommand\Sym[1][n]{{\mathfrak{S}_{#1}}}
\newcommand\hell{\affine{\ell}}
\NewDocumentCommand\Parts{ D(){\ell} O{n} }{\mathrm{P}_{#1,#2}}
\NewDocumentCommand\Partsgood{ D(){\ell} O{n} }{\mathrm{P}_{#1,#2}^{g}}
\newcommand\hParts[1][n]{\affine{\mathrm{P}}_{\hell,#1}}
\newcommand\hPartsgood[1][n]{\affine{\mathrm{P}}_{\hell,#1}^{g}}
\newcommand\Nodes[1][n]{\mathscr{N}_{\ell,#1}}
\newcommand\hcoord{\mathtt{x}_{\affine{\charge}}^{A}}

\newcommand\hcoordc{\mathtt{x}_{\affine{\charge}}^{C}}

\NewDocumentCommand\res{sO{}}{\mathop{\rm res}\nolimits_{#2\rho}}
\newcommand\Affch[1][\blam]{A^{\affine\charge}(#1)}
\newcommand\BX[1][{\mathscr{W}}]{\mathrm{B}_{#1}}

\newcommand\re{\mathsf{r}}

\newcommand\sand[1][\blam]{S_{#1}}

\usepackage{etoolbox}
\newcommand{\DeclareMyOperator}[1]{%
\expandafter\DeclareMathOperator\csname #1\endcsname{#1}
}
\forcsvlist{\DeclareMyOperator}{Shape,defect,gr,End}

\def\1{\mathbf{1}}
\newcommand\onealg{\1_{\mathscr{W}}}
\def\height{\mathop{\rm ht}\nolimits}

\def\Item(#1){\item\textbf{\upshape(#1\upshape)}}

\tikzset{
anchorbase/.style={baseline={([yshift=#1]current bounding box.center)}},
anchorbase/.default={-0.5ex},
dot colour/.initial=black,
dot colour/.default=black,
tinynodes/.style={font=\tiny,text height=0.25ex,text depth=0.05ex},
smallnodes/.style={font=\scriptsize,text height=0.75ex,text depth=0.15ex},
dots/.style={line width=1pt,line cap=round, gray, dash pattern=on 0pt off 2\pgflinewidth},
redstring/.style = {draw=red!50,fill=none,line width=0.35mm,preaction={draw=red,line width=2.5pt,-},nodes={color=red}},
affine/.style= {draw=citron!50,fill=none,
line width=0.35mm,preaction={draw=citron,line width=2.5pt,-},nodes={color=citron}},
solid/.style = {draw=blue,fill=none,dot colour=blue,line width=0.4mm,nodes={color=blue}},
ghost/.style = {draw=darkgray,fill=none,dot colour=darkgray,
densely dashed,line width=0.4mm,nodes={color=darkgray}},
crossline/.style={preaction={draw=white,line width=4.75pt,-},preaction={draw=black,line width=0.9pt,-}},
dot/.style = {
decoration={markings,
post length=0.25mm,
pre length=0.25mm,
mark=at position #1 with {\node[circle,radius=0.3cm,inner sep=-2.0pt,color=\pgfkeysvalueof{/tikz/dot colour},fill=\pgfkeysvalueof{/tikz/dot colour}]{};}
},
postaction={decorate}
},
dot/.default=0.5,
}

\tikzstyle directed=[postaction={decorate,decoration={markings,
mark=at position #1 with {\arrow[line width=0.25mm, black]{>}}}}]

\NewDocumentCommand\crossing{ O{1ex} mmO{} mmO{} } {
\tikz[centered=#1]{
\draw[#2, dot/.list={#4}](0,0)node[below]{$#3$}--++(1,1);
\draw[#5, dot/.list={#7}](1,0)node[below]{$#6$}--++(-1,1);
}%
}

\NewDocumentCommand\DottedIdempotent{ D(){3} omm} {
\begin{tikzpicture}[scale=1.2,anchorbase]
\def\residues{{#4}}
\foreach \res [count=\c,
evaluate=\res as \pos using {\res<=#3 ? #1*\res-\c*#1/15 : #1*(\res-2)-\c*#1/15}
] in {#4} {
\coordinate (\c) at (\pos,0);
\draw[solid](\c)node[below]{$\res$}--++(0,1);
\ifnum\res=#3\relax\else
\draw[ghost](\pos+#1,0)--++(0,1)node[above]{$\res$};
\fi
\ifnum\c=1\relax
\draw[redstring](\pos+#1/15,0)node[below]{$\res$}--++(0,1);
\fi
}
\foreach \pt [evaluate=\pt as \good using {\residues[\pt-1]==#3 ? 0 : 1}
] in {#2} {
\draw[solid,dot](\pt)--++(0,1);
\ifnum\good=1
\draw[ghost,dot]([shift={(#1,0)}]\pt)--++(0,1);
\fi
}
\end{tikzpicture}
}

\NewDocumentCommand\DoubleCrossing{ O{1ex} mmO{} mmO{} } {
\begin{tikzpicture}[scale=1.2,centered=#1]
\draw[#2, dot/.list={#4}] (0,0)node[below]{$#3$} .. controls (1,0.5) .. (0,1);
\draw[#5, dot/.list={#7}] (1,0)node[below]{$#6$} .. controls (0,0.5) .. (1,1);
\end{tikzpicture}\bigskip
}

\newcommand\Ddots[1]{%
\tikz[scale=0.9,centered,anchorbase]{\foreach \y in {1,...,#1} {
\draw[directed=\y/#1] (0,0) to (0.925,0);}
\foreach \x in {0,...,#1} {
\node[circle,inner sep=1.8pt,fill=DarkBlue] at (\x/#1,0){};}
}}

\newcommand\Dldots[2]{%
\tikz[scale=0.9,centered]{\foreach \y in {1,...,#2} {
\draw[directed=\y/#2] (0,0) to (0.925,0);}
\foreach \x [count=\c from 0] in {#1} {
\node[circle,inner sep=1.8pt,fill=DarkBlue,label={above:\scriptsize$\x$}] at (\c/#2,0){};}
}}

\NewDocumentCommand\dotstring{ O{1ex} mm}{%
\tikz[centered=#1]{\draw[#2,dot](0,0)node[below]{$#3$}--++(0,1);}}

\NewDocumentCommand\stringdot{ O{1ex} mmmm}{%
\tikz[centered=#1]{
\draw[#2](0,0)node[below]{$#3$}--++(0,1);
\draw[#4,dot](1,0)node[below]{$#5$}--++(0,1);
}%
}%

\DeclarePairedDelimiterX{\set}[1]{\{}{\}}{\setargs{#1}}
\NewDocumentCommand{\setargs}{>{\SplitArgument{1}{|}}m}{\setargsaux#1}
\NewDocumentCommand{\setargsaux}{mm}
{\IfNoValueTF{#2}{#1} {#1\,\delimsize|\,\mathopen{}#2}}

\def\NewTheorem#1{%
\newaliascnt{#1}{equation}%
\newtheorem{#1}[#1]{#1}%
\aliascntresetthe{#1}%
\expandafter\def\csname #1autorefname\endcsname{#1}%
}
\def\equationautorefname~#1\null{(#1)\null}

\numberwithin{equation}{subsection}

\NewTheorem{Proposition}
\NewTheorem{Theorem}
\NewTheorem{Corollary}
\NewTheorem{Lemma}
\theoremstyle{definition}
\NewTheorem{Definition}
\NewTheorem{Notation}
\NewTheorem{Example}
\AtEndEnvironment{Example}{\null\hfill$\Diamond$}%
\NewTheorem{Examples}
\AtEndEnvironment{Examples}{\vskip-10mm\null\hfill$\Diamond$}%

\theoremstyle{remark}
\NewTheorem{Remark}
\NewTheorem{Assumption}

\def\Item(#1){\item\textbf{\upshape(#1\upshape)}}
\newcounter{relation}

\newcommand\TestRelation[1]{\smallskip\noindent\textbf{Relations (#1):}\space}
\def\relationautorefname~#1\null{\upshape(W$_{#1}$\upshape)}
\let\ref\autoref
\let\eqref\autoref

\newcommand\Std{\mathop{\rm Std}\nolimits}
\newcommand\SStd[1][\charge]{\mathop{\rm SStd}\nolimits_{#1}}
\newcommand\hSStd[1][\hat\charge]{\mathop{\rm SStd}\nolimits_{#1}}

\def\N{\mathbb{Z}_{\geq 0}}

\def\Z{\mathbb{Z}}

\def\y{\mathbf{y}}

\let\<=\langle
\let\>=\rangle

\NewDocumentCommand\mI{ sd() }{\mathtt{m}_{\IfBooleanT{#1}{\Sub}\Gamma}\IfNoValueF{#2}{(#2)}}
\def\pmod#1{\space(\text{mod }#1)}

\newcommand\bi{\mathbf{i}}
\newcommand\bj{\mathbf{j}}
\newcommand\bk{\mathbf{k}}

\newcommand\R{\mathbb{R}}

\def\map#1#2{\,{:}\,#1\!\longrightarrow\!#2}

\let\gedom=\trianglerighteq
\let\gdom=\vartriangleright
\let\ledom=\trianglelefteq
\let\ldom=\vartriangleleft

\usepackage[hypertexnames=false]{hyperref}
\usepackage{bookmark}
\hypersetup{
pdftoolbar=true,
pdfmenubar=true,
pdffitwindow=false,
pdfstartview={FitH},
pdftitle={Cellularity and subdivision of KLR and weighted KLRW algebras},
pdfauthor={Andrew Mathas and Daniel Tubbenhauer},
pdfsubject={},
pdfcreator={Andrew Mathas and Daniel Tubbenhauer},
pdfproducer={Andrew Mathas and Daniel Tubbenhauer},
pdfkeywords={},
pdfnewwindow=true,
colorlinks=true,
linkcolor=mydarkblue,
citecolor=teal,
filecolor=magenta,
urlcolor=orchid,
linkbordercolor=lava,
citebordercolor=teal,
urlbordercolor=orchid,
linktocpage=true
}
\usepackage{xltabular}
\usepackage{tabulary}
\usepackage{booktabs}
\newcommand{\mystrut}{\rule[-0.2\baselineskip]{0pt}{1.1\baselineskip}}

\newcommand{\nnfootnote}[1]{%
\begin{NoHyper}
\renewcommand\thefootnote{}\footnote{#1}%
\addtocounter{footnote}{-1}%
\end{NoHyper}
}

\def\makeautorefname#1#2{\csdef{#1autorefname}{#2}}
\makeautorefname{section}{Section}%
\makeautorefname{subsection}{Section}%
\makeautorefname{subsubsection}{Section}%

\begin{document}
\title[KLR and weighted KLRW algebras]{Cellularity and subdivision of KLR and weighted KLRW algebras}
\author[A. Mathas and D. Tubbenhauer]{Andrew Mathas and Daniel Tubbenhauer}

\address{A.M.: The University of Sydney, School of Mathematics and Statistics F07, Office Carslaw 718, NSW 2006, Australia,
\href{http://www.maths.usyd.edu.au/u/mathas/}{www.maths.usyd.edu.au/u/mathas/},
\href{https://orcid.org/0000-0001-7565-5798}{ORCID 0000-0001-7565-5798}}
\email{andrew.mathas@sydney.edu.au}

\address{D.T.: The University of Sydney, School of Mathematics and Statistics F07, Office Carslaw 827, NSW 2006, Australia, \href{http://www.dtubbenhauer.com}{www.dtubbenhauer.com}, \href{https://orcid.org/0000-0001-7265-5047}{ORCID 0000-0001-7265-5047}}
\email{daniel.tubbenhauer@sydney.edu.au}

\begin{abstract}
Weighted KLRW algebras are diagram algebras generalizing KLR
algebras. This paper undertakes a systematic study of these algebras culminating in the construction of homogeneous
affine cellular bases in affine types $A$ and $C$, which immediately gives cellular bases for the cyclotomic quotients of these algebras. In addition, we construct subdivision homomorphisms that relate weighted KLRW algebras for different quivers. As an application we obtain new results about the (cyclotomic) KLR algebras of affine type, including (re)proving that the cyclotomic KLR algebras of type $A^{(1)}_{e}$ and $C^{(1)}_{e}$ are graded cellular algebras.
\end{abstract}

\nnfootnote{\textit{Mathematics Subject Classification 2020.} Primary:
20C08, 20G43; Secondary: 18M30, 18N25.}
\nnfootnote{\textit{Keywords.} KLR algebras, weighted KLRW algebras, cellular bases, Hecke and Schur algebras.}

\addtocontents{toc}{\protect\setcounter{tocdepth}{1}}

\maketitle

\tableofcontents


\section{Introduction}


Khovanov--Lauda \cite{KhLa-cat-quantum-sln-first}, \cite{KhLa-cat-quantum-sln-second} and  Rouquier \cite{Ro-2-kac-moody}, \cite{Ro-quiver-hecke} independently introduced the
\emph{KLR algebras}, or quiver Hecke algebras, motivated by questions in categorification. Webster \cite{We-knot-invariants} further generalized these algebras to \emph{KLRW algebras}, which give categorifications of tensor products of simple highest weight modules. All of these algebras admit finite dimensional quotients, called \emph{cyclotomic KLR(W) algebras}. These algebras are graded algebras that
can be defined diagrammatically using generators and relations and they play a crucial role in categorification and representation theory.

The discovery of the KLR algebras and their properties initiated a major paradigm shift in representation theory. As important special cases these algebras include the group algebras of the symmetric groups, their Hecke algebras, and the cyclotomic Hecke algebras of type $A$, see for example \cite{BrKl-hecke-klr}. The KLR algebras are important because they reveal deep new structures in the module categories of these algebras, for example, their representation theory is related to Lusztig's geometric construction of canonical bases \cite{VaVa-canonical-bases-klr}.

This paper focuses on the
\emph{weighted KLRW algebras} from \cite{We-weighted-klr} and \cite{We-rouquier-dia-algebra}. As a consequences of our main results we obtain new information about the KLR algebras, which is difficult to obtain by working in the KLR setting. Webster's definition of the weighted KLRW algebras embellishes the KLR algebras by introducing a weighting and additional strings to the diagrammatic presentation of the KLR algebras (see \autoref{D:RationalCherednik}), together with some subtle changes in the relations. This paper undertakes a systematic study of these algebras with the following results:

\begin{enumerate}[label=\emph{\upshape(\Alph*)}]

\item We give a self-contained treatment of the basic properties of the weighted KLRW algebras. Many of these results are known to experts but do not appear in the literature.

\item We show that subdividing the underlying quiver induces an isomorphism between the corresponding weighted KLRW algebras; \autoref{T:SubDiv}. Subdivision gives a way to
relate simple modules for weighted KLRW algebras and KLR algebras with different quivers,
such as $A^{(1)}_{e-1}$ and $A^{(1)}_{e}$. These results generalize \cite{Ma-catrep-klr}.

\item The enhanced diagrammatics of weighted KLRW algebras allow us to construct homogeneous affine cellular bases of the weighted KLRW algebras of affine types $A$ and $C$.
Quite strikingly, we obtain cellular bases for the cyclotomic quotients of these algebras%
, essentially for free. As a consequence, we also get graded cellular bases for the corresponding KLR and cyclotomic KLR algebras by idempotent truncation. In cyclotomic type $A$, analogous bases were constructed by Bowman \cite{Bo-many-cellular-structures} but our approach is easier because we first construct bases for the weighted KLRW algebras and then deduce corresponding results for the cyclotomic quotients.

\end{enumerate}

Under certain conditions, the KLR algebras are isomorphic to idempotent subalgebras of the weighted KLRW algebras, see \autoref{P:WebAlg}. In particular, the
weighted KLRW algebras are a much larger class of algebras
(for example, they crucially depend on the choice
of coordinates for the endpoints of the strings) and these algebras usually have more simple modules than the corresponding KLR algebras. This is also true for the finite dimensional quotients of these algebras. Although they look more difficult at first sight, in many respects the weighted KLRW algebras are easier to work with, which makes them a useful tool for understanding the KLR algebras. Our results show that the weighted KLRW algebras are interesting algebras in and of themselves.


\subsection{Future directions}


The first part of this paper holds (almost) without restriction
on the underlying quiver $\Gamma$. In the last two sections
about cellularity we restrict our attention to the quivers
in \autoref{E:Quivers}. Let us further comment why these
quivers and related quivers are special.

\begin{enumerate}[label=$\bullet$]

  \item In a sequel to this paper \cite{MaTu-klrw-algebras-bad}, which is partly motivated by the combinatorics from \cite{Ar-cyclotomic-klr-affine-type}, we generalize many of the result in this paper to the quivers in \autoref{Eq:FunnyQuivers}. The next natural candidates to consider are the weighted KLRW algebras for the quivers of finite type, which would generalize results from \cite{KlLoMi-KLR-affine-cellular-typea} and
\cite{KlLo-klr-affine-cellular-finite-type}. Quivers of finite types $A$ and $C$ are not considered in this paper.

\item As we will see, the classification of simple modules for
the weighted KLRW algebras of affine type $A$ is straightforward.
In fact, these algebras are quasi-hereditary in the appropriate sense.
However, in affine type $C$ we only get partial results. We do not
know how to obtain the classification of simple modules
of the corresponding KLR algebras
(that can, for example, be deduced from \cite{ArMa-simples-complex-reflection}) in any nice way using the weighted KLRW framework.
However, we think that these are interesting questions that may have a nice diagrammatic description.

\item The quivers in \autoref{E:Quivers} seem to be the easiest
from the Fock space point of view.
They are also the easiest from the viewpoint of
the associated braid groups, as summarized in \cite{Al-braids-abcd}.
Moreover, the cellular bases in type $A$ admit a knot-theoretical interpretation \cite{Tu-gln-bases} (strictly speaking these are the cellular bases from \cite{HuMa-klr-basis}, but the ones discussed in this paper have a similar interpretation). A natural question is whether there are any knot-theoretical interpretations (for example, in the spirit of \cite{RoTu-homflypt-handlebody}), of the cellular bases of this paper.

\end{enumerate}

\noindent\textbf{Acknowledgments.}
This project started when A.M. visited
D.T. in Louvain-la-Neuve mid 2015, albeit in a very different form. This paper
was finally finished late 2021. We received a lot of support along the way.

Both authors were supported, in part, by the Australian Research Council.
We were financially supported by the
Universit{\'e} catholique de Louvain, the University of Sydney, the Hausdorff Center of Mathematics (HCM) and the Hausdorff Research Institute for Mathematics (HIM), both in Bonn,
and the Universit{\"a}t Z{\"u}rich, who sponsored several research visits. All of this is gratefully acknowledged.
During the years we were mathematically supported by many people, too many
to be named here. We especially thank
Chris Bowman, Joe Chuang, Huang Lin, Tao Qin, Salim Rostam, Liron Speyer, Catharina Stroppel, Pedro Vaz and Ben Webster for comments, helpful discussions, patiently answering our questions and freely sharing ideas.
Finally, we thank the referee for their
careful reading of the manuscript and their helpful suggestions for improving the paper.

Sadly D.T.s memory sucks: D.T. cannot remember all the support and help
that deserves to be acknowledged here -- omissions are not intentional. D.T. bears the blame.


\section{Weighted KLRW algebras}\label{S:WebsterAlgebras}


In this section we slightly rephrase Webster's definitions of
weighted KLR algebras \cite{We-weighted-klr} and diagrammatic Cherednik
algebras \cite{We-rouquier-dia-algebra}. These
algebras are defined using string diagrams in the usual sense, but they also take into account the positions of the strings. We prove some basic properties of these algebras in \autoref{S:Properties}, some of which appear to be new, some are reformulations of {\loccit}


\subsection{Quiver combinatorics}\label{SS:Quiver}


We start by fixing notation from the classical theory of
Kac--Moody algebras, which can be found in \cite{Ka-infdim-lie}. A matrix
$(a_{ij})_{i,j\in I}$ is a \emph{symmetrizable generalized Cartan
matrix} if $I\subset\Z$, $a_{ii}=2$, $a_{ij}\in\Z_{\leq 0}$ for $i\neq
j$, $a_{ij}=0$ if and only if $a_{ji}=0$, and there is a
symmetrizer $\bd=(d_{i})_{i\in I}$ such that $(d_{i}a_{ij})_{i,j\in I}$
is symmetric. In this paper we only use such matrices with $a_{ij}a_{ji}<4$.

In \cite[\S 4.7]{Ka-infdim-lie} it is explained how
$(a_{ij})_{i,j\in I}$ gives rise to a quiver $\Gamma=(I,E)$, with vertex set
$I$ and edge set $E$ with $a_{ij}$ edges from $i$ to
$j$, for $i,j\in I$. That is, if $|a_{ij}|\geq|a_{ji}|$, then
$\Gamma$ has $|a_{ij}|$ edges from $i$ to $j$, which are oriented
from $i$ to $j$ unless $|a_{ij}|=1$.

\begin{Definition}\label{D:Quiver}
An oriented quiver $\Gamma=(I,E)$ with countable vertex set $I$ and
countable edge set $E$ is \emph{symmetrizable} if it arises from a
symmetrizable generalized Cartan matrix for $a_{ij}a_{ji}<4$
by choosing an orientation on
the simply laced edges.
\end{Definition}

We will see in \autoref{P:FixedOrientation} that the choice of orientation is not important  in this paper.

\begin{Notation}\label{N:Quiver}
If not stated otherwise, we fix an oriented symmetrizable quiver $\Gamma$ and $n,\ell\in\N$. (By convention, whenever $n$ or $\ell$ are zero then the notions involving them are vacuous.)
Let $e=\#I$ and $\#E$ be the sizes of the vertex and edge sets, respectively. We allow $e$ and $\#E$ to be infinite. A \emph{residue} is an element $i\in I$ and an $n$-tuple $\bi\in I^{n}$ is a \emph{residue sequence}.
\end{Notation}

\begin{Notation}
We allow three types of edges $i\to j$, $i\Rightarrow j$ and $i\Rrightarrow j$. We write $i\rightsquigarrow j$ when the multiplicity is unimportant. All of these count as one edge. In particular, the quivers arising in this way include all Dynkin quivers, except affine type $A_{1}$. We include the affine $A_{1}$ quiver by using the convention is that there are two arrows, $0\to 1$ and $1\to 0$, between the vertices $0$ and $1$.
\end{Notation}

\begin{Example}
The earlier sections of this paper apply with
almost no restrictions on the choice of quiver but for the some results of
the paper we restrict our attention to the quivers:
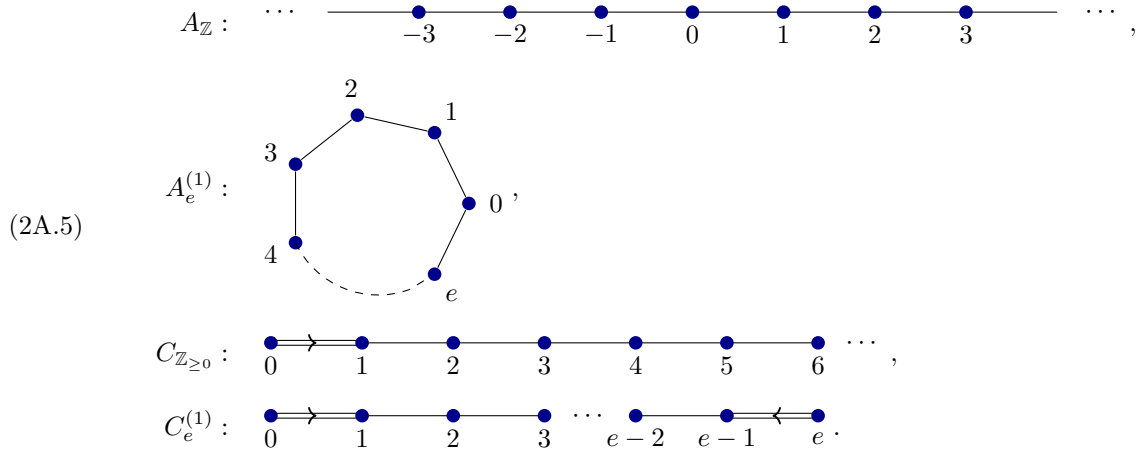
\begin{gather}\label{E:Quivers}
\begin{aligned}
A_{\Z}:\quad&
\begin{tikzpicture}[scale=1.2,anchorbase]
\draw(-4,0)--(4,0);
\foreach \x in {-3,...,3} {
\node[circle,inner sep=1.8pt,fill=DarkBlue] (\x) at (\x,0){};
\node at (\x,-0.25){$\x$};
}
\node at (-4.5,0){$\cdots$};
\node at (4.5,0){$\cdots$};
\end{tikzpicture},
\\[1mm]
A^{(1)}_{e}:\quad&
\begin{tikzpicture}[scale=1.2,anchorbase]
\foreach \r [remember=\r as \rr] in {0,...,4} {
\node[circle,inner sep=1.8pt,fill=DarkBlue] (\r) at (360/7*\r:1){};
\node at (360/7*\r:1.3){$\r$};
\ifnum\r>0\draw(\rr)--(\r);\fi
}
\node[circle,inner sep=1.8pt,fill=DarkBlue] (6) at (360/7*6:1){};
\node at (360/7*6:1.3){$e$};
\draw[dashed](4) arc [start angle=205, end angle=290, radius=1] -- (6);
\draw(6)--(0);
\end{tikzpicture},
\\[1mm]
C_{\N}:\quad&
\begin{tikzpicture}[scale=1.2,anchorbase]
\draw[directed=0.5,double,double distance=0.5mm](0,0)--(1,0);
\draw(1,0)--(6,0);
\node at (6.5,0){$\cdots$};
\foreach \x in {0,...,6} {
\node[circle,inner sep=1.8pt,fill=DarkBlue] (\x) at (\x,0){};
\node at (\x,-0.25){$\x$};
}
\end{tikzpicture},
\\[1mm]
C^{(1)}_{e}:\quad&
\begin{tikzpicture}[scale=1.2,anchorbase]
\draw[directed=0.5,double,double distance=0.5mm](0,0)--(1,0);
\draw[directed=0.5,double,double distance=0.5mm](6,0)--(5,0);
\draw(1,0)--(3,0);
\draw(4,0)--(5,0);
\node at (3.5,0){$\cdots$};
\foreach \x in {0,...,3} {
\node[circle,inner sep=1.8pt,fill=DarkBlue] (\x) at (\x,0){};
\node at (\x,-0.25){$\x$};
}
\foreach \x [evaluate=\x as \c using {int(6-\x)}] in {4,5} {
\node[circle,inner sep=1.8pt,fill=DarkBlue] (\x) at (\x,0){};
\node at (\x,-0.25){$e-\c$};
}
\node[circle,inner sep=1.8pt,fill=DarkBlue] (6) at (6,0){};
\node at (6,-0.25){$e$};
\end{tikzpicture}.
\end{aligned}
\end{gather}
When we refer to these quivers then we fix an orientation on the simply laced edges. Throughout this paper, affine type $A$ and $C$ will always refer
to one of these quivers even though, strictly speaking, $A_{\Z}$ and $C_{\N}$ are not affine quivers.

An explicit example of an oriented quiver is:
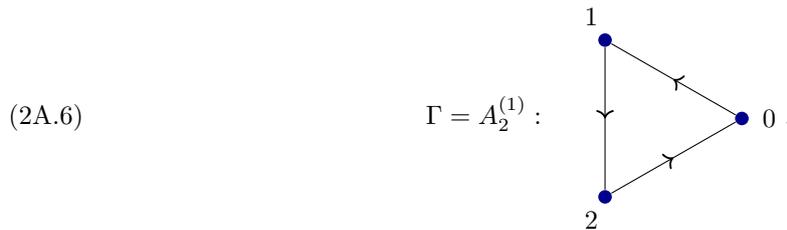
\begin{gather}\label{E:ExampleQuiver}
\Gamma=A^{(1)}_{2}:\quad
\begin{tikzpicture}[scale=1.2,anchorbase]
\node[circle,inner sep=1.8pt,fill=DarkBlue] (0) at (360/3*0:1){};
\node at (360/3*0:1.3){$0$};
\node[circle,inner sep=1.8pt,fill=DarkBlue] (1) at (360/3*1:1){};
\node at (360/3*1:1.3){$1$};
\node[circle,inner sep=1.8pt,fill=DarkBlue] (2) at (360/3*2:1){};
\node at (360/3*2:1.3){$2$};
\draw[directed=0.5](0)--(1);
\draw[directed=0.5](1)--(2);
\draw[directed=0.5](2)--(0);
\end{tikzpicture}
.
\end{gather}
We will use this quiver in several examples below.
\end{Example}

\begin{Remark}
The following quivers are also included in the general theory developed
in this paper:
\begin{gather}\label{Eq:FunnyQuivers}
\begin{aligned}
B_{\N}:\quad&
\begin{tikzpicture}[scale=1.2,anchorbase]
\draw[directed=0.5,double,double distance=0.5mm](1,0)--(0,0);
\draw(1,0)--(6.5,0);
\node at (6.8,0){$\cdots$};
\foreach \x in {0,...,6} {
\node[circle,inner sep=1.8pt,fill=DarkBlue] (\x) at (\x,0){};
\node at (\x,-0.25){$\x$};
}
\end{tikzpicture},
\\[1mm]
A^{(2)}_{2\cdot e}:\quad&
\begin{tikzpicture}[scale=1.2,anchorbase]
\draw[directed=0.5,double,double distance=0.5mm](0,0)--(1,0);
\draw[directed=0.5,double,double distance=0.5mm](5,0)--(6,0);
\draw(1,0)--(3,0);
\draw(4,0)--(5,0);
\node at (3.5,0){$\cdots$};
\foreach \x in {0,...,3} {
\node[circle,inner sep=1.8pt,fill=DarkBlue] (\x) at (\x,0){};
\node at (\x,-0.25){$\x$};
}
\foreach \x [evaluate=\x as \c using {int(6-\x)}] in {4,5} {
\node[circle,inner sep=1.8pt,fill=DarkBlue] (\x) at (\x,0){};
\node at (\x,-0.25){$e-\c$};
}
\node[circle,inner sep=1.8pt,fill=DarkBlue] (6) at (6,0){};
\node at (6,-0.25){$e$};
\end{tikzpicture},
\\[1mm]
D^{(2)}_{e+1}:\quad&
\begin{tikzpicture}[scale=1.2,anchorbase]
\draw[directed=0.5,double,double distance=0.5mm](1,0)--(0,0);
\draw[directed=0.5,double,double distance=0.5mm](5,0)--(6,0);
\draw(1,0)--(3,0);
\draw(4,0)--(5,0);
\node at (3.5,0){$\cdots$};
\foreach \x in {0,...,3} {
\node[circle,inner sep=1.8pt,fill=DarkBlue] (\x) at (\x,0){};
\node at (\x,-0.25){$\x$};
}
\foreach \x [evaluate=\x as \c using {int(6-\x)}] in {4,5} {
\node[circle,inner sep=1.8pt,fill=DarkBlue] (\x) at (\x,0){};
\node at (\x,-0.25){$e-\c$};
}
\node[circle,inner sep=1.8pt,fill=DarkBlue] (6) at (6,0){};
\node at (6,-0.25){$e$};
\end{tikzpicture}
.
\end{aligned}
\end{gather}
Cellular bases for the weighted KLRW algebras associated to these quivers are constructed in \cite{MaTu-klrw-algebras-bad}.
\end{Remark}

Let $\Sym$ be the symmetric group on $\set{1,2,\dots,n}$, viewed as a Coxeter group via the presentation
\begin{gather*}
\Sym=\<s_{1},\dots,s_{n-1}|s_{j}^{2}=1,
s_{j}s_{k}=s_{k}s_{j}\text{ if }|j-k|>1,
s_{j}s_{j+1}s_{j}=s_{j+1}s_{j}s_{j+1}\>,
\end{gather*}
for all admissible $j,k$. Let $s_{j}=(j,j+1)$ and let
$(k,l)\in\Sym$ be the transposition that swaps $k$ and $l$.

Let $Q^{+}=\bigoplus_{i\in I}\N\alpha_{i}$ be the \emph{positive root lattice} of the Kac--Moody algebra determined by $\Gamma$, where $\set{\alpha_{i}|i\in I}$ are the \emph{simple roots}.
The symmetric group $\Sym$ acts on $I^{n}$ by place permutations. The
\emph{height} of
$\beta=\sum_{i\in I}b_{i}\alpha_{i}\in Q^{+}$ is $\height\beta=\sum_{i\in I}b_{i}\geq 0$. Let
$Q^{+}_{n}=\set{\beta\in Q^{+}|\height\beta=n}$. Each
$\beta\in Q^{+}_{n}$ determines the $\Sym$-orbit
$I^\beta=\set{\bi\in I^{n}|\beta=\sum_{k=1}^{n}\alpha_{i_{k}}}$ and
$I^{n}=\bigcup_{\beta\in Q^{+}_{n}}I^{\beta}$. Finally, let $\<{}_{-},{}_{-}\>$ be the Cartan pairing
associated to the quiver $\Gamma$.

\begin{Example}\label{Ex:Beta}
Consider the quiver $A_{I}$ for $I=\Z$. Let $\set{e_{i}|i\in I}$ be the standard
basis of $\R^{I}$, considered as a vector space with inner product
determined by $\<e_{i},e_{j}\>=\delta_{ij}$. Then the simple roots
$\set{\alpha_{i}|i\in I}$ can be defined as $\alpha_{i}=e_{i}-e_{i+1}$ for
$i\in I$. Taking $n=2$, we have
$Q^{+}_{2}=\{\beta_{ij}=\alpha_{i}+\alpha_{j},\beta_{i}=2\alpha_{i}|i,j\in I,i\neq j\}$,
$I^{\beta_{ij}}=\set{(i,j),(j,i)|i,j\in I,i\neq j}$ and
$I^{\beta_{i}}=\set{(i,i)|i\in I}$.
\end{Example}


\subsection{Weighted KLRW diagrams}\label{SS:Diagrams}


The main algebras considered in this paper are defined in terms of
weighted KLRW diagrams, which are the subject of this section.

\begin{Notation}\label{N:ReadingDiagrams}
In illustrations we read diagrams from bottom to top, and the concatenation $E\circ D$ of two diagrams will be viewed as stacking $E$ on top of $D$:
\begin{gather*}
E\circ D
=
\begin{tikzpicture}[scale=1.2,anchorbase,smallnodes,rounded corners]
\node[rectangle,draw,minimum width=0.56cm,minimum height=0.56cm,ultra thick] at(0,0){\raisebox{-0.05cm}{$D$}};
\node[rectangle,draw,minimum width=0.56cm,minimum height=0.56cm,ultra thick] at(0,0.5){\raisebox{-0.05cm}{$E$}};
\end{tikzpicture}
.
\end{gather*}
In particular, left actions and left modules are given by acting from the top.
\end{Notation}

A \emph{string} is a smooth embedding
$\str\map{[0,1]}\R\times[0,1]$ such that $\str(t)\in\R\times\set{t}$, for $t\in[0,1]$. For $i\in I$, an \emph{$i$-string} is a string labeled by $i$. A labeled \emph{string
diagram} is an embedding of finitely many $i$-strings, for possibly
different $i\in I$, in $\R^{2}$ such that each point on these strings has
a local neighborhood that is of one of the following two forms:
\begin{gather}\label{E:GoodCrossings}
\begin{tikzpicture}[scale=1.2,anchorbase,smallnodes,rounded corners]
\draw[solid] (0,0)node[below]{$i$}node[below]{$\phantom{j}$} to (0,0.5)node[above,yshift=-1pt]{$\phantom{i}$};
\end{tikzpicture}
\qquad\text{or}\qquad
\begin{tikzpicture}[scale=1.2,anchorbase,smallnodes,rounded corners]
\draw[solid] (0,0)node[below]{$i$} to (0.5,0.5);
\draw[solid] (0.5,0)node[below]{$j$} to (0,0.5)node[above,yshift=-1pt]{$\phantom{i}$};
\end{tikzpicture}
\qquad\text{for }i,j\in I.
\end{gather}
A \emph{crossing} in a diagram is a point where two strings
intersect. The right-hand diagram in \autoref{E:GoodCrossings} shows how we draw crossings.

A \emph{dot} on a string is a distinguished point on the
string that is not on any crossing or on either endpoint of the string.
We illustrate dots on strings as follows:
\begin{gather*}
\begin{tikzpicture}[scale=1.2,anchorbase,smallnodes,rounded corners]
\draw[solid,dot] (0,0)node[below]{$i$} to (0,0.5)node[above,yshift=-1pt]{$\phantom{i}$};
\end{tikzpicture}.
\end{gather*}
If $\str\map{[0,1]}\R\times[0,1]$ is a string with
$\str(t)=\bigl(\str^{\prime}(t),t\bigr)$ and $\sigma\in\R$, then the
\emph{$\sigma$-shift} of $\str$ is the
string $\str+\sigma\map{[0,1]}\R\times[0,1]$ given
by $(\str+\sigma)(t)=\bigl(\str^{\prime}(t)+\sigma,t\bigr)$, for $t\in[0,1]$.

We apply this terminology below to solid, ghost and red strings, which
we now define. Before we can do this we fix notation to ensure that the
boundary points of our strings are distinct.

\begin{Definition}\label{D:DataShifts}
Recall that we have fixed $n$, $\ell$ and a quiver $\Gamma=(I,E)$.
\begin{enumerate}

\item A \emph{solid positioning} is an $n$-tuple
$\bx=(x_{1},\dots,x_{n})\in\R^{n}$.

\item A \emph{ghost shift} for $\Gamma$ is a function
$\bsig\map{E}{\R_{\neq 0}},\epsilon\mapsto\sigma_{\epsilon}$.

\item A \emph{charge}, or \emph{red positioning}, is a tuple
$\charge=(\kappa_{1},\dots,\kappa_{\ell})\in\R^{\ell}$ such that
$\kappa_{1}<\dots<\kappa_{\ell}$.

\item A \emph{loading} for $(\Gamma,\bsig)$ is a pair $(\charge,\bx)$ where
$\charge$ is a charge, $\bx$ is a solid positioning and
the numbers $x_{i}$, $x_{j}+|\sigma_{\epsilon}|$ and $\kappa_{k}$ are pairwise
distinct, where $1\leq i,j\leq n$, $\epsilon\in E$, $1\leq k\leq\ell$.

\end{enumerate}
\end{Definition}

As we will see shortly, we use $\bx$, $\bsig$ and $\charge$ to determine
the boundary points of solid, ghost and red strings in the diagrams that
we consider.

There are two extremal cases of ghost shifts that play important roles: the \emph{infinitesimal case}, where $\bsig_{\epsilon}=\varepsilon$, and the \emph{asymptotic case}, where $\bsig_{\epsilon}=1/\varepsilon$, where $0<\varepsilon\ll 1$ and all $\epsilon\in E$.

\begin{Remark}\label{R:Weighting}
In \cite{We-weighted-klr} the ghost shift $\bsig$ is called a
\emph{weighting}, similar to the corresponding terminology from graph theory. We draw weighted graphs where the weights are the $\sigma_{\epsilon}$.
\end{Remark}

\begin{Notation}\label{N:Rho}
Unless otherwise stated, we fix $\bsig$ and $\brho=(\rho_{1},\dots,\rho_{\ell})\in I^{\ell}$. (Set $\rho=\rho_{1}$.) Everything below depends on these choices.
\end{Notation}

There is a lot of notation in the following definition of weighted KLRW diagrams. The key point is that we have solid, ghost and red strings with the coordinates of the endpoints being specified by the loadings, which we will fix below.

\begin{Definition}\label{D:WeightedKLRW}({See \cite[Definition 2.3]{We-weighted-klr}, \cite[Definition 4.1]{We-rouquier-dia-algebra}}.)
\label{D:WebsterDiagram}
Fix a pair $(\Gamma,\bsig)$, where $\Gamma$ is a quiver and $\bsig$ is a ghost shift. Suppose
that $\bi\in I^{n}$ and that $(\charge,\bx)$ and
$(\chargetwo,\by)$ are loadings for $(\Gamma,\bsig)$, with $(\charge,\bx)$ and $(\chargetwo,\by)$
being the loadings for the bottom and top, respectively. A
\emph{weighted KLRW diagram} $D$ of
type $(\charge,\bx)\text{-}(\chargetwo,\by)$ and residue $\bi$
is a string diagram consisting of:

\begin{enumerate}

\item \emph{Solid} strings $\str_{1},\dots,\str_{n}$ such that
$\str_{k}$ is an $i_{k}$-string with $\str_k(0)=(x_{k},0)$ and
$\str_{k}(1)=(y_{w(k)},1)$, for some $w\in\Sym$ and $1\leq k\leq n$.

\item For each edge $\epsilon:i\rightsquigarrow j\in E$ with $\sigma_{\epsilon}>0$, every
solid $i$-string $\str$ has a \emph{ghost} $i$-string
$\str[g]_{\epsilon}=\str+\sigma_{\epsilon}$. Similarly, for each edge $\epsilon^{\prime}:k\rightsquigarrow i\in E$ with $\sigma_{\epsilon^{\prime}}<0$, every
solid $i$-string $\str$ has a \emph{ghost} $i$-string and $\str[g]_{\epsilon^{\prime}}=\str-\sigma_{\epsilon^{\prime}}$.

\item \emph{Red} strings $\str[r]_{1},\dots,\str[r]_{\ell}$ such that
$\str[r]_{k}$ is a $\rho_{k}$-string with
$\str[r]_{k}(t)=(t\nu_{k}+(1-t)\kappa_{k},t)$,
for $t\in[0,1]$ and $1\leq k\leq\ell$.

\item Solid strings can be decorated with finitely many dots, and
ghost strings with finitely many ghost dots, such that a dot appears
at position $\str(t)$ if and only if a ghost dot appear at position $\str[g]_{\epsilon}(t)$ on the corresponding ghost string, for each relevant edge $\epsilon\in E$.

\end{enumerate}
\end{Definition}

We will usually simply call a weighted KLRW diagram a \emph{diagram}.
We warn the reader that diagrams have a left-right bias in that ghost are always shifted to the right.

\begin{Remark}
Recall that in this paper $i\Rightarrow j$ and $i\Rrightarrow j$ count as a single edge,
so \autoref{D:WeightedKLRW} gives only one ghost string for such edges.
\end{Remark}

\begin{Remark}
\autoref{D:DataShifts} ensures that the endpoints of the solid, ghost and red strings in a diagram are distinct. In particular, ghost shifts like the following are not allowed:
\begin{gather*}
\begin{tikzpicture}[scale=1.2,anchorbase]
\node[circle,inner sep=1.8pt,fill=DarkBlue] (0) at (0,0){};
\node[circle,inner sep=1.8pt,fill=DarkBlue] (1) at (1,0){};
\node[circle,inner sep=1.8pt,fill=DarkBlue] (2) at (-1,0){};
\draw[directed=0.5](0) to node[above,xshift=0.05cm,yshift=-1pt]{$\sigma>0$} (1);
\draw[directed=0.5](0) to node[above,yshift=-1pt]{$\sigma>0$} (2);
\end{tikzpicture}
.
\end{gather*}
Unlike \cite{We-weighted-klr}, we do not allow ghost shifts to be zero because this would mean that solid strings overlap with their ghost strings, contrary to \autoref{D:DataShifts}.
This said, the zero ghost shift case is captured by the infinitesimal case.
\end{Remark}

To help distinguish between the different types of strings in diagrams, we draw solid strings as in \autoref{E:GoodCrossings}, ghost strings as dashed gray strings (with their labels illustrated at the top)
and red strings as thick red strings, {\cf} \autoref{E:DrawingGhost}. By
\autoref{D:WeightedKLRW}, red strings do not cross each other because,
in contrast to solid and ghost strings, we do not allow a permutation of
their endpoints. Consequently, locally, a diagram is always of one of the following forms.
\begin{gather}\label{E:DrawingGhost}
\begin{tikzpicture}[scale=1.2,anchorbase,smallnodes,rounded corners]
\draw[solid] (0,0)node[below]{$i$}node[below]{$\phantom{j}$} to (0,0.5)node[above,yshift=-1pt]{$\phantom{i}$};
\end{tikzpicture}
,\,
\begin{tikzpicture}[scale=1.2,anchorbase,smallnodes,rounded corners]
\draw[ghost] (0,0)node[below]{$\phantom{i}$}node[below]{$\phantom{j}$} to (0,0.5)node[above,yshift=-1pt]{$i$};
\end{tikzpicture}
,\,
\begin{tikzpicture}[scale=1.2,anchorbase,smallnodes,rounded corners]
\draw[solid,dot] (0,0)node[below]{$i$}node[below]{$\phantom{j}$} to (0,0.5)node[above,yshift=-1pt]{$\phantom{i}$};
\end{tikzpicture}
,\,
\begin{tikzpicture}[scale=1.2,anchorbase,smallnodes,rounded corners]
\draw[ghost,dot] (0,0)node[below]{$\phantom{i}$}node[below]{$\phantom{j}$} to (0,0.5)node[above,yshift=-1pt]{$i$};
\end{tikzpicture}
,\,
\begin{tikzpicture}[scale=1.2,anchorbase,smallnodes,rounded corners]
\draw[redstring] (0,0)node[below]{$\rho$}node[below]{$\phantom{j}$} to (0,0.5)node[above,yshift=-1pt]{$\phantom{i}$};
\end{tikzpicture}
,\,
\begin{tikzpicture}[scale=1.2,anchorbase,smallnodes,rounded corners]
\draw[solid] (0,0)node[below]{$i$} to (0.5,0.5);
\draw[solid] (0.5,0)node[below]{$j$} to (0,0.5)node[above,yshift=-1pt]{$\phantom{i}$};
\end{tikzpicture}
,\,
\begin{tikzpicture}[scale=1.2,anchorbase,smallnodes,rounded corners]
\draw[ghost] (0,0)node[below]{$\phantom{i}$} to (0.5,0.5)node[above,yshift=-1pt]{$i$};
\draw[ghost] (0.5,0)node[below]{$\phantom{i}$} to (0,0.5)node[above,yshift=-1pt]{$j$};
\end{tikzpicture}
,\,
\begin{tikzpicture}[scale=1.2,anchorbase,smallnodes,rounded corners]
\draw[ghost] (0.5,0)node[below]{\phantom{i}} to (0,0.5)node[above,yshift=-1pt]{$j$};
\draw[solid] (0,0)node[below]{$i$} to (0.5,0.5);
\end{tikzpicture}
,\,
\begin{tikzpicture}[scale=1.2,anchorbase,smallnodes,rounded corners]
\draw[ghost] (0,0)node[below]{\phantom{i}} to (0.5,0.5)node[above,yshift=-1pt]{$i$};
\draw[solid] (0.5,0)node[below]{$j$} to (0,0.5)node[above,yshift=-1pt]{$\phantom{i}$};
\end{tikzpicture}
,\,
\begin{tikzpicture}[scale=1.2,anchorbase,smallnodes,rounded corners]
\draw[solid] (0,0)node[below]{$i$} to (0.5,0.5);
\draw[redstring] (0.5,0)node[below]{$\rho$} to (0,0.5)node[above,yshift=-1pt]{$\phantom{i}$};
\end{tikzpicture}
,\,
\begin{tikzpicture}[scale=1.2,anchorbase,smallnodes,rounded corners]
\draw[solid] (0.5,0)node[below]{$i$} to (0,0.5)node[above,yshift=-1pt]{$\phantom{i}$};
\draw[redstring] (0,0)node[below]{$\rho$} to (0.5,0.5);
\end{tikzpicture}
,\,
\begin{tikzpicture}[scale=1.2,anchorbase,smallnodes,rounded corners]
\draw[ghost] (0,0)node[below]{$\phantom{i}$} to (0.5,0.5)node[above,yshift=-1pt]{$i$};
\draw[redstring] (0.5,0)node[below]{$\rho$} to (0,0.5)node[above,yshift=-1pt]{$\phantom{i}$};
\end{tikzpicture}
,\,
\begin{tikzpicture}[scale=1.2,anchorbase,smallnodes,rounded corners]
\draw[ghost] (0.5,0)node[below]{$\phantom{i}$} to (0,0.5)node[above,yshift=-1pt]{$i$};
\draw[redstring] (0,0)node[below]{$\rho$} to (0.5,0.5);
\end{tikzpicture}
.
\end{gather}
Here, and throughout the paper, we put the labels for ghost strings above the string, and the labels for the red and solid strings below the string, to help distinguish these strings.

\begin{Remark}
For those readers who are reading a black-and-white version
of the paper, all strings are distinguishable by their thickness and if
they are solid or dashed. The colors
used in this paper are not essential
and are only used as a visual aid.
\end{Remark}

\begin{Example}
Take $n=1=\ell$ and $\charge=\chargetwo$ and let $0<\varepsilon\ll 1$.
Let $\Gamma$ be the quiver
\begin{tikzpicture}[anchorbase]
\node[circle,inner sep=1.8pt,fill=DarkBlue] (0) at (0,0){};
\node[circle,inner sep=1.8pt,fill=DarkBlue] (1) at (1,0){};
\draw[directed=0.5](0) to (1);
\end{tikzpicture}. Here are three diagrams for $\Gamma$ using the specified choice of ghost shift:
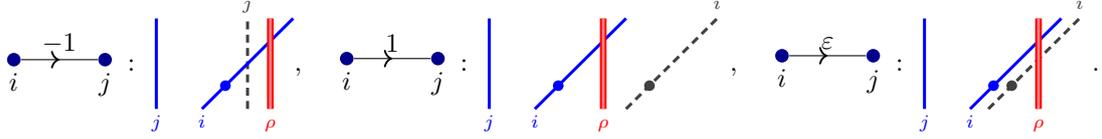
\begin{gather*}
\begin{tikzpicture}[scale=1.2,anchorbase]
\node[circle,inner sep=1.8pt,fill=DarkBlue] (0) at (0,0){};
\node at (0,-0.25){$i$};
\node[circle,inner sep=1.8pt,fill=DarkBlue] (1) at (1,0){};
\node at (1,-0.28){$j$};
\draw[directed=0.5](0) to node[above,yshift=-1pt]{$-1$} node[below]{$\phantom{-1}$} (1);
\end{tikzpicture}
:
\begin{tikzpicture}[scale=1.2,anchorbase,smallnodes,rounded corners]
\draw[ghost](0.5,0)node[below]{$\phantom{i}$}--++(0,1)node[above,yshift=-1pt]{$j$};
\draw[solid,dot=0.25](0,0)node[below]{$i$}--++(1,1);
\draw[solid](-0.5,0)node[below]{$j$}--++(0,1);
\draw[redstring](0.75,0)node[below]{$\rho$}--++(0,1);
\end{tikzpicture}
,\quad
\begin{tikzpicture}[scale=1.2,anchorbase]
\node[circle,inner sep=1.8pt,fill=DarkBlue] (0) at (0,0){};
\node at (0,-0.25){$i$};
\node[circle,inner sep=1.8pt,fill=DarkBlue] (1) at (1,0){};
\node at (1,-0.28){$j$};
\draw[directed=0.5](0) to node[above,yshift=-1pt]{$1$} node[below]{$\phantom{1}$} (1);
\end{tikzpicture}
:
\begin{tikzpicture}[scale=1.2,anchorbase,smallnodes,rounded corners]
\draw[ghost,dot=0.25](1,0)node[below]{$\phantom{i}$}--++ (1,1)node[above,yshift=-1pt]{$i$};
\draw[solid,dot=0.25](0,0)node[below]{$i$}--++(1,1);
\draw[solid](-0.5,0)node[below]{$j$}--++(0,1);
\draw[redstring](0.75,0)node[below]{$\rho$}--++(0,1);
\end{tikzpicture}
,\quad
\begin{tikzpicture}[scale=1.2,anchorbase]
\node[circle,inner sep=1.8pt,fill=DarkBlue] (0) at (0,0){};
\node at (0,-0.25){$i$};
\node[circle,inner sep=1.8pt,fill=DarkBlue] (1) at (1,0){};
\node at (1,-0.28){$j$};
\draw[directed=0.5](0) to node[above,yshift=-1pt]{$\varepsilon$} node[below]{$\phantom{\epsilon}$} (1);
\end{tikzpicture}
:
\begin{tikzpicture}[scale=1.2,anchorbase,smallnodes,rounded corners]
\draw[ghost,dot=0.25](0.2,0)node[below]{$\phantom{i}$}--++ (1,1)node[above,yshift=-1pt]{$i$};
\draw[solid,dot=0.25](0,0)node[below]{$i$}--++(1,1);
\draw[solid](-0.5,0)node[below]{$j$}--++(0,1);
\draw[redstring](0.75,0)node[below]{$\rho$}--++(0,1);
\end{tikzpicture}
.
\end{gather*}
Reading from left to right,
$\sigma_{\epsilon}=-1$, $\sigma_{\epsilon}=1$ and $\sigma_{\epsilon}=\varepsilon$ are the ghost shifts for the edge $\epsilon:i\to j\in E$. In particular, note that the $i$-strings have a ghost when $\sigma>0$ and the $j$-strings have ghosts when $\sigma<0$.
When a vertex in $\Gamma$ is the tail of more than one edge we obtain more ghost strings.
For example,
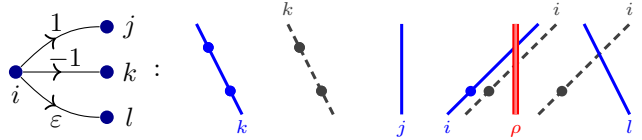
\begin{gather*}
\begin{tikzpicture}[scale=1.2,anchorbase]
\node[circle,inner sep=1.8pt,fill=DarkBlue] (0) at (0,0){};
\node at (0,-0.25){$i$};
\node[circle,inner sep=1.8pt,fill=DarkBlue] (1) at (1,0.5){};
\node at (1.25,0.5){$j$};
\node[circle,inner sep=1.8pt,fill=DarkBlue] (2) at (1,0){};
\node at (1.25,0){$k$};
\node[circle,inner sep=1.8pt,fill=DarkBlue] (3) at (1,-0.5){};
\node at (1.25,-0.5){$l$};
\draw[directed=0.5](0) to[out=45,in=180] node[above,yshift=-1pt]{$1$} (1);
\draw[directed=0.5](0) to node[above,xshift=0.05cm,yshift=-0.075cm]{$-1$} (2);
\draw[directed=0.5](0) to[out=315,in=180] node[below]{$\varepsilon$} (3);
\end{tikzpicture}
:\quad
\begin{tikzpicture}[scale=1.2,anchorbase,smallnodes,rounded corners]
\draw[ghost,dot=0.25](0.2,0)node[below]{$\phantom{i}$}--++ (1,1)node[above,yshift=-1pt]{$i$};
\draw[ghost,dot=0.25](1,0)node[below]{$\phantom{i}$}--++ (1,1)node[above,yshift=-1pt]{$i$};
\draw[ghost,dot=0.25,dot=0.75](-1.25,0)node[below]{$\phantom{i}$}--++(-0.5,1)node[above,yshift=-1pt]{$k$};
\draw[solid,dot=0.25](0,0)node[below]{$i$}--++(1,1);
\draw[solid](-0.5,0)node[below]{$j$}--++(0,1);
\draw[solid,dot=0.25,dot=0.75](-2.25,0)node[below]{$k$}--++(-0.5,1);
\draw[solid](2,0)node[below]{$l$}--++(-0.5,1);
\draw[redstring](0.75,0)node[below]{$\rho$}--++(0,1);
\end{tikzpicture}
,
\end{gather*}
is a diagram for the illustrated pair
$(\Gamma,\bsig)$ and the corresponding loadings.
\end{Example}

\begin{Notation}\label{N:ExampleGhost}
Unless otherwise stated, in examples we usually take
$\sigma_{\epsilon}=1$.
\end{Notation}

The following simple, yet important, classes of diagrams are used throughout this paper.

\begin{Definition}
A \emph{idempotent diagram} is any diagram
with no dots and no crossings and where the $x$-coordinates of the points on each string are constant. A \emph{straight line diagram} is any diagram
with no dots and no crossings.
\end{Definition}

Less formally, the strings in idempotent diagrams are always vertical. The strings in straight line diagrams are not necessarily straight, however, up to isotopy we can always assume that they are straight.

\begin{Example}\label{E:StraightLine}
Prototypical examples of idempotent straight line diagrams are:
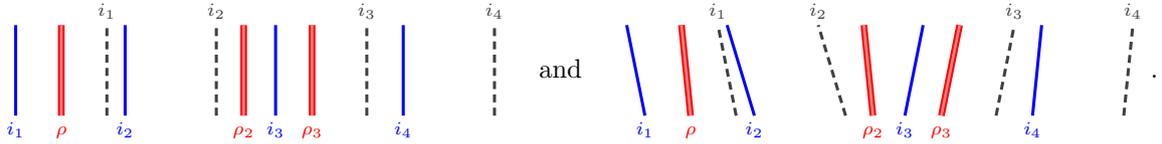
\begin{gather*}
\begin{tikzpicture}[scale=1.2,anchorbase,smallnodes,rounded corners]
\draw[ghost](0,0)node[below]{$\phantom{i}$}--++(0,1)node[above,yshift=-1pt]{$i_{1}$};
\draw[ghost](1.2,0)node[below]{$\phantom{i}$}--++(0,1)node[above,yshift=-1pt]{$i_{2}$};
\draw[ghost](2.85,0)node[below]{$\phantom{i}$}--++(0,1)node[above,yshift=-1pt]{$i_{3}$};
\draw[ghost](4.25,0)node[below]{$\phantom{i}$}--++(0,1)node[above,yshift=-1pt]{$i_{4}$};
\draw[solid](-1,0)node[below]{$i_{1}$}--++(0,1) node[above,yshift=-1pt]{$\phantom{i_{1}}$};
\draw[solid](0.2,0)node[below]{$i_{2}$}--++(0,1);
\draw[solid](1.85,0)node[below]{$i_{3}$}--++(0,1);
\draw[solid](3.25,0)node[below]{$i_{4}$}--++(0,1);
\draw[redstring](-0.5,0)node[below]{$\rho$}--++(0,1);
\draw[redstring](1.5,0)node[below]{$\rho_{2}$}--++(0,1);
\draw[redstring](2.25,0)node[below]{$\rho_{3}$}--++(0,1);
\end{tikzpicture}
\quad\text{and}\quad
\begin{tikzpicture}[scale=1.2,anchorbase,smallnodes,rounded corners]
\draw[ghost](0,0)node[below]{$\phantom{i}$}--++(-0.2,1)node[above,yshift=-1pt]{$i_{1}$};
\draw[ghost](1.2,0)node[below]{$\phantom{i}$}--++(-0.3,1)node[above,yshift=-1pt]{$i_{2}$};
\draw[ghost](2.85,0)node[below]{$\phantom{i}$}--++(0.2,1)node[above,yshift=-1pt]{$i_{3}$};
\draw[ghost](4.25,0)node[below]{$\phantom{i}$}--++(0.1,1)node[above,yshift=-1pt]{$i_{4}$};
\draw[solid](-1,0)node[below]{$i_{1}$}--++(-0.2,1) node[above,yshift=-1pt]{$\phantom{i_{1}}$};
\draw[solid](0.2,0)node[below]{$i_{2}$}--++(-0.3,1);
\draw[solid](1.85,0)node[below]{$i_{3}$}--++(0.2,1);
\draw[solid](3.25,0)node[below]{$i_{4}$}--++(0.1,1);
\draw[redstring](-0.5,0)node[below]{$\rho$}--++(-0.1,1);
\draw[redstring](1.5,0)node[below]{$\rho_{2}$}--++(-0.1,1);
\draw[redstring](2.25,0)node[below]{$\rho_{3}$}--++(0.2,1);
\end{tikzpicture}
.
\end{gather*}
The left diagram is a both an idempotent and a straight line diagram. The right diagram is a straight line diagram but not an idempotent diagram.
\end{Example}

\begin{Definition}\label{D:Bjstrings}
Let $\Webab$ be the set of diagrams of
type $(\charge,\bx)\text{-}(\chargetwo,\by)$ and residue $\bi$ such
that the residues of the
strings at the top of the diagrams,
when read in order from left to right, are given by $\bj$.
Whenever $D\in\Webab$ we
assume that $\bj$ is the permutation of $\bi$ determined by $D$.
\end{Definition}

\begin{Definition}
Two diagrams $D$ and $D^{\prime}$ in $\Webab$ are
equivalent if they differ by an isotopy, which is a smooth
deformation $\theta\map{[0,1]}\Webab$ such that
$\theta(0)=D$ and $\theta(1)=D^{\prime}$.
\end{Definition}

Note that isotopies cannot separate
strings, create crossings, change the number of dots on any string or change the residue of any string. We abuse notation and write
$\Webab$ for the corresponding set of equivalence
classes under isotopy. Note that, up to isotopy, there is a unique idempotent diagram
$\1_{(\charge,\bx),\bi}\in\Webab[(\charge,\bx),\bi][(\charge,\bx),\bi]$
that has no dots and no crossings.

\begin{Example}
Let $\Gamma$ be a quiver with an edge $i\to j$ and take $n=1=\ell$. Then we have one solid string (and also one ghost since $\sigma=1$ by \autoref{N:ExampleGhost}), and one red string, which we set to be at position $\charge=\chargetwo=(0)$. Let $\bx=(-1)$ and $\brho=(i)=\bi$. Then
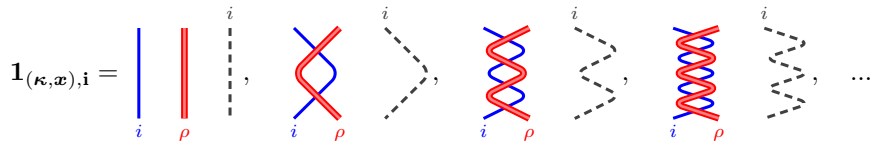
\begin{gather*}
\1_{(\charge,\bx),\bi}
=
\begin{tikzpicture}[scale=1.2,anchorbase,smallnodes,rounded corners]
\draw[ghost](0.5,1)node[above,yshift=-1pt]{$i$}--++(0,-1)node[below]{$\phantom{i}$};
\draw[solid](-0.5,1)--++(0,-1)node[below]{$i$};
\draw[redstring](0,0)node[below]{$\rho$}--++(0,1);
\end{tikzpicture}
,\quad
\begin{tikzpicture}[scale=1.2,anchorbase,smallnodes,rounded corners]
\draw[ghost](0.5,1)node[above,yshift=-1pt]{$i$}--++(0.5,-0.5)--++(-0.5,-0.5)node[below]{$\phantom{i}$};
\draw[solid](-0.5,1)--++(0.5,-0.5)--++(-0.5,-0.5)node[below]{$i$};
\draw[redstring](0,1)--++(-0.5,-0.5)--++(0.5,-0.5)node[below]{$\rho$};
\end{tikzpicture}
,\quad
\begin{tikzpicture}[scale=1.2,anchorbase,smallnodes,rounded corners]
\draw[ghost](0.5,1)node[above,yshift=-1pt]{$i$}--++(0.5,-0.25)--++(-0.5,-0.25)--++(0.5,-0.25)--++(-0.5,-0.25)node[below]{$\phantom{i}$};
\draw[solid](-0.5,1)--++(0.5,-0.25)--++(-0.5,-0.25)--++(0.5,-0.25)--++(-0.5,-0.25)node[below]{$i$};
\draw[redstring](0,1)--++(-0.5,-0.25)--++(0.5,-0.25)--++(-0.5,-0.25)--++(0.5,-0.25)node[below]{$\rho$};
\end{tikzpicture}
,\quad
\begin{tikzpicture}[scale=1.2,anchorbase,smallnodes,rounded corners]
\draw[ghost](0.5,1)node[above,yshift=-1pt]{$i$}--++(0.5,-0.166)--++(-0.5,-0.166)--++(0.5,-0.166)--++(-0.5,-0.166)--++(0.5,-0.166)--++(-0.5,-0.166)node[below]{$\phantom{i}$};
\draw[solid](-0.5,1)--++(0.5,-0.166)--++(-0.5,-0.166)--++(0.5,-0.166)--++(-0.5,-0.166)--++(0.5,-0.166)--++(-0.5,-0.166)node[below]{$i$};
\draw[redstring](0,1)--++(-0.5,-0.166)--++(0.5,-0.166)--++(-0.5,-0.166)--++(0.5,-0.166)--++(-0.5,-0.166)--++(0.5,-0.166)node[below]{$\rho$};
\end{tikzpicture}
,\quad
\dots
\end{gather*}
are examples of diagrams in $\Webab[(\chargetwo,\bx),\bi][(\chargetwo,\bx),\bi]$ that are not equivalent.
\end{Example}

Let $D\in\Webab[(\charge,\bx),\bi][(\chargetwo,\by),\bj]$ and $E\in\Webab[(\chargetwo,\by),\bj][(\chargethree,\bz),\bk]$
be diagrams. Then $E\circ D\in\Webab[(\charge,\bx),\bi][(\chargethree,\bz),\bk]$
is obtained by gluing $D$ under $E$ (see \autoref{N:ReadingDiagrams}) and then rescaling. If $D\in\Webab[(\charge,\bx),\bi][(\chargetwo,\by),\bj]$, then $D=\1_{(\chargetwo,\by),\bj}\circ D\circ\1_{(\charge,\bx),\bi}$.

\begin{Notation}
Straight line diagrams will allow us to show that we can fix a charge $\charge$, {\cf} \autoref{P:FixedCharge}. In order to simplify the notation, from \autoref{SS:WebsterAlgebras} onward we will simply write $\bx$ for a loading $(\charge,\bx)$.
\end{Notation}


\subsection{Weighted KLRW algebras}\label{SS:WebsterAlgebras}


Recall that we have fixed $\bsig$ and $\brho$ as in \autoref{N:Rho} (the ghost shift and the labels for the red strings). We additionally need:

\begin{Notation}
Fix a commutative integral domain $R$, for example $R=\Z$.
Throughout the rest of this section
we also fix $\beta\in Q^{+}_{n}$ of height $n$ as in \autoref{SS:Quiver} (the labels for the solid and ghost strings), and a finite nonempty set $X$ set of loadings, called the \emph{positioning}.
\end{Notation}

We define
\begin{gather*}
\Web_{\beta}^{\brho}(X)
=
\bigcup_{\bx,\by\in X}
\bigcup_{\bi,\bj\in I^{\beta}}
\Webaa.
\end{gather*}
In particular, $\1_{\bx,\bi}\in\Web_{\beta}^{\brho}(X)$, whenever
$\bx\in X$ and $\bi\in I^{\beta}$.

For $D\in\Webaa$ define $y_{r}D$ to be the diagram obtained from $D$
by concatenating with a dotted idempotent on top that has a dot on
the $r$th solid string and a ghost dot on the $r$th ghost string.
We extend this notation so that $f(y_{1},\dots,y_{n})D$ is
the evident linear combination of diagrams for any polynomial
$f(u_{1},\dots,u_{n})\in R[u_{1},\dots,u_{n}]$.

\begin{Example}
For example,
\begin{gather}\label{E:Unsteady}
y_{3}y_{4}^{2}\1_{\bx,\bi}
=
\begin{tikzpicture}[scale=1.2,anchorbase,smallnodes,rounded corners]
\draw[ghost](0,1)node[above,yshift=-1pt]{$i_{1}$}--++(0,-1) node[below]{$\phantom{i}$};
\draw[ghost](1.2,1)node[above,yshift=-1pt]{$i_{2}$}--++(0,-1) node[below]{$\phantom{i}$};
\draw[ghost,dot](2.85,1)node[above,yshift=-1pt]{$i_{3}$}--++(0,-1) node[below]{$\phantom{i}$};
\draw[ghost,dot=0.25,dot](4.25,1)node[above,yshift=-1pt]{$i_{4}$}--++(0,-1) node[below]{$\phantom{i}$};
\draw[solid](-1,1)node[above,yshift=-1pt]{$\phantom{i}$}--++(0,-1) node[below]{$i_{1}$};
\draw[solid](0.2,1)--++(0,-1) node[below]{$i_{2}$};
\draw[solid,dot](1.85,1)--++(0,-1) node[below]{$i_{3}$};
\draw[solid,dot=0.25,dot](3.25,1)--++(0,-1) node[below]{$i_{4}$};
\draw[redstring](-0.5,1)--++(0,-1) node[below]{$\rho_{1}$};
\draw[redstring](1.5,1)--++(0,-1) node[below]{$\rho_{2}$};
\draw[redstring](2.25,1)--++(0,-1) node[below]{$\rho_{3}$};
\end{tikzpicture}
,
\end{gather}
where $\bx$ and $\bi$ can be read-off from the illustration \autoref{E:Unsteady}.
\end{Example}

Let $\mathbf{d}\in(\N)^{e}$ be the symmetrizer of the Kac--Moody data associated to $\Gamma$.
For $i,j\in I$ fix polynomials $Q_{ij}(u,v)\in R[u,v]$, called \emph{$Q$-polynomials}, where $u$ and
$v$ are indeterminates (all of our variables appearing in polynomial rings will be indeterminates) of degrees
$2d_{i}$ and $2d_{j}$, such that:

\begin{enumerate}

\item We assume that $Q_{ii}(u,v)=0$ and that $Q_{ij}(u,v)$ is invertible if $i\ne j$ are not connected by an edge in $\Gamma$.

\item For $i\neq j$ we assume that $Q_{ij}(u,v)$
is homogeneous of degree $2\<\alpha_{i},\alpha_{j}\>$ and the coefficients of all monomials are units.

\item For $i\neq j$ we assume that $Q_{ij}(u,v)=Q_{ji}(v,u)$.

\end{enumerate}

Similar conditions appear in \cite[Section 3.2.3]{Ro-2-kac-moody}, or \cite{Ro-quiver-hecke}, and \cite[Section 2.1]{We-weighted-klr}.

\begin{Example}\label{Ex:QPoly}
For the quivers in \autoref{E:Quivers} standard choices for $i\neq j$ are
\begin{gather}\label{E:QPoly}
Q_{ij}(u,v)=
\left\{
\begin{array}{CC@{\quad}|@{\quad}CC}
$u-v$ & if $i\rightarrow j$, & $v-u$ & if $i\leftarrow j$,
\\
$u-v^{2}$ & if $i\Rightarrow j$, & $v-u^{2}$ & if $i\Leftarrow j$,
\\
$0$ & if $i=j$, & $1$ & otherwise.
\\
\end{array}
\right.
\end{gather}
We will use these choices when constructing
homogeneous (affine) cellular bases in \autoref{S:AffineCellular} and \autoref{S:TypeC}.
\end{Example}

Further, define polynomials $Q_{ijk}(u,v,w)\in R[u,v,w]$ by
\begin{gather*}
Q_{ijk}(u,v,w)=
\begin{cases*}
\frac{Q_{ij}(u,v)-Q_{kj}(u,w)}{w-v}& if $i=k$,\\
0& otherwise.
\end{cases*}
\end{gather*}
Below we abuse notation and write
$Q_{ij}(\y)D=Q_{ij}(y_{r},y_{s})D$ and
$Q_{ijk}(\y)D=Q_{ijk}(y_{r},y_{s},y_{t})D$ for the linear combination of diagrams obtained by putting dots on the corresponding strings $r,s$ and $t$ of residues $i_{r}=i$, $i_{s}=j$ and $i_{t}=k$, respectively.

We are almost ready to define the algebras that we are interested in.
Since each solid string can have several ghosts, the relations that we
use are \emph{multilocal} in the following sense:
We need to simultaneously
apply the relations both in local neighborhoods around the solid strings
and in the corresponding shifted local neighborhoods around the ghost strings.
We can only apply the relations if they can be simultaneously applied in
all neighborhoods to give new diagrams.

To ease notation we sometimes
omit solid strings or ghosts strings from relations or diagrams. In such cases the missing strings are implicit because solid and ghost strings always occur together.

\begin{Notation}\label{N:KLRWdata}
As in \autoref{D:DataShifts}, fix a quiver $\Gamma$, a ghost shift $\bsig$ and a charge $\charge$. In addition, fix $\beta\in Q^{+}_{n}$, polynomials $Q_{ij}(u,v)$ and a positioning set $X$. Even though our notation does not reflect this, the following algebras depend on all of these choices.
\end{Notation}

The following is our formulation of
\cite[Definition 2.4]{We-weighted-klr}, \cite[Definition 4.2]{We-rouquier-dia-algebra}. Recall that we write $i\rightsquigarrow j$ if $\Gamma$ has an edge from $i$ to $j$, of any multiplicity.

\begin{Definition}\label{D:RationalCherednik}
The \emph{weighted KLRW algebra} $\WA(X)$
is the unital associative $R$-algebra generated (as an algebra) by the diagrams
in $\Webab[\beta][\brho](X)$ with multiplication given by
\begin{gather*}
ED=
\begin{cases*}
E\circ D & if $D\in\Webaa$
and $E\in\Webaa[\by,\bj][\bz,\bk]$,
\\
0 & otherwise,
\end{cases*}
\end{gather*}
and subject to the following multilocal relations.
\begin{enumerate}

\item \label{I:Cross}
The \emph{dot sliding relations} hold, that is,
solid and ghost dots can pass through any crossing except:
\begin{gather}\label{R:DotCrossing}
\begin{tikzpicture}[scale=1.2,anchorbase,smallnodes,rounded corners]
\draw[solid](0.5,0.5)node[above,yshift=-1pt]{$\phantom{i}$}--(0,0) node[below]{$i$};
\draw[solid,dot=0.25](0,0.5)--(0.5,0) node[below]{$i$};
\end{tikzpicture}
-
\begin{tikzpicture}[scale=1.2,anchorbase,smallnodes,rounded corners]
\draw[solid](0.5,0.5)node[above,yshift=-1pt]{$\phantom{i}$}--(0,0) node[below]{$i$};
\draw[solid,dot=0.75](0,0.5)--(0.5,0) node[below]{$i$};
\end{tikzpicture}
=
\begin{tikzpicture}[scale=1.2,anchorbase,smallnodes,rounded corners]
\draw[solid](0,0.5)node[above,yshift=-1pt]{$\phantom{i}$}--(0,0) node[below]{$i$};
\draw[solid](0.5,0.5)--(0.5,0) node[below]{$i$};
\end{tikzpicture}
=
\begin{tikzpicture}[scale=1.2,anchorbase,smallnodes,rounded corners]
\draw[solid,dot=0.75](0.5,0.5)node[above,yshift=-1pt]{$\phantom{i}$}--(0,0) node[below]{$i$};
\draw[solid](0,0.5)--(0.5,0) node[below]{$i$};
\end{tikzpicture}
-
\begin{tikzpicture}[scale=1.2,anchorbase,smallnodes,rounded corners]
\draw[solid,dot=0.25](0.5,0.5)node[above,yshift=-1pt]{$\phantom{i}$}--(0,0) node[below]{$i$};
\draw[solid](0,0.5)--(0.5,0) node[below]{$i$};
\end{tikzpicture}
.
\end{gather}

\item \label{I:DoubleCross}
The \emph{Reidemeister II relations} hold except in the following cases:
\begin{gather}\label{R:SolidSolid}
\begin{tikzpicture}[scale=1.2,anchorbase,smallnodes,rounded corners]
\draw[solid](0,1)--++(0.5,-0.5)--++(-0.5,-0.5) node[below]{$i$};
\draw[solid](0.5,1)node[above,yshift=-1pt]{$\phantom{i}$}--++(-0.5,-0.5)--++(0.5,-0.5) node[below]{$i$};
\end{tikzpicture}
=0
.
\end{gather}
\begin{gather}\label{R:GhostSolid}
\begin{tikzpicture}[scale=1.2,anchorbase,smallnodes,rounded corners]
\draw[ghost](0,1)node[above,yshift=-1pt]{$i$}--++(0.5,-0.5)--++(-0.5,-0.5) node[below]{$\phantom{i}$};
\draw[solid](0.5,1)--++(-0.5,-0.5)--++(0.5,-0.5) node[below]{$j$};
\end{tikzpicture}
=Q_{ij}(\y)
\begin{tikzpicture}[scale=1.2,anchorbase,smallnodes,rounded corners]
\draw[ghost](0,1)node[above,yshift=-1pt]{$i$}--++(0,-1)node[below]{$\phantom{i}$};
\draw[solid](0.5,1)--++(0,-1)node[below]{$j$};
\end{tikzpicture}
,\quad
\begin{tikzpicture}[scale=1.2,anchorbase,smallnodes,rounded corners]
\draw[ghost](0.5,1)node[above,yshift=-1pt]{$i$}--++(-0.5,-0.5)--++(0.5,-0.5) node[below]{$\phantom{i}$};
\draw[solid](0,1)--++(0.5,-0.5)--++(-0.5,-0.5) node[below]{$j$};
\end{tikzpicture}
=Q_{ji}(\y)
\begin{tikzpicture}[scale=1.2,anchorbase,smallnodes,rounded corners]
\draw[ghost](0.5,1)node[above,yshift=-1pt]{$i$}--++(0,-1)node[below]{$\phantom{i}$};
\draw[solid](0,1)--++(0,-1)node[below]{$j$};
\end{tikzpicture}
\quad
\text{if $i\rightsquigarrow j$}
.
\end{gather}
\begin{gather}\label{R:RedSolid}
\begin{tikzpicture}[scale=1.2,anchorbase,smallnodes,rounded corners]
\draw[solid](0.5,1)node[above,yshift=-1pt]{$\phantom{i}$}--++(-0.5,-0.5)--++(0.5,-0.5) node[below]{$i$};
\draw[redstring](0,1)--++(0.5,-0.5)--++(-0.5,-0.5) node[below]{$i$};
\end{tikzpicture}
=
\begin{tikzpicture}[scale=1.2,anchorbase,smallnodes,rounded corners]
\draw[solid,dot](0.5,0)node[below]{$i$}--++(0,1)node[above,yshift=-1pt]{$\phantom{i}$};
\draw[redstring](0,0)node[below]{$i$}--++(0,1);
\end{tikzpicture}
,\quad
\begin{tikzpicture}[scale=1.2,anchorbase,smallnodes,rounded corners]
\draw[solid](0,1)node[above,yshift=-1pt]{$\phantom{i}$}--++(0.5,-0.5)--++(-0.5,-0.5) node[below]{$i$};
\draw[redstring](0.5,1)--++(-0.5,-0.5)--++(0.5,-0.5) node[below]{$i$};
\end{tikzpicture}
=
\begin{tikzpicture}[scale=1.2,anchorbase,smallnodes,rounded corners]
\draw[solid,dot](0,0)node[below]{$i$}--++(0,1)node[above,yshift=-1pt]{$\phantom{i}$};
\draw[redstring](0.5,0)node[below]{$i$}--++(0,1);
\end{tikzpicture}.
\end{gather}

\item \label{I:TripleCrossings}
The \emph{Reidemeister III relations} hold except in the following cases:
\begin{gather}\label{R:BraidGSG}
\begin{aligned}
\begin{tikzpicture}[scale=1.2,anchorbase,smallnodes,rounded corners]
\draw[ghost](1,1)node[above,yshift=-1pt]{$i$}--++(1,-1)node[below]{$\phantom{i}$};
\draw[ghost](2,1)node[above,yshift=-1pt]{$i$}--++(-1,-1)node[below]{$\phantom{i}$};
\draw[solid,smallnodes,rounded corners](1.5,1)--++(-0.5,-0.5)--++(0.5,-0.5)node[below]{$j$};
\end{tikzpicture}
&=
\begin{tikzpicture}[scale=1.2,anchorbase,smallnodes,rounded corners]
\draw[ghost](3,1)node[above,yshift=-1pt]{$\phantom{i}$}--++(1,-1)node[below]{$\phantom{i}$};
\draw[ghost](4,1)node[above,yshift=-1pt]{$i$}--++(-1,-1)node[below]{$\phantom{i}$};
\draw[solid,smallnodes,rounded corners](3.5,1)--++(0.5,-0.5)--++(-0.5,-0.5)node[below]{$j$};
\end{tikzpicture}
-Q_{iji}(\y)
\begin{tikzpicture}[scale=1.2,anchorbase,smallnodes,rounded corners]
\draw[ghost](6.2,1)node[above,yshift=-1pt]{$i$}--++(0,-1)node[below]{$\phantom{i}$};
\draw[ghost](7.2,1)node[above,yshift=-1pt]{$i$}--++(0,-1)node[below]{$\phantom{i}$};
\draw[solid](6.7,1)--++(0,-1)node[below]{$j$};
\end{tikzpicture}
\quad
\text{if $i\rightsquigarrow j$}
,\\
\begin{tikzpicture}[scale=1.2,anchorbase,smallnodes,rounded corners]
\draw[solid](1,1)node[above,yshift=-1pt]{$\phantom{i}$}--++(1,-1)node[below]{$j$};
\draw[solid](2,1)--++(-1,-1)node[below]{$j$};
\draw[ghost,smallnodes,rounded corners](1.5,1)node[above,yshift=-1pt]{$i$}--++(-0.5,-0.5)--++(0.5,-0.5)node[below]{$\phantom{i}$};
\end{tikzpicture}
&=
\begin{tikzpicture}[scale=1.2,anchorbase,smallnodes,rounded corners]
\draw[solid](3,1)node[above,yshift=-1pt]{$\phantom{i}$}--++(1,-1)node[below]{$j$};
\draw[solid](4,1)--++(-1,-1)node[below]{$j$};
\draw[ghost,smallnodes,rounded corners](3.5,1)node[above,yshift=-1pt]{$i$}--++(0.5,-0.5)--++(-0.5,-0.5)node[below]{$\phantom{i}$};
\end{tikzpicture}
+Q_{jij}(\y)
\begin{tikzpicture}[scale=1.2,anchorbase,smallnodes,rounded corners]
\draw[solid](7.2,1)node[above,yshift=-1pt]{$\phantom{i}$}--++(0,-1)node[below]{$j$};
\draw[solid](8.2,1)--++(0,-1)node[below]{$j$};
\draw[ghost](7.7,1)node[above,yshift=-1pt]{$i$}--++(0,-1)node[below]{$\phantom{i}$};
\end{tikzpicture}
\quad
\text{if $i\rightsquigarrow j$}
.
\end{aligned}
\end{gather}
\begin{gather}\label{R:BraidSRS}
\begin{tikzpicture}[scale=1.2,anchorbase,smallnodes,rounded corners]
\draw[solid](1,1)node[above,yshift=-1pt]{$\phantom{i}$}--++(1,-1)node[below]{$i$};
\draw[solid](2,1)--++(-1,-1)node[below]{$i$};
\draw[redstring](1.5,1)--++(-0.5,-0.5)--++(0.5,-0.5)node[below]{$i$};
\end{tikzpicture}
=
\begin{tikzpicture}[scale=1.2,anchorbase,smallnodes,rounded corners]
\draw[solid](3,1)node[above,yshift=-1pt]{$\phantom{i}$}--++(1,-1)node[below]{$i$};
\draw[solid](4,1)--++(-1,-1)node[below]{$i$};
\draw[redstring](3.5,1)--++(0.5,-0.5)--++(-0.5,-0.5)node[below]{$i$};
\end{tikzpicture}
-
\begin{tikzpicture}[scale=1.2,anchorbase,smallnodes,rounded corners]
\draw[solid](5,1)node[above,yshift=-1pt]{$\phantom{i}$}--++(0,-1)node[below]{$i$};
\draw[solid](6,1)--++(0,-1)node[below]{$i$};
\draw[redstring](5.5,1)--++(0,-1)node[below]{$i$};
\end{tikzpicture}
.
\end{gather}
\end{enumerate}
The \emph{honest Reidemeister relations} are those relations in (b) and (c) where the strings satisfy the Reidemeister relations, without error terms or dots.
\end{Definition}

There is an asymmetry in \autoref{R:GhostSolid} and \autoref{R:BraidGSG}
that plays an important role in \autoref{SS:KLRW}.

\begin{Example}
Multilocal relations can be tricky to apply. For example, it looks as if we can apply \autoref{R:DotCrossing} to the left-hand side of the following diagram:
\begin{gather*}
\begin{tikzpicture}[scale=1.2,anchorbase,smallnodes,rounded corners]
\draw[ghost](1,0.5)node[above,yshift=-1pt]{$i$}--(1,0) node[below]{$\phantom{i}$};
\draw[ghost](1.5,0.5)node[above,yshift=-1pt]{$i$}--(1.5,0) node[below]{$\phantom{i}$};
\draw[solid](0,0.5)node[above,yshift=-1pt]{$\phantom{i}$}--(0,0) node[below]{$i$};
\draw[solid](0.5,0.5)--(0.5,0) node[below]{$i$};
\draw[solid](1.25,0.5)--(1.25,0) node[below]{$j$};
\end{tikzpicture}
.
\end{gather*}
Relation \autoref{R:DotCrossing} changes both the solid strings and their ghost strings, since the ghost strings are implicitly part of this relation. This relation cannot be applied here because every local neighborhood that contains the two ghost $i$-strings also contains the solid $j$-string.
\end{Example}

\begin{Remark}\label{R:BetaVsN}
Note that we fixed $\beta\in Q^{+}_{n}$
of height $\mathrm{ht}(\beta)=n$, which
amounts to fixing the labels of the
solid strings. We sometimes write
\begin{gather*}
\WA[n](X)=
\bigoplus_{\beta\in Q^{+}_{n}}\WA(X),
\end{gather*}
which is a decomposition of $\WA(X)$ into a direct sum of two-sided
ideals.
Of course, anything we say about $\WA(X)$ has a corresponding statement for $\WA[n](X)$, which we will not usually state explicitly.
\end{Remark}


\subsection{The grading on \texorpdfstring{$\WA(X)$}{W}}\label{SS:WebsterDegrees}


\begin{Notation}
In this paper a \emph{graded} algebra or module
will always mean a $\Z$-graded algebra or module.
\end{Notation}

\begin{Definition}\label{D:Grading}
We endow the algebra $\WA(X)$ with a grading as follows.
The grading is defined
on the diagrams by summing over the contributions from each dot and crossing in the diagram according to the following local (not multilocal) rules:
\begin{gather*}
\deg\begin{tikzpicture}[scale=1.2,anchorbase,smallnodes,rounded corners]
\draw[solid,dot] (0,0)node[below]{$i$} to (0,0.5)node[above,yshift=-1pt]{$\phantom{i}$};
\end{tikzpicture}
=2d_{i}
,\quad
\deg\begin{tikzpicture}[scale=1.2,anchorbase,smallnodes,rounded corners]
\draw[ghost,dot] (0,0)node[below]{$\phantom{i}$} to (0,0.5)node[above,yshift=-1pt]{$i$};
\end{tikzpicture}=0,
\quad
\deg\begin{tikzpicture}[scale=1.2,anchorbase,smallnodes,rounded corners]
\draw[solid] (0,0)node[below]{$i$} to (0.5,0.5);
\draw[solid] (0.5,0)node[below]{$j$} to (0,0.5)node[above,yshift=-1pt]{$\phantom{i}$};
\end{tikzpicture}
=-\delta_{i,j}2d_{i}
,\quad
\deg\begin{tikzpicture}[scale=1.2,anchorbase,smallnodes,rounded corners]
\draw[ghost] (0,0)node[below]{$\phantom{i}$} to (0.5,0.5)node[above,yshift=-1pt]{$i$};
\draw[ghost] (0.5,0)node[below]{$\phantom{i}$} to (0,0.5)node[above,yshift=-1pt]{$j$};
\end{tikzpicture}
=0
,\\
\deg\begin{tikzpicture}[scale=1.2,anchorbase,smallnodes,rounded corners]
\draw[ghost] (0,0)node[below]{$\phantom{i}$} to (0.5,0.5)node[above,yshift=-1pt]{$i$};
\draw[solid] (0.5,0)node[below]{$j$} to (0,0.5)node[above,yshift=-1pt]{$\phantom{i}$};
\end{tikzpicture}
=
\deg\begin{tikzpicture}[scale=1.2,anchorbase,smallnodes,rounded corners]
\draw[ghost] (0.5,0)node[below]{$\phantom{i}$} to (0,0.5)node[above,yshift=-1pt]{$i$};
\draw[solid] (0,0)node[below]{$j$} to (0.5,0.5)node[above,yshift=-1pt]{$\phantom{i}$};
\end{tikzpicture}
=
\begin{cases}
\<\alpha_{i},\alpha_{j}\>
&\text{if $i\rightsquigarrow j$},
\\
0&\text{else},
\end{cases}
\\
\deg\begin{tikzpicture}[scale=1.2,anchorbase,smallnodes,rounded corners]
\draw[solid] (0,0)node[below]{$i$} to (0.5,0.5)node[above,yshift=-1pt]{$\phantom{i}$};
\draw[redstring] (0.5,0)node[below]{$j$} to (0,0.5);
\end{tikzpicture}
=
\deg\begin{tikzpicture}[scale=1.2,anchorbase,smallnodes,rounded corners]
\draw[solid] (0.5,0)node[below]{$i$} to (0,0.5)node[above,yshift=-1pt]{$\phantom{i}$};
\draw[redstring] (0,0)node[below]{$j$} to (0.5,0.5);
\end{tikzpicture}
=\tfrac{1}{2}\delta_{i,j}\<\alpha_{i},\alpha_{i}\>
,\quad
\deg\begin{tikzpicture}[scale=1.2,anchorbase,smallnodes,rounded corners]
\draw[ghost] (0,0)node[below]{$\phantom{i}$} to (0.5,0.5)node[above,yshift=-1pt]{$i$};
\draw[redstring] (0.5,0)node[below]{$j$} to (0,0.5)node[above,yshift=-1pt]{$\phantom{i}$};
\end{tikzpicture}
=
\deg\begin{tikzpicture}[scale=1.2,anchorbase,smallnodes,rounded corners]
\draw[ghost] (0.5,0)node[below]{$\phantom{i}$} to (0,0.5)node[above,yshift=-1pt]{$i$};
\draw[redstring] (0,0)node[below]{$j$} to (0.5,0.5)node[above,yshift=-1pt]{$\phantom{i}$};
\end{tikzpicture}
=0.
\end{gather*}
\end{Definition}

\begin{Lemma}
\autoref{D:Grading} endows $\WA(X)$ with the structure of a graded algebra.
\end{Lemma}

\begin{proof}
The algebra $\WA(X)$ is graded using
these degrees because all of the relations in
\autoref{D:RationalCherednik} are homogeneous with respect to the degree function from \autoref{D:Grading}.
\end{proof}

From now on, we consider $\WA(X)$ as a graded algebra using \autoref{D:Grading}.

\begin{Notation}
Let $\WA(X)\text{-}\mathbf{Mod}_{\Z}$ be the \emph{category of graded $\WA(X)$-modules}, where a \emph{graded category} means a category with hom-spaces enriched in the category of
graded $R$-modules. The corresponding categories of right $\WA(X)$-modules is
$\mathbf{Mod}_{\Z}\text{-}\WA(X)$.
\end{Notation}


\section{Some properties of weighted KLRW algebras}\label{S:Properties}


As in \autoref{N:KLRWdata}, we fix all of the data $\Gamma,\bsig,\beta,X$ {\etc} needed to define the weighted KLRW algebra $\WA(X)$.


\subsection{First properties}\label{SS:FirstProperties}


\begin{Proposition}\label{P:FixedOrientation}
Consider pairs $(\Gamma,\bsig)$ and $(\Gamma^{\prime},\bsig^{\prime})$ where the quivers $\Gamma$ and $\Gamma^{\prime}$ differ by precisely one
edge $\epsilon$, which is oriented $\epsilon:i\rightsquigarrow j$ in $\Gamma$ and $\epsilon^{\prime}:j\rightsquigarrow i$ in $\Gamma^{\prime}$, such that $Q_{ij}(u,v)=-Q_{ji}^{\prime}(v,u)$. Let $\WA(X)$ and $\WA(X)^{\prime}$, respectively, be the corresponding weighted KLRW algebras.
If $\sigma_{\epsilon}=-\sigma_{\epsilon^{\prime}}$, then there is an isomorphism of graded algebras $\WA(X)\cong\WA(X)^{\prime}$.
\end{Proposition}

\begin{proof}
The identity map on diagrams gives the required isomorphism.
\end{proof}

Define an $R$-linear map $({}_{-}){}^{\star}\map{\WA(X)}{\WA(X)}$
by reflecting diagrams in the line $y=\frac{1}{2}$:
\begin{gather}\label{E:StarMap}
\left(
\begin{tikzpicture}[scale=1.2,anchorbase,smallnodes,rounded corners]
\draw[ghost](1.5,1)node[above,yshift=-1pt]{$j$}--++(1,-1);
\draw[ghost](2.5,1)node[above,yshift=-1pt]{$i$}--++(-1,-1);
\draw[ghost](4,1)node[above,yshift=-1pt]{$k$}--++(0.5,-1);
\draw[solid](1,1)node[above,yshift=-1pt]{$\phantom{i}$}--++(1,-1)node[below]{$j$};
\draw[solid](2,1)--++(-1,-1)node[below]{$i$};
\draw[solid](3,1)--++(0.5,-1)node[below]{$k$};
\draw[redstring](1.25,1)--++(0,-1)node[below]{$\rho$};
\end{tikzpicture}
\right)^{\star}
=
\begin{tikzpicture}[scale=1.2,anchorbase,smallnodes,rounded corners]
\draw[ghost](1.5,1)node[above,yshift=-1pt]{$i$}--++(1,-1);
\draw[ghost](2.5,1)node[above,yshift=-1pt]{$j$}--++(-1,-1);
\draw[ghost](4,0)--++(0.5,1)node[above,yshift=-1pt]{$k$};
\draw[solid](1,1)node[above,yshift=-1pt]{$\phantom{i}$}--++(1,-1)node[below]{$i$};
\draw[solid](2,1)--++(-1,-1)node[below]{$j$};
\draw[solid](3,0)node[below]{$k$}--++(0.5,1);
\draw[redstring](1.25,1)--++(0,-1)node[below]{$\rho$};
\end{tikzpicture}
.
\end{gather}

The next two easy, but important, properties are straightforward.

\begin{Lemma}\label{L:Involutions}
\leavevmode
\begin{enumerate}

\item The map $({}_{-}){}^{\star}\map{\WA(X)}{\WA(X)}$ is a homogeneous algebra antiinvolution.

\item The map $({}_{-}){}^{\star}\map{\WA(X)}{\WA(X)}$ gives rise to equivalences of graded $R$-linear abelian categories
\begin{gather*}
({}_{-}){}^{\star}\map{\WA(X)\text{-}\mathbf{Mod}_{\Z}}{\mathbf{Mod}_{\Z}\text{-}\WA(X)},
\quad
({}_{-}){}^{\star}\map{\mathbf{Mod}_{\Z}\text{-}\WA(X)}{\WA(X)\text{-}\mathbf{Mod}_{\Z}},
\end{gather*}
both of which we denote by the same symbol.\qed

\end{enumerate}
\end{Lemma}

\begin{Lemma}\label{L:SLine}
Let $D\in\Webab[\bx,\bi][\by,\bj]$ be a straight line diagram. Then $D$ is invertible in $\WA(X)$, for appropriately chosen $X$, with inverse $D^{\star}$.\qed
\end{Lemma}

The set of positions $X$ is a set of loadings $(\charge,\bx)$, where $\charge$ is a charge and $\bx$ is a solid positioning. In particular, the loadings in $X$ do not necessarily have the same charge, which means that the endpoints of the red strings can vary.

A set of positions $Y$ has \emph{constant charge} $\charge$ if every loading in $Y$ has charge $\charge$.
The next result says that we can assume that the set of positions has constant charge.

\begin{Proposition}\label{P:FixedCharge}
Suppose that $X$ is a set of positions. Then there exists a set of positions $Y$ of constant charge such that $\WA(X)\cong\WA(Y)$ as graded algebras.
\end{Proposition}

\begin{proof}
Conjugation by straight line diagrams defines the required isomorphism by \autoref{L:SLine}:
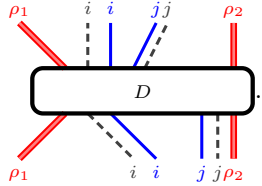
\begin{gather*}
\begin{tikzpicture}[scale=1.2,anchorbase,smallnodes,rounded corners]
\draw[ghost](-0.5,-0.25)node[above,yshift=-1pt]{$\phantom{i}$}--++(0.5,-0.5)node[below]{$i$};
\draw[ghost](-0.5,0.25)node[above,yshift=-1pt]{$\phantom{i}$}--++(0,0.5)node[above,yshift=-1pt]{$i$};
\draw[ghost](0.925,-0.25)node[above,yshift=-1pt]{$\phantom{i}$}--++(0,-0.5)node[below]{$j$};
\draw[ghost](0.125,0.25)node[above,yshift=-1pt]{$\phantom{i}$}--++(0.25,0.5)node[above,yshift=-1pt]{$j$};
\draw[solid](-0.25,-0.25)node[above,yshift=-1pt]{$\phantom{i}$}--++(0.5,-0.5)node[below]{$i$};
\draw[solid](-0.25,0.25)node[above,yshift=-1pt]{$\phantom{i}$}--++(0,0.5)node[above,yshift=-1pt]{$i$};
\draw[solid](0.75,-0.25)node[above,yshift=-1pt]{$\phantom{i}$}--++(0,-0.5)node[below]{$j$};
\draw[solid](0,0.25)node[above,yshift=-1pt]{$\phantom{i}$}--++(0.25,0.5)node[above,yshift=-1pt]{$j$};
\draw[redstring](-0.75,-0.25)node[above,yshift=-1pt]{$\phantom{i}$}--++(-0.5,-0.5)node[below]{$\rho_{1}$};
\draw[redstring](1.1,-0.25)node[above,yshift=-1pt]{$\phantom{i}$}--++(0,-0.5)node[below]{$\rho_{2}$};
\draw[redstring](-0.75,0.25)node[above,yshift=-1pt]{$\phantom{i}$}--++(-0.5,0.5)node[above,yshift=-1pt]{$\rho_{1}$};
\draw[redstring](1.1,0.25)node[above,yshift=-1pt]{$\phantom{i}$}--++(0,0.5)node[above,yshift=-1pt]{$\rho_{2}$};
\node[rectangle,draw,minimum width=2.9cm,minimum height=0.6cm,ultra thick] at(0.1,0){\raisebox{-0.05cm}{$D$}};
\end{tikzpicture}
.
\end{gather*}
Here $D$ is a diagram in $\WA(X)$.
\end{proof}

\begin{Notation}
Without loss of generality, by \autoref{P:FixedCharge}, we can, and will, assume that $X$ has constant charge $\charge$.
\end{Notation}

Let us further explain how one can vary the size of $X$.

\begin{Definition}
Let $\simeq$ be the equivalence relation on the set of loadings where
$\bx\simeq\by$, if for each $\bi\in I^{n}$ there exists a straight line diagram $D\in\Webaa[\bx,\bi][\by,\bi]$.
Extend $\simeq$ to the set of positionings by defining
$X\simeq Y$, if for any $\bx\in X$ there exists $\by\in Y$ such that $\bx\simeq\by$.
\end{Definition}

\begin{Example}
Let $\bx$ and $\by$ be the bottom and top positionings of a straight line diagram, as in \autoref{E:StraightLine}. Then
$\{\bx\}\simeq\{\bx,\by\}\simeq\{\by\}$, since $\bx\simeq\by$. Then \autoref{P:VaryingSize} below shows that from the perspective of
weighted KLRW algebras one of these loadings is redundant.
\end{Example}

\begin{Proposition}\label{P:VaryingSize}
Suppose that $X\simeq Y$.
\begin{enumerate}

\item If $\#X=\#Y$, then $\WA(X)\cong\WA(Y)$ as graded algebras.

\item The algebras $\WA(X)$ and $\WA(Y)$ are graded Morita equivalent.

\end{enumerate}
\end{Proposition}

\begin{proof}
(a). Note that we can restrict to the case where $X$ and $Y$ differ by only one element.
Then the straight line diagrams provide the corresponding isomorphism.

(b). We can restrict to the case where
$\#X+1=\#Y$ and, by (a), that $X\subset Y$. Note that, by assumption, the additional element in $Y$ is related to another element in $X$ by a straight line diagram, showing that the corresponding graded projective $\WA(Y)$-modules are graded equivalent. This together with (a) proves the claim.
\end{proof}

Observe that the identity element of $\WA(X)$ is
\begin{gather*}
\onealg=\sum_{\bx\in X,\bi\in I^{\beta}}\1_{\bx,\bi},
\end{gather*}
and that $\WA(X)\1_{\bx,\bi}$ is a projective $\WA(X)$-module, for $\bx\in X$ and $\bi\in I^{\beta}$.
Note that every projective $\WA(X)$-module appears as a direct summand of some projective $\WA(X)$-module of the form $\WA(X)\1_{\bx,\bi}$. Moreover, for any $\bx\in X$ let $\1_{\bx}=\sum_{\bi\in I^{\beta}}\1_{\bx,\bi}$.
Then $\1_{\bx}$ is an idempotent in $\WA(X)$ and
$\1_{\bx}\WA(X)\1_{\bx}$ and $\1_{\bx,\bi}\WA(X)\1_{\bx,\bi}$ are idempotent subalgebras, see \autoref{SS:KLRW}.


\subsection{A basis of \texorpdfstring{$\WA(X)$}{W}}


If $w\in\Sym$, then a reduced expression for $w$ is any word $w=s_{i_{1}}\dots s_{i_{k}}$ of minimal length, for $1\leq i_{j}<n$. If $w=s_{i_{1}}\dots s_{i_{k}}$ is a reduced expression for $w$, then $w$ has length $k$.
The following \emph{permutation diagrams} play an important role throughout this paper.

\begin{Definition}\label{D:Dw}
Fix $\bx,\by\in X$ and $\bi,\bj\in I^{\beta}$.
For each $w\in\Sym$ such that $\bj=w\cdot\bi$
define a diagram $D_{\bx,\bi}^{\by,\bj}(w)\in\Webab[\bx,\bi][\by,\bj]$
as follows. First, assume as an intermediate step that $\charge=\emptyset$, so that there are no red strings. Then
connect the strings in $D_{\balp,\bi}^{\bbet,\bj}(w)$ from bottom to top by adding solid crossings (and their ghosts)
corresponding to the fixed choice of
reduced expression. Now add the red strings vertically to the diagram and, if necessary, pull each red string to the right so that it is to the right of any solid-solid, ghost-ghost and solid-ghost
crossings on the strings that the red string crosses.
\end{Definition}

The diagram $D_{\bx,\bi}^{\by,\bj}(w)$ depends on the choice of reduced expression, and not just on $w$, and on $\charge$, but we suppress this from our notation.

\begin{Example}
In the situation of the configuration as in
\autoref{R:BraidSRS} we prefer the
leftmost diagram, not the middle. Precisely, let $n=2$, $l=1$, $\charge=(\tfrac{3}{2})$, $\charge^{\prime}=(\tfrac{1}{2})$,
$\charge^{\prime\prime}=(-\tfrac{1}{2})$ and $\balp=\bbet=(-1,1)$. Then, for $\bi=(i,j)$, $\bj=(j,i)$ with $i,j\in I$,
and $s$ being the unique simple reflection in $\Sym$
(with $s$ being the reduced expression) we have
\begin{gather*}
D_{\bx,\bi}^{\by,\bj}(s)
=
\begin{tikzpicture}[scale=1.2,anchorbase,smallnodes,rounded corners]
\draw[ghost](1.5,1)node[above,yshift=-1pt]{$j$}--++(1,-1);
\draw[ghost](2.5,1)node[above,yshift=-1pt]{$i$}--++(-1,-1);
\draw[solid](1,1)node[above,yshift=-1pt]{$\phantom{i}$}--++(1,-1)node[below]{$j$};
\draw[solid](2,1)--++(-1,-1)node[below]{$i$};
\draw[redstring](2.25,1)--++(0,-1)node[below]{$\rho$};
\end{tikzpicture}
\text{for }\charge
,\quad
D_{\bx,\bi}^{\by,\bj}(s)
=
\begin{tikzpicture}[scale=1.2,anchorbase,smallnodes,rounded corners]
\draw[ghost](1.5,1)node[above,yshift=-1pt]{$j$}--++(1,-1);
\draw[ghost](2.5,1)node[above,yshift=-1pt]{$i$}--++(-1,-1);
\draw[solid](1,1)node[above,yshift=-1pt]{$\phantom{i}$}--++(1,-1)node[below]{$j$};
\draw[solid](2,1)--++(-1,-1)node[below]{$i$};
\draw[redstring](1.75,1)--++(0.5,-0.35)--++(0,-0.3)--++(-0.5,-0.35)node[below]{$\rho$};
\end{tikzpicture}
\text{ for }\charge^{\prime}
,\quad
D_{\bx,\bi}^{\by,\bj}(s)
=
\begin{tikzpicture}[scale=1.2,anchorbase,smallnodes,rounded corners]
\draw[ghost](1.5,1)node[above,yshift=-1pt]{$j$}--++(1,-1);
\draw[ghost](2.5,1)node[above,yshift=-1pt]{$i$}--++(-1,-1);
\draw[solid](1,1)--++(1,-1)node[below]{$j$};
\draw[solid](2,1)--++(-1,-1)node[below]{$i$};
\draw[redstring](1.25,1)--++(0.4,-0.35)--++(0,-0.3)--++(-0.4,-0.35)node[below]{$\rho$};
\end{tikzpicture}
\text{for }\charge^{\prime\prime}
.
\end{gather*}
By construction, all crossings on the non-red strings that cross a red string are to the left of that red string.
\end{Example}

As in \autoref{D:Bjstrings}, whenever we write
$D_{\bx,\bi}^{\by,\bj}(w)$ we assume that $\bj=w\cdot\bi$. When
$\bx$, $\by$, $\bi$ and $\bj$ are clear from context, write
$D(w)=D_{\bx,\bi}^{\by,\bj}(w)$. By definition
$D(w)=D(w)\1_{\bx,\bi}=\1_{\by,\bj}D(w)$.

\begin{Lemma}\label{L:FiniteGeneration}
The algebra $\WA(X)$ is a finitely generated $R$-module and spanned by
\begin{gather}\label{E:AffineBasis}
\WABasis=
\set{\1_{\by,\bj}D(w)y_{1}^{a_{1}}\dots y_{n}^{a_{n}}\1_{\bx,\bi}|\bx,\by\in X,\bi,\bj\in I^{\beta},w\in\Sym,a_{1},\dots,a_{n}\in\N}.
\end{gather}
\end{Lemma}

\autoref{P:WABasis} below shows that \autoref{E:AffineBasis} is
a basis of $\WA(X)$.

\begin{proof}
The relations imply that the algebra $\WA(X)$ is filtered with
$\WA(X)=\bigcup_{k\in\Z_{\geq 0}}\WA(X)_{k}$, where
$\WA(X)_{k}$ is spanned by the diagrams with at most
$k$ crossings, for $k\in\Z_{\geq 0}$. In the associated
graded algebra,
\begin{gather*}
\gr\WA(X)
=\bigoplus_{k>0}\WA(X)_{k}/\WA(X)_{k-1},
\end{gather*}
dots and crossing commute, double crossings satisfy
the Reidemeister II relations or are zero, and the crossings
satisfy the Reidemeister III relations. Therefore, the image of $D(w)$ in
$\gr\WA(X)$ is nonzero and depends only on $w$. Moreover, if $D(w)$
and $D(w)^{\prime}$ are two such diagrams for $w\in\Sym$ that are defined using
different reduced expressions, then $D(w)\equiv D(w)^{\prime}\pmod{\WA(X)_{k}}$,
where $w$ has length $k$, which is
defined to be the minimal number of crossings in a permutation diagram for $w$.
Thus, $\gr\WA(X)$ is spanned by the image
of the set $\WABasis$ in $\gr\WA(X)$, implying the result.
\end{proof}

\begin{Remark}
Let $w\in\Sym$. The proof of \autoref{L:FiniteGeneration} shows that, if $D(w)$ and $D(w)^{\prime}$ are diagrams defined using different reduced expressions for $w$, then $D(w)\equiv D(w)^{\prime}\pmod{\WA(X)_{k-1}}$, where $k=\ell(w)$. By the same argument, if $D(w)^{\prime}$ is a diagram defined in the same way as $D(w)$ except that the red strings are pulled to the left in \autoref{D:Dw}, then $D(w)$ and $D(w)^{\prime}$ agree modulo $\WA(X)_{k-1}$.
\end{Remark}

\begin{Example}\label{Ex:Basis}
Consider the quiver $\Gamma$ as in \autoref{E:ExampleQuiver}, and let $n=2$ and $\ell=1$. Let the solid strings have positions $\bx=(0,0.75)$ and place the red string at $1.5$. Recall from \autoref{Ex:Beta} that there are six choices of $\beta$. Let $\beta=\alpha_{0}+\alpha_{1}$. Then
\begin{gather*}
\WABasis=\set[\Bigg]{
\begin{tikzpicture}[scale=1.2,anchorbase,smallnodes,rounded corners]
\draw[ghost,dot=0.15](1,0)--++(0,1)node[above,yshift=-1pt]{$0$};
\draw[ghost,dot=0.15](1.75,0)--++(0,1)node[above,yshift=-1pt]{$1$};
\draw[solid,dot=0.15](0,0)node[below]{$0$}node[left,xshift=0.06cm,yshift=0.23cm]{$a$}--++(0,1)node[above,yshift=-1pt]{$\phantom{i}$};
\draw[solid,dot=0.15](0.75,0)node[below]{$1$}node[right,xshift=-0.06cm,yshift=0.23cm]{$b$}--++(0,1)node[above,yshift=-1pt]{$\phantom{i}$};
\draw[redstring](1.5,0)node[below]{$\rho$}--++(0,1)node[above,yshift=-1pt]{$\phantom{i}$};
\end{tikzpicture}
,\quad
\begin{tikzpicture}[scale=1.2,anchorbase,smallnodes,rounded corners]
\draw[ghost,dot=0.15](1,0)--++(0.75,1)node[above,yshift=-1pt]{$0$};
\draw[ghost,dot=0.15](1.75,0)--++(-0.75,1)node[above,yshift=-1pt]{$1$};
\draw[solid,dot=0.15](0,0)node[below]{$0$}node[left,xshift=0.06cm,yshift=0.23cm]{$a$}--++(0.75,1)node[above,yshift=-1pt]{$\phantom{i}$};
\draw[solid,dot=0.15](0.75,0)node[below]{$1$}node[right,xshift=-0.06cm,yshift=0.23cm]{$b$}--++(-0.75,1)node[above,yshift=-1pt]{$\phantom{i}$};
\draw[redstring](1.5,0)node[below]{$\rho$}--++(0,1)node[above,yshift=-1pt]{$\phantom{i}$};
\end{tikzpicture}
,\quad
\begin{tikzpicture}[scale=1.2,anchorbase,smallnodes,rounded corners]
\draw[ghost,dot=0.15](1,0)--++(0,1)node[above,yshift=-1pt]{$1$};
\draw[ghost,dot=0.15](1.75,0)--++(0,1)node[above,yshift=-1pt]{$0$};
\draw[solid,dot=0.15](0,0)node[below]{$1$}node[left,xshift=0.06cm,yshift=0.23cm]{$a$}--++(0,1)node[above,yshift=-1pt]{$\phantom{i}$};
\draw[solid,dot=0.15](0.75,0)node[below]{$0$}node[right,xshift=-0.06cm,yshift=0.23cm]{$b$}--++(0,1)node[above,yshift=-1pt]{$\phantom{i}$};
\draw[redstring](1.5,0)node[below]{$\rho$}--++(0,1)node[above,yshift=-1pt]{$\phantom{i}$};
\end{tikzpicture}
,\quad
\begin{tikzpicture}[scale=1.2,anchorbase,smallnodes,rounded corners]
\draw[ghost,dot=0.15](1,0)--++(0.75,1)node[above,yshift=-1pt]{$1$};
\draw[ghost,dot=0.15](1.75,0)--++(-0.75,1)node[above,yshift=-1pt]{$0$};
\draw[solid,dot=0.15](0,0)node[below]{$1$}node[left,xshift=0.06cm,yshift=0.23cm]{$a$}--++(0.75,1)node[above,yshift=-1pt]{$\phantom{i}$};
\draw[solid,dot=0.15](0.75,0)node[below]{$0$}node[right,xshift=-0.06cm,yshift=0.23cm]{$b$}--++(-0.75,1)node[above,yshift=-1pt]{$\phantom{i}$};
\draw[redstring](1.5,0)node[below]{$\rho$}--++(0,1)node[above,yshift=-1pt]{$\phantom{i}$};
\end{tikzpicture}
|a,b\in\N},
\end{gather*}
where $a$ and $b$ specify the number of dots. (The ghost strings have the same number of dots.)
\end{Example}

To show that $\WABasis$ is a basis of $\WA(X)$ we introduce a faithful polynomial
module of $\WA(X)$, following the standard approach in this setting, for example see
\cite[\S3.2.2]{Ro-2-kac-moody}, \cite[\S2.3]{KhLa-cat-quantum-sln-first} or
\cite[\S2.2]{We-weighted-klr}. Fix indeterminates
$y_{1},\dots,y_{n}$ over $R$ and let
\begin{gather*}
P_{\beta}(X)=
\bigoplus_{\bx\in X,\bi\in I^{\beta}}
R[y_{1},\dots,y_{n}]\1_{\bx,\bi}
\end{gather*}
be $\#(X\times I^{\beta})$ copies of the polynomial ring $R[y_{1},\dots,y_{n}]$, where $y_{r}\1_{\bx,\bi}$ has degree $2d_{i_{r}}$, for $1\leq r\leq n$; {\cf} \autoref{D:Grading}. The
symmetric group $\Sym$ acts on $P_{\beta}(X)$ by place permutations. In
particular, if $1\leq r<s<n$ let $(r,s)$ be the operator that interchanges
$y_{r}\1_{\bx,\bi}$ and $y_{s}\1_{\by,\bi}$ and fixes
$y_{t}\1_{\bx,\bi}$ for $t\neq r,s$.

\begin{Definition}\label{D:PolynomialAction}
Define an assignment, using local (not multilocal) rules,
from $\WA(X)$ to $P_{\beta}(X)$ by specifying that
\begin{gather}\label{E:Action1}
\1_{\by,\bj}\cdot f(y)\1_{\bx,\bi}=\delta_{\bx\by}\delta_{\bi\bj}f(y)\1_{\bx,\bi}
,\quad
\begin{tikzpicture}[scale=1.2,anchorbase,smallnodes,rounded corners]
\draw[solid,dot] (0,0)node[below]{$i_{r}$} to (0,0.5)node[above,yshift=-1pt]{$\phantom{i}$};
\end{tikzpicture}
\mapsto y_{r}
,\quad
\begin{tikzpicture}[scale=1.2,anchorbase,smallnodes,rounded corners]
\draw[ghost,dot] (0,0)node[below]{$\phantom{i}$} to (0,0.5)node[above,yshift=-1pt]{$i_{r}$};
\end{tikzpicture}
\mapsto 1,
\end{gather}
and all crossings act as the identity except that
\begin{gather}\label{E:Action2}
\begin{tikzpicture}[scale=1.2,anchorbase,smallnodes,rounded corners]
\draw[solid] (0,0)node[below]{$i_{r}$} to (0.5,0.5);
\draw[solid] (0.5,0)node[below]{$i_{s}$} to (0,0.5)node[above,yshift=-1pt]{$\phantom{i}$};
\end{tikzpicture}
\mapsto
\begin{cases*}
\partial_{r,s} & if $i_{r}=i_{s}$,
\\
(r,s) & if $i_{r}\neq i_{s}$,
\end{cases*}
\!
\begin{tikzpicture}[scale=1.2,anchorbase,smallnodes,rounded corners]
\draw[solid] (0,0)node[below]{$i_{r}$} to (0.5,0.5);
\draw[redstring] (0.5,0)node[below]{$i_{s}$} to (0,0.5)node[above,yshift=-1pt]{$\phantom{i}$};
\end{tikzpicture}
\mapsto
\begin{cases*}
x_r & if $i_{r}=i_{s}$,
\\
1 & if $i_{r}\neq i_{s}$,
\end{cases*}
\!
\begin{tikzpicture}[scale=1.2,anchorbase,smallnodes,rounded corners]
\draw[ghost] (0,0) to (0.5,0.5)node[above,yshift=-1pt]{$i_{r}$};
\draw[solid] (0.5,0)node[below]{$i_{s}$} to (0,0.5)node[above,yshift=-1pt]{$\phantom{i}$};
\end{tikzpicture}
\mapsto
\begin{cases}
Q_{i_{r},i_{s}}(y_{r},y_{s})
&\text{if $i\rightsquigarrow j$},
\\
1&\text{else},
\end{cases}
\end{gather}
where $\partial_{r,s}=\frac{(r,s)-1}{y_{s}-y_{r}}$
is the \emph{Demazure operator}.
\end{Definition}

We call $P_{\beta}(X)$ the \emph{polynomial module} of $\WA(X)$.
Note that the action of the crossings maps
$R[y_{1},\dots,y_{n}]\1_{\bx,\bi}$ to $R[y_{1},\dots,y_{n}]\1_{\by,\bj}$ for
related $\bi$ and $\bj$.
Similarly we define:

\begin{Definition}\label{D:PolynomialActionGr}
Define an assignment, using local rules,
from $\gr\WA(X)$ to $P_{\beta}(X)$ by
using the assignment in \autoref{E:Action1}, but instead of
\autoref{E:Action2} we let
\begin{gather*}
\begin{tikzpicture}[scale=1.2,anchorbase,smallnodes,rounded corners]
\draw[solid] (0,0)node[below]{$i_{r}$} to (0.5,0.5);
\draw[solid] (0.5,0)node[below]{$i_{s}$} to (0,0.5)node[above,yshift=-1pt]{$\phantom{i}$};
\end{tikzpicture}
\mapsto
\begin{cases*}
0 & if $i_{r}=i_{s}$,
\\
(r,s) & if $i_{r}\neq i_{s}$,
\end{cases*}
\quad
\begin{tikzpicture}[scale=1.2,anchorbase,smallnodes,rounded corners]
\draw[solid] (0,0)node[below]{$i_{r}$} to (0.5,0.5);
\draw[redstring] (0.5,0)node[below]{$i_{s}$} to (0,0.5)node[above,yshift=-1pt]{$\phantom{i}$};
\end{tikzpicture}
\mapsto
\begin{cases*}
x_r & if $i_{r}=i_{s}$,
\\
1 & if $i_{r}\neq i_{s}$,
\end{cases*}
\quad
\begin{tikzpicture}[scale=1.2,anchorbase,smallnodes,rounded corners]
\draw[ghost] (0,0) to (0.5,0.5)node[above,yshift=-1pt]{$i_{r}$};
\draw[solid] (0.5,0)node[below]{$i_{s}$} to (0,0.5);
\end{tikzpicture}
\mapsto
\begin{cases}
0
&\text{if $i\rightsquigarrow j$},
\\
1&\text{else},
\end{cases}
\end{gather*}
and all other crossings act as the identity.
\end{Definition}

\begin{Lemma}
There is an action of $\WA(X)$ on $P_{\beta}(X)$ given by
\begin{gather*}
Df(y)
=
\begin{tikzpicture}[scale=1.2,anchorbase,smallnodes,rounded corners]
\node[rectangle,draw,minimum width=1cm,minimum height=0.56cm,ultra thick] at(0,0){\raisebox{-0.05cm}{$f(y)$}};
\node[rectangle,draw,minimum width=1cm,minimum height=0.56cm,ultra thick] at(0,0.5){\raisebox{-0.05cm}{$D$}};
\end{tikzpicture}
,
\end{gather*}
where we apply the local
rules from \autoref{D:PolynomialAction} to
$P_{\beta}(X)$ by reading the crossings in $D$ in order from bottom to top. Similarly, $\gr\WA(X)$ acts on $P_{\beta}(X)$ via \autoref{D:PolynomialActionGr}.
\end{Lemma}

\begin{proof}
Both claims follow by what is
now a standard check of the relations, such as
in \cite[\S2.3]{KhLa-cat-quantum-sln-first}, with the slight caveat that solid strings might have associated ghost strings. For example, whenever a ghost crossing appears in one of the relations then so does its solid crossing. Explicitly, for an appropriate choice of ghost shifts and $i,j\in I^{\beta}$, \autoref{R:BraidGSG} becomes
\begin{gather*}
\begin{tikzpicture}[scale=1.2,anchorbase,smallnodes,rounded corners]
\draw[ghost](1,1)node[above,yshift=-1pt]{$i$}--++(1,-1);
\draw[ghost](2,1)node[above,yshift=-1pt]{$i$}--++(-1,-1);
\draw[ghost,smallnodes,rounded corners](2.5,1)node[above,yshift=-1pt]{$j$}--++(-0.5,-0.5)--++(0.5,-0.5);
\draw[solid](-0.5,1)node[above,yshift=-1pt]{$\phantom{i}$}--++(1,-1)node[below]{$i$};
\draw[solid](0.5,1)--++(-1,-1)node[below]{$i$};
\draw[solid,smallnodes,rounded corners](1.5,1)--++(-0.5,-0.5)--++(0.5,-0.5)node[below]{$j$};
\end{tikzpicture}
=
\begin{tikzpicture}[scale=1.2,anchorbase,smallnodes,rounded corners]
\draw[ghost](3,1)node[above,yshift=-1pt]{$i$}--++(1,-1);
\draw[ghost](4,1)node[above,yshift=-1pt]{$i$}--++(-1,-1);
\draw[ghost,smallnodes,rounded corners](4.5,1)node[above,yshift=-1pt]{$j$}--++(0.5,-0.5)--++(-0.5,-0.5);
\draw[solid](1.5,1)node[above,yshift=-1pt]{$\phantom{i}$}--++(1,-1)node[below]{$i$};
\draw[solid](2.5,1)--++(-1,-1)node[below]{$i$};
\draw[solid,smallnodes,rounded corners](3.5,1)--++(0.5,-0.5)--++(-0.5,-0.5)node[below]{$j$};
\end{tikzpicture}
-Q_{iji}(\y)
\begin{tikzpicture}[scale=1.2,anchorbase,smallnodes,rounded corners]
\draw[ghost](6.2,1)node[above,yshift=-1pt]{$i$}--++(0,-1);
\draw[ghost](7.2,1)node[above,yshift=-1pt]{$i$}--++(0,-1);
\draw[ghost](7.7,1)node[above,yshift=-1pt]{$j$}--++(0,-1);
\draw[solid](4.7,1)--++(0,-1)node[below]{$i$};
\draw[solid](5.7,1)--++(0,-1)node[below]{$i$};
\draw[solid](6.7,1)--++(0,-1)node[below]{$j$};
\end{tikzpicture}
.
\end{gather*}
Hence, the action of the two sides of this relation coincide because
\begin{gather*}
\partial_{r,s}\big(Q_{ij}(r,t)f\big)
=
Q_{ij}(s,t)\partial_{r,s}(f)
-Q_{iji}(r,t,s),
\end{gather*}
where we have named the solid strings $r$, $s$, and $t$ in order. Note that this requires the Leibniz rule for $\partial_{r,s}$. That is, $\partial_{r,s}\big(Q_{ij}(r,t)f\big)=
Q_{ij}(r,t)\partial_{r,s}(f)+\big((r,s)\centerdot f\big)\partial_{r,s}\big(Q_{ij}(r,t)\big)$.
\end{proof}

Recall the set $\WABasis$ from \autoref{E:AffineBasis}.

\begin{Proposition}\label{P:WABasis}
The algebra $\WA(X)$ is free as an $R$-module with homogeneous basis $\WABasis$.
\end{Proposition}

\begin{proof}
By \autoref{L:FiniteGeneration} the elements in $\WABasis(X)$ span $\WA(X)$. Hence, it remains
to show that $\WABasis(X)$ is linearly
independent and for this it is
sufficient to show that the images of
$\WABasis(X)$ in $\gr\WA(X)$ are
linearly independent. Using the action of $\gr\WA(X)$ on
$P_{\beta}(X)$, {\cf} \autoref{D:PolynomialActionGr}, we see that $D(w)$ acts by
sending $f(y_{1},\dots,y_{n})$ to $f(y_{w(1)},\dots,y_{w(n)})$. Hence, the elements of $\WABasis(X)$ map to
linearly independent automorphisms of
$R[y_{1},\dots,y_{n}]$ because $R[y_{1},\dots,y_{n}]$ is free as a module over itself.
\end{proof}

As an immediate corollary of \autoref{P:WABasis} we obtain:

\begin{Corollary}\label{C:FaithfulPoly}
The polynomial module $P_{\beta}(X)$ is a faithful
$\WA(X)$-module.
\end{Corollary}

\begin{Remark}
Note that the basis of \autoref{P:WABasis} is a homogeneous basis. On the other hand, $P_{\beta}(X)$ is not a graded $\WA(X)$-module. However, $P_{\beta}(X)$ is a graded $\gr\WA(X)$-module.
\end{Remark}


\subsection{The center of \texorpdfstring{$\WA(X)$}{W}}


Let $R[y_{1},\dots,y_{a}]^{\Sym[a]}$ be the ring of symmetric
polynomials in $y_{1},\dots,y_{a}$, where $a\in\Z_{\geq 1}$. For $\beta=\sum_{i\in I}a_{i}\alpha_{i}$, let
$R[\beta]^{\Sym[\beta]}=\bigotimes_{i\in I}R[y_{1},\dots,y_{a_{i}}]^{\Sym[a_{i}]}$.
By \autoref{P:WABasis}, we can view $R[\beta]^{\Sym[\beta]}$ as a graded subalgebra of $\WA(X)$ via the map $f\mapsto f\onealg$.

\begin{Example}
Recall from \autoref{Ex:Beta} that $\beta_{ij}=\alpha_{i}+\alpha_{j}$ and $\beta_{i}=2\alpha_{i}$, for $i,j\in I$. Then the generators for $R[\beta_{ij}]^{\Sym[\beta_{ij}]}$ are
\begin{gather*}
\onealg=
\begin{tikzpicture}[scale=1.2,anchorbase,smallnodes,rounded corners]
\draw[ghost](1,0)--++(0,1)node[above,yshift=-1pt]{$i$};
\draw[ghost](1.5,0)--++(0,1)node[above,yshift=-1pt]{$j$};
\draw[solid](0,0)node[below]{$i$}--++(0,1)node[above,yshift=-1pt]{$\phantom{i}$};
\draw[solid](0.5,0)node[below]{$j$}--++(0,1)node[above,yshift=-1pt]{$\phantom{i}$};
\draw[redstring](0.25,0)node[below]{$\rho$}--++(0,1)node[above,yshift=-1pt]{$\phantom{i}$};
\end{tikzpicture}
\!+\!
\begin{tikzpicture}[scale=1.2,anchorbase,smallnodes,rounded corners]
\draw[ghost](1,0)--++(0,1)node[above,yshift=-1pt]{$j$};
\draw[ghost](1.5,0)--++(0,1)node[above,yshift=-1pt]{$i$};
\draw[solid](0,0)node[below]{$j$}--++(0,1)node[above,yshift=-1pt]{$\phantom{i}$};
\draw[solid](0.5,0)node[below]{$i$}--++(0,1)node[above,yshift=-1pt]{$\phantom{i}$};
\draw[redstring](0.25,0)node[below]{$\rho$}--++(0,1)node[above,yshift=-1pt]{$\phantom{i}$};
\end{tikzpicture}
\!,\!
\begin{tikzpicture}[scale=1.2,anchorbase,smallnodes,rounded corners]
\draw[ghost,dot](1,0)--++(0,1)node[above,yshift=-1pt]{$i$};
\draw[ghost](1.5,0)--++(0,1)node[above,yshift=-1pt]{$j$};
\draw[solid,dot](0,0)node[below]{$i$}--++(0,1)node[above,yshift=-1pt]{$\phantom{i}$};
\draw[solid](0.5,0)node[below]{$j$}--++(0,1)node[above,yshift=-1pt]{$\phantom{i}$};
\draw[redstring](0.25,0)node[below]{$\rho$}--++(0,1)node[above,yshift=-1pt]{$\phantom{i}$};
\end{tikzpicture}
\!+\!
\begin{tikzpicture}[scale=1.2,anchorbase,smallnodes,rounded corners]
\draw[ghost](1,0)--++(0,1)node[above,yshift=-1pt]{$j$};
\draw[ghost,dot](1.5,0)--++(0,1)node[above,yshift=-1pt]{$i$};
\draw[solid](0,0)node[below]{$j$}--++(0,1)node[above,yshift=-1pt]{$\phantom{i}$};
\draw[solid,dot](0.5,0)node[below]{$i$}--++(0,1)node[above,yshift=-1pt]{$\phantom{i}$};
\draw[redstring](0.25,0)node[below]{$\rho$}--++(0,1)node[above,yshift=-1pt]{$\phantom{i}$};
\end{tikzpicture}
\!,\!
\begin{tikzpicture}[scale=1.2,anchorbase,smallnodes,rounded corners]
\draw[ghost](1,0)--++(0,1)node[above,yshift=-1pt]{$i$};
\draw[ghost,dot](1.5,0)--++(0,1)node[above,yshift=-1pt]{$j$};
\draw[solid](0,0)node[below]{$i$}--++(0,1)node[above,yshift=-1pt]{$\phantom{i}$};
\draw[solid,dot](0.5,0)node[below]{$j$}--++(0,1)node[above,yshift=-1pt]{$\phantom{i}$};
\draw[redstring](0.25,0)node[below]{$\rho$}--++(0,1)node[above,yshift=-1pt]{$\phantom{i}$};
\end{tikzpicture}
\!+\!
\begin{tikzpicture}[scale=1.2,anchorbase,smallnodes,rounded corners]
\draw[ghost,dot](1,0)--++(0,1)node[above,yshift=-1pt]{$j$};
\draw[ghost](1.5,0)--++(0,1)node[above,yshift=-1pt]{$i$};
\draw[solid,dot](0,0)node[below]{$j$}--++(0,1)node[above,yshift=-1pt]{$\phantom{i}$};
\draw[solid](0.5,0)node[below]{$i$}--++(0,1)node[above,yshift=-1pt]{$\phantom{i}$};
\draw[redstring](0.25,0)node[below]{$\rho$}--++(0,1)node[above,yshift=-1pt]{$\phantom{i}$};
\end{tikzpicture}
\!,
\end{gather*}
and the generators for $R[\beta_{i}]^{\Sym[\beta_{i}]}$ are
\begin{gather*}
\onealg=
\begin{tikzpicture}[scale=1.2,anchorbase,smallnodes,rounded corners]
\draw[ghost](1,0)--++(0,1)node[above,yshift=-1pt]{$i$};
\draw[ghost](1.5,0)--++(0,1)node[above,yshift=-1pt]{$i$};
\draw[solid](0,0)node[below]{$i$}--++(0,1)node[above,yshift=-1pt]{$\phantom{i}$};
\draw[solid](0.5,0)node[below]{$i$}--++(0,1)node[above,yshift=-1pt]{$\phantom{i}$};
\draw[redstring](0.25,0)node[below]{$\rho$}--++(0,1)node[above,yshift=-1pt]{$\phantom{i}$};
\end{tikzpicture}
,\quad
\begin{tikzpicture}[scale=1.2,anchorbase,smallnodes,rounded corners]
\draw[ghost,dot](1,0)--++(0,1)node[above,yshift=-1pt]{$i$};
\draw[ghost](1.5,0)--++(0,1)node[above,yshift=-1pt]{$i$};
\draw[solid,dot](0,0)node[below]{$i$}--++(0,1)node[above,yshift=-1pt]{$\phantom{i}$};
\draw[solid](0.5,0)node[below]{$i$}--++(0,1)node[above,yshift=-1pt]{$\phantom{i}$};
\draw[redstring](0.25,0)node[below]{$\rho$}--++(0,1)node[above,yshift=-1pt]{$\phantom{i}$};
\end{tikzpicture}
+
\begin{tikzpicture}[scale=1.2,anchorbase,smallnodes,rounded corners]
\draw[ghost](1,0)--++(0,1)node[above,yshift=-1pt]{$i$};
\draw[ghost,dot](1.5,0)--++(0,1)node[above,yshift=-1pt]{$i$};
\draw[solid](0,0)node[below]{$i$}--++(0,1)node[above,yshift=-1pt]{$\phantom{i}$};
\draw[solid,dot](0.5,0)node[below]{$i$}--++(0,1)node[above,yshift=-1pt]{$\phantom{i}$};
\draw[redstring](0.25,0)node[below]{$\rho$}--++(0,1)node[above,yshift=-1pt]{$\phantom{i}$};
\end{tikzpicture}
,\quad
\begin{tikzpicture}[scale=1.2,anchorbase,smallnodes,rounded corners]
\draw[ghost,dot](1,0)--++(0,1)node[above,yshift=-1pt]{$i$};
\draw[ghost,dot](1.5,0)--++(0,1)node[above,yshift=-1pt]{$i$};
\draw[solid,dot](0,0)node[below]{$i$}--++(0,1)node[above,yshift=-1pt]{$\phantom{i}$};
\draw[solid,dot](0.5,0)node[below]{$i$}--++(0,1)node[above,yshift=-1pt]{$\phantom{i}$};
\draw[redstring](0.25,0)node[below]{$\rho$}--++(0,1)node[above,yshift=-1pt]{$\phantom{i}$};
\end{tikzpicture}
.
\end{gather*}
for suitable choices of ghost shifts etc.
\end{Example}

Let $Z({}_{-})$ be the center.

\begin{Lemma}\label{L:A1Case}
Let $\Gamma$ be of type $A^{(1)}_{1}$ and $\charge=\emptyset$. Then there is an isomorphism of graded algebras
\begin{gather*}
Z\bigl(\WA(X)\bigr)
\cong
R[y_{1},\dots,y_{n}]^{\Sym[n]}.
\end{gather*}
\end{Lemma}

\begin{proof}
By \autoref{P:VaryingSize}, if $\charge=\emptyset$, then $\WA(X)$ is Morita equivalent to the nil Hecke algebra for any choice of positions for $X$. Hence, we can prove this statement {\muta} as
for the nil Hecke algebra, see \cite[Section 2]{Ma-symmetric-schubert-loci} and \cite[Proposition 3.5]{La-categorification-sl2}.
\end{proof}

Almost exactly as in \cite[Theorem 2.9]{KhLa-cat-quantum-sln-first} we have the following.

\begin{Proposition}
There is an isomorphism
$Z\bigl(\WA(X)\bigr)\cong R[\beta]^{\Sym[\beta]}$ of graded
algebras given by
\begin{gather*}
\nu\colon
R[\beta]^{\Sym[\beta]}\to Z\bigl(\WA(X)\bigr),
f\mapsto f\onealg.
\end{gather*}
Moreover, $\WA(X)$ is free and of finite rank over its center.
\end{Proposition}

\begin{proof}
Using the relations in
\autoref{D:RationalCherednik} it is easy to see that for any symmetric polynomial $f\in R[\beta]^{\Sym[\beta]}$ the element
$f\onealg$ is central. The only relation for which commutativity is not immediate is \autoref{R:DotCrossing} and, in
this situation, symmetrically placed dots slide freely. Hence, $\nu$ is a well-defined algebra homomorphism and it is injective by \autoref{P:WABasis}.

Let $\1_{X,\bi}=\sum_{\bx\in X}\1_{\bx,\bi}$. Using the same arguments as in \cite[Theorem 2.9]{KhLa-cat-quantum-sln-first} to prove that $\nu$ is surjective it suffices to show that
the composition
\begin{gather*}
R[\beta]^{\Sym[\beta]}
\xrightarrow{\nu}
Z\bigl(\WA(X)\bigr)
\xrightarrow{\1_{X,\bi}({}_{-})\1_{X,\bi}}
Z\bigl(\1_{X,\bi}\WA(X)\1_{X,\bi}\bigr)
\end{gather*}
is an isomorphism for all $\bi\in I^{\beta}$, where
$\bi=(\underbrace{i_{1},\dots,i_{1}}_{a_{1}},\dots,\underbrace{i_{r},\dots,i_{r}}_{a_{r}})$.
For such $\bi$ we have
\begin{gather*}
\1_{X,\bi}\WA(X)\1_{X,\bi}
\cong
\bigotimes_{k=1}^{r}
\WA[a_{k}](X),
\end{gather*}
where each algebra appearing on the right-hand side is the subalgebra of $\1_{X,\bi}\WA(X)\1_{X,\bi}$ having only identities outside of the indicated region.
Hence, we have identified the center since,
by \autoref{L:A1Case}, the center of $\WA[a_{k}](X)$ is $R[y_{1},\dots,y_{a_{k}}]^{\Sym[a_{k}]}$.

The final claim now follows from \autoref{P:WABasis}.
\end{proof}

Standard arguments now yield:

\begin{Proposition}
\leavevmode
\begin{enumerate}

\item The algebra $\WA(X)$ is left and right Noetherian.

\item The algebra $\WA(X)$ is indecomposable.

\item Suppose that $R$ is a field. Then every simple
$\WA(X)$-module is finite dimensional.\qed

\end{enumerate}
\end{Proposition}


\subsection{Cyclotomic quotients}\label{SS:CyclotomicQuotients}


We now define \emph{cyclotomic or steadied quotients} of $\WA(X)$. As we will see in \autoref{P:WebAlg}, these quotients generalize cyclotomic KLR(W) algebras.

Recall that $X$ is the set of allowable endpoints of the solid and red strings. A string is \emph{bounded} by $X$ if the $x$-coordinates of its points are bounded, on the left and right, by the $x$-coordinates of the points in $X$.

\begin{Definition}
A diagram $\1_{\bx,\bi}$ is (right) \emph{unsteady}
if it contains a solid string that can be pulled arbitrarily far to the right when the red strings are bounded by $X$. A diagram is \emph{unsteady} if it contains an unsteady string and otherwise it is \emph{steady}.
\end{Definition}

\begin{Definition}
The \emph{cyclotomic weighted KLRW algebra} $\WAc(X)$
is the quotient of $\WA(X)$ by the two-sided ideal generated by all
diagrams that factor through an unsteady
idempotent diagram.
\end{Definition}

Similarly, we can define left unsteady diagrams. Reflecting diagrams shows that a quotient algebra that is defined by factoring out by the two-sided ideal of diagrams that factor through some left unsteady idempotent diagram is isomorphic to a cyclotomic weighted KLRW algebra and {\vive}.
We work with right unsteady diagrams because we already have a left-right bias for the ghosts strings in \autoref{D:WeightedKLRW}.

\begin{Example}
If $\rho=i$, then following diagrams are unsteady and steady, respectively.
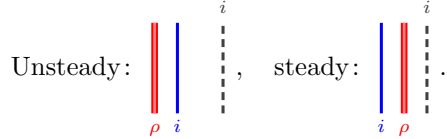
\begin{gather*}
\text{Unsteady}\colon
\begin{tikzpicture}[scale=1.2,anchorbase,smallnodes,rounded corners]
\draw[ghost](1,0)--++(0,1)node[above,yshift=-1pt]{$i$};
\draw[solid](0.5,0)node[below]{$i$}--++(0,1)node[above,yshift=-1pt]{$\phantom{i}$};
\draw[redstring](0.25,0)node[below]{$\rho$}--++(0,1)node[above,yshift=-1pt]{$\phantom{i}$};
\end{tikzpicture}
,\quad
\text{steady}\colon
\begin{tikzpicture}[scale=1.2,anchorbase,smallnodes,rounded corners]
\draw[ghost](0.5,0)--++(0,1)node[above,yshift=-1pt]{$i$};
\draw[solid](0,0)node[below]{$i$}--++(0,1)node[above,yshift=-1pt]{$\phantom{i}$};
\draw[redstring](0.25,0)node[below]{$\rho$}--++(0,1)node[above,yshift=-1pt]{$\phantom{i}$};
\end{tikzpicture}
.
\end{gather*}
In the right-hand diagram the solid string cannot be pulled further
to the right because it cannot be pulled past the red string, so the diagram is steady. However, if $\rho\neq i$
we can use \autoref{R:RedSolid}:
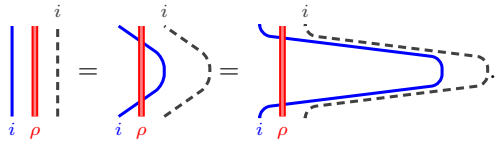
\begin{gather*}
\begin{tikzpicture}[scale=1.2,anchorbase,smallnodes,rounded corners]
\draw[ghost](0.5,0)--++(0,1)node[above,yshift=-1pt]{$i$};
\draw[solid](0,0)node[below]{$i$}--++(0,1)node[above,yshift=-1pt]{$\phantom{i}$};
\draw[redstring](0.25,0)node[below]{$\rho$}--++(0,1)node[above,yshift=-1pt]{$\phantom{i}$};
\end{tikzpicture}
=
\begin{tikzpicture}[scale=1.2,anchorbase,smallnodes,rounded corners]
\draw[ghost](0.5,0)--++(0.5,0.35)--++(0,0.3)--++(-0.5,0.35)node[above,yshift=-1pt]{$i$};
\draw[solid](0,0)node[below]{$i$}--++(0.5,0.35)--++(0,0.3)--++(-0.5,0.35)node[above,yshift=-1pt]{$\phantom{i}$};
\draw[redstring](0.25,0)node[below]{$\rho$}--++(0,1)node[above,yshift=-1pt]{$\phantom{i}$};
\end{tikzpicture}
=
\begin{tikzpicture}[scale=1.2,anchorbase,smallnodes,rounded corners]
\draw[ghost](0.5,0)--++(0,0.1)--++(2,0.25)--++(0,0.3)--++(-2,0.25)--++(0,0.1)node[above,yshift=-1pt]{$i$};
\draw[solid](0,0)node[below]{$i$}--++(0,0.1)--++(2,0.25)--++(0,0.3)--++(-2,0.25)--++(0,0.1)node[above,yshift=-1pt]{$\phantom{i}$};
\draw[redstring](0.25,0)node[below]{$\rho$}--++(0,1)node[above,yshift=-1pt]{$\phantom{i}$};
\end{tikzpicture}.
\end{gather*}
Hence, if $\rho\neq i$, then both of the diagrams above are unsteady, and so zero in $\WAc(X)$.
\end{Example}

We abuse notation and identify the elements of $\WA(X)$ with their images in $\WAc(X)$.

\begin{Proposition}
The algebra $\WAc(X)$ is finite dimensional.
\end{Proposition}

\begin{proof}
By \autoref{P:WABasis}, the algebra $\WAc(X)$ is spanned by the image of $\WABasis$ in $\WAc(X)$. Hence, it is enough to show that the dotted idempotent
$y_{1}^{a_{1}}\dots y_{n}^{a_{n}}\1_{\bx,\bi}$ is zero when $a_{1},\dots,a_{n}\in\N$ are big enough. We prove this by induction on the number $\ell$ of red strings.
If $\ell=0$, then
any diagram is unsteady and the claim follows. If $\ell>0$, then
take the leftmost red string and
record the minimal number $m$ of dots that we need to add to the solid
strings so that \autoref{R:RedSolid} can be applied to
pull the red string to the left of all of the solid strings. Now remove this red string from the diagram. By induction there exist $a_{1}^{\prime},\dots,a_{n}^{\prime}$
annihilating the dotted idempotent. Setting $a_{r}=a_{r}^{\prime}+m$ now implies that
$y_{1}^{a_{1}}\dots y_{n}^{a_{n}}\1_{\bx,\bi}=0$, completing the proof.
\end{proof}

As we will see below,
the unsteady condition, which crucially depends on the
choice of red strings, corresponds to the cyclotomic ideal in the KLR world.
It is difficult to tell what elements belong to the cyclotomic ideal. In contrast, it is very easy to tell when a diagram is unsteady, which is one of the main advantages of the weighted KLRW framework.


\subsection{Induction and restriction}


We now discuss the analog of \cite[Section 2.6]{KhLa-cat-quantum-sln-first} and \cite[Section 2.4]{We-weighted-klr} for weighted KLRW algebras, which apply almost without change.
To this end, recall that $X$ is a set of positionings. For each such set
and each $z\in\R$ we let $X\set{z}=\set{\bx+z|\bx\in X}$ be the set of
shifted positionings.
More generally, set $XY\{z\}=X\cup(Y\{z\})$ (we use similar notations below) and let $\mathrm{d}(X,Y)=\max(X)-\min(Y)$, where $\min$ and $\max$ have the obvious meanings.

The following is clear from \autoref{P:WABasis}.

\begin{Proposition}\label{P:IndRes}
For $\beta,\beta^{\prime}\in Q^{+}$ and $z>\mathrm{d}(X,X^{\prime})\in\R$, we have an embedding of graded algebras
\begin{gather*}
\iota_{\beta,\beta^{\prime}}\map{\WA(X)\otimes\WA[\beta^{\prime}][\brho^{\prime}](X^{\prime})}{\WA[\beta+\beta^{\prime}][\brho+\brho^{\prime}](XX^{\prime}\{z\})},
\quad
\begin{tikzpicture}[scale=1.2,anchorbase,smallnodes,rounded corners]
\node[rectangle,draw,minimum width=0.5cm,minimum height=0.5cm,ultra thick] at(0,0){\raisebox{-0.05cm}{$D$}};
\end{tikzpicture}
\otimes
\begin{tikzpicture}[scale=1.2,anchorbase,smallnodes,rounded corners]
\node[rectangle,draw,minimum width=0.5cm,minimum height=0.5cm,ultra thick] at(0,0){\raisebox{-0.05cm}{$E$}};
\end{tikzpicture}
\mapsto
\begin{tikzpicture}[scale=1.2,anchorbase,smallnodes,rounded corners]
\node[rectangle,draw,minimum width=0.5cm,minimum height=0.5cm,ultra thick] at(0.5,0){\raisebox{-0.05cm}{$E$}};
\node[rectangle,draw,minimum width=0.5cm,minimum height=0.5cm,ultra thick] at(0.05,0){\raisebox{-0.05cm}{$D$}};
\end{tikzpicture}
.
\end{gather*}
Moreover, $\Ind{\1}_{\beta,\beta^{\prime}}=\iota_{\beta,\beta^{\prime}}(\onealg\otimes\1_{\WA[\beta^{\prime}][\brho^{\prime}](X^{\prime})})$ is an idempotent.\qed
\end{Proposition}

For any $z>\mathrm{d}(X,X^{\prime})\in\R$ we have a restriction functor
\begin{gather*}
\Res{F}^{\beta,\beta^{\prime}}\map{\WA[\beta+\beta^{\prime}](XX^{\prime}\{z\})\text{-}\mathbf{Mod}_{\Z}}{\bigl(\WA(X)\otimes\WA[\beta^{\prime}][\brho^{\prime}](X^{\prime})\bigr)\text{-}\mathbf{Mod}_{\Z}},
\end{gather*}
that sends $M$ to $\Ind{\1}_{\beta,\beta^{\prime}}M$. There are also induction and coinduction functors
\begin{gather*}
\Ind{F}_{\beta,\beta^{\prime}}\map{\bigl(\WA(X)\otimes\WA[\beta^{\prime}][\brho^{\prime}](X^{\prime})\bigr)\text{-}\mathbf{Mod}_{\Z}}{\WA[\beta+\beta^{\prime}](XX^{\prime}\{z\})\text{-}\mathbf{Mod}_{\Z}},
\\
\Ind{F}_{\beta,\beta^{\prime}}=\WA(XX^{\prime}\{z\})\Ind{\1}_{\beta,\beta^{\prime}}
\otimes_{\WA(X)\otimes\WA(X^{\prime})}{}_{-},
\\
\CoInd{F}_{\beta,\beta^{\prime}}\map{\bigl(\WA(X)\otimes\WA[\beta^{\prime}][\brho^{\prime}](X^{\prime})\bigr)\text{-}\mathbf{Mod}_{\Z}}{\WA(XX^{\prime}\{z\})\text{-}\mathbf{Mod}_{\Z}},
\\
\CoInd{F}_{\beta,\beta^{\prime}}=\mathrm{Hom}_{(\WA(X)\otimes\WA[\beta^{\prime}][\brho^{\prime}](X^{\prime}))\text{-}\mathbf{Mod}_{\Z}}\big(\WA[\beta+\beta^{\prime}](XX^{\prime}\{z\}),{}_{-}\big)
.
\end{gather*}

\begin{Remark}
A straight line diagram argument implies that
the functors $\Res{F}^{\beta,\beta^{\prime}}$, $\Ind{F}_{\beta,\beta^{\prime}}$ and $\CoInd{F}_{\beta,\beta^{\prime}}$
do not depend on $z$, which is why we do not include $z$ in the notation.
\end{Remark}

Let $P_{\bx,\bi}=\WA(X)\1_{\bx,\bi}$ be the projective
$\WA(X)$-module generated by the idempotent $\1_{\bx,\bi}$.

\begin{Proposition}
\leavevmode
\begin{enumerate}

\item The functor $\Ind{F}_{\beta,\beta^{\prime}}$ is the left adjoint of $\Res{F}^{\beta,\beta^{\prime}}$.

\item The functor $\CoInd{F}_{\beta,\beta^{\prime}}$ is the right adjoint of $\Res{F}^{\beta,\beta^{\prime}}$.

\item The functors $\Res{F}^{\beta,\beta^{\prime}}$, $\Ind{F}_{\beta,\beta^{\prime}}$ and $\CoInd{F}_{\beta,\beta^{\prime}}$ are $R$-linear additive and homogeneous.

\item The functor $\Res{F}^{\beta,\beta^{\prime}}$ is exact, $\Ind{F}_{\beta,\beta^{\prime}}$ is right exact and
$\CoInd{F}_{\beta,\beta^{\prime}}$ is left exact.

\item We have
$\Ind{F}_{\beta,\beta^{\prime}}(P_{\bx,\bi}\otimes P_{\by,\bj})\cong P_{\bx\by\{z\},\bi\bj}$
as graded $\WA[\beta+\beta^{\prime}](XX^{\prime}\{z\})$-modules.

\item For $\bz=\bx\by\{z\}$ we have $\Res{F}^{\beta,\beta^{\prime}}(P_{\bz,\bk})\cong\bigoplus_{\bi\ast\bj=\bk}(P_{\bx,\bi}\otimes P_{\by,\bj})$
as graded $\WA(X)\otimes\WA[\beta^{\prime}](X^{\prime})$-modules, where the direct sum runs over all shuffles of $\bk$ (see \cite[Section 2.6]{KhLa-cat-quantum-sln-first}).

\end{enumerate}
\end{Proposition}

\begin{proof}
The first four statements follow from the usual Yoga, while
the final two claims can be proven {\muta} as in \cite[Section 2.6]{KhLa-cat-quantum-sln-first}.
\end{proof}


\subsection{Relationship to (non-weighted) KLRW algebras}\label{SS:KLRW}


Let $\TA$ be Webster's tensor product
algebra attached to the datum $\beta$ for solid, and $\charge$ and $\brho$ for red strings. We will not recall the definition of
$\TA$, which is given in \cite[Chapter 4]{We-knot-invariants} via string diagrams. We note that
$\TA[n][\brho]$ is the algebra $\tilde{T}_{n}^{\lambda}$ in Webster's notation;
see \autoref{R:BetaVsN}.
By \cite[Theorem 4.18]{We-knot-invariants}, the algebra $\TA$ is a generalization of the
\emph{KLR algebra} attached to $\beta$ \cite{KhLa-cat-quantum-sln-first}, \cite{Ro-2-kac-moody}, so we call $\TA$ a \emph{KLRW algebra}. Let $\TAc$ be the \emph{cyclotomic quotient} of $\TA$, as defined in \cite[Chapter 4]{We-knot-invariants}.
Finally, the algebra $\TA$ has a relation that is a
symmetric version of \autoref{R:GhostSolid}, which will be
crucial in the proof of \autoref{P:WebAlg} below. Namely:
\begin{gather*}
\begin{tikzpicture}[scale=1.2,anchorbase,smallnodes,rounded corners]
\draw[solid](0,1)node[above,yshift=-1pt]{$\phantom{i}$}--++(0.5,-0.5)--++(-0.5,-0.5) node[below]{$i$};
\draw[solid](0.5,1)--++(-0.5,-0.5)--++(0.5,-0.5) node[below]{$j$};
\end{tikzpicture}
=Q_{ij}(\y)
\begin{tikzpicture}[scale=1.2,anchorbase,smallnodes,rounded corners]
\draw[solid](0,1)node[above,yshift=-1pt]{$\phantom{i}$}--++(0,-1)node[below]{$i$};
\draw[solid](0.5,1)--++(0,-1)node[below]{$j$};
\end{tikzpicture}
,\quad
\begin{tikzpicture}[scale=1.2,anchorbase,smallnodes,rounded corners]
\draw[solid](0.5,1)node[above,yshift=-1pt]{$\phantom{i}$}--++(-0.5,-0.5)--++(0.5,-0.5) node[below]{$i$};
\draw[solid](0,1)--++(0.5,-0.5)--++(-0.5,-0.5) node[below]{$j$};
\end{tikzpicture}
=Q_{ji}(\y)
\begin{tikzpicture}[scale=1.2,anchorbase,smallnodes,rounded corners]
\draw[solid](0.5,1)node[above,yshift=-1pt]{$\phantom{i}$}--++(0,-1)node[below]{$i$};
\draw[solid](0,1)--++(0,-1)node[below]{$j$};
\end{tikzpicture},
\end{gather*}
regardless of the orientation of the underlying quiver.
The definition of $\TA$ and $\TAc$ also involves a choice of positions $(\bx,\charge)$ for the strands in \cite[Chapter 4]{We-knot-invariants}. By conjugating by straight line diagrams, we can assume that $\bx=(x_{1},\dots,x_{n})$ are nonintegral points and that $\charge$ consists of integral points. We call this a \emph{KLRW positioning}.

The following should be compared with \cite[Proposition 2.14]{We-weighted-klr}.

\begin{Proposition}\label{P:WebAlg}
Let $X\simeq\set{\bx}$ where $\bx$ is a KLRW positioning. Suppose that we are in the infinitesimal-case and $\Gamma$ has no parallel edges. We also fix one choice of $Q$-polynomials.
\begin{enumerate}

\item The algebras $\TA$ and $\WA(X)$ are graded Morita equivalent.

\item Furthermore, for any $\by\in X$ with $\bx\simeq\by$, we have an isomorphism of graded algebras
\begin{gather*}
\begin{gathered}
\mu_{T}\map{\TA}{1_{\by}\WA(X)1_{\by}},
\quad
\begin{tikzpicture}[scale=1.2,anchorbase,smallnodes,rounded corners]
\draw[solid] (0,0)node[below]{$i$}node[below]{$\phantom{j}$} to (0,1)node[above,yshift=-1pt]{$\phantom{i}$};
\end{tikzpicture}
\mapsto
\begin{tikzpicture}[scale=1.2,anchorbase,smallnodes,rounded corners]
\draw[ghost] (0.15,0) to (0.15,1)node[above,yshift=-1pt]{$i$};
\draw[solid] (0,0)node[below]{$i$}node[below]{$\phantom{j}$} to (0,1)node[above,yshift=-1pt]{$\phantom{i}$};
\end{tikzpicture}
.
\end{gathered}
\end{gather*}

\item The isomorphism $\mu_{T}$ descends to an isomorphism of graded algebras
$\tilde{\mu}_{T}\map{\TAc}{1_{\by}\WAc(X)1_{\by}}$.

\end{enumerate}
\end{Proposition}

\begin{proof}[Sketch of proof.]
The assignment $\mu_{T}$ is clearly homogeneous. Moreover, we can use the
polynomial module to show that $\mu_{T}$ is injective, and
thus bijective, by \autoref{P:WABasis} and \cite[Proposition 4.16]{We-knot-invariants}, which is the corresponding statement for
$\TA$. Hence, it suffices to prove that $\mu_{T}$
is well-defined. This follows from the combinatorics of how solid and ghost
strings interact. For example, we have
\begin{gather*}
\begin{tikzpicture}[scale=1.2,anchorbase]
\node[circle,inner sep=1.8pt,fill=DarkBlue] (0) at (0,0){};
\node at (0,-0.25){$i$};
\node[circle,inner sep=1.8pt,fill=DarkBlue] (1) at (1,0){};
\node at (1,-0.25){$j$};
\draw[->](0) to node[above,yshift=-1pt]{$\varepsilon$} (1);
\end{tikzpicture}
\text{ or }
\begin{tikzpicture}[scale=1.2,anchorbase]
\node[circle,inner sep=1.8pt,fill=DarkBlue] (0) at (0,0){};
\node at (0,-0.25){$i$};
\node[circle,inner sep=1.8pt,fill=DarkBlue] (1) at (1,0){};
\node at (1,-0.25){$j$};
\draw[->](1) to node[above,yshift=-1pt]{$\varepsilon$} (0);
\end{tikzpicture}
:\qquad
\begin{tikzpicture}[scale=1.2,anchorbase,smallnodes,rounded corners]
\draw[solid](0,1)node[above,yshift=-1pt]{$\phantom{i}$}--++(0.5,-0.5)--++(-0.5,-0.5) node[below]{$i$};
\draw[solid](0.5,1)--++(-0.5,-0.5)--++(0.5,-0.5) node[below]{$j$};
\end{tikzpicture}
=Q_{ij}(\by)
\begin{tikzpicture}[scale=1.2,anchorbase,smallnodes,rounded corners]
\draw[solid](0,1)node[above,yshift=-1pt]{$\phantom{i}$}--++(0,-1)node[below]{$i$};
\draw[solid](0.5,1)--++(0,-1)node[below]{$j$};
\end{tikzpicture}
\mapsto
\begin{tikzpicture}[scale=1.2,anchorbase,smallnodes,rounded corners]
\draw[ghost](0.1,1)node[above,yshift=-1pt]{$i$}--++(0.5,-0.5)--++(-0.5,-0.5);
\draw[ghost](0.6,1)node[above,yshift=-1pt]{$j$}--++(-0.5,-0.5)--++(0.5,-0.5);
\draw[solid](0,1)node[above,yshift=-1pt]{$\phantom{i}$}--++(0.5,-0.5)--++(-0.5,-0.5)node[below]{$i$};
\draw[solid](0.5,1)--++(-0.5,-0.5)--++(0.5,-0.5)node[below]{$j$};
\end{tikzpicture}
=Q_{ij}(\y)
\begin{tikzpicture}[scale=1.2,anchorbase,smallnodes,rounded corners]
\draw[ghost](0.1,1)node[above,yshift=-1pt]{$i$}--++(0,-1);
\draw[ghost](0.6,1)node[above,yshift=-1pt]{$j$}--++(0,-1);
\draw[solid](0,1)node[above,yshift=-1pt]{$\phantom{i}$}--++(0,-1)node[below]{$i$};
\draw[solid](0.5,1)--++(0,-1)node[below]{$j$};
\end{tikzpicture}.
\end{gather*}
for the two orientations given above.
The left-hand side is a defining relation in $\TA$, the right-hand side
is a relation in $\WA(X)$, which holds by
\autoref{R:GhostSolid} and the fact that solid-solid and ghost-ghost strings
satisfy the Reidemeister II relation. In more detail, for the two orientations considered
above we have:
\begin{gather*}
\begin{tikzpicture}[scale=1.2,anchorbase,smallnodes,rounded corners]
\draw[ghost](0.1,1)node[above,yshift=-1pt]{$i$}--++(0.5,-0.5)--++(-0.5,-0.5);
\draw[ghost](0.6,1)node[above,yshift=-1pt]{$j$}--++(-0.5,-0.5)--++(0.5,-0.5)node[below]{$j$};
\draw[solid](0,1)node[above,yshift=-1pt]{$\phantom{i}$}--++(0.5,-0.5)--++(-0.5,-0.5)node[below]{$i$};
\draw[solid](0.5,1)--++(-0.5,-0.5)--++(0.5,-0.5)node[below]{$j$};
\end{tikzpicture}
=Q_{ij}(\y)
\begin{tikzpicture}[scale=1.2,anchorbase,smallnodes,rounded corners]
\draw[ghost](0.1,1)node[above,yshift=-1pt]{$i$}--++(0.5,-0.5)--++(-0.5,-0.5);
\draw[ghost](0.6,1)node[above,yshift=-1pt]{$j$}--++(-0.25,-0.5)--++(0.25,-0.5);
\draw[solid](0,1)node[above,yshift=-1pt]{$\phantom{i}$}--++(0.25,-0.5)--++(-0.25,-0.5)node[below]{$i$};
\draw[solid](0.5,1)--++(-0.5,-0.5)--++(0.5,-0.5)node[below]{$j$};
\end{tikzpicture}
=Q_{ij}(\y)
\begin{tikzpicture}[scale=1.2,anchorbase,smallnodes,rounded corners]
\draw[ghost](0.1,1)node[above,yshift=-1pt]{$i$}--++(0.3,-0.5)--++(-0.3,-0.5);
\draw[ghost](0.6,1)node[above,yshift=-1pt]{$j$}--++(-0,-0.5)--++(0,-0.5);
\draw[solid](0,1)node[above,yshift=-1pt]{$\phantom{i}$}--++(0,-0.5)--++(-0,-0.5)node[below]{$i$};
\draw[solid](0.5,1)--++(-0.3,-0.5)--++(0.3,-0.5)node[below]{$j$};
\end{tikzpicture}
=Q_{ij}(\y)
\begin{tikzpicture}[scale=1.2,anchorbase,smallnodes,rounded corners]
\draw[ghost](0.1,1)node[above,yshift=-1pt]{$i$}--++(0,-1);
\draw[ghost](0.6,1)node[above,yshift=-1pt]{$j$}--++(0,-1);
\draw[solid](0,1)node[above,yshift=-1pt]{$\phantom{i}$}--++(0,-1)node[below]{$i$};
\draw[solid](0.5,1)--++(0,-1)node[below]{$j$};
\end{tikzpicture}
,
\begin{tikzpicture}[scale=1.2,anchorbase,smallnodes,rounded corners]
\draw[ghost](0.1,1)node[above,yshift=-1pt]{$i$}--++(0.5,-0.5)--++(-0.5,-0.5);
\draw[ghost](0.6,1)node[above,yshift=-1pt]{$j$}--++(-0.5,-0.5)--++(0.5,-0.5)node[below]{$j$};
\draw[solid](0,1)node[above,yshift=-1pt]{$\phantom{i}$}--++(0.5,-0.5)--++(-0.5,-0.5)node[below]{$i$};
\draw[solid](0.5,1)--++(-0.5,-0.5)--++(0.5,-0.5)node[below]{$j$};
\end{tikzpicture}
=
\begin{tikzpicture}[scale=1.2,anchorbase,smallnodes,rounded corners]
\draw[ghost](0.1,1)node[above,yshift=-1pt]{$i$}--++(0.5,-0.5)--++(-0.5,-0.5);
\draw[ghost](0.6,1)node[above,yshift=-1pt]{$j$}--++(-0.25,-0.5)--++(0.25,-0.5);
\draw[solid](0,1)node[above,yshift=-1pt]{$\phantom{i}$}--++(0.25,-0.5)--++(-0.25,-0.5)node[below]{$i$};
\draw[solid](0.5,1)--++(-0.5,-0.5)--++(0.5,-0.5)node[below]{$j$};
\end{tikzpicture}
=
\begin{tikzpicture}[scale=1.2,anchorbase,smallnodes,rounded corners]
\draw[ghost](0.1,1)node[above,yshift=-1pt]{$i$}--++(0.3,-0.5)--++(-0.3,-0.5);
\draw[ghost](0.6,1)node[above,yshift=-1pt]{$j$}--++(-0,-0.5)--++(0,-0.5);
\draw[solid](0,1)node[above,yshift=-1pt]{$\phantom{i}$}--++(0,-0.5)--++(-0,-0.5)node[below]{$i$};
\draw[solid](0.5,1)--++(-0.3,-0.5)--++(0.3,-0.5)node[below]{$j$};
\end{tikzpicture}
=Q_{ij}(\y)
\begin{tikzpicture}[scale=1.2,anchorbase,smallnodes,rounded corners]
\draw[ghost](0.1,1)node[above,yshift=-1pt]{$i$}--++(0,-1);
\draw[ghost](0.6,1)node[above,yshift=-1pt]{$j$}--++(0,-1);
\draw[solid](0,1)node[above,yshift=-1pt]{$\phantom{i}$}--++(0,-1)node[below]{$i$};
\draw[solid](0.5,1)--++(0,-1)node[below]{$j$};
\end{tikzpicture}
,
\end{gather*}
where we have straightened some of the $i$-strings to improve readability.
Note that the appearance of the $Q$-polynomial in the $\WA(X)$ relations
depends on the orientation. Moreover, the ghost $j$-string and solid $i$-string,
do not play a role in these calculations, so this argument
also works when these ghost strings do not appear. Note that multiple ghosts $i$ or $j$-strings do not affect the argument as the relations of $\WA(X)$ depend on edges in the quiver and not on the vertices. (This calculation requires the assumption that $\Gamma$ does not have parallel edges.) All
other relations can be checked {\muta},
proving (b). The claim in (a) now follows from (b) using \autoref{P:VaryingSize}.
For the final claim, observe that $\mu_{T}$ gives a bijection between unsteady diagrams for the two algebras.
\end{proof}

\begin{Example}
Let $\Gamma$ be $A_{\Z}$ or $A_{e}^{(1)}$ in \autoref{E:Quivers}. Let $\ell=1$, $\bx=0.9\cdot(1,\dots,n)$ and $\kappa=(n)$. Let
$\bsig$ be constant $0.01$ and define $Q_{ij}(\y)$ as in \autoref{E:QPoly}. Then $\1_{\bx}\WA(X)\1_{\bx}$ is graded isomorphic to the KLR algebra from \cite{KhLa-cat-quantum-sln-first} or \cite{Ro-2-kac-moody} of the corresponding type.
\end{Example}

For the rest of the paper we will identify $\1_{\bx}\WA(X)\1_{\bx}$ and the
KLR(W) algebras of \cite{KhLa-cat-quantum-sln-first}, \cite{Ro-2-kac-moody},
\cite{We-knot-invariants} using the choices in \autoref{P:WebAlg}. Similarly, we also identify the cyclotomic quotients of these algebras.

\begin{Remark}
For suitable choices of $X$, the KLR algebra is an idempotent subalgebra of $\WA(X)$ and the cyclotomic KLR algebra is an idempotent subalgebra of $\WAc(X)$. In particular, the weighted KLRW algebras and KLR algebras usually have different numbers of simple modules.
\end{Remark}


\section{Varying the quiver}\label{S:Subdivision}


The aim of this section is to define
isomorphisms between weighted KLRW algebras
attached to different quivers.
Our main definition is
\autoref{D:DiagramSubdivision} and partially inspired by
\cite{ChMi-runner-removal} and
\cite{Ma-catrep-klr}. The isomorphism in
\autoref{T:SubDiv} makes it possible to compare the representation theories of
weighted KLRW algebras for different quivers, for example
for quivers of types $A^{(1)}_{e-1}$ and $A^{(1)}_{e}$.
Before coming to the generalization of these results, for completeness, we briefly discuss a much easier (but less interesting), way to vary the quiver.


\subsection{Another induction and restriction}


We keep using the conventions in \autoref{N:Quiver}.

\begin{Definition}
Let $\Gamma=(\Gamma,\bsig)$ and
$\Indz{\Gamma}=(\Indz{\Gamma},\Indz{\bsig})$ be two quivers as in
\autoref{SS:Quiver} together with a choice of ghost shifts.
Then $\Indz{\Gamma}$ is an
\emph{induction} of $\Gamma$ and $\Gamma$ is a
\emph{restriction} of $\Indz{\Gamma}$
if $\Gamma$ is a weighted subgraph ({\cf} \autoref{R:Weighting}) of
$\Indz{\Gamma}$.
\end{Definition}

\begin{Example}\label{Ex:EasyQuiverVary}
If $d<e$, then $(A_{d},\bsig)$ is a restriction of $(A^{(1)}_{e},\Indz{\bsig})$.
In contrast, the quiver $A^{(1)}_{d}$ is an induction or a restriction of the quiver $A^{(1)}_{e}$ if and only if $e=d$ (and the ghost shifts match).
\end{Example}

In the following we write $\Indz{{}_{-}}$ for everything related to $\Indz{\Gamma}$.

\begin{Proposition}
For any induction $\Indz{\Gamma}$
we have an embedding and a projection of graded algebras
\begin{gather*}
\Indz{f}\colon\WA[n][\rho](X)\hookrightarrow\Indz{\WA[n][\rho](X)},
\quad
\Resz{f}\colon\Indz{\WA[n][\rho](X)}\twoheadrightarrow\WA[n][\rho](X),
\end{gather*}
given by sending each generator of $\WA[n][\rho](X)$
to the same named element in $\Indz{\WA[n][\rho](X)}$ respectively
by annihilation $\1_{\bx,\bi}$ for $\bi$ with $i_{k}\in\Indz{I}\setminus I$ for some $k=1,\dots,n$.
Moreover, $\Resz{f}\circ\Indz{f}=\mathrm{Id}_{\WA[n][\rho](X)}$.
\end{Proposition}

\begin{proof}
The only claim that is not immediate by construction is that $\Indz{f}$ is an embedding. This however follows from \autoref{P:WABasis}.
\end{proof}

We thus get the associated restriction, induction and coinduction functors, and the analog of \autoref{P:IndRes}, all of which we leave to the reader to spell out.


\subsection{Subdividing quivers}\label{SS:Subdivision}


Now we come to one of our main definitions.

\begin{Definition}\label{D:Subdvisionedge}
We call \Ddots{2} a \emph{subdivision} of the simply laced edge \Ddots{1}.
\end{Definition}

\begin{Definition}\label{D:Subdvision}
Let $\Gamma$ and $\Sub{\Gamma}$ be two quivers as in
\autoref{SS:Quiver}. Then $\Sub{\Gamma}$ is a
\emph{(simply laced) subdivision} of $\Gamma$ if $\Sub{\Gamma}$ is obtained from
$\Gamma$ by subdividing a finite number of simply laced edges and, potentially, relabeling the vertices.
\end{Definition}

By definition, any subdivision of $\Gamma$ is obtained by successively replacing edges \Ddots{1} in a quiver with \Ddots{2}. In particular, the orientations of the subdivided edges are compatible with the original orientation of the subdivided edge.

\begin{Examples}\label{Ex:Subdvision}
With respect to the quivers in \autoref{E:Quivers} we have:
\begin{enumerate}

\item The quiver $A_{\Z}$ is a subdivision of itself.

\item If $e\geq d$, then the quiver $A^{(1)}_{e}$ is a subdivision of $A^{(1)}_{d}$.
(Compare with \autoref{Ex:EasyQuiverVary}.)

\item The quiver $C_{\N}$ is a subdivision of itself since we only
subdivide simply laced edges.

\item If $e\geq d$, then the quiver $C^{(1)}_{e}$ is a subdivision of $C^{(1)}_{d}$.

\item Conversely, any subdivision of the quivers in (a)--(d) is a quiver of the same kind.
\end{enumerate}
\end{Examples}

\begin{Notation}
As a general rule, we place a bar above all of the associated Cartan data for the quiver
$\Sub{\Gamma}$. For example,
$\Sub{Q}^{+}=\bigoplus_{i\in\Sub{I}}\N\Sub{\alpha}_{i}$, and so on.
\end{Notation}

For $\Gamma$ we fix a subdivision $\Sub{\Gamma}$.
Subdivision determines two injective maps, denoted by the same symbol, $\Submap\map{I}{\Sub{I}}$
and $\Submap\map{E}{\Sub{E}}$, such
that if $r\colon i\to j$ is in $E$, then
\begin{gather*}
\begin{tikzcd}[ampersand replacement=\&,column sep=1.5cm]
\Submap(i)=i_{0}\ar[r,"\Submap(r)=r_{0}"]
\&
i_{1}\ar[r,"r_{1}"]
\&
\dots\ar[r,"r_{k-2}"]
\&
i_{k-1}\ar[r,"r_{k-1}"]
\&
i_{k}
=
\Submap(j)
\end{tikzcd}
\end{gather*}
is the subdivided edge in $\Sub{E}$. In particular, the edge $r\colon i\to j$ in $\Gamma$ is replaced with $k$ edges in $\Sub{\Gamma}$.

\begin{Definition}
Let $\beta=\sum_{i\in I_{\Gamma}}b_{i}\alpha_{i}\in Q^{+}_{\Gamma}$.
Define
\begin{gather*}
\Sub{\beta}=\sum_{i\in I_{\Gamma}}b_{i}\Sub{\alpha}_{S(i)}
\in Q^{+}_{\Sub{\Gamma}}.
\end{gather*}
\end{Definition}

\begin{Example}\label{Ex:SubdivQuiver}
Consider the quiver $\Gamma$ from \autoref{E:ExampleQuiver}, and the two subdivisions given by
\begin{gather*}
\Gamma\colon\quad
\begin{tikzpicture}[scale=1.2,anchorbase]
\node[circle,inner sep=1.8pt,fill=DarkBlue] (0) at (360/3*0:1){};
\node at (360/3*0:1.3){$0$};
\node[circle,inner sep=1.8pt,fill=DarkBlue] (1) at (360/3*1:1){};
\node at (360/3*1:1.3){$1$};
\node[circle,inner sep=1.8pt,fill=DarkBlue] (2) at (360/3*2:1){};
\node at (360/3*2:1.3){$2$};
\draw[->](0)--(1);
\draw[->](1)--(2);
\draw[->](2)--(0);
\end{tikzpicture}
,\quad
\Sub{\Gamma}\colon\quad
\begin{tikzpicture}[scale=1.2,anchorbase]
\node[circle,inner sep=1.8pt,fill=DarkBlue] (0) at (360/3*0:1){};
\node at (360/3*0:1.3){$0$};
\node[circle,inner sep=1.8pt,fill=DarkBlue] (1) at (360/3*1:1){};
\node at (360/3*1:1.3){$1$};
\node[circle,inner sep=1.8pt,fill=DarkBlue] (2) at (360/3*2:1){};
\node at (360/3*2:1.3){$2$};
\node[circle,inner sep=1.8pt,fill=DarkBlue] (4) at (360/3*0.5:0.5){};
\node at (360/3*0.5:0.8){$3$};
\draw[->](0)--(4);
\draw[->](4)--(1);
\draw[->](1)--(2);
\draw[->](2)--(0);
\end{tikzpicture}
,\quad
\Sub{\Gamma}^{\prime}\colon\quad
\begin{tikzpicture}[scale=1.2,anchorbase]
\node[circle,inner sep=1.8pt,fill=DarkBlue] (0) at (360/3*0:1){};
\node at (360/3*0:1.3){$0$};
\node[circle,inner sep=1.8pt,fill=DarkBlue] (1) at (360/3*1:1){};
\node at (360/3*1:1.3){$2$};
\node[circle,inner sep=1.8pt,fill=DarkBlue] (2) at (360/3*2:1){};
\node at (360/3*2:1.3){$3$};
\node[circle,inner sep=1.8pt,fill=DarkBlue] (4) at (360/3*0.5:0.5){};
\node at (360/3*0.5:0.8){$1$};
\draw[->](0)--(4);
\draw[->](4)--(1);
\draw[->](1)--(2);
\draw[->](2)--(0);
\end{tikzpicture}
.
\end{gather*}
The subdivision map $\Submap$ sends $i$ to $i$, so all of the original edges keep their name, and $\Sub{\beta}=\beta$ using the Kac--Moody data for $\Sub{\Gamma}$ rather than $\Gamma$.
For $\Sub{\Gamma}^{\prime}$, the map $\Submap'$ sends $0\mapsto 0$, $1\mapsto 2$, and $2\mapsto 3$, so that the edges of $\Gamma$, and the subscripts of $\beta$, change accordingly.
\end{Example}

To apply subdivision to weighted KLRW algebras we need to add a weighting to the subdivided quiver.

\begin{Definition}\label{D:WeithedSubdvisionedge}
Let $\Gamma$ be an admissible quiver and let $\bsig$ be a ghost shift for $\Gamma$.
Given a weighted simply laced edge
$\tikz[scale=0.9,centered,anchorbase]{
\draw[directed=1] (0,0) to node[above,yshift=-1pt]{$\sigma$}node[below]{$\phantom{\sigma}$} (0.925,0);
\node[circle,inner sep=1.8pt,fill=DarkBlue] at (0,0){};
\node[circle,inner sep=1.8pt,fill=DarkBlue] at (1,0){};
}$
in $\Gamma$, let
$\tikz[scale=0.9,centered,anchorbase]{
\draw[directed=1] (0,0) to node[above,yshift=-1pt]{$\sigma$}node[below]{$\phantom{\Sub{\sigma}}$} (0.425,0);
\draw[directed=1] (0.5,0) to node[above,yshift=-1pt]{$\Sub{\sigma}$}node[below]{$\phantom{\sigma}$} (0.925,0);
\node[circle,inner sep=1.8pt,fill=DarkBlue] at (0,0){};
\node[circle,inner sep=1.8pt,fill=DarkBlue] at (0.5,0){};
\node[circle,inner sep=1.8pt,fill=DarkBlue] at (1,0){};
}$
be weighted subdivision of this edge, where $\Sub{\sigma}\in\R_{\neq 0}$.
\end{Definition}

We define \emph{weighted (simply laced) subdivisions} of $\Gamma$
accordingly, and we use the terminology of subdivision
and weighted subdivision interchangeably.
We note also that we will make a choice of $\Gamma$ and a subdivision $\Sub{\Gamma}$, but omit these from the notation.


\subsection{Subdividing weighted KLRW diagrams}


In view of the paragraph after \autoref{D:Subdvision}, throughout this section we restrict to the case
where the quiver $\Sub{\Gamma}$ is obtained from $\Gamma$ by
subdividing a single edge $r\colon i\to j$ to $\Submap(i)\to\Sub{i}\to\Submap(j)$, for some $i,j\in I$.
Fix $\beta\in Q^{+}$ and extend $\Submap\map{I}{\Sub{I}}$ to a map $I^{\beta}\to\Sub{I}^{\Sub{\beta}}$ by replacing all occurrences of $i$ with $\Submap(i)$ and $\Sub{i}$, and leaving the others untouched.
Abusing notation, if $\bi\in I^{\beta}$ let $\Submap(\bi)$ be the resulting sequence in $\Sub{I}^{\Sub{\beta}}$.

Informally, the subdivision map on diagrams is given by ``fattening'' the
ghost $i$-strings using the rule
\begin{gather}\label{E:Colors}
\begin{tikzpicture}[scale=1.2,anchorbase,smallnodes,rounded corners]
\draw[ghost](1,0)--++(0,1)node[above,yshift=-1pt]{$i$};
\draw[solid](0,0)node[below]{$i$}--++(0,1)node[above,yshift=-1pt]{$\phantom{i}$};
\draw[redstring](0.5,0)node[below]{$\rho$}--++(0,1);
\end{tikzpicture}
\mapsto
\begin{tikzpicture}[scale=1.2,anchorbase,smallnodes,rounded corners]
\draw[ghost](4.0,0)--++(0,1)node[above,yshift=-1pt]{$s(i)$};
\draw[ghost,spinach](5.1,0)--++(0,1)node[above,yshift=-1pt,spinach]{$\Sub{i}$};
\draw[solid](3.0,0)node[below]{$s(i)$}--++(0,1);
\draw[solid,spinach](4.1,0)node[below,spinach]{$\Sub{i}$}--++(0,1);
\draw[redstring](3.5,0)node[below]{$\rho$}--++(0,1);
\draw[->] (4.1,-0.5)node[below]{added solid} to (4.1,-0.3);
\draw[->] (5.1,1.5)node[above]{added ghost} to (5.1,1.3);
\end{tikzpicture}
\end{gather}
and keeping all other strings as they are, except that the $j$-strings are relabeled as $\Submap(j)$-strings.

\begin{Notation}
For simplicity of notation, unless stated otherwise,
we assume from now on that $I$ is a subset of $\Sub{I}$,
$\Submap$ is the identity map on vertices and edges that are not subdivided, and subdividing replaces the edge $r\colon i\to j\in E$ with $r\colon i\to\Sub{i}$ and $\Sub{r}\colon\Sub{i}\to j$. In particular, the weighting of $r\colon i\to j\in E$ will be the weighting of $r\colon i\to\Sub{i}\in\Sub{E}$, while
$\Sub{r}$ has ghost shift $\Sub{\sigma}_{\Sub{r}}\in\R_{\neq 0}$.
As in \autoref{E:Colors}, we use (green) colors in diagrams to highlight the $\Sub{i}$-strings.
\end{Notation}

\begin{Example}
Consider the subdivision $\Sub{\Gamma}$ of $\Gamma$ as in \autoref{Ex:SubdivQuiver}. Then
\begin{gather*}
\begin{tikzpicture}[scale=1.2,anchorbase,smallnodes,rounded corners]
\draw[ghost](0,1)node[above,yshift=-1pt]{$1$}--++(-1,-1);
\draw[ghost](-1,1)node[above,yshift=-1pt]{$0$}--++(1,-1);
\draw[ghost](1.5,1)node[above,yshift=-1pt]{$2$}--++(1,-1);
\draw[ghost](2.5,1)node[above,yshift=-1pt]{$0$}--++(-1,-1);
\draw[solid](-0.5,1)--++(-1,-1)node[below]{$1$};
\draw[solid](-1.5,1)--++(1,-1)node[below]{$0$};
\draw[solid](1,1)node[above,yshift=-1pt]{$\phantom{i}$}--++(1,-1)node[below]{$2$};
\draw[solid](2,1)--++(-1,-1)node[below]{$0$};
\draw[redstring](1.25,1)--++(0,-1)node[below]{$\rho_{1}$};
\draw[redstring](3,1)--++(0,-1)node[below]{$\rho_{2}$};
\end{tikzpicture}
\mapsto
\begin{tikzpicture}[scale=1.2,anchorbase,smallnodes,rounded corners]
\draw[ghost](0,1)node[above,yshift=-1pt]{$1$}--++(-1,-1);
\draw[ghost](-1,1)node[above,yshift=-1pt]{$0$}--++(1,-1);
\draw[ghost,spinach](-0.25,1)node[above,yshift=-1pt,spinach]{$3$}--++(1,-1);
\draw[ghost](1.5,1)node[above,yshift=-1pt]{$2$}--++(1,-1);
\draw[ghost](2.5,1)node[above,yshift=-1pt]{$0$}--++(-1,-1);
\draw[ghost,spinach](3.25,1)node[above,yshift=-1pt,spinach]{$3$}--++(-1,-1);
\draw[solid](-0.5,1)--++(-1,-1)node[below]{$1$};
\draw[solid](-1.5,1)--++(1,-1)node[below]{$0$};
\draw[solid](1,1)node[above,yshift=-1pt]{$\phantom{i}$}--++(1,-1)node[below]{$2$};
\draw[solid](2,1)--++(-1,-1)node[below]{$0$};
\draw[solid,spinach](-0.75,1)--++(1,-1)node[below,spinach]{$3$};
\draw[solid,spinach](2.75,1)--++(-1,-1)node[below,spinach]{$3$};
\draw[redstring](1.25,1)--++(0,-1)node[below]{$\rho_{1}$};
\draw[redstring](3,1)--++(0,-1)node[below]{$\rho_{2}$};
\end{tikzpicture}
\end{gather*}
is an example of what the subdivision map, which we are about to define, actually does.
\end{Example}

As the distances between strings are crucial, we need to clearly specify the positions of all of the strings in the subdivided diagram. This makes the following formal definition of subdivision look more complicated than it actually is.

\begin{Definition}\label{D:DiagramSubdivision}
Let $D\in\Webab[\bx,\bi][\by,\bj]$, for $\bx,\by\in X$ and $\bi,\bj\in I^\beta$. Fix a \emph{translation factor} $t\in\R_{\neq 0}$.
The \emph{$t$-subdivision} of $D$ is the diagram $D_{t}$ obtained by adding new solid $\Sub{i}$-strings by translating each ghost $i$-string by $t$-units and then adding new ghost $\Sub{i}$-strings by translating the new solid $\Sub{i}$-strings by $\Sub{\sigma}_{\Sub{r}}$ units. All other strings and dots are untouched and no dots are added to the new strings.
\end{Definition}

We can always ensure that the diagram $D_{t}$ satisfies the conditions \autoref{Ex:Subdvision} by replacing $t$ with $t+\varepsilon$, for $\varepsilon$ sufficiently small. Hence, we can always ensure that $D_{t}$ is a diagram in the sense of \autoref{D:WebsterDiagram}.

Note that subdivision involves choices related to $\Gamma$, such as
ghost shifts $\bsig$, the choice of ghost shift $\Sub{\sigma}_{\,\Sub{i}}$ for $\Sub{\Gamma}$, and an additional choice of translation factor $t$.

The following is immediate:

\begin{Lemma}\label{L:Small}
Suppose that $0<|\varepsilon|\ll 1$.
Then $\Sub{D}_{t}$ and
$\Sub{D}_{t+\varepsilon}$ are the same up to isotopy and
conjugation by straight line diagrams.\qed
\end{Lemma}

Straight line diagrams and \autoref{L:Small} imply
that $\Sub{D}_{t}$ is essentially unique for small variations of $t$ and we call $\Sub{D}=\Sub{D}_{t}$ a \emph{subdivision} of $D$ if no confusion regarding $t$
can arise. Define
$\Sub{\balp}$ and $\Sub{\bbet}$ to be the endpoints of $\Sub{D}$ so that
$\Sub{D}\in\Webabs$.
Let
\begin{gather*}
\Sub{X}=\set{\Sub{\bx}|\bx\in X}
\quad\text{so that}\quad
\Sub{\Web}_{\Sub{\beta}}(\Sub{X})=\bigcup_{\bx,\by\in\Sub{X}}
\bigcup_{\bi,\bj\in I^{\beta}}\Webabs.
\end{gather*}
Note that $\Sub{\Web}_{\beta}(\Sub{X})$ is a set of (isotopy classes of) diagrams.
Moreover, by construction, subdivision gives a well-defined map:
\begin{gather}\label{E:subdivision}
\sg\map{\Web_{\beta}(X)}\Web_{\Sub{\beta}}(\Sub{X}),
D\mapsto\Sub{D}.
\end{gather}

\begin{Example}\label{Ex:DivisionMap}
In the setup of \autoref{Ex:Basis}, with the subdivision $\Sub{\Gamma}$ from \autoref{Ex:SubdivQuiver}, we have $\Sub{\beta}=\alpha_{0}+\alpha_{1}+\alpha_{3}$ and
\begin{gather*}
\sg(\WABasis)
\!=\!
\Bigg\{
\begin{tikzpicture}[scale=1.2,anchorbase,smallnodes,rounded corners]
\draw[ghost,dot=0.1](1,0)--++(0,1)node[above,yshift=-1pt]{$0$};
\draw[ghost,dot=0.1](1.75,0)--++(0,1)node[above,yshift=-1pt]{$1$};
\draw[ghost,spinach](1.35,0)--++(0,1)node[above,yshift=-1pt,spinach]{$3$};
\draw[solid,dot=0.1](0,0)node[below]{$0$}node[left,xshift=0.08cm,yshift=0.2cm]{$a$}--++(0,1);
\draw[solid,dot=0.1](0.75,0)node[below]{$1$}node[right,xshift=-0.08cm,yshift=0.2cm]{$b$}--++(0,1);
\draw[solid,spinach](1.2,0)node[below,spinach]{$3$}--++(0,1);
\draw[redstring](1.5,0)node[below]{$\rho$}--++(0,1);
\end{tikzpicture}
,\!\!
\begin{tikzpicture}[scale=1.2,anchorbase,smallnodes,rounded corners]
\draw[ghost,dot=0.1](1,0)--++(0.75,1)node[above,yshift=-1pt]{$0$};
\draw[ghost,dot=0.1](1.75,0)--++(-0.75,1)node[above,yshift=-1pt]{$1$};
\draw[ghost,spinach](1.35,0)--++(0.75,1)node[above,yshift=-1pt,spinach]{$3$};
\draw[solid,dot=0.1](0,0)node[below]{$0$}node[left,xshift=0.08cm,yshift=0.2cm]{$a$}--++(0.75,1);
\draw[solid,dot=0.1](0.75,0)node[below]{$1$}node[right,xshift=-0.08cm,yshift=0.2cm]{$b$}--++(-0.75,1);
\draw[solid,spinach](1.2,0)node[below,spinach]{$3$}--++(0.75,1);
\draw[redstring](1.5,0)node[below]{$\rho$}--++(0,1);
\end{tikzpicture}
,\!\!
\begin{tikzpicture}[scale=1.2,anchorbase,smallnodes,rounded corners]
\draw[ghost,dot=0.1](1,0)--++(0,1)node[above,yshift=-1pt]{$1$};
\draw[ghost,dot=0.1](1.75,0)--++(0,1)node[above,yshift=-1pt]{$0$};
\draw[ghost,spinach](2.1,0)--++(0,1)node[above,yshift=-1pt,spinach]{$3$};
\draw[solid,dot=0.1](0,0)node[below]{$1$}node[left,xshift=0.08cm,yshift=0.2cm]{$a$}--++(0,1);
\draw[solid,dot=0.1](0.75,0)node[below]{$0$}node[right,xshift=-0.08cm,yshift=0.2cm]{$b$}--++(0,1);
\draw[solid,spinach](1.95,0)node[below,spinach]{$3$}--++(0,1);
\draw[redstring](1.5,0)node[below]{$\rho$}--++(0,1);
\end{tikzpicture}
,\!\!
\begin{tikzpicture}[scale=1.2,anchorbase,smallnodes,rounded corners]
\draw[ghost,dot=0.1](1,0)--++(0.75,1)node[above,yshift=-1pt]{$1$};
\draw[ghost,dot=0.1](1.75,0)--++(-0.75,1)node[above,yshift=-1pt]{$0$};
\draw[ghost,spinach](2.1,0)--++(-0.75,1)node[above,yshift=-1pt,spinach]{$3$};
\draw[solid,dot=0.1](0,0)node[below]{$1$}node[left,xshift=0.08cm,yshift=0.2cm]{$a$}--++(0.75,1);
\draw[solid,dot=0.1](0.75,0)node[below]{$0$}node[right,xshift=-0.08cm,yshift=0.2cm]{$b$}--++(-0.75,1);
\draw[solid,spinach](1.95,0)node[below,spinach]{$3$}--++(-0.75,1);
\draw[redstring](1.5,0)node[below]{$\rho$}--++(0,1);
\end{tikzpicture}
\Bigg\vert a,b\in\N
\Bigg\}
,
\end{gather*}
for $0<t\ll 1$ and $0<\Sub{\sigma}_{\Sub{i}}\ll 1$.
In general, $\sg$ is injective but not surjective. For example, the
diagram
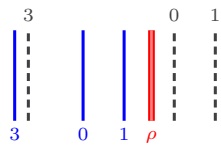
\begin{gather*}
\begin{tikzpicture}[scale=1.2,anchorbase,smallnodes,rounded corners]
\draw[ghost](0.15,0)--++(0,1)node[above,yshift=-1pt]{$3$};
\draw[ghost](1.75,0)--++(0,1)node[above,yshift=-1pt]{$0$};
\draw[ghost](2.2,0)--++(0,1)node[above,yshift=-1pt]{$1$};
\draw[solid](0,0)node[below]{$3$}--++(0,1);
\draw[solid](0.75,0)node[below]{$0$}--++(0,1);
\draw[solid](1.2,0)node[below]{$1$}--++(0,1);
\draw[redstring](1.5,0)node[below]{$\rho$}--++(0,1);
\end{tikzpicture}
\end{gather*}
does not belong to $\sg(\WABasis)$.
\end{Example}

The next definition ensures that the map $\sg$ preserves the degrees
of diagrams.

\begin{Definition}\label{D:admissible}
  A subdivision $(\Gamma,\Sub{\Gamma})$ is \emph{homogeneous} if
  $|t|\ll1$ and $|\Sub{\sigma}_{\Sub{r}}|\ll 1$ are sufficiently small.
\end{Definition}

\begin{Lemma}\label{L:DegreePreserving}
   Suppose that $(\Gamma,\Sub{\Gamma})$ is \emph{homogeneous}. Then
   $\deg D=\deg\Sub{D}$, for $D\in\Web_{\beta}(X)$.
\end{Lemma}

\begin{proof}
  It is enough to consider the case when subdivision sends $\Dldots{i,j}{1}$ to
  $\Dldots{i,\Sub{i},j}{2}$. We need to check that subdivision respects the degrees of
the diagrams in \autoref{SS:WebsterDegrees} that contain a
ghost $i$-string. The most interesting case is the $(i,j)$-ghost-ghost
crossing, where the result is immediate if $\<\alpha_{\,\Sub{i}},\alpha_{k}\>=0$. Moreover,
$\<\alpha_{\,\Sub{i}},\alpha_{k}\>\neq 0$ only if
$k\in\set{i,j}$. If $k=i$ then, locally,
\begin{gather*}
\deg
\begin{tikzpicture}[scale=1.2,anchorbase,smallnodes,rounded corners]
\draw[ghost](0,0)--++(1,1)node[above,yshift=-1pt]{$i$};
\draw[ghost](1,0)node[below]{$\phantom{i}$}--++(-1,1)node[above,yshift=-1pt]{$i$};
\end{tikzpicture}
=0
\mapsto
\deg
\begin{tikzpicture}[scale=1.2,anchorbase,smallnodes,rounded corners]
\draw[ghost](0,0)--++(1,1)node[above,yshift=-1pt]{$i$};
\draw[ghost](1,0)--++(-1,1)node[above,yshift=-1pt]{$i$};
\draw[ghost,spinach](0.4,0)--++(1,1)node[above,yshift=-1pt,spinach]{$\Sub{i}$};
\draw[ghost,spinach](1.4,0)--++(-1,1)node[above,yshift=-1pt,spinach]{$\Sub{i}$};
\draw[solid,spinach](0.25,0)node[below,spinach]{$\Sub{i}$}--++(1,1);
\draw[solid,spinach](1.25,0)node[below,spinach]{$\Sub{i}$}--++(-1,1);
\end{tikzpicture}
=
\deg
\begin{tikzpicture}[scale=1.2,anchorbase,smallnodes,rounded corners]
\draw[ghost](0,0)--++(1,1)node[above,yshift=-1pt]{$i$};
\draw[solid,spinach](1.25,0)node[below,spinach]{$\Sub{i}$}--++(-1,1);
\end{tikzpicture}
+
\deg
\begin{tikzpicture}[scale=1.2,anchorbase,smallnodes,rounded corners]
\draw[ghost](1,0)--++(-1,1)node[above,yshift=-1pt]{$i$};
\draw[solid,spinach](0.25,0)node[below,spinach]{$\Sub{i}$}--++(1,1);
\end{tikzpicture}
+
\deg
\begin{tikzpicture}[scale=1.2,anchorbase,smallnodes,rounded corners]
\draw[solid,spinach](0.25,0)node[below,spinach]{$\Sub{i}$}--++(1,1)node[above,yshift=-1pt]{$\phantom{i}$};
\draw[solid,spinach](1.25,0)node[below,spinach]{$\Sub{i}$}--++(-1,1);
\end{tikzpicture}
=
0
,
\end{gather*}
where we only illustrate nonzero contributions on the right-hand side.
If $k=j$, then, noting that
$(\Sub{i},j)$-solid-ghost crossings are of degree zero because of the choice of
orientation, we have locally
\begin{gather*}
\deg
\begin{tikzpicture}[scale=1.2,anchorbase,smallnodes,rounded corners]
\draw[ghost](0,0)node[below]{$\phantom{i}$}--++(1,1)node[above,yshift=-1pt]{$i$};
\draw[ghost](1,0)--++(-1,1)node[above,yshift=-1pt]{$j$};
\end{tikzpicture}
=
0
\mapsto
\deg
\begin{tikzpicture}[scale=1.2,anchorbase,smallnodes,rounded corners]
\draw[ghost](0,0)--++(1,1)node[above,yshift=-1pt]{$i$};
\draw[ghost](1,0)--++(-1,1)node[above,yshift=-1pt]{$j$};
\draw[ghost,spinach](0.4,0)--++(1,1)node[above,yshift=-1pt,spinach]{$\Sub{i}$};
\draw[solid,spinach](0.25,0)node[below,spinach]{$\Sub{i}$}--++(1,1);
\end{tikzpicture}
=
0.
\end{gather*}
Another crucial case, which requires the assumption that $0<|\Sub{\sigma}_{\Sub{r}}|\ll 1$, is
\begin{gather*}
\deg
\begin{tikzpicture}[scale=1.2,anchorbase,smallnodes,rounded corners]
\draw[ghost](0,0)--++(1,1)node[above,yshift=-1pt]{$i$};
\draw[solid](1,0)node[below]{$j$}--++(-1,1);
\end{tikzpicture}
=
1
\mapsto
\deg
\begin{tikzpicture}[scale=1.2,anchorbase,smallnodes,rounded corners]
\draw[ghost](0,0)--++(1,1)node[above,yshift=-1pt]{$i$};
\draw[ghost,spinach](0.5,0)--++(1,1)node[above,yshift=-1pt,spinach]{$\Sub{i}$};
\draw[solid](1,0)node[below]{$j$}--++(-1,1);
\draw[solid,spinach](0.25,0)node[below,spinach]{$\Sub{i}$}--++(1,1);
\end{tikzpicture}
=
\deg
\begin{tikzpicture}[scale=1.2,anchorbase,smallnodes,rounded corners]
\draw[solid](1,0)node[below,yshift=-1pt]{$j$}--++(-1,1);
\draw[ghost,spinach](0.5,0)--++(1,1)node[above,spinach,yshift=-1pt]{$\Sub{i}$};
\end{tikzpicture}
=1,
\end{gather*}
where we use the fact that the $(i,j)$-ghost-solid and the
$(\Sub{i},j)$-ghost-solid crossings are of degrees $0$ and $1$, respectively, after subdivision. Checking that degrees are preserved in the remaining
diagrams is similar.
\end{proof}

It is easy to construct examples where \autoref{L:DegreePreserving} fails by taking $|t|$ or $|\Sub{\sigma}_{\Sub{i}}|$ sufficiently large.

We assume from now on that $(\Gamma,\Sub{\Gamma})$ is homogeneous. Before we can identify the image of $\sg$, we need more definitions.

\begin{Definition}
Let $\Sub{I}_{\text{bad}}=\set{\bi\in\Sub{I}^{n}|i_{r}=\Sub{i}=i_{r+1} \text{ for some }1\leq r<n}$, be the set of \emph{bad residue sequences}. A \emph{bad idempotent diagram} is an idempotent $\1_{\Sub{\bx},\Sub{\bi}}$, for some $\Sub{\bi}\in\Sub{I}_{\text{bad}}$ and $\Sub{\bx}\in\Sub{X}$. A \emph{bad diagram} is a diagram that factors through a bad idempotent diagram. Define
\begin{gather*}
\1_{\text{bad}}=\sum_{\Sub{\bi}\in\Sub{I}_{\text{bad}},\Sub{\bx}\in\Sub{X}}\1_{\Sub{\bx},\Sub{\bi}},
\end{gather*}
to be the sum of all bad idempotent diagrams.
\end{Definition}

To detect bad diagrams, define
a \emph{horizontal cut} through a diagram $D$ to be a
horizontal line $H(D)=\set[\big]{(x,a)|x\in\R}\subset\R\times[0,1]$,
for $a\in(0,1)$, such that the intersections $H(D)\cap D$ are locally of the form
\begin{gather*}
\begin{tikzpicture}[scale=1.2,anchorbase,smallnodes,rounded corners]
\draw[solid] (0,0)node[below]{$i$} to (0,1)node[above,yshift=-1pt]{$\phantom{i}$};
\draw[densely dotted,thick] (-0.5,0.5)node[left]{$H(D)$} to (0.5,0.5);
\end{tikzpicture}
,\quad
\begin{tikzpicture}[scale=1.2,anchorbase,smallnodes,rounded corners]
\draw[ghost] (0,0)node[below]{$\phantom{i}$} to (0,1)node[above,yshift=-1pt]{$i$};
\draw[densely dotted,thick] (-0.5,0.5)node[left]{$H(D)$} to (0.5,0.5);
\end{tikzpicture}
,\quad
\begin{tikzpicture}[scale=1.2,anchorbase,smallnodes,rounded corners]
\draw[redstring] (0,0)node[below]{$\rho_{i}$} to (0,1)node[above,yshift=-1pt]{$\phantom{i}$};
\draw[densely dotted,thick] (-0.5,0.5)node[left]{$H(D)$} to (0.5,0.5);
\end{tikzpicture}
.
\end{gather*}
In other words, a horizontal cut intersects $D$ generically, in the sense that it avoids crossings and dots.
By definition, a diagram is bad
if there is some neighborhood of a horizontal cut that is locally of the form
\begin{gather}\label{E:LocallyBad}
\begin{tikzpicture}[scale=1.2,anchorbase,smallnodes,rounded corners]
\draw[solid,spinach](0,0)node[below,spinach]{$\Sub{i}$}--++(0,1)node[above,yshift=-1pt]{$\phantom{\Sub{i}}$};
\draw[solid,spinach](0.5,0)node[below,spinach]{$\Sub{i}$}--++(0,1);
\draw[densely dotted,thick] (-0.5,0.5)node[left]{$H(D)$} to (1,0.5);
\end{tikzpicture}
,\quad
\begin{tikzpicture}[scale=1.2,anchorbase,smallnodes,rounded corners]
\draw[ghost,spinach](0,0)node[below,spinach]{$\phantom{\Sub{i}}$}--++(0,1)node[above,yshift=-1pt,spinach]{$\Sub{i}$};
\draw[ghost,spinach](0.5,0)--++(0,1)node[above,yshift=-1pt,spinach]{$\Sub{i}$};
\draw[densely dotted,thick] (-0.5,0.5)node[left]{$H(D)$} to (1,0.5);
\end{tikzpicture}
.
\end{gather}

\begin{Definition}\label{D:SubWebster}
If $\Sub{\Gamma}$ is a subdivision of $\Gamma$, then define
$\Oneg=\sum_{\bx\in X}\sum_{\bi\in I^\beta}\1_{\Sub{\bx},\Sub{\bi}}$ and
\begin{gather}\label{E:IdemTrunc}
  \SAX=\Oneg\WAs(\Sub{X})\Oneg/
\Oneg\WAs(\Sub{X})\1_{\text{bad}}\WAs(\Sub{X})\Oneg.
\end{gather}
\end{Definition}

We identify a diagram in $\WA(\Sub{X})$ with its image in $\SAX$.

\begin{Lemma}\label{L:KillDiagrams}
  Let $B\in\SAX$ be a bad diagram. Then $B=0$.
\end{Lemma}

\begin{proof}
Every bad diagram $B$ has a horizontal cut of the form of \autoref{E:LocallyBad} and therefore factors through a bad idempotent diagram in the sense that $B=D^{\prime}\1_{bad}D^{\prime\prime}$, for some diagrams
$D^{\prime}$ and $D^{\prime\prime}$.
The claim then follows by the definition of $\SAX$.
\end{proof}


\subsection{The subdivision isomorphism}


\begin{Notation}
From now on we assume that
$Q_{ij}(u,v)=au+bv$, for units $a,b\in R$ such that $a=-b$, is the polynomial associated to the edge \Ddots{1} that we are subdividing.
Moreover, we assume that the two polynomials $Q_{i,\Sub{i}}(u,v)$ and $Q_{\Sub{i},j}(u,v)$ for the edges \Ddots{2} are both equal to $Q_{ij}(u,v)$.
\end{Notation}

The following theorem extends the subdivision map $\sg$ on diagrams from \autoref{E:subdivision} to an algebra map.

\begin{Theorem}\label{T:SubDiv}
Suppose that $\Sub{\Gamma}$ is a homogeneous
subdivision of $\Gamma$.
Then there is an isomorphism of graded algebras $\sg\map{\WA(X)}{\SAX}$.
\end{Theorem}

\begin{proof}
By construction, if
$D\in\Web_{\beta}(X)$, then $\sg(D)\in\SAX$. So, by
\autoref{P:WABasis}, we can view the map $\sg\map{\WA(X)}{\SAX}$ as a morphism of $R$-modules. Moreover, by \autoref{L:DegreePreserving},
$\sg\map{\WA(X)}{\SAX}$ is homogeneous,
so it remains to show that $\sg$ is a bijective morphism of algebras.

We first show that $\sg$ is an algebra homomorphism by checking the
relations from \autoref{D:RationalCherednik}. We only need to
consider those relations that involve ghost $i$-strings because all other relations are unchanged by subdivision.

We consider each of the relations (a)--(c) from
\autoref{D:RationalCherednik} by checking some of the exceptional
relations \autoref{R:DotCrossing}--\autoref{R:BraidSRS} together
with some of the non-exceptional relations. Moreover, we will only need the condition from \autoref{D:admissible} that $\Sub{\Gamma}$ is a homogeneous subdivision for a few relations, with all other relations being preserved without this condition.

\TestRelation{a}
We start by checking \autoref{R:DotCrossing}. As we will see, this
is one of the relations that forces the appearance of $\Oneg$ in the
definition of $\SAX$, which in turn invokes \autoref{L:KillDiagrams}. Applying $\sg$ to the left-hand side of
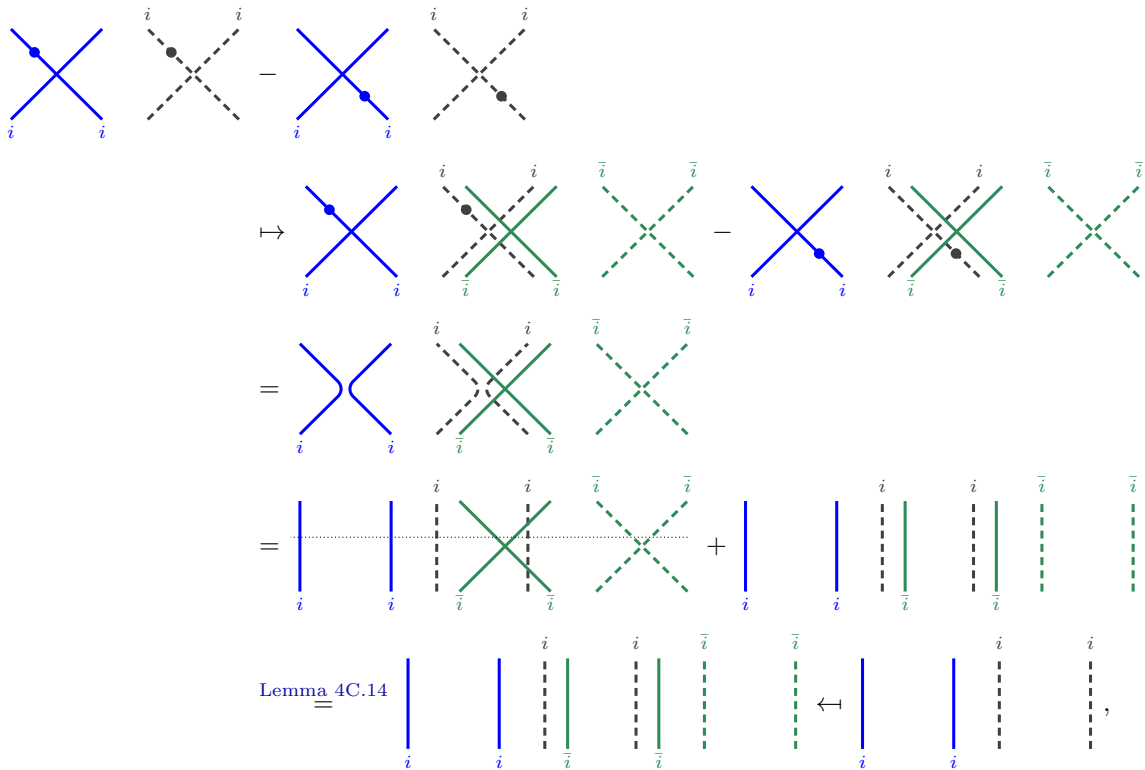
\autoref{R:DotCrossing} gives
\begin{align*}
\begin{tikzpicture}[scale=1.2,anchorbase,smallnodes,rounded corners]
\draw[ghost](1.5,0)--++(1,1)node[above,yshift=-1pt]{$i$};
\draw[ghost,dot=0.75](2.5,0)--++(-1,1)node[above,yshift=-1pt]{$i$};
\draw[solid](0,0)node[below]{$i$}--++(1,1);
\draw[solid,dot=0.75](1,0)node[below]{$i$}--++(-1,1);
\end{tikzpicture}
&-
\begin{tikzpicture}[scale=1.2,anchorbase,smallnodes,rounded corners]
\draw[ghost](1.5,0)--++(1,1)node[above,yshift=-1pt]{$i$};
\draw[ghost,dot=0.25](2.5,0)--++(-1,1)node[above,yshift=-1pt]{$i$};
\draw[solid](0,0)node[below]{$i$}--++(1,1);
\draw[solid,dot=0.25](1,0)node[below]{$i$}--++(-1,1);
\end{tikzpicture}
\\
&\mapsto
\begin{tikzpicture}[scale=1.2,anchorbase,smallnodes,rounded corners]
\draw[ghost](1.5,0)--++(1,1)node[above,yshift=-1pt]{$i$};
\draw[ghost,dot=0.75](2.5,0)--++(-1,1)node[above,yshift=-1pt]{$i$};
\draw[ghost,spinach](3.25,0)--++(1,1)node[above,yshift=-1pt,spinach]{$\Sub{i}$};
\draw[ghost,spinach](4.25,0)--++(-1,1)node[above,yshift=-1pt,spinach]{$\Sub{i}$};
\draw[solid](0,0)node[below]{$i$}--++(1,1);
\draw[solid,dot=0.75](1,0)node[below]{$i$}--++(-1,1);
\draw[solid,spinach](1.75,0)node[below,spinach]{$\Sub{i}$}--++(1,1);
\draw[solid,spinach](2.75,0)node[below,spinach]{$\Sub{i}$}--++(-1,1);
\end{tikzpicture}
-
\begin{tikzpicture}[scale=1.2,anchorbase,smallnodes,rounded corners]
\draw[ghost](1.5,0)--++(1,1)node[above,yshift=-1pt]{$i$};
\draw[ghost,dot=0.25](2.5,0)--++(-1,1)node[above,yshift=-1pt]{$i$};
\draw[ghost,spinach](3.25,0)--++(1,1)node[above,yshift=-1pt,spinach]{$\Sub{i}$};
\draw[ghost,spinach](4.25,0)--++(-1,1)node[above,yshift=-1pt,spinach]{$\Sub{i}$};
\draw[solid](0,0)node[below]{$i$}--++(1,1);
\draw[solid,dot=0.25](1,0)node[below]{$i$}--++(-1,1);
\draw[solid,spinach](1.75,0)node[below,spinach]{$\Sub{i}$}--++(1,1);
\draw[solid,spinach](2.75,0)node[below,spinach]{$\Sub{i}$}--++(-1,1);
\end{tikzpicture}
\\&=
\begin{tikzpicture}[scale=1.2,anchorbase,smallnodes,rounded corners]
\draw[ghost](1.5,0)--++(0.5,0.5)--++(-0.5,0.5)node[above,yshift=-1pt]{$i$};
\draw[ghost](2.5,0)--++(-0.5,0.5)--++(0.5,0.5)node[above,yshift=-1pt]{$i$};
\draw[ghost,spinach](3.25,0)--++(1,1)node[above,yshift=-1pt,spinach]{$\Sub{i}$};
\draw[ghost,spinach](4.25,0)--++(-1,1)node[above,yshift=-1pt,spinach]{$\Sub{i}$};
\draw[solid](0,0)node[below]{$i$}--++(0.5,0.5)--++(-0.5,0.5);
\draw[solid](1,0)node[below]{$i$}--++(-0.5,0.5)--++(0.5,0.5);
\draw[solid,spinach](1.75,0)node[below,spinach]{$\Sub{i}$}--++(1,1);
\draw[solid,spinach](2.75,0)node[below,spinach]{$\Sub{i}$}--++(-1,1);
\end{tikzpicture}
\\
&=
\begin{tikzpicture}[scale=1.2,anchorbase,smallnodes,rounded corners]
\draw[ghost](1.5,0)--++(0,1)node[above,yshift=-1pt]{$i$};
\draw[ghost](2.5,0)--++(0,1)node[above,yshift=-1pt]{$i$};
\draw[ghost,spinach](3.25,0)--++(1,1)node[above,yshift=-1pt,spinach]{$\Sub{i}$};
\draw[ghost,spinach](4.25,0)--++(-1,1)node[above,yshift=-1pt,spinach]{$\Sub{i}$};
\draw[solid](0,0)node[below]{$i$}--++(0,1);
\draw[solid](1,0)node[below]{$i$}--++(0,1);
\draw[solid,spinach](1.75,0)node[below,spinach]{$\Sub{i}$}--++(1,1);
\draw[solid,spinach](2.75,0)node[below,spinach]{$\Sub{i}$}--++(-1,1);
\draw[densely dotted] (4.25,0.6)--++(-4.35,0);
\end{tikzpicture}
+
\begin{tikzpicture}[scale=1.2,anchorbase,smallnodes,rounded corners]
\draw[ghost](1.5,0)--++(0,1)node[above,yshift=-1pt]{$i$};
\draw[ghost](2.5,0)--++(0,1)node[above,yshift=-1pt]{$i$};
\draw[ghost,spinach](3.25,0)--++(0,1)node[above,yshift=-1pt,spinach]{$\Sub{i}$};
\draw[ghost,spinach](4.25,0)--++(0,1)node[above,yshift=-1pt,spinach]{$\Sub{i}$};
\draw[solid](0,0)node[below]{$i$}--++(0,1);
\draw[solid](1,0)node[below]{$i$}--++(0,1);
\draw[solid,spinach](1.75,0)node[below,spinach]{$\Sub{i}$}--++(0,1);
\draw[solid,spinach](2.75,0)node[below,spinach]{$\Sub{i}$}--++(0,1);
\end{tikzpicture}
\\
&\stackrel{\text{\autoref{L:KillDiagrams}}}{=}
\begin{tikzpicture}[scale=1.2,anchorbase,smallnodes,rounded corners]
\draw[ghost](1.5,0)--++(0,1)node[above,yshift=-1pt]{$i$};
\draw[ghost](2.5,0)--++(0,1)node[above,yshift=-1pt]{$i$};
\draw[ghost,spinach](3.25,0)--++(0,1)node[above,yshift=-1pt,spinach]{$\Sub{i}$};
\draw[ghost,spinach](4.25,0)--++(0,1)node[above,yshift=-1pt,spinach]{$\Sub{i}$};
\draw[solid](0,0)node[below]{$i$}--++(0,1);
\draw[solid](1,0)node[below]{$i$}--++(0,1);
\draw[solid,spinach](1.75,0)node[below,spinach]{$\Sub{i}$}--++(0,1);
\draw[solid,spinach](2.75,0)node[below,spinach]{$\Sub{i}$}--++(0,1);
\end{tikzpicture}
\mapsfrom
\begin{tikzpicture}[scale=1.2,anchorbase,smallnodes,rounded corners]
\draw[ghost](1.5,0)--++(0,1)node[above,yshift=-1pt]{$i$};
\draw[ghost](2.5,0)--++(0,1)node[above,yshift=-1pt]{$i$};
\draw[solid](0,0)node[below]{$i$}--++(0,1);
\draw[solid](1,0)node[below]{$i$}--++(0,1);
\end{tikzpicture}
,
\end{align*}
where the second equality uses \autoref{R:BraidGSG}, and the fact that $Q_{iji}(u,v,w)=1$, and the last equality follows by drawing a horizontal cut just above the middle of
the left-hand diagram, as indicated above. Hence, by the symmetry of \autoref{R:DotCrossing}, $\sg$ respects \autoref{R:DotCrossing} as required.

Note that it is not necessary to check what happens to
the solid and ghost strings outside
of the ``central'' region of a
diagram that contains the ghost $i$-strings because the behavior of these strings is controlled by the strings in the central region.
Consequently, we will omit these strings from the diagrams when
checking the remaining relations and focus on the part of the diagram that contains the ghost $i$-strings.

Continuing with our check of the
relations in (a), if $k\neq i$, then
\begin{gather*}
\begin{tikzpicture}[scale=1.2,anchorbase,smallnodes,rounded corners]
\draw[ghost](0,0)--++(1,1)node[above,yshift=-1pt]{$i$};
\draw[solid,dot=0.75](1,0)node[below]{$k$}--++(-1,1);
\end{tikzpicture}
\mapsto
\begin{tikzpicture}[scale=1.2,anchorbase,smallnodes,rounded corners]
\draw[ghost](0,0)--++(1,1)node[above,yshift=-1pt]{$i$};
\draw[solid,dot=0.75](1,0)node[below]{$k$}--++(-1,1);
\draw[solid,spinach](0.25,0)node[below,spinach]{$\Sub{i}$}--++(1,1);
\end{tikzpicture}
=
\begin{tikzpicture}[scale=1.2,anchorbase,smallnodes,rounded corners]
\draw[ghost](0,0)--++(1,1)node[above,yshift=-1pt]{$i$};
\draw[solid,dot=0.25](1,0)node[below]{$k$}--++(-1,1);
\draw[solid,spinach](0.25,0)node[below,spinach]{$\Sub{i}$}--++(1,1);
\end{tikzpicture}
\mapsfrom
\begin{tikzpicture}[scale=1.2,anchorbase,smallnodes,rounded corners]
\draw[ghost](0,0)--++(1,1)node[above,yshift=-1pt]{$i$};
\draw[solid,dot=0.25](1,0)node[below]{$k$}--++(-1,1);
\end{tikzpicture},
\end{gather*}
since the dot slides freely in this case.
All other dot sliding rules can be proven {\muta}.

\TestRelation{b}
Recall that we are assuming that $Q_{ij}(u,v)=au+bv$, where $a=-b$ is a
unit.  Consider the double crossing relation \autoref{R:SolidSolid}, where we do not
draw the solid-solid crossings:
\begin{gather*}
\begin{tikzpicture}[scale=1.2,anchorbase,smallnodes,rounded corners]
\draw[ghost](0,0)--++(0.5,0.5)--++(-0.5,0.5)node[above,yshift=-1pt]{$i$};
\draw[ghost](0.5,0)node[below]{$\phantom{i}$}--++(-0.5,0.5)--++(0.5,0.5)node[above,yshift=-1pt]{$i$};
\end{tikzpicture}
\mapsto
\begin{tikzpicture}[scale=1.2,anchorbase,smallnodes,rounded corners]
\draw[ghost](0,0)node[below]{$i$}--++(0.5,0.5)--++(-0.5,0.5)node[above,yshift=-1pt]{$i$};
\draw[ghost](0.5,0)node[below]{$i$}--++(-0.5,0.5)--++(0.5,0.5)node[above,yshift=-1pt]{$i$};
\draw[solid,spinach](0.2,0)node[below,spinach]{$\Sub{i}$}--++(0.5,0.5)--++(-0.5,0.5);
\draw[solid,spinach](0.7,0)node[below,spinach]{$\Sub{i}$}--++(-0.5,0.5)--++(0.5,0.5);
\end{tikzpicture}
=
a
\begin{tikzpicture}[scale=1.2,anchorbase,smallnodes,rounded corners]
\draw[ghost,dot](0,0)node[below]{$i$}--++(0.25,0.5)--++(-0.25,0.5)node[above,yshift=-1pt]{$i$};
\draw[ghost](0.5,0)node[below]{$i$}--++(-0.5,0.5)--++(0.5,0.5)node[above,yshift=-1pt]{$i$};
\draw[solid,spinach](0.2,0)node[below,spinach]{$\Sub{i}$}--++(0.5,0.5)--++(-0.5,0.5);
\draw[solid,spinach](0.7,0)node[below,spinach]{$\Sub{i}$}--++(-0.25,0.5)--++(0.25,0.5);
\end{tikzpicture}
+b
\begin{tikzpicture}[scale=1.2,anchorbase,smallnodes,rounded corners]
\draw[ghost](0,0)node[below]{$i$}--++(0.25,0.5)--++(-0.25,0.5)node[above,yshift=-1pt]{$i$};
\draw[ghost](0.5,0)node[below]{$i$}--++(-0.5,0.5)--++(0.5,0.5)node[above,yshift=-1pt]{$i$};
\draw[solid,spinach](0.2,0)node[below,spinach]{$\Sub{i}$}--++(0.5,0.5)--++(-0.5,0.5);
\draw[solid,spinach,dot](0.7,0)node[below,spinach]{$\Sub{i}$}--++(-0.25,0.5)--++(0.25,0.5);
\end{tikzpicture}
=0,
\end{gather*}
as required. It is automatic that $\sg$ respects relation
\autoref{R:RedSolid} since there are no ghost strings. For \autoref{R:GhostSolid} observe that, because we only
ever subdivide simply laced edges, it is enough to consider the cases when there is an edge $i\to k$, in which case in $\WA(X)$ we have:
\begin{gather*}
\begin{tikzpicture}[scale=1.2,anchorbase,smallnodes,rounded corners]
\draw[ghost](0,0)--++(0.5,0.5)--++(-0.5,0.5)node[above,yshift=-1pt]{$i$};
\draw[solid](0.5,0)node[below]{$k$}--++(-0.5,0.5)--++(0.5,0.5)node[above,yshift=-1pt]{$\phantom{i}$};
\end{tikzpicture}
=
a
\begin{tikzpicture}[scale=1.2,anchorbase,smallnodes,rounded corners]
\draw[ghost,dot](0,0)--++(0,1)node[above,yshift=-1pt]{$i$};
\draw[solid](0.5,0)node[below]{$k$}--++(0,1)node[above,yshift=-1pt]{$\phantom{i}$};
\end{tikzpicture}
+b
\begin{tikzpicture}[scale=1.2,anchorbase,smallnodes,rounded corners]
\draw[ghost](0,0)--++(0,1)node[above,yshift=-1pt]{$i$};
\draw[solid,dot](0.5,0)node[below]{$k$}--++(0,1)node[above,yshift=-1pt]{$\phantom{i}$};
\end{tikzpicture}
,\quad\text{and}\quad
\begin{tikzpicture}[scale=1.2,anchorbase,smallnodes,rounded corners]
\draw[ghost](0.5,0)--++(-0.5,0.5)--++(0.5,0.5)node[above,yshift=-1pt]{$i$};
\draw[solid](0,0)node[below]{$k$}--++(0.5,0.5)--++(-0.5,0.5)node[above,yshift=-1pt]{$\phantom{i}$};
\end{tikzpicture}
=
a
\begin{tikzpicture}[scale=1.2,anchorbase,smallnodes,rounded corners]
\draw[ghost,dot](0.5,0)--++(0,1)node[above,yshift=-1pt]{$i$};
\draw[solid](0,0)node[below]{$k$}--++(0,1)node[above,yshift=-1pt]{$\phantom{i}$};
\end{tikzpicture}
+b
\begin{tikzpicture}[scale=1.2,anchorbase,smallnodes,rounded corners]
\draw[ghost](0.5,0)--++(0,1)node[above,yshift=-1pt]{$i$};
\draw[solid,dot](0,0)node[below]{$k$}--++(0,1)node[above,yshift=-1pt]{$\phantom{i}$};
\end{tikzpicture}
.
\end{gather*}
The most
interesting case for the left-hand relation is when $k=j$
(the other cases are automatic since we subdivide the edge $i\to j$), where we need $0<|\Sub{\sigma}_{\,\Sub{r}}|\ll 1$:
\begin{gather*}
\begin{tikzpicture}[scale=1.2,anchorbase,smallnodes,rounded corners]
\draw[ghost](0,0)--++(0.5,0.5)--++(-0.5,0.5)node[above,yshift=-1pt]{$i$};
\draw[solid](0.5,0)node[below]{$j$}--++(-0.5,0.5)--++(0.5,0.5);
\end{tikzpicture}
\mapsto
\begin{tikzpicture}[scale=1.2,anchorbase,smallnodes,rounded corners]
\draw[ghost](0,0)--++(0.5,0.5)--++(-0.5,0.5)node[above,yshift=-1pt]{$i$};
\draw[ghost,spinach](0.4,0)--++(0.5,0.5)--++(-0.5,0.5)node[above,yshift=-1pt,spinach]{$\Sub{i}$};
\draw[solid](0.5,0)node[below,xshift=0.2cm]{$j$}--++(-0.5,0.5)--++(0.5,0.5);
\draw[solid,spinach](0.25,0)node[below,spinach]{$\Sub{i}$}--++(0.5,0.5)--++(-0.5,0.5);
\end{tikzpicture}
=
\begin{tikzpicture}[scale=1.2,anchorbase,smallnodes,rounded corners]
\draw[ghost](0,0)--++(0,1)node[above,yshift=-1pt]{$i$};
\draw[ghost,spinach](0.4,0)--++(0.25,0.5)--++(-0.25,0.5)node[above,yshift=-1pt,spinach]{$\Sub{i}$};
\draw[solid](0.5,0)node[below,xshift=0.2cm]{$j$}--++(-0.25,0.5)--++(0.25,0.5);
\draw[solid,spinach](0.25,0)node[below,spinach]{$\Sub{i}$}--++(0.25,0.5)--++(-0.25,0.5);
\end{tikzpicture}
=
a
\begin{tikzpicture}[scale=1.2,anchorbase,smallnodes,rounded corners]
\draw[ghost](0,0)--++(0,1)node[above,yshift=-1pt]{$i$};
\draw[ghost,spinach,dot](0.4,0)--++(0,1)node[above,yshift=-1pt,spinach]{$\Sub{i}$};
\draw[solid](0.5,0)node[below,xshift=0.2cm]{$j$}--++(0,1);
\draw[solid,spinach,dot](0.25,0)node[below,spinach]{$\Sub{i}$}--++(0,1);
\end{tikzpicture}
+b
\begin{tikzpicture}[scale=1.2,anchorbase,smallnodes,rounded corners]
\draw[ghost](0,0)--++(0,1)node[above,yshift=-1pt]{$i$};
\draw[ghost,spinach](0.4,0)--++(0,1)node[above,yshift=-1pt,spinach]{$\Sub{i}$};
\draw[solid,dot](0.5,0)node[below,xshift=0.2cm]{$j$}--++(0,1);
\draw[solid,spinach](0.25,0)node[below,spinach]{$\Sub{i}$}--++(0,1);
\end{tikzpicture}
=
-
\begin{tikzpicture}[scale=1.2,anchorbase,smallnodes,rounded corners]
\draw[ghost](0,0)node[below]{$i$}--++(0.4,0.5)--++(-0.4,0.5)node[above,yshift=-1pt]{$i$};
\draw[ghost,spinach](0.4,0)--++(-0.18,0.5)--++(0.18,0.5)node[above,yshift=-1pt,spinach]{$\Sub{i}$};
\draw[solid](0.5,0)node[below,xshift=0.2cm]{$j$}--++(0,1);
\draw[solid,spinach](0.25,0)node[below,spinach]{$\Sub{i}$}--++(-0.18,0.5)--++(0.18,0.5);
\draw[densely dotted] (0.7,0.5)--++(-0.7,0);
\end{tikzpicture}
+
a
\begin{tikzpicture}[scale=1.2,anchorbase,smallnodes,rounded corners]
\draw[ghost,dot](0,0)--++(0,1)node[above,yshift=-1pt]{$i$};
\draw[ghost,spinach](0.4,0)--++(0,1)node[above,yshift=-1pt,spinach]{$\Sub{i}$};
\draw[solid](0.5,0)node[below,xshift=0.2cm]{$j$}--++(0,1);
\draw[solid,spinach](0.25,0)node[below,spinach]{$\Sub{i}$}--++(0,1);
\end{tikzpicture}
+b
\begin{tikzpicture}[scale=1.2,anchorbase,smallnodes,rounded corners]
\draw[ghost](0,0)--++(0,1)node[above,yshift=-1pt]{$i$};
\draw[ghost,spinach](0.4,0)--++(0,1)node[above,yshift=-1pt,spinach]{$\Sub{i}$};
\draw[solid,dot](0.5,0)node[below,xshift=0.2cm]{$j$}--++(0,1);
\draw[solid,spinach](0.25,0)node[below,spinach]{$\Sub{i}$}--++(0,1);
\end{tikzpicture}
\\
\stackrel{\text{\autoref{L:KillDiagrams}}}{=}
a
\begin{tikzpicture}[scale=1.2,anchorbase,smallnodes,rounded corners]
\draw[ghost,dot](0,0)--++(0,1)node[above,yshift=-1pt]{$i$};
\draw[ghost,spinach](0.4,0)--++(0,1)node[above,yshift=-1pt,spinach]{$\Sub{i}$};
\draw[solid](0.5,0)node[below,xshift=0.2cm]{$j$}--++(0,1);
\draw[solid,spinach](0.25,0)node[below,spinach]{$\Sub{i}$}--++(0,1);
\end{tikzpicture}
+b
\begin{tikzpicture}[scale=1.2,anchorbase,smallnodes,rounded corners]
\draw[ghost](0,0)--++(0,1)node[above,yshift=-1pt]{$i$};
\draw[ghost,spinach](0.4,0)--++(0,1)node[above,yshift=-1pt,spinach]{$\Sub{i}$};
\draw[solid,dot](0.5,0)node[below,xshift=0.2cm]{$j$}--++(0,1);
\draw[solid,spinach](0.25,0)node[below,spinach]{$\Sub{i}$}--++(0,1);
\end{tikzpicture}
\mapsfrom
a
\begin{tikzpicture}[scale=1.2,anchorbase,smallnodes,rounded corners]
\draw[ghost,dot](0,0)--++(0,1)node[above,yshift=-1pt]{$i$};
\draw[solid](0.5,0)node[below]{$j$}--++(0,1);
\end{tikzpicture}
+b
\begin{tikzpicture}[scale=1.2,anchorbase,smallnodes,rounded corners]
\draw[ghost](0,0)--++(0,1)node[above,yshift=-1pt]{$i$};
\draw[solid,dot](0.5,0)node[below]{$j$}--++(0,1);
\end{tikzpicture}
.
\end{gather*}
We finish by showing that $\sg$ respects one version of the
Reidemeister II relations for solid-solid and ghost-ghost crossing, all other cases being similar. If $k\neq i$, then
the Reidemeister II relations hold in the image of $\sg$:
\begin{gather*}
\begin{tikzpicture}[scale=1.2,anchorbase,smallnodes,rounded corners]
\draw[ghost](0,0)node[below]{$\phantom{i}$}--++(0.5,0.5)--++(-0.5,0.5)node[above,yshift=-1pt]{$i$};
\draw[ghost](0.5,0)--++(-0.5,0.5)--++(0.5,0.5)node[above,yshift=-1pt]{$k$};
\end{tikzpicture}
\mapsto
\begin{tikzpicture}[scale=1.2,anchorbase,smallnodes,rounded corners]
\draw[ghost](0,0)--++(0.5,0.5)--++(-0.5,0.5)node[above,yshift=-1pt]{$i$};
\draw[ghost](0.5,0)--++(-0.5,0.5)--++(0.5,0.5)node[above,yshift=-1pt]{$k$};
\draw[solid,spinach](0.2,0)node[below,spinach]{$\Sub{i}$}--++(0.5,0.5)--++(-0.5,0.5);
\end{tikzpicture}
=
\begin{tikzpicture}[scale=1.2,anchorbase,smallnodes,rounded corners]
\draw[ghost](0,0)--++(0,1)node[above,yshift=-1pt]{$i$};
\draw[ghost](0.5,0)--++(0,1)node[above,yshift=-1pt]{$k$};
\draw[solid,spinach](0.2,0)node[below,spinach]{$\Sub{i}$}--++(0,1);
\end{tikzpicture}
\mapsfrom
\begin{tikzpicture}[scale=1.2,anchorbase,smallnodes,rounded corners]
\draw[ghost](0,0)node[below]{$\phantom{i}$}--++(0,1)node[above,yshift=-1pt]{$i$};
\draw[ghost](0.5,0)--++(0,1)node[above,yshift=-1pt]{$k$};
\end{tikzpicture}
.
\end{gather*}

\TestRelation{c}
First, we consider \autoref{R:BraidGSG}. The most interesting case is when we subdivide $i$ and $k=j$, where we again need the assumption that $|\Sub{\sigma}_{\,\Sub{r}}|\ll 1$:
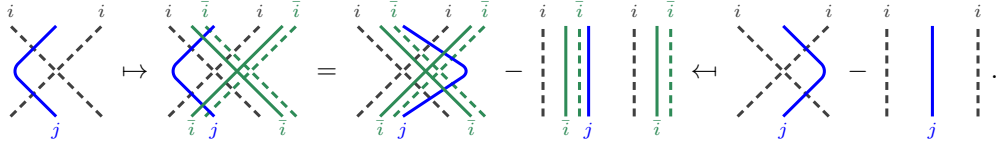
\begin{gather*}
\begin{tikzpicture}[scale=1.2,anchorbase,smallnodes,rounded corners]
\draw[ghost](1,0)--++(1,1)node[above,yshift=-1pt]{$i$};
\draw[ghost](2,0)--++(-1,1)node[above,yshift=-1pt]{$i$};
\draw[solid](1.5,0)node[below]{$j$}--++(-0.5,0.5)--++(0.5,0.5);
\end{tikzpicture}
\mapsto
\begin{tikzpicture}[scale=1.2,anchorbase,smallnodes,rounded corners]
\draw[ghost](1,0)--++(1,1)node[above,yshift=-1pt]{$i$};
\draw[ghost](2,0)--++(-1,1)node[above,yshift=-1pt]{$i$};
\draw[ghost,spinach](1.4,0)--++(1,1)node[above,yshift=-1pt,spinach]{$\Sub{i}$};
\draw[ghost,spinach](2.4,0)--++(-1,1)node[above,yshift=-1pt,spinach]{$\Sub{i}$};
\draw[solid](1.5,0)node[below]{$j$}--++(-0.5,0.5)--++(0.5,0.5);
\draw[solid,spinach](1.25,0)node[below,spinach]{$\Sub{i}$}--++(1,1);
\draw[solid,spinach](2.25,0)node[below,spinach]{$\Sub{i}$}--++(-1,1);
\end{tikzpicture}
=
\begin{tikzpicture}[scale=1.2,anchorbase,smallnodes,rounded corners]
\draw[ghost](1,0)--++(1,1)node[above,yshift=-1pt]{$i$};
\draw[ghost](2,0)--++(-1,1)node[above,yshift=-1pt]{$i$};
\draw[ghost,spinach](1.4,0)--++(1,1)node[above,yshift=-1pt,spinach]{$\Sub{i}$};
\draw[ghost,spinach](2.4,0)--++(-1,1)node[above,yshift=-1pt,spinach]{$\Sub{i}$};
\draw[solid](1.5,0)node[below]{$j$}--++(0.75,0.5)--++(-0.75,0.5);
\draw[solid,spinach](1.25,0)node[below,spinach]{$\Sub{i}$}--++(1,1);
\draw[solid,spinach](2.25,0)node[below,spinach]{$\Sub{i}$}--++(-1,1);
\end{tikzpicture}
-
\begin{tikzpicture}[scale=1.2,anchorbase,smallnodes,rounded corners]
\draw[ghost](1,0)--++(0,1)node[above,yshift=-1pt]{$i$};
\draw[ghost](2,0)--++(0,1)node[above,yshift=-1pt]{$i$};
\draw[ghost,spinach](1.4,0)--++(0,1)node[above,yshift=-1pt,spinach]{$\Sub{i}$};
\draw[ghost,spinach](2.4,0)--++(0,1)node[above,yshift=-1pt,spinach]{$\Sub{i}$};
\draw[solid](1.5,0)node[below]{$j$}--++(0,1);
\draw[solid,spinach](1.25,0)node[below,spinach]{$\Sub{i}$}--++(0,1);
\draw[solid,spinach](2.25,0)node[below,spinach]{$\Sub{i}$}--++(0,1);
\end{tikzpicture}
\mapsfrom
\begin{tikzpicture}[scale=1.2,anchorbase,smallnodes,rounded corners]
\draw[ghost](1,0)--++(1,1)node[above,yshift=-1pt]{$i$};
\draw[ghost](2,0)--++(-1,1)node[above,yshift=-1pt]{$i$};
\draw[solid](1.5,0)node[below]{$j$}--++(0.5,0.5)--++(-0.5,0.5);
\end{tikzpicture}
-
\begin{tikzpicture}[scale=1.2,anchorbase,smallnodes,rounded corners]
\draw[ghost](1,0)--++(0,1)node[above,yshift=-1pt]{$i$};
\draw[ghost](2,0)--++(0,1)node[above,yshift=-1pt]{$i$};
\draw[solid](1.5,0)node[below]{$j$}--++(0,1);
\end{tikzpicture}
.
\end{gather*}
All other variants of this relation can be checked {\ver}.
It is immediate that \autoref{R:RedSolid} is preserved under $\sg$, so it remains to check the Reidemeister III relations. To this end, one of the most
interesting cases is
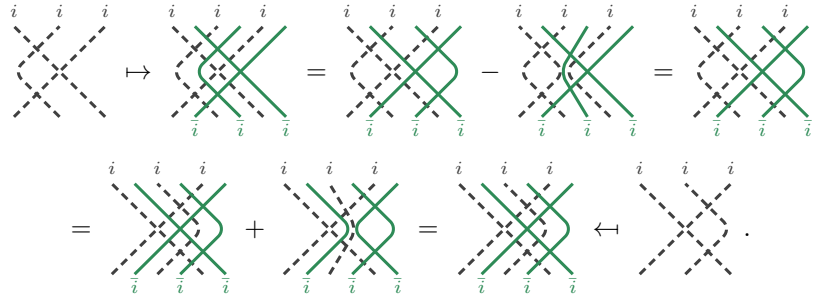
\begin{gather*}
\begin{tikzpicture}[scale=1.2,anchorbase,smallnodes,rounded corners]
\draw[ghost](1,0)node[below]{$\phantom{i}$}--++(1,1)node[above,yshift=-1pt]{$i$};
\draw[ghost](1.5,0)--++(-0.5,0.5)--++(0.5,0.5)node[above,yshift=-1pt]{$i$};
\draw[ghost](2,0)--++(-1,1)node[above,yshift=-1pt]{$i$};
\end{tikzpicture}
\mapsto
\begin{tikzpicture}[scale=1.2,anchorbase,smallnodes,rounded corners]
\draw[ghost](1,0)--++(1,1)node[above,yshift=-1pt]{$i$};
\draw[ghost](1.5,0)--++(-0.5,0.5)--++(0.5,0.5)node[above,yshift=-1pt]{$i$};
\draw[ghost](2,0)--++(-1,1)node[above,yshift=-1pt]{$i$};
\draw[solid,spinach](1.25,0)node[below,spinach]{$\Sub{i}$}--++(1,1);
\draw[solid,spinach](1.75,0)node[below,spinach]{$\Sub{i}$}--++(-0.5,0.5)--++(0.5,0.5);
\draw[solid,spinach](2.25,0)node[below,spinach]{$\Sub{i}$}--++(-1,1);
\end{tikzpicture}
=
\begin{tikzpicture}[scale=1.2,anchorbase,smallnodes,rounded corners]
\draw[ghost](1,0)--++(1,1)node[above,yshift=-1pt]{$i$};
\draw[ghost](1.5,0)--++(-0.5,0.5)--++(0.5,0.5)node[above,yshift=-1pt]{$i$};
\draw[ghost](2,0)--++(-1,1)node[above,yshift=-1pt]{$i$};
\draw[solid,spinach](1.25,0)node[below,spinach]{$\Sub{i}$}--++(1,1);
\draw[solid,spinach](1.75,0)node[below,spinach]{$\Sub{i}$}--++(0.5,0.5)--++(-0.5,0.5);
\draw[solid,spinach](2.25,0)node[below,spinach]{$\Sub{i}$}--++(-1,1);
\end{tikzpicture}
-
\begin{tikzpicture}[scale=1.2,anchorbase,smallnodes,rounded corners]
\draw[ghost](1,0)--++(0.5,0.5)--++(-0.5,0.5)node[above,yshift=-1pt]{$i$};
\draw[ghost](1.5,0)--++(-0.5,0.5)--++(0.5,0.5)node[above,yshift=-1pt]{$i$};
\draw[ghost](2,0)--++(-0.5,0.5)--++(0.5,0.5)node[above,yshift=-1pt]{$i$};
\draw[solid,spinach](1.25,0)node[below,spinach]{$\Sub{i}$}--++(1,1);
\draw[solid,spinach](1.75,0)node[below,spinach]{$\Sub{i}$}--++(-0.3,0.5)--++(+0.3,0.5);
\draw[solid,spinach](2.25,0)node[below,spinach]{$\Sub{i}$}--++(-1,1);
\end{tikzpicture}
=
\begin{tikzpicture}[scale=1.2,anchorbase,smallnodes,rounded corners]
\draw[ghost](1,0)--++(1,1)node[above,yshift=-1pt]{$i$};
\draw[ghost](1.5,0)--++(-0.5,0.5)--++(0.5,0.5)node[above,yshift=-1pt]{$i$};
\draw[ghost](2,0)--++(-1,1)node[above,yshift=-1pt]{$i$};
\draw[solid,spinach](1.25,0)node[below,spinach]{$\Sub{i}$}--++(1,1);
\draw[solid,spinach](1.75,0)node[below,spinach]{$\Sub{i}$}--++(0.5,0.5)--++(-0.5,0.5);
\draw[solid,spinach](2.25,0)node[below,spinach]{$\Sub{i}$}--++(-1,1);
\end{tikzpicture}
\\
=
\begin{tikzpicture}[scale=1.2,anchorbase,smallnodes,rounded corners]
\draw[ghost](1,0)--++(1,1)node[above,yshift=-1pt]{$i$};
\draw[ghost](1.5,0)--++(0.5,0.5)--++(-0.5,0.5)node[above,yshift=-1pt]{$i$};
\draw[ghost](2,0)--++(-1,1)node[above,yshift=-1pt]{$i$};
\draw[solid,spinach](1.25,0)node[below,spinach]{$\Sub{i}$}--++(1,1);
\draw[solid,spinach](1.75,0)node[below,spinach]{$\Sub{i}$}--++(0.5,0.5)--++(-0.5,0.5);
\draw[solid,spinach](2.25,0)node[below,spinach]{$\Sub{i}$}--++(-1,1);
\end{tikzpicture}
+
\begin{tikzpicture}[scale=1.2,anchorbase,smallnodes,rounded corners]
\draw[ghost](1,0)--++(1,1)node[above,yshift=-1pt]{$i$};
\draw[ghost](1.5,0)--++(0.3,0.5)--++(-0.3,0.5)node[above,yshift=-1pt]{$i$};
\draw[ghost](2,0)--++(-1,1)node[above,yshift=-1pt]{$i$};
\draw[solid,spinach](1.25,0)node[below,spinach]{$\Sub{i}$}--++(0.5,0.5)--++(-0.5,0.5);
\draw[solid,spinach](1.75,0)node[below,spinach]{$\Sub{i}$}--++(0.5,0.5)--++(-0.5,0.5);
\draw[solid,spinach](2.25,0)node[below,spinach]{$\Sub{i}$}--++(-0.5,0.5)--++(0.5,0.5);
\end{tikzpicture}
=
\begin{tikzpicture}[scale=1.2,anchorbase,smallnodes,rounded corners]
\draw[ghost](1,0)--++(1,1)node[above,yshift=-1pt]{$i$};
\draw[ghost](1.5,0)--++(0.5,0.5)--++(-0.5,0.5)node[above,yshift=-1pt]{$i$};
\draw[ghost](2,0)--++(-1,1)node[above,yshift=-1pt]{$i$};
\draw[solid,spinach](1.25,0)node[below,spinach]{$\Sub{i}$}--++(1,1);
\draw[solid,spinach](1.75,0)node[below,spinach]{$\Sub{i}$}--++(0.5,0.5)--++(-0.5,0.5);
\draw[solid,spinach](2.25,0)node[below,spinach]{$\Sub{i}$}--++(-1,1);
\end{tikzpicture}
\mapsfrom
\begin{tikzpicture}[scale=1.2,anchorbase,smallnodes,rounded corners]
\draw[ghost](1,0)node[below]{$\phantom{i}$}--++(1,1)node[above,yshift=-1pt]{$i$};
\draw[ghost](1.5,0)--++(0.5,0.5)--++(-0.5,0.5)node[above,yshift=-1pt]{$i$};
\draw[ghost](2,0)--++(-1,1)node[above,yshift=-1pt]{$i$};
\end{tikzpicture}.
\end{gather*}
The arguments for the other Reidemeister III relations are easier.

It remains to show that $\sg$ is bijective. Recall the basis $\WABasis(X)$ of $\WA(X)$ from \autoref{E:AffineBasis} and the basis
$\WAsBasis(\Sub{X})$ of $\WAs(\Sub{X})$. To see that $\sg(\WABasis)\subset\WAsBasis$ note that every element of $\WABasis$ is mapped to an element of $\WAsBasis$ because a permutation diagram $D(w)$ is sent to a permutation diagram $D(\Sub{w})$ because cabling preserves the property of being reduced so, for example,
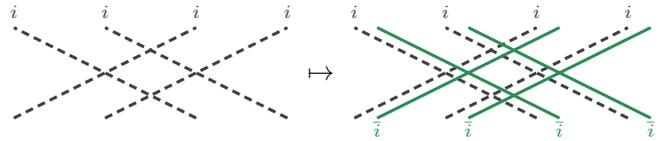
\begin{gather*}
\begin{tikzpicture}[scale=1.2,anchorbase,smallnodes,rounded corners]
\draw[ghost](0,0)node[below]{$\phantom{i}$}--++(2,1)node[above,yshift=-1pt]{$i$};
\draw[ghost](1,0)--++(2,1)node[above,yshift=-1pt]{$i$};
\draw[ghost](2,0)--++(-2,1)node[above,yshift=-1pt]{$i$};
\draw[ghost](3,0)--++(-2,1)node[above,yshift=-1pt]{$i$};
\end{tikzpicture}
\mapsto
\begin{tikzpicture}[scale=1.2,anchorbase,smallnodes,rounded corners]
\draw[ghost](0,0)--++(2,1)node[above,yshift=-1pt]{$i$};
\draw[ghost](1,0)--++(2,1)node[above,yshift=-1pt]{$i$};
\draw[ghost](2,0)--++(-2,1)node[above,yshift=-1pt]{$i$};
\draw[ghost](3,0)--++(-2,1)node[above,yshift=-1pt]{$i$};
\draw[solid,spinach](0.25,0)node[below,spinach]{$\Sub{i}$}--++(2,1);
\draw[solid,spinach](1.25,0)node[below,spinach]{$\Sub{i}$}--++(2,1);
\draw[solid,spinach](2.25,0)node[below,spinach]{$\Sub{i}$}--++(-2,1);
\draw[solid,spinach](3.25,0)node[below,spinach]{$\Sub{i}$}--++(-2,1);
\end{tikzpicture}
,
\end{gather*}
while dots on strings are untouched.
Moreover, using the faithfulness of the polynomial module from \autoref{C:FaithfulPoly}, we see that $\sg$ is an inclusion of $\WABasis$ into $\WAsBasis$. Note that $\sg$ does not surject $\WABasis$ onto $\WAsBasis$ because, for example, $\Sub{i}$-strings in the image of $\sg$ never carry dots.
However, ignoring dots for the time being, cabling arguments
show that every permutation diagram $D(\Sub{w})$ which
begins and ends with an idempotent as in
\autoref{E:IdemTrunc} is in the image of $\sg$.
Here we use that the permutation of the $i$-strings
and of the $\Sub{i}$-strings cannot be different since
otherwise \autoref{E:LocallyBad} annihilates the diagram. For example,
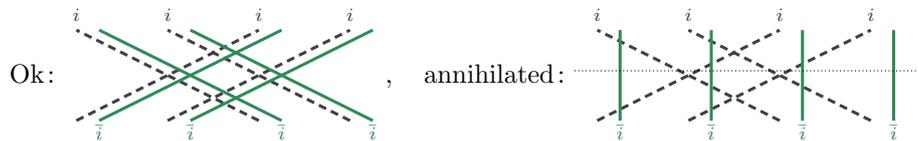
\begin{gather*}
\text{Ok}\colon
\begin{tikzpicture}[scale=1.2,anchorbase,smallnodes,rounded corners]
\draw[ghost](0,0)--++(2,1)node[above,yshift=-1pt]{$i$};
\draw[ghost](1,0)--++(2,1)node[above,yshift=-1pt]{$i$};
\draw[ghost](2,0)--++(-2,1)node[above,yshift=-1pt]{$i$};
\draw[ghost](3,0)--++(-2,1)node[above,yshift=-1pt]{$i$};
\draw[solid,spinach](0.25,0)node[below,spinach]{$\Sub{i}$}--++(2,1);
\draw[solid,spinach](1.25,0)node[below,spinach]{$\Sub{i}$}--++(2,1);
\draw[solid,spinach](2.25,0)node[below,spinach]{$\Sub{i}$}--++(-2,1);
\draw[solid,spinach](3.25,0)node[below,spinach]{$\Sub{i}$}--++(-2,1);
\end{tikzpicture}
,\quad
\text{annihilated}\colon
\begin{tikzpicture}[scale=1.2,anchorbase,smallnodes,rounded corners]
\draw[ghost](0,0)--++(2,1)node[above,yshift=-1pt]{$i$};
\draw[ghost](1,0)--++(2,1)node[above,yshift=-1pt]{$i$};
\draw[ghost](2,0)--++(-2,1)node[above,yshift=-1pt]{$i$};
\draw[ghost](3,0)--++(-2,1)node[above,yshift=-1pt]{$i$};
\draw[solid,spinach](0.25,0)node[below,spinach]{$\Sub{i}$}--++(0,1);
\draw[solid,spinach](1.25,0)node[below,spinach]{$\Sub{i}$}--++(0,1);
\draw[solid,spinach](2.25,0)node[below,spinach]{$\Sub{i}$}--++(0,1);
\draw[solid,spinach](3.25,0)node[below,spinach]{$\Sub{i}$}--++(0,1);
\draw[densely dotted] (-0.25,0.55) to (3.5,0.55);
\end{tikzpicture}
.
\end{gather*}
The right-hand diagram is annihilated because of the illustrated horizontal cut.

To deal with the dots we can use \autoref{R:GhostSolid} as follows: Take any solid $\Sub{i}$-string with a dot and apply
\begin{gather*}
a
\begin{tikzpicture}[scale=1.2,anchorbase,smallnodes,rounded corners]
\draw[ghost](0,0)--++(0,1)node[above,yshift=-1pt]{$i$};
\draw[solid,dot,spinach](0.5,0)node[below,spinach]{$\Sub{i}$}--++(0,1);
\end{tikzpicture}
=
-
\begin{tikzpicture}[scale=1.2,anchorbase,smallnodes,rounded corners]
\draw[ghost](0,0)--++(0.5,0.5)--++(-0.5,0.5)node[above,yshift=-1pt]{$i$};
\draw[solid,spinach](0.5,0)node[below,spinach]{$\Sub{i}$}--++(-0.5,0.5)--++(0.5,0.5);
\draw[densely dotted] (0.55,0.5)--++(-0.6,0);
\end{tikzpicture}
+b
\begin{tikzpicture}[scale=1.2,anchorbase,smallnodes,rounded corners]
\draw[ghost,dot](0,0)--++(0,1)node[above,yshift=-1pt]{$i$};
\draw[solid,spinach](0.5,0)node[below,spinach]{$\Sub{i}$}--++(0,1);
\end{tikzpicture}
=b
\begin{tikzpicture}[scale=1.2,anchorbase,smallnodes,rounded corners]
\draw[ghost,dot](0,0)--++(0,1)node[above,yshift=-1pt]{$i$};
\draw[solid,spinach](0.5,0)node[below,spinach]{$\Sub{i}$}--++(0,1);
\end{tikzpicture}.
\end{gather*}
Hereby we observe that the double crossing is zero
since the idempotent we get at the illustrated
horizontal cut is not compatible with the idempotents in \autoref{E:IdemTrunc}. Thus, we can get rid of all dots on solid $\Sub{i}$-strings
(recall that $a,b$ are units), concluding the proof.
\end{proof}

\begin{Example}
Back to \autoref{Ex:DivisionMap}. For the following element it is not immediate that it lies in the image of $\sg$, but it does):
\begin{gather*}
a
\begin{tikzpicture}[scale=1.2,anchorbase,smallnodes,rounded corners]
\draw[ghost](1,0)--++(0,1)node[above,yshift=-1pt]{$0$};
\draw[ghost](1.75,0)--++(0,1)node[above,yshift=-1pt]{$1$};
\draw[ghost,spinach,dot=0.1](1.35,0)--++(0,1)node[above,yshift=-1pt,spinach]{$3$};
\draw[solid](0,0)node[below]{$0$}--++(0,1);
\draw[solid](0.75,0)node[below]{$1$}--++(0,1);
\draw[solid,spinach,dot=0.1](1.2,0)node[below,spinach]{$3$}--++(0,1);
\draw[redstring](1.5,0)node[below]{$\rho$}--++(0,1);
\end{tikzpicture}
=
-
\begin{tikzpicture}[scale=1.2,anchorbase,smallnodes,rounded corners]
\draw[ghost](1,0)--++(0.15,0.5)--++(-0.15,0.5)node[above,yshift=-1pt]{$0$};
\draw[ghost](1.75,0)--++(0,1)node[above,yshift=-1pt]{$1$};
\draw[ghost,spinach](1.35,0)--++(-0.15,0.5)--++(0.15,0.5)node[above,yshift=-1pt,spinach]{$3$};
\draw[solid](0,0)node[below]{$0$}--++(0.15,0.5)--++(-0.15,0.5);
\draw[solid](0.75,0)node[below]{$1$}--++(0,1);
\draw[solid,spinach](1.2,0)node[below,spinach]{$3$}--++(-0.15,0.5)--++(0.15,0.5);
\draw[redstring](1.5,0)node[below]{$\rho$}--++(0,1);
\draw[densely dotted] (1.8,0.5)--++(-1.85,0);
\end{tikzpicture}
+b
\begin{tikzpicture}[scale=1.2,anchorbase,smallnodes,rounded corners]
\draw[ghost,dot=0.1](1,0)--++(0,1)node[above,yshift=-1pt]{$0$};
\draw[ghost](1.75,0)--++(0,1)node[above,yshift=-1pt]{$1$};
\draw[ghost,spinach](1.35,0)--++(0,1)node[above,yshift=-1pt,spinach]{$3$};
\draw[solid,dot=0.1](0,0)node[below]{$0$}--++(0,1);
\draw[solid](0.75,0)node[below]{$1$}--++(0,1);
\draw[solid,spinach](1.2,0)node[below,spinach]{$3$}--++(0,1);
\draw[redstring](1.5,0)node[below]{$\rho$}--++(0,1);
\end{tikzpicture}
=b
\begin{tikzpicture}[scale=1.2,anchorbase,smallnodes,rounded corners]
\draw[ghost,dot=0.1](1,0)--++(0,1)node[above,yshift=-1pt]{$0$};
\draw[ghost](1.75,0)--++(0,1)node[above,yshift=-1pt]{$1$};
\draw[ghost,spinach](1.35,0)--++(0,1)node[above,yshift=-1pt,spinach]{$3$};
\draw[solid,dot=0.1](0,0)node[below]{$0$}--++(0,1);
\draw[solid](0.75,0)node[below]{$1$}--++(0,1);
\draw[solid,spinach](1.2,0)node[below,spinach]{$3$}--++(0,1);
\draw[redstring](1.5,0)node[below]{$\rho$}--++(0,1);
\end{tikzpicture}
\in\mathrm{im}(\sg),
\end{gather*}
because we can apply \autoref{L:KillDiagrams} to the second diagram from the left.
\end{Example}

There is an analogous subdivision map $\sg^{\prime}\map{\WA(X)}{\SAX}$ given by putting the solid $i$-string to the left of the ghost $s(i)$-string, so that \autoref{E:Colors} becomes:
\begin{gather*}
\begin{tikzpicture}[scale=1.2,anchorbase,smallnodes,rounded corners]
\draw[ghost](1,0)--++(0,1)node[above,yshift=-1pt]{$i$};
\draw[solid](0,0)node[below]{$i$}--++(0,1)node[above,yshift=-1pt]{$\phantom{i}$};
\draw[redstring](0.5,0)node[below]{$\rho$}--++(0,1);
\end{tikzpicture}
\mapsto
\begin{tikzpicture}[scale=1.2,anchorbase,smallnodes,rounded corners]
\draw[ghost](4.0,0)--++(0,1)node[above,yshift=-1pt]{$s(i)$};
\draw[ghost,spinach](4.9,0)--++(0,1)node[above,yshift=-1pt,spinach]{$\Sub{i}$};
\draw[solid](3.0,0)node[below]{$s(i)$}--++(0,1);
\draw[solid,spinach](3.9,0)node[below,spinach]{$\Sub{i}$}--++(0,1);
\draw[redstring](3.5,0)node[below]{$\rho$}--++(0,1);
\draw[->] (3.9,-0.5)node[below]{added solid} to (3.9,-0.3);
\draw[->] (4.9,1.5)node[above]{added ghost} to (4.9,1.3);
\end{tikzpicture}
.
\end{gather*}
Everything above works verbatim for this subdivision. The map $\sg^{\prime}\map{\WA(X)}{\SAX}$ is better adjusted to the cyclotomic quotients:

\begin{Proposition}
The isomorphism $\sg^{\prime}\map{\WA(X)}{\SAX}$ of \autoref{T:SubDiv} descends to an isomorphism of the cyclotomic quotients.
\end{Proposition}

\begin{proof}
This follows immediately from \autoref{T:SubDiv} since $\sg^{\prime}$ gives a bijection between unsteady diagrams of $\WA(X)$ and the unsteady diagrams of $\SAX$.
\end{proof}

In the infinitesimal-case both \autoref{P:WebAlg} and \autoref{T:SubDiv} apply. Hence, we obtain the following result for the (cyclotomic) KLRW algebra:

\begin{Corollary}
We have an
isomorphism of graded algebras
$\sg^{\prime}\map{\TA}{\Sub{\TA}}$,
where $\Sub{\TA}$ is defined as in \autoref{E:IdemTrunc}. There are similar isomorphisms for the cyclotomic quotients.
\end{Corollary}

This corollary generalizes \cite[Theorem 2.12]{Ma-catrep-klr} to arbitrary quivers and to weighted KLRW algebras. In the special case of affine type~$A$, this result should be related to \cite[Theorem 2.2]{JaMa-equating-decomposition-numbers} and the Morita equivalences of \cite[Theorem 2]{ChMi-runner-removal}.


\section{Homogeneous (affine) cellular bases in type \texorpdfstring{$A$}{A}}\label{S:AffineCellular}


In the next two sections we analyze special cases where
the weighted KLRW algebra is a graded (affine) cellular
in the sense of \cite{GrLe-cellular} and \cite{KoXi-affine-cellular}.
Our construction of the homogeneous (affine) cellular basis
works for the quivers in \autoref{E:Quivers} and
is motivated by \cite{We-rouquier-dia-algebra}
and \cite{Bo-many-cellular-structures}, which only consider cyclotomic algebras
and type $A$. A special case of our constructions
gives homogeneous (affine) cellular basis
for the corresponding KLR(W) algebras and their cyclotomic quotients.

In this section we discuss the quivers of type $A$. In type $C$ the
arguments and results are essentially
the same, although they look quite different at first sight because
we use different recipes for constructing the diagrammatic bases. We explain the differences carefully in \autoref{S:TypeC}.

\begin{Remark}\label{R:Strategy}
In both cases, type $A$ and $C$, the crucial ingredient is the definition of
certain idempotent diagrams $\1_{\blam}$, for each
$\ell$-partition $\blam$ (and a bit more general). This idempotent is constructed roughly as follows:
\begin{enumerate}

\item Choose an ordering of the nodes of $\blam$. (We only do this implicitly
in type $A$, ordering by diagonal and height.)

\item This ordering gives an ordering of solid strings, which are in bijection with the nodes.

\item Place each string as far to the right as possible until it is held in
check by the relations and one of the earlier strings. For example, the relations
do not let us pull a ghost $i$-string through a solid $j$-string if $i\rightsquigarrow j$.
Similarly, a solid $i$-string cannot be pulled through a red $i$-string.

\item The cellular algebra partial order is given by looking how far strings are to the right.

\end{enumerate}
Given $\1_{\blam}$, the cellular bases are constructed by putting permutation
diagrams, indexed by semistandard $\blam$-tableaux, above and below $\1_{\blam}$.
\end{Remark}

The following example explains how we will construct the cellular bases for the weighted KLRW algebras. We have not yet introduced all of the notation used in the example because we intend that the reader will use this example to understand the definitions as they read through the paper.

\begin{Example}\label{E:Strategy}
In both, affine type $A$ and $C$, the strategy to construct the idempotent diagrams $\1_{\blam}$ is roughly as follows. For definiteness, we take the quiver of type $A_{\Z}$ with orientation $i\to(i+1)$. Combinatorially, and diagrammatically, this coincides with the quivers of type $A^{(1)}_e$ and $C^{(1)}_e$ with $e$ large enough and $i$ ``far from the boundary''. Let $\blam=(3,2)$, a partition of~$5$. We think of $\blam$ as being constructed inductively by adding nodes:
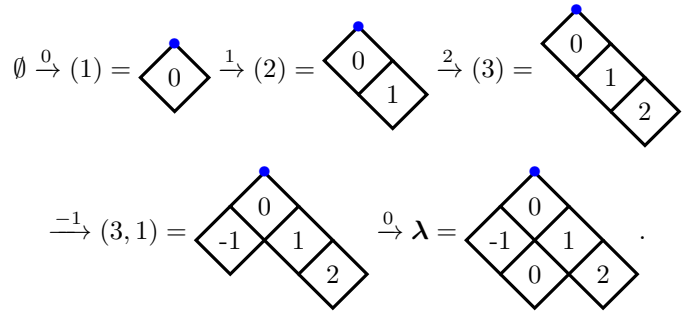
\begin{gather*}
\emptyset\xrightarrow{0}
(1)=\begin{tikzpicture}[scale=1.2,anchorbase,scale=0.75]
\draw[very thick] (0,-1) to (-0.5,-1.5) to (0,-2) to (0.5,-1.5) to (0,-1);
\node at (0,-1.5){$0$};
\node[blue] at (0,-1){$\bullet$};
\end{tikzpicture}
\xrightarrow{1}
(2)=\begin{tikzpicture}[scale=1.2,anchorbase,scale=0.75]
\draw[very thick] (0,-1) to (-0.5,-1.5) to (0,-2) to (0.5,-1.5) to (0,-1);
\draw[very thick] (0.5,-1.5) to (1,-2) to (0.5,-2.5) to (0,-2);
\node at (0,-1.5){$0$};
\node at (0.5,-2){$1$};
\node[blue] at (0,-1){$\bullet$};
\end{tikzpicture}
\xrightarrow{2}
(3)=\begin{tikzpicture}[scale=1.2,anchorbase,scale=0.75]
\draw[very thick] (0,-1) to (-0.5,-1.5) to (0,-2) to (0.5,-1.5) to (0,-1);
\draw[very thick] (0.5,-1.5) to (1,-2) to (0.5,-2.5) to (0,-2);
\draw[very thick] (1,-2) to (1.5,-2.5) to (1,-3) to (0.5,-2.5);
\node at (0,-1.5){$0$};
\node at (0.5,-2){$1$};
\node at (1,-2.5){$2$};
\node[blue] at (0,-1){$\bullet$};
\end{tikzpicture}
\\
\xrightarrow{-1}
(3,1)=\begin{tikzpicture}[scale=1.2,anchorbase,scale=0.75]
\draw[very thick] (-0.5,-1.5) to (-1,-2) to (-0.5,-2.5) to (0,-2);
\draw[very thick] (0,-1) to (-0.5,-1.5) to (0,-2) to (0.5,-1.5) to (0,-1);
\draw[very thick] (0.5,-1.5) to (1,-2) to (0.5,-2.5) to (0,-2);
\draw[very thick] (1,-2) to (1.5,-2.5) to (1,-3) to (0.5,-2.5);
\node at (0,-1.5){$0$};
\node at (-0.5,-2){$\text{-}1$};
\node at (0.5,-2){$1$};
\node at (1,-2.5){$2$};
\node[blue] at (0,-1){$\bullet$};
\end{tikzpicture}
\xrightarrow{0}
\blam=
\begin{tikzpicture}[scale=1.2,anchorbase,scale=0.75]
\draw[very thick] (-0.5,-1.5) to (-1,-2) to (-0.5,-2.5) to (0,-2);
\draw[very thick] (0,-1) to (-0.5,-1.5) to (0,-2) to (0.5,-1.5) to (0,-1);
\draw[very thick] (0.5,-1.5) to (1,-2) to (0.5,-2.5) to (0,-2);
\draw[very thick] (1,-2) to (1.5,-2.5) to (1,-3) to (0.5,-2.5);
\draw[very thick] (-0.5,-2.5) to (0,-3) to (0.5,-2.5);
\node at (0,-1.5){$0$};
\node at (0,-2.5){$0$};
\node at (-0.5,-2){$\text{-}1$};
\node at (0.5,-2){$1$};
\node at (1,-2.5){$2$};
\node[blue] at (0,-1){$\bullet$};
\end{tikzpicture}
.
\end{gather*}
As we will explain, we display the diagrams using Russian notation.  We have filled the nodes with their residues. By analogy, we want to construct the
idempotent diagram $\1_{\blam}$ by inductively adding solid strings with labels $0$, $0$, $-1$, $1$ and $2$, together with their ghosts. We construct $\1_{\blam}$ by putting each string as far to the right as possible:
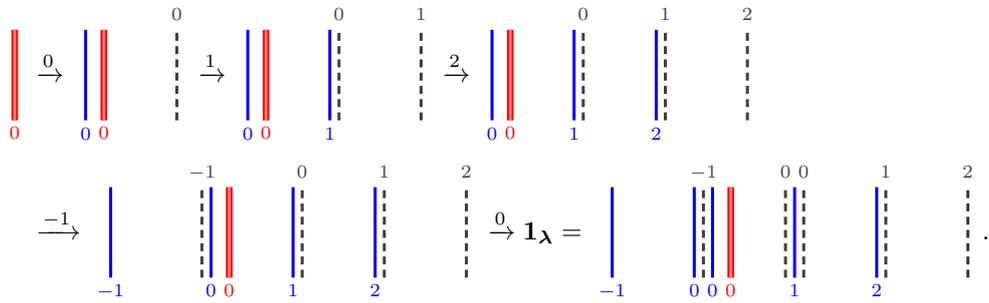
\begin{gather*}
\begin{aligned}
\begin{tikzpicture}[scale=1.2,anchorbase,smallnodes,rounded corners]
\draw[redstring](0.3,0)node[below]{$0$}--++(0,1)node[above,yshift=-1pt]{$\phantom{i}$};
\end{tikzpicture}
&\xrightarrow{0}
\begin{tikzpicture}[scale=1.2,anchorbase,smallnodes,rounded corners]
\draw[ghost](1.1,0)node[below]{$\phantom{i}$}--++(0,1)node[above,yshift=-1pt]{$0$};
\draw[solid](0.1,0)node[below]{$0$}--++(0,1)node[above,yshift=-1pt]{$\phantom{i}$};
\draw[redstring](0.3,0)node[below]{$0$}--++(0,1)node[above,yshift=-1pt]{$\phantom{i}$};
\end{tikzpicture}
\xrightarrow{1}
\begin{tikzpicture}[scale=1.2,anchorbase,smallnodes,rounded corners]
\draw[ghost](1.1,0)node[below]{$\phantom{i}$}--++(0,1)node[above,yshift=-1pt]{$0$};
\draw[ghost](2.0,0)node[below]{$\phantom{i}$}--++(0,1)node[above,yshift=-1pt]{$1$};
\draw[solid](0.1,0)node[below]{$0$}--++(0,1)node[above,yshift=-1pt]{$\phantom{i}$};
\draw[solid](1.0,0)node[below]{$1$}--++(0,1)node[above,yshift=-1pt]{$\phantom{i}$};
\draw[redstring](0.3,0)node[below]{$0$}--++(0,1)node[above,yshift=-1pt]{$\phantom{i}$};
\end{tikzpicture}
\xrightarrow{2}
\begin{tikzpicture}[scale=1.2,anchorbase,smallnodes,rounded corners]
\draw[ghost](1.1,0)node[below]{$\phantom{i}$}--++(0,1)node[above,yshift=-1pt]{$0$};
\draw[ghost](2.0,0)node[below]{$\phantom{i}$}--++(0,1)node[above,yshift=-1pt]{$1$};
\draw[ghost](2.9,0)node[below]{$\phantom{i}$}--++(0,1)node[above,yshift=-1pt]{$2$};
\draw[solid](0.1,0)node[below]{$0$}--++(0,1)node[above,yshift=-1pt]{$\phantom{i}$};
\draw[solid](1.0,0)node[below]{$1$}--++(0,1)node[above,yshift=-1pt]{$\phantom{i}$};
\draw[solid](1.9,0)node[below]{$2$}--++(0,1)node[above,yshift=-1pt]{$\phantom{i}$};
\draw[redstring](0.3,0)node[below]{$0$}--++(0,1)node[above,yshift=-1pt]{$\phantom{i}$};
\end{tikzpicture}
\\
&\xrightarrow{-1}
\begin{tikzpicture}[scale=1.2,anchorbase,smallnodes,rounded corners]
\draw[ghost](0,0)node[below]{$\phantom{i}$}--++(0,1)node[above,yshift=-1pt]{$-1$};
\draw[ghost](1.1,0)node[below]{$\phantom{i}$}--++(0,1)node[above,yshift=-1pt]{$0$};
\draw[ghost](2.0,0)node[below]{$\phantom{i}$}--++(0,1)node[above,yshift=-1pt]{$1$};
\draw[ghost](2.9,0)node[below]{$\phantom{i}$}--++(0,1)node[above,yshift=-1pt]{$2$};
\draw[solid](-1,0)node[below]{$-1$}--++(0,1)node[above,yshift=-1pt]{$\phantom{i}$};
\draw[solid](0.1,0)node[below]{$0$}--++(0,1)node[above,yshift=-1pt]{$\phantom{i}$};
\draw[solid](1.0,0)node[below]{$1$}--++(0,1)node[above,yshift=-1pt]{$\phantom{i}$};
\draw[solid](1.9,0)node[below]{$2$}--++(0,1)node[above,yshift=-1pt]{$\phantom{i}$};
\draw[redstring](0.3,0)node[below]{$0$}--++(0,1)node[above,yshift=-1pt]{$\phantom{i}$};
\end{tikzpicture}
\xrightarrow{0}
\1_{\blam}=
\begin{tikzpicture}[scale=1.2,anchorbase,smallnodes,rounded corners]
\draw[ghost](0,0)node[below]{$\phantom{i}$}--++(0,1)node[above,yshift=-1pt]{$-1$};
\draw[ghost](0.9,0)node[below]{$\phantom{i}$}--++(0,1)node[above,yshift=-1pt]{$0$};
\draw[ghost](1.1,0)node[below]{$\phantom{i}$}--++(0,1)node[above,yshift=-1pt]{$0$};
\draw[ghost](2.0,0)node[below]{$\phantom{i}$}--++(0,1)node[above,yshift=-1pt]{$1$};
\draw[ghost](2.9,0)node[below]{$\phantom{i}$}--++(0,1)node[above,yshift=-1pt]{$2$};
\draw[solid](-1,0)node[below]{$-1$}--++(0,1)node[above,yshift=-1pt]{$\phantom{i}$};
\draw[solid](-0.1,0)node[below]{$0$}--++(0,1)node[above,yshift=-1pt]{$\phantom{i}$};
\draw[solid](0.1,0)node[below]{$0$}--++(0,1)node[above,yshift=-1pt]{$\phantom{i}$};
\draw[solid](1.0,0)node[below]{$1$}--++(0,1)node[above,yshift=-1pt]{$\phantom{i}$};
\draw[solid](1.9,0)node[below]{$2$}--++(0,1)node[above,yshift=-1pt]{$\phantom{i}$};
\draw[redstring](0.3,0)node[below]{$0$}--++(0,1)node[above,yshift=-1pt]{$\phantom{i}$};
\end{tikzpicture}
.
\end{aligned}
\end{gather*}
To explain what is happening here, recall from \autoref{D:WeightedKLRW} that the relations of the form
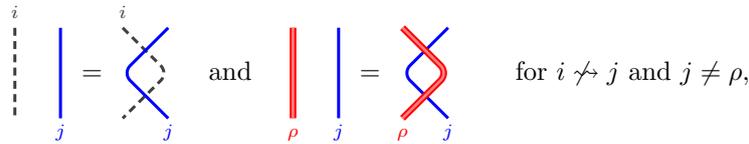
\begin{gather*}
\begin{tikzpicture}[scale=1.2,anchorbase,smallnodes,rounded corners]
\draw[ghost](0,1)node[above,yshift=-1pt]{$i$}--++(0,-1)node[below]{$\phantom{i}$};
\draw[solid](0.5,1)--++(0,-1)node[below]{$j$};
\end{tikzpicture}
=
\begin{tikzpicture}[scale=1.2,anchorbase,smallnodes,rounded corners]
\draw[ghost](0,1)node[above,yshift=-1pt]{$i$}--++(0.5,-0.5)--++(-0.5,-0.5) node[below]{$\phantom{i}$};
\draw[solid](0.5,1)--++(-0.5,-0.5)--++(0.5,-0.5) node[below]{$j$};
\end{tikzpicture}
\quad\text{and}\quad
\begin{tikzpicture}[scale=1.2,anchorbase,smallnodes,rounded corners]
\draw[solid](0.5,1)--++(0,-1)node[below]{$j$};
\draw[redstring](0,0)node[below]{$\rho$}--++(0,1)node[above,yshift=-1pt]{$\phantom{i}$};
\end{tikzpicture}
=\begin{tikzpicture}[scale=1.2,anchorbase,smallnodes,rounded corners]
\draw[solid](0.5,1)--++(-0.5,-0.5)--++(0.5,-0.5) node[below]{$j$};
\draw[redstring](0,0)node[below]{$\rho$}--++(0.5,0.5)--++(-0.5,0.5)node[above,yshift=-1pt]{$\phantom{i}$};
\end{tikzpicture}
\qquad\text{for $i\not\rightsquigarrow j$ and $j\neq\rho$,}
\end{gather*}
are \emph{honest Reidemeister II relations}. If we add strings to the left of the red string, then we can pull them to the right of the other strings whenever the honest Reidemeister II relations apply. When the honest relations do not apply we think of the strings as being stuck.

In the diagrams above, we start with the red $0$-string and then add the $0$-string and its ghost. The honest Reidemeister relations II allow us to pull the ghost $0$-string through the red string at which point the solid $0$-string is stuck behind the red string. Notice that because the solid $0$-string is stuck we cannot pull the ghost $0$-string further to the right. Next, we add the solid $1$-string and its ghost. The honest Reidemeister relations allow the ghost $1$-string to be pulled past all of the strings in the diagram but the solid $1$-string is stuck behind the $0$-ghost. Continuing in this way we obtain the diagrams above where, in each case, the strings are as far to the right as possible before becoming stuck.
This process ensures that the diagram that has all of the solid and ghost strings to the left of the red string, in row reading order, factors through the idempotent diagram $\1_{\blam}$.

The motivation for this construction is that a diagram is unsteady if a solid string can be pulled to the right of all of the other strings in the diagram. It is not obvious that the diagrams above are steady because it might be possible to pull a stuck string further to the right using the other relations in \autoref{D:WeightedKLRW}. That stuck strings cannot be pulled further to the right, or equivalently that these diagrams are steady, is a consequence of our main results.
\end{Example}

To describe how strings become stuck behind other strings we will work with explicit coordinates, or positioning functions, for the idempotent diagrams $\1_{\blam}$. Even though we employ the same strategy in affine types $A$ and $C$ the diagram combinatorics is different, so we treat the two cases separately. The next section introduces the positioning functions and tableaux combinatorics in affine type~$A$.


\subsection{Positioning and tableaux}\label{SS:Tableaux}


Let $\Gamma=(I,E)$ be a cyclic quiver of type $A^{(1)}_{e}$ and we identify $I$ with $\Z/(e+1)\Z$.
Without loss of generality, {\cf}
\autoref{P:FixedOrientation},
we will choose the orientation to be $i\to i+1$.
We fix ghost shifts
$\bsig\map{E}{\R_{\neq 0}}$ with $\sigma_{\epsilon}=1$ for all
$\epsilon\in E$.

As in
previous sections fix $n,\ell\in\Z_{\geq 0}$. Below we will work with $\WA[n](X)$
instead of $\WA(X)$ as the notation is a bit simpler.

Choose $\charge=(\kappa_{1},\dots,\kappa_{\ell})\in\Z^{\ell}$
with $\kappa_{1}<\dots<\kappa_{\ell}$
and $\brho=(\rho_{1},\dots,\rho_{\ell})\in I^{\ell}$.

Although we do not need them now, for the algebra $\WA[n](X)$ we will choose the
$Q$-polynomials in accordance with \autoref{E:QPoly}.

\begin{Remark}
The choices of the underlying quiver,
the $Q$-polynomials and the ghost shifts being $1$
while the $\kappa_{m}$ are integer units apart will play crucial roles.
\end{Remark}

\begin{Remark}
We could allow $e=\infty$, corresponding to
the quiver of type $A_{\Z}$. As this case is captured by taking $e$ to be sufficiently large, we assume that $e$ is finite because this simplifies the exposition below.
\end{Remark}

Let $\Parts$ be the set of \emph{$\ell$-partitions of $n$}. Identify
$\blam=(\lambda^{(1)}|\dots|\lambda^{(\ell)})\in\Parts$ with its
\emph{diagram}
\begin{gather*}
\set[\big]{(m,r,c)|1\leq k\leq\ell,1\leq r\text{ and }1\leq c\leq\lambda^{(m)}_{r}}.
\end{gather*}
A \emph{node} is any ordered triple
$(m,r,c)\in\Nodes=\bigcup_{\blam\in\Parts}\blam$, where $m$ is the component index, $r$ the row index and $c$ the column index.

\begin{Notation}
The reader unfamiliar with the combinatorics of partitions and tableaux used in this paper is referred to \cite[Section
3.3]{HuMa-klr-basis}. Compared with that paper, the only significant difference
is that we illustrate $\ell$-partitions and $\ell$-tableaux
using the Russian convention. For example, we draw the partition $(3,1)\in\Parts(1)[4]$ as:
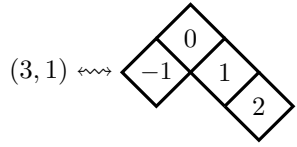
\begin{gather*}
(3,1)\leftrightsquigarrow
\begin{tikzpicture}[scale=1.2,anchorbase,scale=0.75]
\draw[very thick] (0,2) to (-0.5,1.5) to (-1,2) to (-0.5,2.5) to (0,2);
\draw[very thick] (0.5,1.5) to (1,2) to (0.5,2.5) to (0,2) to (0.5,1.5);
\draw[very thick] (0.5,1.5) to (1,1) to (1.5,1.5) to (1,2);
\draw[very thick] (-0.5,2.5) to (0,3) to (0.5,2.5);
\node at (0,2.5){$0$};
\node at (-0.5,2){\scalebox{0.95}{${-}1$}};
\node at (0.5,2){$1$};
\node at (1,1.5){$2$};
\end{tikzpicture}
\end{gather*}
The entries is the nodes are the differences $c-r$ between the column index $c$ and the row index $r$,
which is the \emph{content} of the associated node.
The Russian convention will give us a nice pictorial interpretation
of the coordinates of the strings in $\1_{\blam}$, as described in \autoref{R:Strategy}.
\end{Notation}

As we will see, the set of $\ell$-partitions $\Parts$ indexes the cells
of the cyclotomic weighted KLRW algebras but we need to extend this notation
slightly for $\WA[n](X)$. Set $\hell=\ell+n(e+1)$, and $\hell=1$ for $n=\ell=0$, so that $\Parts(\hell)$ is the
set of $\hell$-partitions of $n$.
If $n>0$, we will see that a proper subset of
$\Parts(\hell)$ indexes the cells of $\WA[n](X)$, not
$\Parts(\hell)$ itself. Identify $\Parts$ with the
left-adjusted subset of $\Parts(\hell)$, given by having empty partitions from $\ell+1$ onward. The nodes with $m>\ell$ are the \emph{affine nodes}. More generally, we attach the adjective \emph{affine} to anything associated to nodes $(m,r,c)$ with $m>\ell$.

Define
$\affine{\charge}=
(\affine{\kappa}_{1},\dots,\affine{\kappa}_{\hell})\in\Z^{\hell}$
and
$\affine{\brho}=(\affine{\rho}_{1},\dots,\affine{\rho}_{\hell})\in I^{\hell}$ by
\begin{gather*}
\affine{\kappa}_{m}=
\begin{cases*}
\kappa_{m} & if $1\leq m\leq\ell$,\\
\kappa_{\ell}+2n(m-\ell) & otherwise,
\end{cases*}
\quad\text{and}\quad
\affine{\rho}_{m}=
\begin{cases*}
\rho_{m} & if $1\leq m\leq\ell$,\\
\floor{m-\ell-1}{n}+(e+1)\Z & otherwise.
\end{cases*}
\end{gather*}
If $\ell=0$, then $\kappa_{\ell}=0$ in the above formula, by convention.
Note that $\affine{\kappa}_{1}<\dots<\affine{\kappa}_{\hell}$.
Define the \emph{($\affine{\brho}$-)residue} of the node $(m,r,c)$ to be
$\res(m,r,c)=c-r+\affine{\rho}_{m}+(e+1)\Z\in I$.

Fix $0<\varepsilon<\tfrac{1}{2n\hell}$.
Motivated by \cite[Section 5]{Bo-many-cellular-structures}, define a
\emph{positioning function}:

\begin{Definition}\label{D:Acoordinates}
Let $\hcoord\map{\Nodes}{\R}$ be the map
\begin{gather*}
\hcoord(m,r,c)=\affine{\kappa}_{m}+(c-r)-\tfrac{m}{\hell}-(c+r)\varepsilon.
\end{gather*}
For $\blam\in\Parts(\hell)$ let
$\hcoord(\blam)=\set{\hcoord(m,r,c)|(m,r,c)\in\blam}$ and write
$\hcoord(\blam)=\set{x^{\blam}_{1}<\dots<x^{\blam}_{n}}$.
The \emph{suspension point} of $\lambda^{(m)}$ is $\hcoord(m,1,1)$,
for $1\leq m\leq\hell$.
\end{Definition}

The following lemma ensures that all of the rightmost
coordinates for nodes in some $\blam\in\Parts(\ell)$ are always to the
left of the coordinates of all affine nodes:

\begin{Lemma}\label{L:RightLeft}
We have $\hcoord(\ell,1,n)<\hcoord(k,r,c)$ for all $k\in\set{\ell+1,\dots,\hell}$ and all $r,c\in\set{1,\dots,n}$.
\end{Lemma}

\begin{proof}
The claim is immediate for $\ell=0$. Otherwise,
the difference between $\hcoord(\ell,1,n)=\kappa_{\ell}+(n-1)-\tfrac{\ell}{\hell}-(n+1)\varepsilon$
and $\hcoord(\ell+1,n,1)=\kappa_{\ell}+2n(\ell+1-\ell)+(1-n)-
\tfrac{\ell+1}{\hell}-(1+n)\varepsilon$ is $\hcoord(\ell+1,n,1)-\hcoord(\ell,1,n)=1-\tfrac{1}{\hell}>0$.
This completes the proof as $\hcoord(\ell+1,n,1)\le\hcoord(\ell+k,r,c)$.
\end{proof}

For completeness, we note that the value $c-r$ is
sometimes called the diagonal or content, and the value $c+r$ is the height
of the component $\lambda^{(m)}$.
In examples, we will often use the
\emph{approximate coordinates}, which are obtained by taking the limit
$\varepsilon\to 0$ in the coordinates $\hcoord(m,r,c)$.

\begin{Lemma}
Let $\alpha$ and $\alpha^{\prime}$ be distinct nodes in $\blam$,
for $\blam\in\Parts$. Then $\hcoord(\alpha)\neq\hcoord(\alpha^{\prime})$.
\end{Lemma}

\begin{proof}
This follows easily by induction on $n$, using the fact that $\varepsilon$ is small.
\end{proof}

The red strings have coordinates $\affine{\kappa}_{m}$, for $1\leq m\leq\hell$, which are slightly to the right of the suspension points $\hcoord(m,1,1)$ of the components $\lambda^{(m)}$ in the Russian diagram of a $\hell$-partition. As in \autoref{E:Strategy}, we think of the $m$th red string as being stuck at position $\affine{\kappa}_{m}$. By convention, the suspension points, and the coordinates of the red strings, move to the right as $m$ increases:
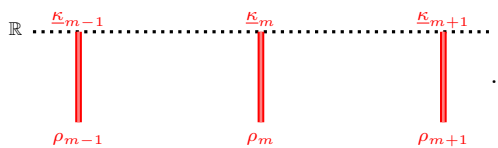
\begin{gather*}
\begin{tikzpicture}[scale=1.2,anchorbase,smallnodes,rounded corners]
\draw[very thick,dotted] (-2.5,1)node[left]{$\R$} to (2.5,1);
\draw[redstring](-2,0)node[below,red]{$\rho_{m-1}$}--++(0,1)node[above,yshift=0pt]{$\affine{\kappa}_{m-1}$};
\draw[redstring](0,0)node[below,red]{$\rho_{m}$}--++(0,1)node[above,yshift=0pt]{$\affine{\kappa}_{m}$};
\draw[redstring](2,0)node[below,red]{$\rho_{m+1}$}--++(0,1)node[above,yshift=0pt]{$\affine{\kappa}_{m+1}$};
\end{tikzpicture}
.
\end{gather*}
For $\ell<m\leq\hell$, the \emph{affine red strings} have $x$-coordinates $\affine{\kappa}_{m}$, which are to the right of all (honest) red strings.
Although we often refer to the affine red strings, and sometimes draw them in
diagrams, these strings are not part of the diagrams and are only used as a visual aid. Our diagrams will therefore contain four types of strings:
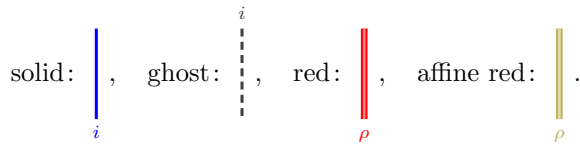
\begin{gather*}
\text{solid}\colon
\begin{tikzpicture}[scale=1.2,anchorbase,smallnodes,rounded corners]
\draw[solid](0,0)node[below]{$i$}--++(0,1)node[above,yshift=-1pt]{$\phantom{i}$};
\end{tikzpicture}
,\quad
\text{ghost}\colon
\begin{tikzpicture}[scale=1.2,anchorbase,smallnodes,rounded corners]
\draw[ghost](0,1)node[above,yshift=-1pt]{$i$}--++(0,-1)node[below]{$\phantom{i}$};
\end{tikzpicture}
,\quad
\text{red}\colon
\begin{tikzpicture}[scale=1.2,anchorbase,smallnodes,rounded corners]
\draw[redstring](0,0)node[below,red]{$\rho$}--++(0,1)node[above,yshift=-1pt]{$\phantom{i}$};
\end{tikzpicture}
,\quad
\text{affine red}\colon
\begin{tikzpicture}[scale=1.2,anchorbase,smallnodes,rounded corners]
\draw[affine](0,0)node[below]{$\rho$}--++(0,1)node[above,yshift=-1pt]{$\phantom{i}$};
\end{tikzpicture}
.
\end{gather*}
We draw affine red strings in a different color to red strings to help distinguish them.

\begin{Remark}
The relations in \autoref{D:RationalCherednik} do not allow us to pull solid
$i$-strings through red $i$-strings. Although affine red strings do not
exist, and do not contribute to the relations, we
think of them as blocking solid $i$-strings in the same way
that (honest) red $i$-strings do.
\end{Remark}

\begin{Definition}
For $\blam\in\Parts(\hell)$ let
$\1_{\blam}$ be the idempotent diagram with positions $\hcoord(\blam)$
such that the solid string with position $\hcoord(m,r,c)$ has residue $\res(m,r,c)$.
\end{Definition}

Note that $\1_{\blam}$ has only $\ell$ red strings,
their coordinates are given by $\kappa$, after which there are $\hell-\ell$ affine red strings.

\begin{Notation}
In illustrations we will indicate the suspension points on the Russian diagrams of $\hell$-partitions with a red dot. If we do not draw a suspension point, such as in \autoref{E:BoxConfiguration} below, then will only be showing only a local component in a diagram.
\end{Notation}

\begin{Remark}\label{R:HangingPictures}
Informally, the positioning function tilts an $\hell$-partition $\blam$ in Russian
notation by a small clockwise rotation, and then projects the top corner of each node
`up to $\R$' so that $\affine{\kappa}_{m}-\tfrac{m}{\hell}-2\varepsilon$ is the position
of the node in $\lambda^{(m)}$ with $c=r=1$. For example, if $\lambda_{5}=(2,2)$ (note
that, strictly speaking, we should write $\lambda_{5}=\bigl((2,2)\bigr)$ but we will
not use this cumbersome notation), then we have in approximate
coordinates:
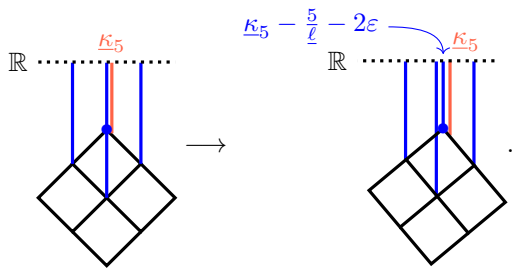
\begin{gather*}
\begin{tikzpicture}[scale=1.2,anchorbase,scale=0.75]
\draw[very thick] (-0.5,-1.5) to (-1,-2) to (-0.5,-2.5) to (0,-2);
\draw[very thick] (0,-1) to (-0.5,-1.5) to (0,-2) to (0.5,-1.5) to (0,-1);
\draw[very thick] (0.5,-1.5) to (1,-2) to (0.5,-2.5) to (0,-2);
\draw[very thick] (-0.5,-2.5) to (0,-3) to (0.5,-2.5);
\draw[very thick,blue] (-0.5,-1.5) to (-0.5,-0.01);
\draw[very thick,blue] (0,-2) to (0,-0.01);
\draw[very thick,tomato] (0.075,-1.05) to (0.075,-0.01)node[above,tomato]{$\affine{\kappa}_{5}$};
\draw[very thick,blue] (0.5,-1.5) to (0.5,-0.01);
\draw[very thick,dotted] (-1,0)node[left]{$\R$} to (1,0);
\node[blue] at (0,-1){$\bullet$};
\end{tikzpicture}
\longrightarrow
\raisebox{0.18cm}{$\begin{tikzpicture}[scale=1.2,anchorbase,scale=0.75]
\draw[very thick,blue] (-0.55,-1.45) to (-0.55,-0.01);
\draw[very thick,blue] (-0.1,-2) to (-0.1,-0.01);
\draw[very thick,blue] (0.0,-1) to (0.0,-0.01);
\draw[very thick,tomato] (0.1,-1.1) to (0.1,-0.01)node[above,tomato,xshift=0.21cm]{$\affine{\kappa}_{5}$};
\draw[very thick,blue] (0.45,-1.55) to (0.45,-0.01);
\begin{scope}[rotate around={-5:(0,-1)}]
\draw[very thick] (0.5,-1.5) to (1,-2) to (0.5,-2.5) to (0,-2);
\draw[very thick] (0,-1) to (0.5,-1.5) to (0,-2) to (-0.5,-1.5) to (0,-1);
\draw[very thick] (-0.5,-1.5) to (-1,-2) to (-0.5,-2.5) to (0,-2);
\draw[very thick] (0.5,-2.5) to (0,-3) to (-0.5,-2.5);
\end{scope}
\draw[very thick,dotted] (0.75,0) to (-1.25,0)node[left]{$\R$};
\node[blue] at (0,-1){$\bullet$};
\draw[->,thin,blue](-0.8,0.5)node[left]{$\affine{\kappa}_{5}-\tfrac{5}{\hell}-2\varepsilon$} to [out=0,in=90](-0.0,0.1);
\end{tikzpicture}$}
.
\end{gather*}
In other words, the further a box is to the southeast, the larger the
associated coordinate in $\R$, the distance between the various components of
$\blam$ is $(\affine{\kappa}_{m}-\affine{\kappa}_{m^{\prime}})-\tfrac{m-m^{\prime}}{\hell}$,
and the distance between nodes in the same component is given by steps of length approximately
$1$. The shift $\varepsilon$ ensures that two different strings are not mapped to the same
point in $\R$ and it is also important for the order
of the strings since it tilts the diagram slightly.
Consider a part of $\blam$ of the following form:
\begin{gather}\label{E:BoxConfiguration}
\begin{tikzpicture}[scale=1.2,anchorbase,scale=0.75]
\draw[very thick] (-0.5,1.5) to (-1,2) to (-0.5,2.5) to (0,2);
\draw[very thick] (0,1) to (-0.5,1.5) to (0,2) to (0.5,1.5) to (0,1);
\draw[very thick] (0.5,1.5) to (1,2) to (0.5,2.5) to (0,2);
\draw[very thick] (-0.5,2.5) to (0,3) to (0.5,2.5);
\node at (0,1.5){$i$};
\node at (0,2.5){$i$};
\node at (-0.5,2){\scalebox{0.85}{$i{-}1$}};
\node at (0.5,2){\scalebox{0.85}{$i{+}1$}};
\end{tikzpicture}
\quad\text{and}\quad
\begin{tikzpicture}[scale=1.2,anchorbase,scale=0.75]
\draw[very thick] (-0.5,1.5) to (-1,2) to (-0.5,2.5) to (0,2);
\draw[very thick] (0,1) to (-0.5,1.5) to (0,2) to (0.5,1.5) to (0,1);
\draw[very thick] (0.5,1.5) to (1,2) to (0.5,2.5) to (0,2);
\draw[very thick] (-0.5,2.5) to (0,3) to (0.5,2.5);
\node at (0,1.5){\scalebox{0.85}{${-}4\varepsilon$}};
\node at (0,2.5){\scalebox{0.85}{${-}2\varepsilon$}};
\node at (-0.5,2){\scalebox{0.85}{${-}3\varepsilon$}};
\node at (0.5,2){\scalebox{0.85}{${-}3\varepsilon$}};
\end{tikzpicture}
\leftrightsquigarrow
\begin{tikzpicture}[scale=1.2,anchorbase,smallnodes,rounded corners]
\draw[ghost](0.1,0)node[below]{$\phantom{i}$}--++(0,1)node[above,yshift=-1pt]{$i{-}1$};
\draw[ghost](1,0)node[below]{$\phantom{i}$}--++(0,1)node[above,yshift=-1pt]{$i$};
\draw[ghost](1.2,0)node[below]{$\phantom{i}$}--++(0,1)node[above,yshift=-1pt]{$i$};
\draw[ghost](2.1,0)node[below]{$\phantom{i}$}--++(0,1)node[above,yshift=-1pt]{$i{+}1$};
\draw[solid](-0.9,0)node[below]{$i{-}1$}--++(0,1)node[above,yshift=-1pt]{$\phantom{i}$};
\draw[solid](0,0)node[below]{$i$}--++(0,1)node[above,yshift=-1pt]{$\phantom{i}$};
\draw[solid](0.2,0)node[below]{$i$}--++(0,1)node[above,yshift=-1pt]{$\phantom{i}$};
\draw[solid](1.1,0)node[below]{$i{+}1$}--++(0,1)node[above,yshift=-1pt]{$\phantom{i}$};
\end{tikzpicture}
,
\end{gather}
We have filled the nodes with their residues and $\varepsilon$-tilts, with the right-hand side
illustrating the corresponding part of $\1_{\blam}$,
ignoring all other nodes and strings. In words, two $i$ nodes on the same diagonal
have a ghost $(i-1)$-string, respectively solid $(i+1)$-string, in between the associated
solid and ghost $i$-strings. Moreover, we have the following configurations:
\begin{gather}\label{E:CasesTypeA}
\begin{aligned}
\begin{tikzpicture}[scale=1.2,anchorbase,scale=0.75]
\draw[very thick] (-0.5,1.5) to (-1,2) to (-0.5,2.5) to (0,2);
\draw[very thick] (0,1) to (-0.5,1.5) to (0,2) to (0.5,1.5) to (0,1);
\node at (-0.5,2){$i$};
\node at (0,1.5){\scalebox{0.85}{$i{+}1$}};
\end{tikzpicture}
\quad\text{and}\quad
\begin{tikzpicture}[scale=1.2,anchorbase,scale=0.75]
\draw[very thick] (-0.5,1.5) to (-1,2) to (-0.5,2.5) to (0,2);
\draw[very thick] (0,1) to (-0.5,1.5) to (0,2) to (0.5,1.5) to (0,1);
\node at (-0.5,2){\scalebox{0.85}{${-}2\varepsilon$}};
\node at (0,1.5){\scalebox{0.85}{${-}3\varepsilon$}};
\end{tikzpicture}
&\leftrightsquigarrow
\begin{tikzpicture}[scale=1.2,anchorbase,smallnodes,rounded corners]
\draw[ghost](1.2,0)node[below]{$\phantom{i}$}--++(0,1)node[above,yshift=-1pt]{$i$};
\draw[ghost](2.1,0)node[below]{$\phantom{i}$}--++(0,1)node[above,yshift=-1pt]{$i{+}1$};
\draw[solid](0.2,0)node[below]{$i$}--++(0,1)node[above,yshift=-1pt]{$\phantom{i}$};
\draw[solid](1.1,0)node[below]{$i{+}1$}--++(0,1)node[above,yshift=-1pt]{$\phantom{i}$};
\end{tikzpicture}
,\\
\begin{tikzpicture}[scale=1.2,anchorbase,scale=0.75]
\draw[very thick] (0.5,1.5) to (1,2) to (0.5,2.5) to (0,2);
\draw[very thick] (0,1) to (0.5,1.5) to (0,2) to (-0.5,1.5) to (0,1);
\node at (0.5,2){$i$};
\node at (0,1.5){\scalebox{0.85}{$i{-}1$}};
\end{tikzpicture}
\quad\text{and}\quad
\begin{tikzpicture}[scale=1.2,anchorbase,scale=0.75]
\draw[very thick] (0.5,1.5) to (1,2) to (0.5,2.5) to (0,2);
\draw[very thick] (0,1) to (0.5,1.5) to (0,2) to (-0.5,1.5) to (0,1);
\node at (0.5,2){\scalebox{0.85}{${-}2\varepsilon$}};
\node at (0,1.5){\scalebox{0.85}{${-}3\varepsilon$}};
\end{tikzpicture}
&\leftrightsquigarrow
\begin{tikzpicture}[scale=1.2,anchorbase,smallnodes,rounded corners]
\draw[ghost](0.1,0)node[below]{$\phantom{i}$}--++(0,1)node[above,yshift=-1pt]{$i{-}1$};
\draw[ghost](1.2,0)node[below]{$\phantom{i}$}--++(0,1)node[above,yshift=-1pt]{$i$};
\draw[solid](-0.9,0)node[below]{$i{-}1$}--++(0,1)node[above,yshift=-1pt]{$\phantom{i}$};
\draw[solid](0.2,0)node[below]{$i$}--++(0,1)node[above,yshift=-1pt]{$\phantom{i}$};
\end{tikzpicture}
.
\end{aligned}
\end{gather}
Once again, we have not illustrated other strings that may appear in this diagram.
\end{Remark}

\begin{Example}\label{Ex:BigPositioning}
Let $n=7$, $\ell=2$ and take the quiver from \autoref{E:ExampleQuiver}
with the indicated labeling. Let $\charge=(1,2)$ and $\rho=(1,2)$. Then
$\hell=23$, $\affine{\charge}=(1,2,14,28,42,\dots,294)$ and
$\affine{\brho}=(1,2,0^{7},1^{7},2^{7})$, where the exponents indicate
repeated entries. Take $(3,2|1^{2})$
and consider it as the $23$-partition
$\blam=(3,1|\emptyset|1^{2}|\emptyset|\dots|\emptyset)$.
Let $\res$ and $\varepsilon\text{-tilt}$ be the sets of all residues
and $\varepsilon$-tilts, respectively.
Using $\varepsilon=0.01$ and framing the data corresponding to suspension points we have:
\begin{gather*}
\begin{tikzpicture}[scale=1,anchorbase]
  \begin{scope}[rotate around={-5:(0,-1)}]
\draw[very thick] (-0.5,-1.5) to (-1,-2) to (-0.5,-2.5) to (0,-2);
\draw[very thick] (0,-1) to (-0.5,-1.5) to (0,-2) to (0.5,-1.5) to (0,-1);
\draw[very thick] (0.5,-1.5) to (1,-2) to (0.5,-2.5) to (0,-2);
\draw[very thick] (1,-2) to (1.5,-2.5) to (1,-3) to (0.5,-2.5);
\draw[very thick] (-0.5,-2.5) to (0,-3) to (0.5,-2.5);
\end{scope}
\begin{scope}[rotate around={-5:(4,-1)}]
\draw[very thick] (3.5,-1.5) to (3,-2) to (3.5,-2.5) to (4,-2);
\draw[very thick] (4,-1) to (3.5,-1.5) to (4,-2) to (4.5,-1.5) to (4,-1);
\end{scope}
\draw[very thick,blue] (-0.55,-1.45) to (-0.55,-0.01);
\draw[very thick,blue] (0,-1) to (0,-0.01);
\draw[very thick,tomato] (0.1,-1.15) to (0.1,-0.01)node[above,yshift=-1pt]{$1$};
\draw[very thick,blue] (-0.1,-2) to (-0.1,-0.01);
\draw[very thick,blue] (0.44,-1.55) to (0.44,-0.01);
\draw[very thick,blue] (0.9,-2.1) to (0.9,-0.01);
\draw[very thick,blue] (3.48,-1.45) to (3.48,-0.01);
\draw[very thick,blue] (4,-1) to (4,-0.01);
\draw[very thick,tomato] (4.1,-1.1) to (4.1,-0.01)node[above,yshift=-1pt]{$14$};
\draw[very thick,dotted] (-1,0)node[left]{$\R$} to (4.5,0);
\node[blue] at (0,-1){$\bullet$};
\node[blue] at (4,-1){$\bullet$};
\end{tikzpicture}
,\quad
\begin{gathered}
\hcoord(\blam)\approx
\set{-0.07,0.92,\framebox{$0.94$},1.93,2.92,12.84,\framebox{$13.85$}},
\\
\affine{\charge}=\bigl(\framebox{$1$},\framebox{$2$},14,28,42,\dots,294\bigr),
\\
\res=(0,1,\framebox{$1$},2,0,2,\framebox{$0$}),
\\
\varepsilon\text{-tilt}=(-0.03,-0.04,\framebox{$-0.02$},-0.03,-0.04,-0.03,\framebox{$0.02$}),
\end{gathered}
\\
\1_{\blam}=
\begin{tikzpicture}[scale=1.2,anchorbase,smallnodes,rounded corners]
\draw[ghost](0.2,0)node[below]{$\phantom{i}$}--++(0,1)node[above,yshift=-1pt]{$0$};
\draw[ghost](1.1,0)node[below]{$\phantom{i}$}--++(0,1)node[above,yshift=-1pt]{$1$};
\draw[ghost](1.3,0)node[below]{$\phantom{i}$}--++(0,1)node[above,yshift=-1pt]{$1$};
\draw[ghost](2.2,0)node[below]{$\phantom{i}$}--++(0,1)node[above,yshift=-1pt]{$2$};
\draw[ghost](3.1,0)node[below]{$\phantom{i}$}--++(0,1)node[above,yshift=-1pt]{$0$};
\draw[ghost](5.0,0)node[below]{$\phantom{i}$}--++(0,1)node[above,yshift=-1pt]{$2$};
\draw[ghost](6.1,0)node[below]{$\phantom{i}$}--++(0,1)node[above,yshift=-1pt]{$0$};
\draw[solid](-0.8,0)node[below]{$0$}--++(0,1)node[above,yshift=-1pt]{$\phantom{i}$};
\draw[solid](0.1,0)node[below]{$1$}--++(0,1)node[above,yshift=-1pt]{$\phantom{i}$};
\draw[solid](0.3,0)node[below]{$1$}--++(0,1)node[above,yshift=-1pt]{$\phantom{i}$};
\draw[solid](1.2,0)node[below]{$2$}--++(0,1)node[above,yshift=-1pt]{$\phantom{i}$};
\draw[solid](2.1,0)node[below]{$0$}--++(0,1)node[above,yshift=-1pt]{$\phantom{i}$};
\draw[solid](4.0,0)node[below]{$2$}--++(0,1)node[above,yshift=-1pt]{$\phantom{i}$};
\draw[solid](5.1,0)node[below]{$0$}--++(0,1)node[above,yshift=-1pt]{$\phantom{i}$};
\draw[redstring](0.45,0)node[below]{$1$}--++(0,1)node[above,yshift=-1pt]{$\phantom{i}$};
\draw[redstring](2.35,0)node[below]{$2$}--++(0,1)node[above,yshift=-1pt]{$\phantom{i}$};
\draw[affine](5.25,0)node[below]{$0$}--++(0,1)node[above,xshift=0.25cm,yshift=-0.02]{$\phantom{i}$};
\end{tikzpicture}
.
\end{gather*}
We have also drawn the corresponding idempotent diagram $\1_{\blam}$
(without precise scaling of the endpoints). We emphasize that the
affine red string, drawn in grayish yellow, does not exist and are drawn only to
show the suspension point for $\affine{\kappa}_{3}$. We have not drawn the remaining
$20$ affine red strings as they do not play a role here. Note that $\varepsilon$
shifts the coordinate of $(m,r+1,c+1)$ slightly to the left of $(m,r,c)$. The
$\varepsilon$-tilts are also illustrated.
\end{Example}

The idempotent diagram $\1_{\blam}$ in \autoref{Ex:BigPositioning}
is prototypical and shows the point of the positioning function:
every solid string in this $\1_{\blam}$ is ``as far to the right as possible''.
That is, one can check that no Reidemeister II relations from \autoref{D:RationalCherednik}
can be used to pull strings in $\1_{\blam}$ further to the right.
(Here we think of affine red strings as they are blocking other strings.)
See also \autoref{E:Strategy} for some additional explanation along the same lines.

\begin{Example}\label{Ex:Ordered}
An important class of examples is given by $\hell$-partitions $\blam$
that have only one non-trivial component $\lambda^{(i)}$, and this
component is equal to $\lambda^{(i)}=(1^{k})$, for some $k\in\N$.
In this case, $\1_{\blam}$ has alternating
solid and ghost strings. For example,
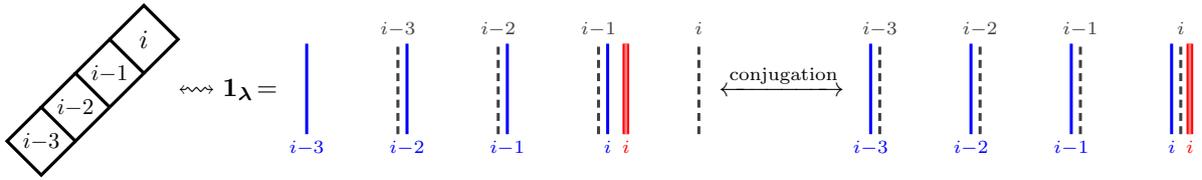
\begin{gather*}
\begin{tikzpicture}[scale=1.2,anchorbase,scale=0.75]
\draw[very thick] (0.5,1.5) to (1,2) to (0.5,2.5) to (0,2);
\draw[very thick] (0,1) to (0.5,1.5) to (0,2) to (-0.5,1.5) to (0,1);
\draw[very thick] (-0.5,1.5) to (-1,1) to (-0.5,0.5) to (0,1);
\draw[very thick] (-1,1) to (-1.5,0.5) to (-1,0) to (0-0.5,0.5);
\node at (0.5,2){$i$};
\node at (0,1.5){\scalebox{0.85}{$i{-}1$}};
\node at (-0.5,1){\scalebox{0.85}{$i{-}2$}};
\node at (-1,0.5){\scalebox{0.85}{$i{-}3$}};
\end{tikzpicture}
\!\!\leftrightsquigarrow
\1_{\blam}\!=\!
\begin{tikzpicture}[scale=1.2,anchorbase,smallnodes,rounded corners]
\draw[ghost](-2.1,0)node[below]{$\phantom{i}$}--++(0,1)node[above,yshift=-1pt]{$i{-}3$};
\draw[ghost](-1,0)node[below]{$\phantom{i}$}--++(0,1)node[above,yshift=-1pt]{$i{-}2$};
\draw[ghost](0.1,0)node[below]{$\phantom{i}$}--++(0,1)node[above,yshift=-1pt]{$i{-}1$};
\draw[ghost](1.2,0)node[below]{$\phantom{i}$}--++(0,1)node[above,yshift=-1pt]{$i$};
\draw[solid](-3.1,0)node[below]{$i{-}3$}--++(0,1)node[above,yshift=-1pt]{$\phantom{i}$};
\draw[solid](-2,0)node[below]{$i{-}2$}--++(0,1)node[above,yshift=-1pt]{$\phantom{i}$};
\draw[solid](-0.9,0)node[below]{$i{-}1$}--++(0,1)node[above,yshift=-1pt]{$\phantom{i}$};
\draw[solid](0.2,0)node[below]{$i$}--++(0,1)node[above,yshift=-1pt]{$\phantom{i}$};
\draw[redstring](0.4,0)node[below]{$i$}--++(0,1)node[above,yshift=-1pt]{$\phantom{i}$};
\end{tikzpicture}
\xleftrightarrow{\text{conjugation}}\hspace*{-0.1cm}
\begin{tikzpicture}[scale=1.2,anchorbase,smallnodes,rounded corners]
\draw[ghost](-3.0,0)node[below]{$\phantom{i}$}--++(0,1)node[above,yshift=-1pt]{$i{-}3$};
\draw[ghost](-1.9,0)node[below]{$\phantom{i}$}--++(0,1)node[above,yshift=-1pt]{$i{-}2$};
\draw[ghost](-0.8,0)node[below]{$\phantom{i}$}--++(0,1)node[above,yshift=-1pt]{$i{-}1$};
\draw[ghost](0.3,0)node[below]{$\phantom{i}$}--++(0,1)node[above,yshift=-1pt]{$i$};
\draw[solid](-3.1,0)node[below]{$i{-}3$}--++(0,1)node[above,yshift=-1pt]{$\phantom{i}$};
\draw[solid](-2,0)node[below]{$i{-}2$}--++(0,1)node[above,yshift=-1pt]{$\phantom{i}$};
\draw[solid](-0.9,0)node[below]{$i{-}1$}--++(0,1)node[above,yshift=-1pt]{$\phantom{i}$};
\draw[solid](0.2,0)node[below]{$i$}--++(0,1)node[above,yshift=-1pt]{$\phantom{i}$};
\draw[redstring](0.4,0)node[below]{$i$}--++(0,1)node[above,yshift=-1pt]{$\phantom{i}$};
\end{tikzpicture}
\end{gather*}
Crossing the ghost $i$-string and the red $i$-string is an invertible degree zero diagram. By adding straight lines for the remaining strings, conjugating by this diagram gives an isomorphisms between two seemingly different weighted KRLW algebras. Therefore, the subalgebra spanned by diagrams with the endpoints on the left-hand side is isomorphic to the weighted KLRW algebra with the endpoints on the right-hand side. Under this isomorphism the diagrams on the left-hand side are sent to infinitesimal diagrams on the right-hand side.
\end{Example}

We can now define the set of $\hell$-partitions that we need to index our basis,
as well as the set of endpoints $X$ that we will consider.

\begin{Definition}\label{D:ThePartitions}
Let $\hParts[0]=\set{{(\emptyset|\dots|\emptyset)}}$, the set containing only
the empty $\hell$-partition, and for $n\geq 1$
define $\hParts$ recursively by the condition that
$\blam\in\hParts$ only if $\blam=\bmu\cup\alpha$,
where $\bmu\in\hParts[n-1]$ and $\alpha$ is an addable $i$-node of
$\bmu$ such that:
\begin{gather}\label{E:Condition}
\text{whenever $\beta$ is an addable $i$-node of $\bmu$
with $\hcoord(\beta)<\hcoord(\alpha)$ then
$\hcoord(\beta)\leq\hcoord(\ell,1,n)$}.
\end{gather}
Finally, set
\begin{gather}\label{E:XSet}
X=\bigcup_{\blam\in\hParts}\hcoord(\blam).
\end{gather}
\end{Definition}

The condition \autoref{E:Condition} on $\alpha$ is
satisfied if $\alpha$ is added to $\lambda^{(m)}$ for $1\leq m\leq\ell$ since
$\hcoord(\ell,1,n)$ is the maximal position an $\ell$-partition of $n$ could have.
In particular, $\Parts\subset\hParts$.
Moreover, by \autoref{L:RightLeft}, an affine node
$\alpha$ with $\hcoord(\ell,1,n)<\hcoord(\alpha)$ is only added
to $\hParts$ at the leftmost available slot.

\begin{Example}\label{Ex:AffineCompare}
Let $n=2$, $\ell=0$, $\charge=\emptyset$ and $e=1$. In this case
we have $\hell=4$, so that $\Parts(\hell)$ is the set of $4$-partitions of $2$,
and has thus has $14$ elements. We also have $\affine{\charge}=(4,8,12,16)$
and $\affine{\brho}=(0,0,1,1)$ and $\hcoord(\ell,1,n)=1$.
Next, we note that all possible nodes are affine,
so \autoref{E:Condition} is always relevant. Then
\begin{align*}
\hParts[0]&=\set[\Big]{{(\emptyset|\dots|\emptyset)}},
\\
\hParts[1]&=\set[\Big]{{\big(1|\emptyset|\emptyset|\emptyset\big),\big(\emptyset|\emptyset|1|\emptyset\big)}},
\\
\hParts[2]&=\set[\Big]{{\big(1|1|\emptyset|\emptyset\big),
\big(1^{2}|\emptyset|\emptyset|\emptyset\big),
\big(\emptyset|\emptyset|1|1\big),\big(1|\emptyset|1|\emptyset\big)}}.
\end{align*}
The four $4$-partitions in $\hParts[2]$ are obtained by
adding leftmost nodes.
We stress that, since the coordinate of the
second node in $\big(1,1|\emptyset|\emptyset|\emptyset\big)$ is
approximately $3-\tfrac{1}{4}$, which is bigger than $\hcoord(\ell,1,n)=1$,
we do not have $\blam=\big(2|\emptyset|\emptyset|\emptyset\big)$ in $\hParts[2]$.
The idempotent diagram $\1_{\blam}$ is
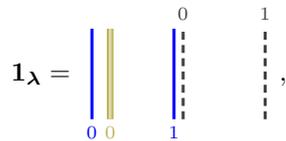
\begin{gather*}
\1_{\blam}
=
\begin{tikzpicture}[scale=1.2,anchorbase,smallnodes,rounded corners]
\draw[ghost](1,0)node[below]{$\phantom{i}$}--++(0,1)node[above,yshift=-1pt]{$0$};
\draw[ghost](1.9,0)node[below]{$\phantom{i}$}--++(0,1)node[above,yshift=-1pt]{$1$};
\draw[solid](0,0)node[below]{$0$}--++(0,1)node[above,yshift=-1pt]{$\phantom{i}$};
\draw[solid](0.9,0)node[below]{$1$}--++(0,1)node[above,yshift=-1pt]{$\phantom{i}$};
\draw[affine](0.2,0)node[below]{$0$}--++(0,1)node[above,yshift=-1pt]{$\phantom{i}$};
\end{tikzpicture}
,
\end{gather*}
while the idempotent diagrams for the elements of $\hParts[2]$
are in the infinitesimal case.

By \cite{Ro-2-kac-moody} and \cite{BrKl-hecke-klr}, the
(affine) KLR algebra of affine type $A$ is isomorphic to the affine Hecke algebra
of the same type, which has simple modules indexed by the set of \emph{aperiodic multisegments}
by \cite[Theorem B]{ArMa-simples-complex-reflection}. In the example above,
$\set[\big]{{[0|0],[0,1],[1|1],[1,0]}}$ is the set of aperiodic multisegments
(see \cite[before Theorem B]{ArMa-simples-complex-reflection} for the definitions).
Under the isomorphism of \autoref{P:WebAlg}, the diagram $\1_{\blam}$ does not
appear in the idempotent subalgebra for the (affine) KLR of this type, which explains
why we do not want $\big(2|\emptyset|\emptyset|\emptyset\big)$ in $\hParts[2]$.
\end{Example}

Our definition of semistandard tableaux is
the same as in \cite[Definition 2.11]{We-rouquier-dia-algebra} or \cite[Section 1.3]{Bo-many-cellular-structures}:

\begin{Definition}
\label{D:Semistandard}
Let $\blam,\bmu\in\hParts$. A \emph{$\blam$-tableau of
type $\bmu$} is a bijection $\bT\map\blam\hcoord(\bmu)$.
Such a tableau is \emph{semistandard} if:
\begin{enumerate}

\item We have $\bT(m,1,1)\leq\kappa_{m}$ for $1\leq m\leq\ell$.

\item We have $\bT(m,r,c)+1<\bT(m,r-1,c)$ for all $(m,r,c),(m,r-1,c)\in\blam$.

\item We have $\bT(m,r,c)<\bT(m,r,c-1)+1$ for all
$(m,r,c),(m,r,c-1)\in\blam$.

\end{enumerate}
Let $\hSStd(\blam,\bmu)$ be the set of semistandard $\blam$-tableaux of
type $\bmu$ and set $\hSStd(\blam)=\bigcup_{\bmu}\hSStd(\blam,\bmu)$.
\end{Definition}

As usual, we think of tableaux as fillings of the nodes:

\begin{Example}\label{Ex:BigPositioning2}
In the setup from \autoref{Ex:BigPositioning},
three $\blam$-tableaux of type $\lambda$ are
\begin{gather*}
\bT=
\begin{tikzpicture}[scale=1.2,anchorbase,scale=0.75]
\draw[very thick] (-0.5,-1.5) to (-1,-2) to (-0.5,-2.5) to (0,-2);
\draw[very thick] (0,-1) to (-0.5,-1.5) to (0,-2) to (0.5,-1.5) to (0,-1);
\draw[very thick] (0.5,-1.5) to (1,-2) to (0.5,-2.5) to (0,-2);
\draw[very thick] (1,-2) to (1.5,-2.5) to (1,-3) to (0.5,-2.5);
\draw[very thick] (-0.5,-2.5) to (0,-3) to (0.5,-2.5);
\draw[very thick] (2.5,-1.5) to (2,-2) to (2.5,-2.5) to (3,-2);
\draw[very thick] (3,-1) to (2.5,-1.5) to (3,-2) to (3.5,-1.5) to (3,-1);
\node at (0,-1.5){\scalebox{0.7}{$0.94$}};
\node at (0,-2.5){\scalebox{0.7}{$0.92$}};
\node at (-0.5,-2){\scalebox{0.7}{$\text{-}0.07$}};
\node at (0.5,-2){\scalebox{0.7}{$1.93$}};
\node at (1,-2.5){\scalebox{0.7}{$2.92$}};
\node at (2.5,-2){\scalebox{0.65}{$12.84$}};
\node at (3,-1.5){\scalebox{0.65}{$13.85$}};
\end{tikzpicture}
,
\bT^{\prime}=
\begin{tikzpicture}[scale=1.2,anchorbase,scale=0.75]
\draw[very thick] (-0.5,-1.5) to (-1,-2) to (-0.5,-2.5) to (0,-2);
\draw[very thick] (0,-1) to (-0.5,-1.5) to (0,-2) to (0.5,-1.5) to (0,-1);
\draw[very thick] (0.5,-1.5) to (1,-2) to (0.5,-2.5) to (0,-2);
\draw[very thick] (1,-2) to (1.5,-2.5) to (1,-3) to (0.5,-2.5);
\draw[very thick] (-0.5,-2.5) to (0,-3) to (0.5,-2.5);
\draw[very thick] (2.5,-1.5) to (2,-2) to (2.5,-2.5) to (3,-2);
\draw[very thick] (3,-1) to (2.5,-1.5) to (3,-2) to (3.5,-1.5) to (3,-1);
\node at (0,-1.5){\scalebox{0.7}{$0.92$}};
\node at (0,-2.5){\scalebox{0.7}{$0.94$}};
\node at (-0.5,-2){\scalebox{0.7}{$\text{-}0.07$}};
\node at (0.5,-2){\scalebox{0.7}{$1.93$}};
\node at (1,-2.5){\scalebox{0.7}{$2.92$}};
\node at (2.5,-2){\scalebox{0.65}{$12.84$}};
\node at (3,-1.5){\scalebox{0.65}{$13.85$}};
\end{tikzpicture}
,
\bT^{\prime\prime}=
\begin{tikzpicture}[scale=1.2,anchorbase,scale=0.75]
\draw[very thick] (-0.5,-1.5) to (-1,-2) to (-0.5,-2.5) to (0,-2);
\draw[very thick] (0,-1) to (-0.5,-1.5) to (0,-2) to (0.5,-1.5) to (0,-1);
\draw[very thick] (0.5,-1.5) to (1,-2) to (0.5,-2.5) to (0,-2);
\draw[very thick] (1,-2) to (1.5,-2.5) to (1,-3) to (0.5,-2.5);
\draw[very thick] (-0.5,-2.5) to (0,-3) to (0.5,-2.5);
\draw[very thick] (2.5,-1.5) to (2,-2) to (2.5,-2.5) to (3,-2);
\draw[very thick] (3,-1) to (2.5,-1.5) to (3,-2) to (3.5,-1.5) to (3,-1);
\node at (0,-1.5){\scalebox{0.7}{$0.94$}};
\node at (0,-2.5){\scalebox{0.7}{$0.92$}};
\node at (-0.5,-2){\scalebox{0.7}{$2.92$}};
\node at (0.5,-2){\scalebox{0.7}{$1.93$}};
\node at (1,-2.5){\scalebox{0.7}{$\text{-}0.07$}};
\node at (2.5,-2){\scalebox{0.65}{$12.84$}};
\node at (3,-1.5){\scalebox{0.65}{$13.85$}};
\end{tikzpicture}
.
\end{gather*}
Note that, for example, $\bT(1,2,1)+1\approx 0.93<\bT(1,1,1)\approx 0.94$
and $\bT(1,1,2)\approx 1.93<\bT(1,1,1)+1\approx 1.94$. Similarly, one checks
that $\bT$ is a semistandard $\blam$-tableau of shape $\lambda$.
In contrast, $\bT^{\prime}$ is not semistandard because, for example,
$\bT^{\prime}(1,2,1)+1\approx 0.93>\bT^{\prime}(1,1,1)\approx 0.92$. The rightmost
tableau $\bT^{\prime\prime}$ is also not semistandard as
$\bT^{\prime\prime}(1,2,1)+1\approx 3.92>\bT^{\prime\prime}(1,1,1)\approx 0.94$.

Take now the $23$-partition of $7$ given by $\bmu=(3,2,1^{2}|\emptyset|\dots|\emptyset)$, which has
approximate coordinates $\hcoord(\bmu)\approx
\set{-2.09,-1.08,-0.07,0.92,\framebox{$0.94$},1.93,2.92}$. Consider
the $\blam$-tableaux of type $\bmu$ given by
\begin{gather*}
\bS=
\begin{tikzpicture}[scale=1.2,anchorbase,scale=0.75]
\draw[very thick] (-0.5,-1.5) to (-1,-2) to (-0.5,-2.5) to (0,-2);
\draw[very thick] (0,-1) to (-0.5,-1.5) to (0,-2) to (0.5,-1.5) to (0,-1);
\draw[very thick] (0.5,-1.5) to (1,-2) to (0.5,-2.5) to (0,-2);
\draw[very thick] (1,-2) to (1.5,-2.5) to (1,-3) to (0.5,-2.5);
\draw[very thick] (-0.5,-2.5) to (0,-3) to (0.5,-2.5);
\draw[very thick] (2.5,-1.5) to (2,-2) to (2.5,-2.5) to (3,-2);
\draw[very thick] (3,-1) to (2.5,-1.5) to (3,-2) to (3.5,-1.5) to (3,-1);
\node at (0,-1.5){\scalebox{0.7}{$0.94$}};
\node at (0,-2.5){\scalebox{0.7}{$0.92$}};
\node at (-0.5,-2){\scalebox{0.7}{$\text{-}0.07$}};
\node at (0.5,-2){\scalebox{0.7}{$1.93$}};
\node at (1,-2.5){\scalebox{0.7}{$2.92$}};
\node at (2.5,-2){\scalebox{0.7}{$\text{-}2.09$}};
\node at (3,-1.5){\scalebox{0.7}{$\text{-}1.08$}};
\end{tikzpicture}
,
\bS^{\prime}=
\begin{tikzpicture}[scale=1.2,anchorbase,scale=0.75]
\draw[very thick] (-0.5,-1.5) to (-1,-2) to (-0.5,-2.5) to (0,-2);
\draw[very thick] (0,-1) to (-0.5,-1.5) to (0,-2) to (0.5,-1.5) to (0,-1);
\draw[very thick] (0.5,-1.5) to (1,-2) to (0.5,-2.5) to (0,-2);
\draw[very thick] (1,-2) to (1.5,-2.5) to (1,-3) to (0.5,-2.5);
\draw[very thick] (-0.5,-2.5) to (0,-3) to (0.5,-2.5);
\draw[very thick] (2.5,-1.5) to (2,-2) to (2.5,-2.5) to (3,-2);
\draw[very thick] (3,-1) to (2.5,-1.5) to (3,-2) to (3.5,-1.5) to (3,-1);
\node at (0,-1.5){\scalebox{0.7}{$\text{-}1.08$}};
\node at (0,-2.5){\scalebox{0.7}{$0.92$}};
\node at (-0.5,-2){\scalebox{0.7}{$\text{-}2.09$}};
\node at (0.5,-2){\scalebox{0.7}{$1.93$}};
\node at (1,-2.5){\scalebox{0.7}{$2.92$}};
\node at (2.5,-2){\scalebox{0.7}{$\text{-}0.07$}};
\node at (3,-1.5){\scalebox{0.7}{$0.94$}};
\end{tikzpicture}
,
\bS^{\prime\prime}=
\begin{tikzpicture}[scale=1.2,anchorbase,scale=0.75]
\draw[very thick] (-0.5,-1.5) to (-1,-2) to (-0.5,-2.5) to (0,-2);
\draw[very thick] (0,-1) to (-0.5,-1.5) to (0,-2) to (0.5,-1.5) to (0,-1);
\draw[very thick] (0.5,-1.5) to (1,-2) to (0.5,-2.5) to (0,-2);
\draw[very thick] (1,-2) to (1.5,-2.5) to (1,-3) to (0.5,-2.5);
\draw[very thick] (-0.5,-2.5) to (0,-3) to (0.5,-2.5);
\draw[very thick] (2.5,-1.5) to (2,-2) to (2.5,-2.5) to (3,-2);
\draw[very thick] (3,-1) to (2.5,-1.5) to (3,-2) to (3.5,-1.5) to (3,-1);
\node at (0,-1.5){\scalebox{0.7}{$\text{-}1.08$}};
\node at (0,-2.5){\scalebox{0.7}{$0.92$}};
\node at (-0.5,-2){\scalebox{0.7}{$\text{-}2.09$}};
\node at (0.5,-2){\scalebox{0.7}{$0.94$}};
\node at (1,-2.5){\scalebox{0.7}{$\text{-}0.07$}};
\node at (2.5,-2){\scalebox{0.7}{$1.93$}};
\node at (3,-1.5){\scalebox{0.7}{$2.92$}};
\end{tikzpicture}
.
\end{gather*}
Of these only $\bS$ is semistandard.
\end{Example}

Recall the map $({}_{-})^{\star}$ from \autoref{E:StarMap} and that we
write $\hcoord(\blam)=\set{x^{\blam}_{1}<\dots<x^{\blam}_{n}}$ in \autoref{D:Acoordinates}.
Then each tableau $\bT$ defines an associated permutation $w_{\bT}$
by tracing the map $\bT\map\blam\hcoord(\bmu)$ minimally. More explicitly,
let $y^{\ba}=y_{1}^{a_{1}}\dots y_{n}^{a_{n}}$ for $\ba=(a_{1},\dots,a_{n})\in\N^{n}$ and define:

\begin{Definition}\label{D:DST}
For $\bT\in\hSStd(\blam,\bmu)$ define the permutation $w_{\bT}\in\Sym$
by requiring that
\begin{gather*}
x^{\bmu}_{w_{\bT}(k)}=\bT(m,r,c)\text{ whenever }
x^{\blam}_{k}=\hcoord(m,r,c),
\quad\text{for $1\leq k\leq n$ and $(m,r,c)\in\blam$}.
\end{gather*}
Let $c=(\hcoord(\blam),\bi^{\blam})$ and $d=(\hcoord(\bmu),\bi^{\bmu})$.
Define
$D_{\bT}=D^{c}_{d}(w_{\bT})$,
a diagram in
$\Web^{b}_{a}$.
Given $\bS\in\hSStd(\blam,\bnu)$ and $\bT\in\hSStd(\blam,\bmu)$ set
\begin{gather*}
D_{\bS\bT}^{\ba}=D_{\bS}^{\star}y^{\ba}\1_{\blam}D_{\bT},\quad\text{ for }
\ba=(a_{1},\dots,a_{n})\in\N^{n}.
\end{gather*}
Set $D_{\bS\bT}=D_{\bS\bT}^{(0,\dots,0)}$.
\end{Definition}

As in \autoref{D:Dw}, the diagrams $D_{\bS\bT}^{\ba}$ are
well-defined only up to the choices of reduced expressions for $w_{\bS}$ and
$w_{\bT}$.

\begin{Example}\label{Ex:T}
Let $\bS$ be the leftmost $\blam$-tableau of type $\bmu$ in \autoref{Ex:BigPositioning2}.
The corresponding diagrams are:
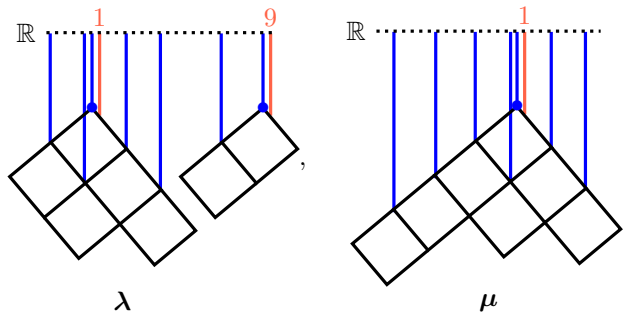
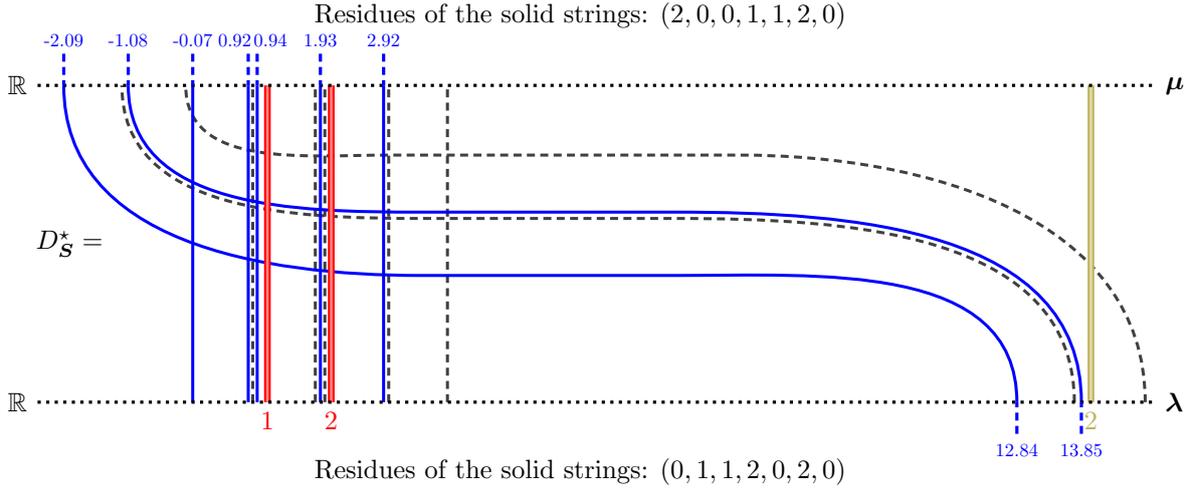
\begin{gather*}
\begin{tikzpicture}[scale=1,anchorbase]
\draw[very thick,blue] (-0.55,-1.45) to (-0.55,-0.01);
\draw[very thick,blue] (0,-1) to (0,-0.01);
\draw[very thick,tomato] (0.1,-1.1) to (0.1,-0.01)node[above,yshift=-1pt]{$1$};
\draw[very thick,blue] (-0.1,-2) to (-0.1,-0.01);
\draw[very thick,blue] (0.45,-1.55) to (0.45,-0.01);
\draw[very thick,blue] (0.9,-2.1) to (0.9,-0.01);
\draw[very thick,blue] (1.7,-1.45) to (1.7,-0.01);
\draw[very thick,blue] (2.25,-1) to (2.25,-0.01);
\draw[very thick,tomato] (2.35,-1.1) to (2.35,-0.01)node[above,yshift=-1pt]{$9$};
\draw[very thick,dotted] (-0.6,0)node[left]{$\R$} to (2.4,0);
\node at (0.4,-3.55) {$\blam$};
\begin{scope}[rotate around={-5:(0,-1)}]
\draw[very thick] (-0.5,-1.5) to (-1,-2) to (-0.5,-2.5) to (0,-2);
\draw[very thick] (0,-1) to (-0.5,-1.5) to (0,-2) to (0.5,-1.5) to (0,-1);
\draw[very thick] (0.5,-1.5) to (1,-2) to (0.5,-2.5) to (0,-2);
\draw[very thick] (1,-2) to (1.5,-2.5) to (1,-3) to (0.5,-2.5);
\draw[very thick] (-0.5,-2.5) to (0,-3) to (0.5,-2.5);
\end{scope}
\begin{scope}[rotate around={-5:(2.25,-1)}]
\draw[very thick] (1.75,-1.5) to (1.25,-2) to (1.75,-2.5) to (2.25,-2);
\draw[very thick] (2.25,-1) to (1.75,-1.5) to (2.25,-2) to (2.75,-1.5) to (2.25,-1);
 \end{scope}
\node[blue] at (0,-1){$\bullet$};
\node[blue] at (2.25,-1){$\bullet$};
\end{tikzpicture}
,\quad
\begin{tikzpicture}[scale=1,anchorbase]
\draw[very thick,blue] (-0.55,-1.45) to (-0.55,-0.01);
\draw[very thick,blue] (0,-1) to (0,-0.01);
\draw[very thick,tomato] (0.1,-1.1) to (0.1,-0.01)node[above,yshift=-1pt]{$1$};
\draw[very thick,blue] (-0.085,-2) to (-0.085,-0.01);
\draw[very thick,blue] (0.45,-1.55) to (0.45,-0.01);
\draw[very thick,blue] (0.9,-2.1) to (0.9,-0.01);
\draw[very thick,blue] (-1.07,-1.95) to (-1.07,-0.01);
\draw[very thick,blue] (-1.62,-2.4) to (-1.62,-0.01);
\begin{scope}[rotate around={-5:(0,-1)}]
\draw[very thick] (-0.5,-1.5) to (-1,-2) to (-0.5,-2.5) to (0,-2);
\draw[very thick] (0,-1) to (-0.5,-1.5) to (0,-2) to (0.5,-1.5) to (0,-1);
\draw[very thick] (0.5,-1.5) to (1,-2) to (0.5,-2.5) to (0,-2);
\draw[very thick] (1,-2) to (1.5,-2.5) to (1,-3) to (0.5,-2.5);
\draw[very thick] (-0.5,-2.5) to (0,-3) to (0.5,-2.5);
\draw[very thick] (-1,-2) to (-1.5,-2.5) to (-1,-3) to (-0.5,-2.5);
\draw[very thick] (-1.5,-2.5) to (-2,-3) to (-1.5,-3.5) to (-1,-3);
\end{scope}
\draw[very thick,dotted] (-1.85,0)node[left]{$\R$} to (1.1,0);
\node at (-0.375,-3.6) {$\bmu$};
\node[blue] at (0,-1){$\bullet$};
\end{tikzpicture}
\\
\begin{tikzpicture}[scale=0.84,anchorbase]
\draw[ghost](0.87,0)to[out=90,in=270](0.87,5);
\draw[ghost](1.85,0)to[out=90,in=270](1.85,5);
\draw[ghost](2,0)to[out=90,in=270](2,5);
\draw[ghost](3,0)to[out=90,in=270](3,5);
\draw[ghost](3.92,0)to[out=90,in=270](3.92,5);
\draw[ghost](13.74,0)to[out=90,in=0](7.5,2.9)to(3.5,2.9)to[out=180,in=270](-1.18,5);
\draw[ghost](14.85,0)to[out=90,in=0](8.5,3.9)to(3.0,3.9)to[out=180,in=270](-0.18,5);
\draw[solid](-0.07,0)to[out=90,in=270](-0.07,5);
\draw[solid](0.8,0)to[out=90,in=270](0.8,5);
\draw[solid](0.94,0)to[out=90,in=270](0.94,5);
\draw[solid](1.93,0)to[out=90,in=270](1.93,5);
\draw[solid](2.92,0)to[out=90,in=270](2.92,5);
\draw[solid](12.84,0)to[out=90,in=0](7.5,2)to(3.5,2)to[out=180,in=270](-2.09,5);
\draw[solid](13.85,0)to[out=90,in=0](7.5,3)to(3.5,3)to[out=180,in=270](-1.08,5);
\draw[redstring](1.1,5)node[above,yshift=-1pt]{$\phantom{i}$}to(1.1,0)node[below]{$1$};	\draw[redstring](2.1,5)node[above,yshift=-1pt]{$\phantom{i}$}to(2.1,0)node[below]{$2$};
\draw[affine](14,5)node[above,yshift=-1pt]{$\phantom{i}$}to(14,0)node[below]{$2$};
\draw[very thick,dotted] (-2.5,0)node[left]{$\R$} to (15,0)node[right]{$\blam$};
\draw[very thick,dotted] (-2.5,5)node[left]{$\R$} to (15,5)node[right]{$\bmu$};
\node at (-2,2.5){$D_{\bS}^{\star}=$};
\draw[very thick,blue,densely dashed] (-2.09,5.5)node[above,yshift=-1pt]{\scalebox{0.7}{$\text{-}2.09$}} to (-2.09,5);
\draw[very thick,blue,densely dashed] (-1.08,5.5)node[above,yshift=-1pt]{\scalebox{0.7}{$\text{-}1.08$}} to (-1.08,5);
\draw[very thick,blue,densely dashed] (-0.07,5.5)node[above,yshift=-1pt]{\scalebox{0.7}{$\text{-}0.07$}} to (-0.07,5);
\draw[very thick,blue,densely dashed] (0.8,5.5)node[above,yshift=-1pt,xshift=-5pt]{\scalebox{0.7}{$0.92$}} to (0.8,5);
\draw[very thick,blue,densely dashed] (0.94,5.5)node[above,yshift=-1pt,xshift=5pt]{\scalebox{0.7}{$0.94$}} to (0.94,5);
\draw[very thick,blue,densely dashed] (1.93,5.5)node[above,yshift=-1pt]{\scalebox{0.7}{$1.93$}} to (1.93,5);
\draw[very thick,blue,densely dashed] (2.92,5.5)node[above,yshift=-1pt]{\scalebox{0.7}{$2.92$}} to (2.92,5);
\draw[very thick,blue,densely dashed] (12.84,-0.5)node[below]{\scalebox{0.7}{$12.84$}} to (12.84,0);
\draw[very thick,blue,densely dashed] (13.85,-0.5)node[below]{\scalebox{0.7}{$13.85$}} to (13.85,0);
\node at (6,-1.1){Residues of the solid strings: $(0,1,1,2,0,2,0)$};
\node at (6,6.1){Residues of the solid strings: $(2,0,0,1,1,2,0)$};
\end{tikzpicture}
\end{gather*}
That is, the diagram $D_{\bS}^{\star}$ traces out the bijection $\bS\map\blam\hcoord(\bmu)$ with a minimal number of crossings. Consequently,
in the setup from \autoref{Ex:BigPositioning}
we obtain
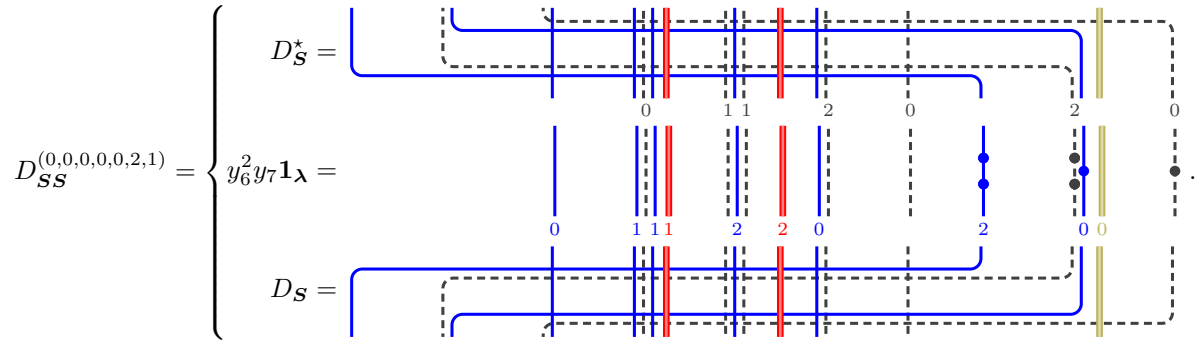
\begin{gather*}
D_{\bS\bS}^{(0,0,0,0,0,2,1)}
=
\left\{
\begin{aligned}
D_{\bS}^{\star}=
&\hspace*{0.17cm}\begin{tikzpicture}[scale=1.2,anchorbase,smallnodes,rounded corners]
\draw[ghost](0.2,0)--++(0,1);
\draw[ghost](1.1,0)--++(0,1);
\draw[ghost](1.3,0)--++(0,1);
\draw[ghost](2.2,0)--++(0,1);
\draw[ghost](3.1,0)--++(0,1);
\draw[ghost](4.9,0)--++(0,0.35)--++(-6.9,0)--++(0,0.65);
\draw[ghost](6,0)--++(0,0.85)--++(-6.9,0)--++(0,0.15);
\draw[solid](-0.8,0)--++(0,1);
\draw[solid](0.1,0)--++(0,1);
\draw[solid](0.3,0)--++(0,1);
\draw[solid](1.2,0)--++(0,1);
\draw[solid](2.1,0)--++(0,1);
\draw[solid](3.9,0)--++(0,0.25)--++(-6.9,0)--++(0,0.75);
\draw[solid](5,0)--++(0,0.75)--++(-6.9,0)--++(0,0.25);
\draw[redstring](0.45,0)--++(0,1);
\draw[redstring](1.7,0)--++(0,1);
\draw[affine](5.2,0)--++(0,1);
\end{tikzpicture}
\\[-5pt]
y_{6}^{2}y_{7}\1_{\blam}=
&\begin{tikzpicture}[scale=1.2,anchorbase,smallnodes,rounded corners]
\draw[white](3.72,0)--++(0,0.25)--++(-6.9,0)--++(0,0.75);
\draw[ghost](0.2,0)node[below]{$\phantom{i}$}--++(0,1)node[above,yshift=-1pt]{$0$};
\draw[ghost](1.1,0)node[below]{$\phantom{i}$}--++(0,1)node[above,yshift=-1pt]{$1$};
\draw[ghost](1.3,0)node[below]{$\phantom{i}$}--++(0,1)node[above,yshift=-1pt]{$1$};
\draw[ghost](2.2,0)node[below]{$\phantom{i}$}--++(0,1)node[above,yshift=-1pt]{$2$};
\draw[ghost](3.1,0)node[below]{$\phantom{i}$}--++(0,1)node[above,yshift=-1pt]{$0$};
\draw[ghost,dot=0.35,dot=0.65](4.9,0)node[below]{$\phantom{i}$}--++(0,1)node[above,yshift=-1pt]{$2$};
\draw[ghost,dot](6,0)node[below]{$\phantom{i}$}--++(0,1)node[above,yshift=-1pt]{$0$};
\draw[solid](-0.8,0)node[below]{$0$}--++(0,1)node[above,yshift=-1pt]{$\phantom{i}$};
\draw[solid](0.1,0)node[below]{$1$}--++(0,1)node[above,yshift=-1pt]{$\phantom{i}$};
\draw[solid](0.3,0)node[below]{$1$}--++(0,1)node[above,yshift=-1pt]{$\phantom{i}$};
\draw[solid](1.2,0)node[below]{$2$}--++(0,1)node[above,yshift=-1pt]{$\phantom{i}$};
\draw[solid](2.1,0)node[below]{$0$}--++(0,1)node[above,yshift=-1pt]{$\phantom{i}$};
\draw[solid,dot=0.35,dot=0.65](3.9,0)node[below]{$2$}--++(0,1)node[above,yshift=-1pt]{$\phantom{i}$};
\draw[solid,dot](5,0)node[below]{$0$}--++(0,1)node[above,yshift=-1pt]{$\phantom{i}$};
\draw[redstring](0.45,0)node[below]{$1$}--++(0,1)node[above,yshift=-1pt]{$\phantom{i}$};
\draw[redstring](1.7,0)node[below]{$2$}--++(0,1)node[above,yshift=-1pt]{$\phantom{i}$};
\draw[affine](5.2,0)node[below]{$0$}--++(0,1)node[above,yshift=-1pt]{$\phantom{i}$};
\end{tikzpicture}
\\[-5pt]
D_{\bS}=
&\hspace*{0.17cm}\begin{tikzpicture}[scale=1.2,anchorbase,smallnodes,rounded corners,yscale=-1]
\draw[ghost](0.2,0)--++(0,1);
\draw[ghost](1.1,0)--++(0,1);
\draw[ghost](1.3,0)--++(0,1);
\draw[ghost](2.2,0)--++(0,1);
\draw[ghost](3.1,0)--++(0,1);
\draw[ghost](4.9,0)--++(0,0.35)--++(-6.9,0)--++(0,0.65);
\draw[ghost](6,0)--++(0,0.85)--++(-6.9,0)--++(0,0.15);
\draw[solid](-0.8,0)--++(0,1);
\draw[solid](0.1,0)--++(0,1);
\draw[solid](0.3,0)--++(0,1);
\draw[solid](1.2,0)--++(0,1);
\draw[solid](2.1,0)--++(0,1);
\draw[solid](3.9,0)--++(0,0.25)--++(-6.9,0)--++(0,0.75);
\draw[solid](5,0)--++(0,0.75)--++(-6.9,0)--++(0,0.25);
\draw[redstring](0.45,0)--++(0,1);
\draw[redstring](1.7,0)--++(0,1);
\draw[affine](5.2,0)--++(0,1);
\end{tikzpicture}
\end{aligned}
.
\right.
\end{gather*}
(We stack these together to obtain $D_{\bS\bS}^{(0,0,0,0,0,2,1)}$.)
Note that the strings in $\1_{\blam}$ are as far to the right as possible, with the affine red string keeping the rightmost solid $0$-string in check.
\end{Example}


\subsection{The affine cellular basis in type $A$}\label{SS:AAffineCellular}


Retain the conventions from the precious section.
We construct an affine cellular basis of the weighted
KLRW algebra $\WA[n](X)$. To this end, we choose the $Q$-polynomials as in \autoref{E:QPoly}. Recall that $\hcoord(\ell,1,n)$
is the maximal position that an $\ell$-partition of $n$ is allowed to have. Define the set
\begin{gather*}
\Affch=\set{\ba=(a_{1},\dots,a_{n})\in\N^{n}|a_{k}=0\text{ whenever }
\hcoord(\blam)_{k}\leq\hcoord(\ell,1,n)},
\end{gather*}
which will index the possible exponents in $D^{\ba}_{\bS\bT}$.
We consider the set
\begin{gather}\label{E:TBasis}
\BX=
\set[\big]{D^{\ba}_{\bS\bT}|\blam\in\hParts,\bS,\bT\in\hSStd(\blam),\ba\in\Affch}.
\end{gather}

\begin{Example}
Note that, if $\ba\in\Affch$, then $a_{k}>0$ only
if $\hcoord(\blam)_{k}$ is to the right of all possible coordinates for
$\ell$-partitions.
In particular, if $\bS$ is as in \autoref{Ex:T}, then $D_{\bS\bS}^{(0,0,0,0,0,2,1)}\in\BX$
whereas $D_{\bS\bS}^{(0,0,0,0,1,2,1)}\notin\BX$.
\end{Example}

\begin{Definition}\label{D:Dominance}
Let $\blam,\bmu\in\hParts$. Then $\blam$ \emph{is dominated by} $\bmu$,
written $\blam\ledom_{A}\bmu$, if there exists a bijection $d\colon\blam\to\bmu$ such that
$\hcoord(\alpha)\leq\hcoord\big(d(\alpha)\big)$, for all $\alpha\in\blam$. Write
$\blam\ldom_{A}\bmu$ if $\blam\ledom_{A}\bmu$ and $\blam\ne\bmu$.
\end{Definition}

Equivalently, if $(x_{1}<\dots<x_{n})$ and $(y_{1}<\dots<y_{n})$ are the
$\hcoord$-coordinates of the solid strings in $\1_{\blam}$ and $\1_{\bmu}$,
respectively, then $\blam\ledom_{A}\bmu$ if and only if
$x_{r}\leq y_{r}$ for all $r=1,\dots,n$. Thus, $\blam\in\hParts$ becomes more
dominant when we move strings, or nodes, to the right.

Note that the above equivalent formulation works for any idempotent diagram,
not just for those of the form $\1_{\blam}$, and we will also use it
in this more general setting.

\begin{Remark}
Note that \autoref{D:Dominance} is different from the ordering used in
\cite[Section 1.1]{Bo-many-cellular-structures}.
It is not clear to us whether the
results below hold using Bowman's ordering.
\end{Remark}

If $S=\1_{\bx,\bi}$ is an idempotent diagram, let
$\hcoord(S)=(x^{S}_{1}<\dots<x^{S}_{n})$ be the corresponding increasing sequence
of coordinates. A diagram $D$ \emph{factors through} $S$ if $D=D^{\prime}SD^{\prime\prime}$,
for some $D^{\prime},D^{\prime\prime}\in\WA[n](X)$.
If $\blam\in\hParts$, then $\blam$ \emph{is dominated by} $S$ if $\hcoord(\blam)\ldom_{A}\hcoord(S)$ and, similarly,
$S$ \emph{is dominated by} $\blam$ if $\hcoord(S)\ldom_{A}\hcoord(\blam)$.

As we will see, if $X$ is as in
\autoref{E:XSet}, then $\BX$ is a
homogeneous affine cellular basis of $\WA[n](X)$ with respect to the dominance
order and with antiinvolution $({}_{-}){}^{\star}$; {\cf} \autoref{D:Grading}. (Everything below is with respect to the choice of
these structures; the dominance order is very important and
sometimes highlighted.)

It is time to recall the definition of a graded affine cellular algebra
from \cite{KoXi-affine-cellular}, incorporating a grading
as in \cite{HuMa-klr-basis}. We slightly adapt \cite{KoXi-affine-cellular}
to our needs by, for example, casting the definition in terms of bases, whereas K{\"o}nig--Xi give a ring theoretic formulation.

\begin{Definition}\label{D:CellularAlgebra}
Let $R$ be a commutative ring with a unit.
Let $A$ be a locally unital graded $R$-algebra.
A \emph{graded affine cell datum} for $A$ is a quintuple $(\Pcal,T,\sand[],C,\deg)$, where:
\begin{itemize}
\item $\Pcal=(\Pcal,\leq)$ is a poset,
\item $T=\bigcup_{\lambda\in\Pcal}T(\lambda)$ is a collection of finite sets,
\item $\sand[]=\bigoplus_{\lambda\in\Pcal}\sand[\lambda]$ is a direct sum of
quotients of polynomial rings such that $B(\lambda)$ is a homogeneous
basis of $\sand[\lambda]$ (we write $\deg$ for the degree function on $\sand[\lambda]$),
\item $C\map{\coprod_{\lambda\in\Pcal}T(\lambda)\times B(\lambda)\times
T(\lambda)}{A};(S,\ba,T)\mapsto C^{\ba}_{ST}$ is an injective map,
\item $\deg\map{\coprod_{\lambda\in\Pcal}T(\lambda)}{\Z}$ is a function,
\end{itemize}
such that:
\begin{enumerate}[label=\upshape(AC${}_{\arabic*}$\upshape)]

\item For $\lambda\in\Pcal$, $S,T\in T(\lambda)$ and $\ba\in B(\lambda)$, $C^{\ba}_{ST}$
is homogeneous of degree $\deg(S)+\deg\ba+\deg(T)$.

\item The set $\set{C^{\ba}_{ST}|\lambda\in\Pcal,S,T\in T(\lambda),\ba\in B(\lambda)}$
is a basis of $A$.

\item For all $x\in A$ there exist scalars $r_{SU}\in R$ that do not depend
on $T$ or on $\ba$, such that
\begin{gather*}
xC^{\ba}_{ST}\equiv
\sum_{U\in T(\lambda)}r_{SU}C^{\ba}_{UT}\pmod{A^{>\lambda}},
\end{gather*}
where $A^{>\lambda}$ is the $R$-submodule of $A$ spanned by
$\set{C^{\bb}_{UV}|\mu\in\Pcal,\mu>\lambda,U,V\in T(\mu),\bb\in B(\mu)}$.

\item Let $A(\lambda)=A^{\geq\lambda}/A^{>\lambda}$, where
$A^{\geq\lambda}$ is the $R$-submodule of $A$ spanned by $\set{C^{\bb}_{UV}|\mu\in\Pcal,
\mu\geq\lambda,U,V\in T(\mu),\bb\in B(\mu)}$. Then $A(\lambda)$ is isomorphic to
$\Delta(\lambda)\otimes_{\sand[\lambda]}\nabla(\lambda)$ for free graded right and left
$\sand[\lambda]$-modules $\Delta(\lambda)$ and $\nabla(\lambda)$, respectively.

\item There is an antiinvolution $({}_{-}){}^{\star}\map{A}{A}$ of~$A$ such that
$(C^{\ba}_{ST})^{\star}\equiv C^{\ba}_{TS}\pmod{A^{>\lambda}}$, for all $S,T\in T(\lambda)$
and $\ba\in B(\lambda)$, for $\lambda\in P$. This antiinvolution identifies $\Delta(\lambda)$
and $\nabla(\lambda)$.

\end{enumerate}
The algebra $A$ is a \emph{graded affine cellular algebra} if it
has a graded affine cell datum and it is an \emph{affine cellular algebra}
if $\deg(S)=0$ for all $S\in T$. The image of $C$ is an \emph{homogeneous affine cellular basis} of $A$.

A graded \emph{cell datum} for $A$ is a graded affine cell datum such
that $\sand[\lambda]\cong R$, for all $\lambda\in\Pcal$. In this case
the image of $C$ is an \emph{homogeneous cellular basis} of $A$.
\end{Definition}

\begin{Remark}
Axiom (AC${}_{5}$) is a slightly weaker assumption than is commonly given in the literature
where it is usually assumed only that the linear map $({}_{-})^{\star}\colon A\to A$ given
by $(C^{\ba}_{ST})^{\star}=C^{\ba}_{TS}\pmod{A^{>\lambda}}$ is an antiinvolution of $A$.
Assumption (AC${}_{5}$) has some advantages in characteristic $2$, see \cite[Remark 2.4]{GoGr-cellularity-jones-basic}.
\end{Remark}

Cellular algebras were introduced by Graham--Lehrer \cite{GrLe-cellular} and
K{\"o}nig--Xi \cite{KoXi-affine-cellular} generalized their definition to the affine
case. It is easy to see that an affine cellular algebra in the sense of \autoref{D:CellularAlgebra}
is an affine cellular algebra in the sense of K{\"o}nig--Xi. For us the picture to keep in mind is
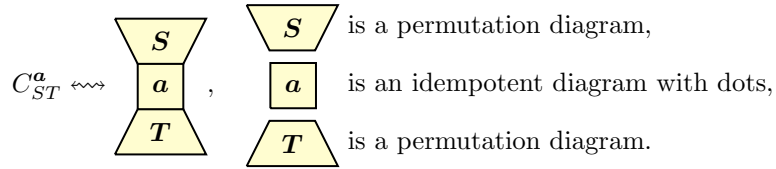
\begin{gather*}
C^{\ba}_{ST}
\leftrightsquigarrow
\begin{tikzpicture}[scale=1.2,anchorbase,scale=1]
\draw[line width=0.75,color=black,fill=cream] (0,-0.5) to (0.25,0) to (0.75,0) to (1,-0.5) to (0,-0.5);
\node at (0.5,-0.25){$\bT$};
\draw[line width=0.75,color=black,fill=cream] (0,1) to (0.25,0.5) to (0.75,0.5) to (1,1) to (0,1);
\node at (0.5,0.75){$\bS$};
\draw[line width=0.75,color=black,fill=cream] (0.25,0) to (0.25,0.5) to (0.75,0.5) to (0.75,0) to (0.25,0);
\node at (0.5,0.25){$\ba$};
\end{tikzpicture},
\quad
\begin{aligned}
\begin{tikzpicture}[scale=1.2,anchorbase,scale=1]
\draw[line width=0.75,color=black,fill=cream] (0,1) to (0.25,0.5) to (0.75,0.5) to (1,1) to (0,1);
\node at (0.5,0.75){$\bS$};
\end{tikzpicture}
&\text{ is a permutation diagram,}
\\
\begin{tikzpicture}[scale=1.2,anchorbase,scale=1]
\draw[white,ultra thin] (0,0) to (1,0);
\draw[line width=0.75,color=black,fill=cream] (0.25,0) to (0.25,0.5) to (0.75,0.5) to (0.75,0) to (0.25,0);
\node at (0.5,0.25){$\ba$};
\end{tikzpicture}
&\text{ is an idempotent diagram with dots,}
\\
\begin{tikzpicture}[scale=1.2,anchorbase,scale=1]
\draw[line width=0.75,color=black,fill=cream] (0,-0.5) to (0.25,0) to (0.75,0) to (1,-0.5) to (0,-0.5);
\node at (0.5,-0.25){$\bT$};
\end{tikzpicture}
&\text{ is a permutation diagram.}
\end{aligned}
\end{gather*}
(To avoid confusion, the cellular basis of the weighted KLRW algebra is denoted by $D^{\ba}_{\bS\bT}$, and not $C^{\ba}_{ST}$.)

Consider $\hParts$ as a poset ordered by $\ledom_{A}$. We stress again that
$\hParts$ is a proper subset of the set of all $\hell$-partitions of $n$ whenever $n>0$. For $D_{\bS\bT}^{\ba}=D_{\bS}^{\star}y^{\ba}\1_{\blam}D_{\bT}$
we define $\deg\bS=\deg D_{\bS}$, $\deg\ba=\deg y^{\ba}\1_{\blam}$
and $\deg\bT=\deg D_{\bT}$.

\begin{Theorem}\label{T:Basis}
The set $\BX$ is a homogeneous affine cellular basis of $\WA[n](X)$ with respect to the poset $(\hParts,\ledom_{A})$.
\end{Theorem}

We defer the proof of \autoref{T:Basis} until
later. The proof itself is not very technical, but
we prefer to first draw the readers attention to
some consequences of \autoref{T:Basis}.

\begin{Corollary}\label{C:TypeAKLR}
The set $\BX[{\mathscr{R}}]=\set[\big]{D_{\bS\bT}|\blam\in\Parts,\bS,\bT\in\SStd(\blam)}$
is a homogeneous cellular basis of
the cyclotomic weighted KLRW algebra $\WAc[n](X)$.
\end{Corollary}

\begin{proof}
By \autoref{T:Basis} it is enough to prove that all of the diagrams in $\BX[{\mathscr{R}}]$ are steady. By way of contradiction, suppose that $D_{\bS\bT}$ is not steady, for some $\bS,\bT\in\SStd(\blam)$ and $\blam\in\Parts$. Then we can pull some strings in the diagram $D_{\bS\bT}$ further to the right and, in particular, past all of the red strings in the diagram. Diagrams with strings to the right of the red strings are more dominant, so \autoref{T:Basis} implies that $D_{\bS\bT}$ can be written as a linear combination of diagrams $D_{\bU\bV}$, for some $\bU,\bV\in\Std(\bmu)$ with $\bmu\in\hParts\setminus\Parts$. This contradicts the linear independence of the basis in \autoref{T:Basis}. Hence, the diagram $D_{\bS\bT}$ is steady, completing the proof.
\end{proof}

A \emph{standard} tableau is a semistandard tableau of type $\bom=(1^{n}|0|\dots|0)$.
Let $\Std(\blam)$ be the set of standard $\blam$-tableaux. Define
\begin{gather}\label{E:1Abom}
\1_{A,n}=\sum_{\bi\in I^{n}}\1_{\hcoord(\bom),\bi}\in\WA[n](X).
\end{gather}
Identify the KLRW algebra $\TA[n]$, and its cyclotomic quotient $\TAc[n]$, with the
idempotent subalgebras $\1_{A,n}\WA[n](X)\1_{A,n}$, and $\1_{A,n}\WAc[n](X)\1_{A,n}$
respectively, using \autoref{P:WebAlg} and \autoref{Ex:Ordered}. We remind the reader that $\TA[n]$ and $\TAc[n]$ include the KLR algebras of \cite{KhLa-cat-quantum-sln-first}, \cite{Ro-2-kac-moody} as special cases.

\begin{Corollary}\label{C:KLRCellular}
The set $\BX[{\mathcal{W}}]{=}\set[\big]{D^{\ba}_{\bs\bt}|\blam\in\hParts,\bs,\bt\in\Std(\blam),\ba\in\Affch}$
is a homogeneous affine cellular basis of
$\TA[n]$ and
$\BX[{\mathcal{R}}]{=}\set[\big]{D_{\bs\bt}|\blam\in\Parts,\bs,\bt\in\Std(\blam)}$
is a homogeneous cellular basis of $\TAc[n]$.
\end{Corollary}

\begin{proof}
This follows by idempotent truncation applied to \autoref{T:Basis} and \autoref{C:TypeAKLR}.
\end{proof}

\begin{Remark}\label{R:RecollectionTensor}
In \cite[Section 4.4]{We-knot-invariants} Webster defined more general (weighted) KLRW algebras where the red strings are decorated by arbitrary dominant weights. As Webster explains, these more general algebras are idempotent truncations of the weighted KLRW algebras considered in this paper, which only have fundamental weights on the red strings. Hence, by for example \cite[Proposition 4.3]{KoXi-cellular}, the results above immediately imply that Webster's more general algebras are cellular. This remark also applies in \autoref{S:TypeC}.
\end{Remark}

We can now read off the ranks of the algebras from
\autoref{C:TypeAKLR} and \autoref{C:KLRCellular}. The non-cyclotomic algebras are free of infinite rank, whereas for the cyclotomic quotients
\begin{gather}\label{E:DimFormulas}
\mathrm{rank}_{R}\big(\WAc[n](X)\big)
=\sum_{\blam\in\hParts}|\SStd(\blam)|^{2}
\quad\text{and}\quad
\mathrm{rank}_{R}\big(\TAc[n]\big)
=\sum_{\blam\in\Parts}|\Std(\blam)|^{2}.
\end{gather}
The latter count matches \cite[Theorem 4.20]{BrKl-graded-decomposition-hecke}, the former matches  \cite[Theorem 6.23]{Bo-many-cellular-structures}. There are similar formulas for the graded ranks and graded dimensions.

\begin{Remark}\label{R:KLRComparision}
\autoref{C:KLRCellular} gives cellular bases for the KLRW algebras of type $A^{(1)}_{e}$
for each choice of loading, which is a choice of positionings for the solid and red strings.
Arguing as in \cite[Proposition 7.3]{Bo-many-cellular-structures}, in the asymptotic case the
cellular basis of \autoref{C:KLRCellular} coincides with the cellular basis constructed in \cite{HuMa-klr-basis} modulo higher terms.
\end{Remark}

\begin{Example}\label{E:Calc}
This example concerns the finite dimensional cyclotomic weighted KLRW algebra of type $A^{(1)}_{2}$ and level $\ell=1$ with $\rho=(0)$. Assume that the red string is placed at $\charge=(0)$.
We have $\hell=10$, and we assume that $\varepsilon=0.001$.

In this case, there are three $1$-partitions, namely
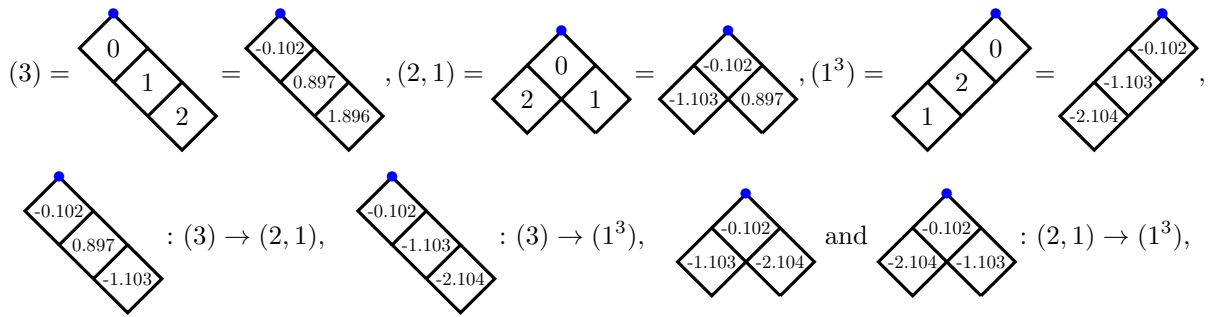
\begin{gather*}
(3)=
\begin{tikzpicture}[scale=1.2,anchorbase,scale=0.75]
\draw[very thick] (0,2) to (-0.5,2.5) to (0,3) to (0.5,2.5);
\draw[very thick] (0.5,1.5) to (1,2) to (0.5,2.5) to (0,2) to (0.5,1.5);
\draw[very thick] (0.5,1.5) to (1,1) to (1.5,1.5) to (1,2);
\node at (0,2.5){$0$};
\node at (0.5,2){$1$};
\node at (1,1.5){$2$};
\node[blue] at (0,3){$\bullet$};
\end{tikzpicture}
\!=\!
\begin{tikzpicture}[scale=1.2,anchorbase,scale=0.75]
\draw[very thick] (0,2) to (-0.5,2.5) to (0,3) to (0.5,2.5);
\draw[very thick] (0.5,1.5) to (1,2) to (0.5,2.5) to (0,2) to (0.5,1.5);
\draw[very thick] (0.5,1.5) to (1,1) to (1.5,1.5) to (1,2);
\node at (0,2.5){\scalebox{0.7}{$\text{-}0.102$}};
\node at (0.5,2){\scalebox{0.7}{$0.897$}};
\node at (1,1.5){\scalebox{0.7}{$1.896$}};
\node[blue] at (0,3){$\bullet$};
\end{tikzpicture}
,
(2,1)=
\begin{tikzpicture}[scale=1.2,anchorbase,scale=0.75]
\draw[very thick] (0,2) to (-0.5,1.5) to (-1,2) to (-0.5,2.5) to (0,2);
\draw[very thick] (0.5,1.5) to (1,2) to (0.5,2.5) to (0,2) to (0.5,1.5);
\draw[very thick] (-0.5,2.5) to (0,3) to (0.5,2.5);
\node at (0,2.5){$0$};
\node at (-0.5,2){$2$};
\node at (0.5,2){$1$};
\node[blue] at (0,3){$\bullet$};
\end{tikzpicture}
\!=\!
\begin{tikzpicture}[scale=1.2,anchorbase,scale=0.75]
\draw[very thick] (0,2) to (-0.5,1.5) to (-1,2) to (-0.5,2.5) to (0,2);
\draw[very thick] (0.5,1.5) to (1,2) to (0.5,2.5) to (0,2) to (0.5,1.5);
\draw[very thick] (-0.5,2.5) to (0,3) to (0.5,2.5);
\node at (0,2.5){\scalebox{0.7}{$\text{-}0.102$}};
\node at (-0.5,2){\scalebox{0.7}{$\text{-}1.103$}};
\node at (0.5,2){\scalebox{0.7}{$0.897$}};
\node[blue] at (0,3){$\bullet$};
\end{tikzpicture}
,
(1^{3})=
\begin{tikzpicture}[scale=1.2,anchorbase,scale=0.75]
\draw[very thick] (0.5,1.5) to (1,2) to (0.5,2.5) to (0,2);
\draw[very thick] (0,1) to (0.5,1.5) to (0,2) to (-0.5,1.5) to (0,1);
\draw[very thick] (-0.5,1.5) to (-1,1) to (-0.5,0.5) to (0,1);
\node at (0.5,2){$0$};
\node at (0,1.5){$2$};
\node at (-0.5,1){$1$};
\node[blue] at (0.5,2.5){$\bullet$};
\end{tikzpicture}
\!=\!
\begin{tikzpicture}[scale=1.2,anchorbase,scale=0.75]
\draw[very thick] (0.5,1.5) to (1,2) to (0.5,2.5) to (0,2);
\draw[very thick] (0,1) to (0.5,1.5) to (0,2) to (-0.5,1.5) to (0,1);
\draw[very thick] (-0.5,1.5) to (-1,1) to (-0.5,0.5) to (0,1);
\node at (0.5,2){\scalebox{0.7}{$\text{-}0.102$}};
\node at (0,1.5){\scalebox{0.7}{$\text{-}1.103$}};
\node at (-0.5,1){\scalebox{0.7}{$\text{-}2.104$}};
\node[blue] at (0.5,2.5){$\bullet$};
\end{tikzpicture}
,\\
\begin{tikzpicture}[scale=1.2,anchorbase,scale=0.75]
\draw[very thick] (0,2) to (-0.5,2.5) to (0,3) to (0.5,2.5);
\draw[very thick] (0.5,1.5) to (1,2) to (0.5,2.5) to (0,2) to (0.5,1.5);
\draw[very thick] (0.5,1.5) to (1,1) to (1.5,1.5) to (1,2);
\node at (0,2.5){\scalebox{0.7}{$\text{-}0.102$}};
\node at (0.5,2){\scalebox{0.7}{$0.897$}};
\node at (1,1.5){\scalebox{0.7}{$\text{-}1.103$}};
\node[blue] at (0,3){$\bullet$};
\end{tikzpicture}
\colon(3)\to (2,1),\quad
\begin{tikzpicture}[scale=1.2,anchorbase,scale=0.75]
\draw[very thick] (0,2) to (-0.5,2.5) to (0,3) to (0.5,2.5);
\draw[very thick] (0.5,1.5) to (1,2) to (0.5,2.5) to (0,2) to (0.5,1.5);
\draw[very thick] (0.5,1.5) to (1,1) to (1.5,1.5) to (1,2);
\node at (0,2.5){\scalebox{0.7}{$\text{-}0.102$}};
\node at (0.5,2){\scalebox{0.7}{$\text{-}1.103$}};
\node at (1,1.5){\scalebox{0.7}{$\text{-}2.104$}};
\node[blue] at (0,3){$\bullet$};
\end{tikzpicture}
\colon(3)\to (1^{3}),\quad
\begin{tikzpicture}[scale=1.2,anchorbase,scale=0.75]
\draw[very thick] (0,2) to (-0.5,1.5) to (-1,2) to (-0.5,2.5) to (0,2);
\draw[very thick] (0.5,1.5) to (1,2) to (0.5,2.5) to (0,2) to (0.5,1.5);
\draw[very thick] (-0.5,2.5) to (0,3) to (0.5,2.5);
\node at (0,2.5){\scalebox{0.7}{$\text{-}0.102$}};
\node at (-0.5,2){\scalebox{0.7}{$\text{-}1.103$}};
\node at (0.5,2){\scalebox{0.7}{$\text{-}2.104$}};
\node[blue] at (0,3){$\bullet$};
\end{tikzpicture}
\text{ and }
\begin{tikzpicture}[scale=1.2,anchorbase,scale=0.75]
\draw[very thick] (0,2) to (-0.5,1.5) to (-1,2) to (-0.5,2.5) to (0,2);
\draw[very thick] (0.5,1.5) to (1,2) to (0.5,2.5) to (0,2) to (0.5,1.5);
\draw[very thick] (-0.5,2.5) to (0,3) to (0.5,2.5);
\node at (0,2.5){\scalebox{0.7}{$\text{-}0.102$}};
\node at (-0.5,2){\scalebox{0.7}{$\text{-}2.104$}};
\node at (0.5,2){\scalebox{0.7}{$\text{-}1.103$}};
\node[blue] at (0,3){$\bullet$};
\end{tikzpicture}
\colon(2,1)\to (1^{3}),
\end{gather*}
that we have filled with their residues and coordinates. We have also illustrated seven semistandard tableaux that we will use momentarily. The three
associated idempotent diagrams are
\begin{gather*}
\1_{(3)}=\!
\begin{tikzpicture}[scale=1.2,anchorbase,smallnodes,rounded corners]
\draw[ghost](1,0)node[below]{$\phantom{i}$}--++(0,1)node[above,yshift=-1pt]{$0$};
\draw[ghost](1.9,0)node[below]{$\phantom{i}$}--++(0,1)node[above,yshift=-1pt]{$1$};
\draw[ghost](2.8,0)node[below]{$\phantom{i}$}--++(0,1)node[above,yshift=-1pt]{$2$};
\draw[solid](0,0)node[below]{$0$}--++(0,1)node[above,yshift=-1pt]{$\phantom{i}$};
\draw[solid](0.9,0)node[below]{$1$}--++(0,1)node[above,yshift=-1pt]{$\phantom{i}$};
\draw[solid](1.8,0)node[below]{$2$}--++(0,1)node[above,yshift=-1pt]{$\phantom{i}$};
\draw[redstring](0.2,0)node[below]{$0$}--++(0,1)node[above,yshift=-1pt]{$\phantom{i}$};
\end{tikzpicture}
,
\1_{(2,1)}=\!
\begin{tikzpicture}[scale=1.2,anchorbase,smallnodes,rounded corners]
\draw[ghost](1,0)node[below]{$\phantom{i}$}--++(0,1)node[above,yshift=-1pt]{$0$};
\draw[ghost](1.9,0)node[below]{$\phantom{i}$}--++(0,1)node[above,yshift=-1pt]{$1$};
\draw[ghost](-0.1,0)node[below]{$\phantom{i}$}--++(0,1)node[above,yshift=-1pt]{$2$};
\draw[solid](0,0)node[below]{$0$}--++(0,1)node[above,yshift=-1pt]{$\phantom{i}$};
\draw[solid](0.9,0)node[below]{$1$}--++(0,1)node[above,yshift=-1pt]{$\phantom{i}$};
\draw[solid](-1.1,0)node[below]{$2$}--++(0,1)node[above,yshift=-1pt]{$\phantom{i}$};
\draw[redstring](0.2,0)node[below]{$0$}--++(0,1)node[above,yshift=-1pt]{$\phantom{i}$};
\end{tikzpicture}
,
\1_{(1^{3})}=\!
\begin{tikzpicture}[scale=1.2,anchorbase,smallnodes,rounded corners]
\draw[ghost](1,0)node[below]{$\phantom{i}$}--++(0,1)node[above,yshift=-1pt]{$0$};
\draw[ghost](-1.2,0)node[below]{$\phantom{i}$}--++(0,1)node[above,yshift=-1pt]{$1$};
\draw[ghost](-0.1,0)node[below]{$\phantom{i}$}--++(0,1)node[above,yshift=-1pt]{$2$};
\draw[solid](0,0)node[below]{$0$}--++(0,1)node[above,yshift=-1pt]{$\phantom{i}$};
\draw[solid](-2.2,0)node[below]{$1$}--++(0,1)node[above,yshift=-1pt]{$\phantom{i}$};
\draw[solid](-1.1,0)node[below]{$2$}--++(0,1)node[above,yshift=-1pt]{$\phantom{i}$};
\draw[redstring](0.2,0)node[below]{$0$}--++(0,1)node[above,yshift=-1pt]{$\phantom{i}$};
\end{tikzpicture}
.
\end{gather*}
We now introduce notation for the various permutation diagrams we have. We use $D_{\bmu}^{\blam}$
for the permutation that has the coordinates of $\1_{\blam}$ at the bottom and goes to the coordinates defined by $\bmu$. For the seven semistandard tableaux above, in order, we then have $D_{(3)}^{(3)}=\1_{(3)}$, $D_{(2,1)}^{(2,1)}=\1_{(2,1)}$, $D_{(1^{3})}^{(1^{3})}=\1_{(1^{3})}$, and four more where we use the letter $E$ for the rightmost tableaux. Then our basis of $\WAc[3](X)$ is
\begin{gather*}
\BX[{\mathscr{R}}]
=
\left\{
\begin{gathered}
\1_{(3)},
\1_{(3)}D^{(3)}_{(2,1)},
\1_{(3)}D^{(3)}_{(1^{3})},
\\
\1_{(2,1)},
\1_{(2,1)}D^{(2,1)}_{(1^{3})},
\1_{(2,1)}E^{(2,1)}_{(1^{3})},
(D^{(3)}_{(2,1)})^{\star}\1_{(3)},
(D^{(3)}_{(2,1)})^{\star}\1_{(3)}D^{(3)}_{(2,1)},
(D^{(3)}_{(2,1)})^{\star}\1_{(3)}D^{(3)}_{(1^{3})},
\\
\1_{(1^{3})},
(D^{(3)}_{(1^{3})})^{\star}\1_{(3)},
(D^{(2,1)}_{(1^{3})})^{\star}\1_{(2,1)},
(E^{(2,1)}_{(1^{3})})^{\star}\1_{(2,1)},
(D^{(3)}_{(1^{3})})^{\star}\1_{(3)}D^{(3)}_{(1^{3})},
(D^{(3)}_{(1^{3})})^{\star}\1_{(3)}D^{(3)}_{(3)},
\\
(D^{(2,1)}_{(1^{3})})^{\star}\1_{(2,1)}D^{(2,1)}_{(1^{3})},
(E^{(2,1)}_{(1^{3})})^{\star}\1_{(2,1)}E^{(2,1)}_{(1^{3})},
(D^{(2,1)}_{(1^{3})})^{\star}\1_{(2,1)}E^{(2,1)}_{(1^{3})},
(E^{(2,1)}_{(1^{3})})^{\star}\1_{(2,1)}D^{(2,1)}_{(1^{3})}
\end{gathered}
\right\}.
\end{gather*}
Thus, the total rank of $\WAc[3](X)$ is nineteen, which is the rank of the Schur algebra $\mathscr{S}(3,3)$~\cite[Theorem 4.13]{Ma-hecke-schur}.

The cyclotomic KLR algebra is isomorphic to the idempotent truncation $\TAc[3]\cong\1_{A,3}\WAc[3](X)\1_{A,3}$, which is $6$ dimensional with basis
\begin{gather*}
\left\{
\begin{gathered}
\1_{(1^{3})},
(D^{(3)}_{(1^{3})})^{\star}\1_{(3)}D^{(3)}_{(1^{3})},
(D^{(2,1)}_{(1^{3})})^{\star}\1_{(2,1)}D^{(2,1)}_{(1^{3})},
(E^{(2,1)}_{(1^{3})})^{\star}\1_{(2,1)}E^{(2,1)}_{(1^{3})},
\\
(D^{(2,1)}_{(1^{3})})^{\star}\1_{(2,1)}E^{(2,1)}_{(1^{3})},
(E^{(2,1)}_{(1^{3})})^{\star}\1_{(2,1)}D^{(2,1)}_{(1^{3})}
\end{gathered}
\right\}.
\end{gather*}

The reader might ask, ``Where are the dots?'' . To answer this we give an example that shows that crossings encode dots. Applying \autoref{R:GhostSolid} and two honest Reidemeister II
relations shows that:
\begin{gather*}
y_{1}\1_{(2,1)}=
\begin{tikzpicture}[scale=1.2,anchorbase,smallnodes,rounded corners]
\draw[ghost](1,0)node[below]{$\phantom{i}$}--++(0,1)node[above,yshift=-1pt]{$0$};
\draw[ghost](1.9,0)node[below]{$\phantom{i}$}--++(0,1)node[above,yshift=-1pt]{$1$};
\draw[ghost,dot](-0.1,0)node[below]{$\phantom{i}$}--++(0,1)node[above,yshift=-1pt]{$2$};
\draw[solid](0,0)node[below]{$0$}--++(0,1)node[above,yshift=-1pt]{$\phantom{i}$};
\draw[solid](0.9,0)node[below]{$1$}--++(0,1)node[above,yshift=-1pt]{$\phantom{i}$};
\draw[solid,dot](-1.1,0)node[below]{$2$}--++(0,1)node[above,yshift=-1pt]{$\phantom{i}$};
\draw[redstring](0.2,0)node[below]{$0$}--++(0,1)node[above,yshift=-1pt]{$\phantom{i}$};
\end{tikzpicture}
=
\begin{tikzpicture}[scale=1.2,anchorbase,smallnodes,rounded corners]
\draw[ghost](1,0)node[below]{$\phantom{i}$}--++(0,1)node[above,yshift=-1pt]{$0$};
\draw[ghost](1.9,0)node[below]{$\phantom{i}$}--++(0,1)node[above,yshift=-1pt]{$1$};
\draw[ghost](-0.1,0)node[below]{$\phantom{i}$}--++(0,0.1)--++(2.9,0.3)--++(0,0.2)--++(-2.9,0.3)--++(0,0.1)node[above,yshift=-1pt]{$2$};
\draw[solid](0,0)node[below]{$0$}--++(0,1)node[above,yshift=-1pt]{$\phantom{i}$};
\draw[solid](0.9,0)node[below]{$1$}--++(0,1)node[above,yshift=-1pt]{$\phantom{i}$};
\draw[solid](-1.1,0)node[below]{$2$}--++(0,0.1)--++(2.9,0.3)--++(0,0.2)--++(-2.9,0.3)--++(0,0.1)node[above,yshift=-1pt]{$\phantom{i}$};
\draw[redstring](0.2,0)node[below]{$0$}--++(0,1)node[above,yshift=-1pt]{$\phantom{i}$};
\end{tikzpicture}
\in\WAc[3](X)
.
\end{gather*}
In particular, the construction of the basis $\BX[{\mathscr{R}}]$ is very different to that of \cite{HuMa-klr-basis}. However, as this is a small example it turns out that the basis $\BX[{\mathscr{R}}]$ coincides with the cellular basis given in \cite{HuMa-klr-basis}. In general, the cellular ideals given by $\BX[{\mathscr{R}}]$ and by \cite{HuMa-klr-basis} coincide in the asymptotic case, which can be proved  following~\cite[Proposition~7.3]{Bo-many-cellular-structures}.
\end{Example}

\begin{Remark}\label{R:Calc}
Hu and Shi~\cite{HuSh-monomial-klr-basis}  have given dimension formulas for the cyclotomic KLR algebras for symmetrizable quivers.
If the reader wants to verify that \autoref{E:Calc} agrees with their formulas we have written some SageMath code that calculates dimensions of the cyclotomic KLR algebras, using for example the online calculator of SageMath \url{https://sagecell.sagemath.org/}.
This code is available, with a self-contained explanation, on GitHub \cite{MaTu-sagemath-finite-type-klrw}.
\end{Remark}


\subsection{Simple modules}\label{SS:SimplesA}


Let us recall parts of the general theory of (graded affine) cellular
algebras from \cite{GrLe-cellular}, \cite{HuMa-klr-basis},
\cite{KoXi-affine-cellular}, \cite{EhTu-relcell}, \cite{TuVa-handlebody}
or \cite{Tu-sandwich-cellular}.

For $\blam\in\hParts$ and $\bT\in\SStd(\blam)$
the graded cellular structure defines a
graded (left) cell module $\Delta(\blam,\bT)$ via
\begin{gather*}
\Delta_{\bT}(\blam)
=\bigl\langle D^{\ba}_{\bS\bT}|\bS\in\hSStd(\blam)\text{ and }\ba\in\Affch\bigr\rangle_{R}
\end{gather*}
with the $\WA(X)$-action defined modulo $\ldom$-higher order terms.
By \autoref{D:CellularAlgebra}, $\Delta_{\bT}(\blam)\cong\Delta_{\bT^{\prime}}(\blam)$ as $\WA(X)$-modules,
so we drop the second superscript. The $\Delta(\blam)$
are the \emph{cell modules} of $\WA(X)$.

We also need more general cell modules. To define them
let $B(\blam)=R[y^{\ba}|\ba\in\Affch]$, a polynomial subring of $R[y_{1},\dots,y_{n}]$.
By convention, we set $B(\blam)=R$ in the cyclotomic case.

Let $K$ be a simple $B(\blam)$-module.
Then the corresponding \emph{affine cell module} is
\begin{gather*}
\Delta(\blam,K)=\Delta(\blam)\otimes_{B(\lambda)}K.
\end{gather*}
Note that $\Delta(\blam,K)\cong\Delta(\blam)$ is a graded $\WA(X)$-module if and only
if $K\cong R$ is the trivial $B(\blam)$-module, since $R$ is the only graded $B(\lambda)$-module.

Each cell module $\Delta(\blam)$ has an associated
\emph{cellular pairing} determined by
\begin{gather*}
\langle
D^{\ba}_{\bS\bT},
D^{\bb}_{\bU\bV}
\rangle
=
r_{\bT\bU}
\leftrightsquigarrow
\langle
D^{\ba}_{\bS\bT},
D^{\bb}_{\bU\bV}
\rangle
=
\text{ coefficient of $1$ of }
\begin{tikzpicture}[scale=1.2,anchorbase,scale=1]
\draw[line width=0.75,color=black,fill=cream] (0,1) to (0.25,0.5) to (0.75,0.5) to (1,1) to (0,1);
\node at (0.5,0.75){$\bU$};
\draw[line width=0.75,color=black,fill=cream] (0,1) to (0.25,1.5) to (0.75,1.5) to (1,1) to (0,1);
\node at (0.5,1.25){$\bT$};
\draw[line width=0.75,color=black,fill=cream] (0.25,0) to (0.25,0.5) to (0.75,0.5) to (0.75,0) to (0.25,0);
\node at (0.5,0.25){$\bb$};
\draw[line width=0.75,color=black,fill=cream] (0.25,1.5) to (0.25,2) to (0.75,2) to (0.75,1.5) to (0.25,1.5);
\node at (0.5,1.75){$\ba$};
\end{tikzpicture}
\text{ in }B(\blam).
\end{gather*}
See \cite[Section 2.2]{KoXi-affine-cellular} for the precise definition.
This can be extended to $\Delta(\blam,K)$ by using the identity
on $K$.
Let $\mathrm{rad}\,\Delta(\blam,K)$ be the
radical of the bilinear form on $\Delta(\blam,K)$, and define
$L(\blam,K)=\Delta(\blam,K)/
\mathrm{rad}\,\Delta(\blam,K)$.

In the following theorem let $q^{s}M$ be the graded module obtained by shifting the grading on $M$ up by $s\in\Z$. Let $S_{B(\blam)}$
be a choice of simple $B(\blam)$-modules, up to isomorphism.

\begin{Theorem}\label{T:SimplesA}
Suppose that $R$ is a field.
\begin{enumerate}

\item The set $\set{L(\blam,K)|\blam\in\hParts,K\in S_{B(\blam)}}$
is a complete and non-redundant set of simple $\WA[n](X)$-modules.

\item The set $\set{q^{s}L(\blam)|\blam\in\hParts,s\in\Z}$ is
a complete and non-redundant set of graded simple $\WA[n](X)$-modules.

\item The set $\set{L(\blam)|\blam\in\Parts}$ is a complete and
non-redundant set of simple $\WAc[n](X)$-modules.

\item The set $\set{q^{s}L(\blam)|\blam\in\Parts,s\in\Z}$ is a
complete and non-redundant set of graded simple $\WAc[n](X)$-modules.

\end{enumerate}

\end{Theorem}

\begin{proof}
The proof of this theorem is easy using the results from the
previous section and the usual arguments in the theory of cellular algebras:

The standard arguments in the theory of cellular algebras imply that each cell
module $\Delta(\blam)$ has a simple head $L(\blam)$, provided that the cellular pairing does not vanish.
The image of the idempotent $\1_{\blam}$ in $\Delta(\blam)$ has nonzero inner
product with itself, so this pairing cannot vanish. This also implies that
$L(\blam)\cong q^{d}L(\blam)$ if and only if $d=0$ and $\bmu=\blam$, although this is also a
consequence of the theory of cellular algebras. For general
$K$ the same arguments work.

It remains to show that every simple $\WA(X)$-module is isomorphic to $L(\blam,K)$,
for some $\blam\in\hParts$. However, this
is also an immediate consequence of standard cellular algebra arguments.

The statement
about the grading then follows from the ungraded cases
together with the observation that
$K=R$ is the only graded simple $B(\blam)$-module.

In the cyclotomic case we use the same arguments but also note that
$B(\blam)\cong R$ in this case.
\end{proof}

The following is an almost immediate consequence
(quasi-hereditary algebras and their affine versions are studied, for example, in
\cite{ClPaSc-h-weight-qh} respectively \cite{Kl-affine-quasi-hereditary}):

\begin{Corollary}\label{C:qhereditary}
Suppose that $R$ is a field.
The algebra $\WA(X)$ is a graded affine quasi-hereditary algebra and $\WAc(X)$ is a graded quasi-hereditary algebra.
\end{Corollary}

\begin{proof}
Note that \autoref{T:SimplesA} shows that every (graded) cell module contributes
a unique associated (graded) simple module. Thus, the claim follows by
the general theory of cellular algebras; see, for example,  \cite[Remark (3.10)]{GrLe-cellular}.
\end{proof}

The corresponding statements
for the (cyclotomic) KLR algebras are more delicate because of the idempotent truncation
in the proof of \autoref{C:KLRCellular}. The simple modules of these algebras are
classified in \cite{Ar-classification-cyclotomic-hecke} (proving a conjecture
of \cite{ArMa-simples-complex-reflection}), and the graded simple modules in
\cite[Theorem 5.13]{BrKl-graded-decomposition-hecke} and
\cite[Theorem 10.5]{Bo-many-cellular-structures}. For example, \cite{KlLoMi-KLR-affine-cellular-typea}
shows that when $n=2$ then, up to shift, the (affine) KLR algebra of type $A^{(1)}_{1}$ has four
graded simple modules whereas the corresponding weighted KLRW algebras have more graded simple
modules, with the precise number depending on the choice of $\brho$. See also \autoref{Ex:AffineCompare}.

\begin{Remark}
\autoref{T:SubDiv} applies to the algebra $\WA[n](X)$ in the infinitesimal
case and hence to the algebras $\TA[n]$ and $\TAc[n]$. It would be interesting
to explicitly describe what this tells us about the decomposition numbers of these
algebras for different quivers.
\end{Remark}


\subsection{Proof of cellularity in type $A$}
\label{SS:TypeA}


The following two lemmas are instances of
\emph{pulling strings and jumping dots to the right}.
In these lemmas we color, in green, the strings
where either the string pulls, or the dot jumps, to the right.

\begin{Lemma}\label{L:GeneralSliding}
For any quiver and any choice of $Q$-polynomials we have
\begin{gather}\label{R:IISliding}
\begin{tikzpicture}[scale=1.2,anchorbase,smallnodes,rounded corners]
\draw[solid,spinach](0,1)node[above,yshift=-1pt]{$\phantom{i}$}--++(0,-1)node[below]{$i$};
\draw[solid](0.5,1)--++(0,-1)node[below]{$i$};
\end{tikzpicture}
=
\begin{tikzpicture}[scale=1.2,anchorbase,smallnodes,rounded corners]
\draw[solid,dot](0,1)node[above,yshift=-1pt]{$\phantom{i}$}--++(0.5,-0.5)--++(-0.5,-0.5) node[below]{$i$};
\draw[solid,dot=0.1](0.5,1)--++(-0.5,-0.5)--++(0.5,-0.5) node[below]{$i$};
\end{tikzpicture}
-
\begin{tikzpicture}[scale=1.2,anchorbase,smallnodes,rounded corners]
\draw[solid,dot,dot=0.9](0,1)node[above,yshift=-1pt]{$\phantom{i}$}--++(0.5,-0.5)--++(-0.5,-0.5) node[below]{$i$};
\draw[solid](0.5,1)--++(-0.5,-0.5)--++(0.5,-0.5) node[below]{$i$};
\end{tikzpicture}
.
\end{gather}
\begin{gather}\label{R:IIISliding}
\begin{tikzpicture}[scale=1.2,anchorbase,smallnodes,rounded corners]
\draw[solid,dot,spinach](0,1)node[above,yshift=-1pt]{$\phantom{i}$}--++(0,-1)node[below]{$i$};
\draw[solid](0.5,1)--++(0,-1)node[below]{$i$};
\end{tikzpicture}
=
\begin{tikzpicture}[scale=1.2,anchorbase,smallnodes,rounded corners]
\draw[solid](0,1)node[above,yshift=-1pt]{$\phantom{i}$}--++(0,-1)node[below]{$i$};
\draw[solid,dot](0.5,1)--++(0,-1)node[below]{$i$};
\end{tikzpicture}
+
\begin{tikzpicture}[scale=1.2,anchorbase,smallnodes,rounded corners]
\draw[solid,dot=0.4,dot=0.6](0,1)node[above,yshift=-1pt]{$\phantom{i}$}--++(0.5,-0.5)--++(-0.5,-0.5) node[below]{$i$};
\draw[solid,dot=0.1](0.5,1)--++(-0.5,-0.5)--++(0.5,-0.5) node[below]{$i$};
\end{tikzpicture}
-
\begin{tikzpicture}[scale=1.2,anchorbase,smallnodes,rounded corners]
\draw[solid,dot=0.4,dot=0.6,dot=0.9](0,1)node[above,yshift=-1pt]{$\phantom{i}$}--++(0.5,-0.5)--++(-0.5,-0.5)node[below]{$i$};
\draw[solid](0.5,1)--++(-0.5,-0.5)--++(0.5,-0.5) node[below]{$i$};
\end{tikzpicture}
.
\end{gather}
\end{Lemma}

In the both equations we are pulling the leftmost string to the right and in the second equation the dot also jumps to the right.

\begin{proof}
Equation \autoref{R:IISliding} follows from \autoref{R:DotCrossing}
and \autoref{R:SolidSolid}. Then \autoref{R:IIISliding} comes from \autoref{R:IISliding} and \autoref{R:DotCrossing}.
\end{proof}

Our choice of $Q$-polynomial implies the next result,
which is similar to \cite[(5.1) and (5.2)]{Bo-many-cellular-structures}.

\begin{Lemma}\label{L:Sliding}
The following hold.
\begin{gather}\label{R:GhostSolidTypeA}
\begin{tikzpicture}[scale=1.2,anchorbase,smallnodes,rounded corners]
\draw[ghost,dot,spinach](0,1)node[above,yshift=-1pt]{$\phantom{i}$}--++(0,-1)node[below]{$i$};
\draw[solid](0.5,1)--++(0,-1)node[below]{$i{+}1$};
\end{tikzpicture}
=
\begin{tikzpicture}[scale=1.2,anchorbase,smallnodes,rounded corners]
\draw[ghost](0,1)node[above,yshift=-1pt]{$\phantom{i}$}--++(0,-1)node[below]{$i$};
\draw[solid,dot](0.5,1)--++(0,-1)node[below]{$i{+}1$};
\end{tikzpicture}
+
\begin{tikzpicture}[scale=1.2,anchorbase,smallnodes,rounded corners]
\draw[ghost](0,1)node[above,yshift=-1pt]{$\phantom{i}$}--++(0.5,-0.5)--++(-0.5,-0.5) node[below]{$i$};
\draw[solid](0.5,1)--++(-0.5,-0.5)--++(0.5,-0.5) node[below]{$i{+}1$};
\end{tikzpicture}
,\quad
\begin{tikzpicture}[scale=1.2,anchorbase,smallnodes,rounded corners]
\draw[ghost](0.5,1)node[above,yshift=-1pt]{$\phantom{i}$}--++(0,-1)node[below]{$i$};
\draw[solid,dot,spinach](0,1)--++(0,-1)node[below]{$i{+}1$};
\end{tikzpicture}
=
\begin{tikzpicture}[scale=1.2,anchorbase,smallnodes,rounded corners]
\draw[ghost,dot](0.5,1)node[above,yshift=-1pt]{$\phantom{i}$}--++(0,-1)node[below]{$i$};
\draw[solid](0,1)--++(0,-1)node[below]{$i{+}1$};
\end{tikzpicture}
-
\begin{tikzpicture}[scale=1.2,anchorbase,smallnodes,rounded corners]
\draw[ghost](0.5,1)node[above,yshift=-1pt]{$\phantom{i}$}--++(-0.5,-0.5)--++(0.5,-0.5) node[below]{$i$};
\draw[solid](0,1)--++(0.5,-0.5)--++(-0.5,-0.5) node[below]{$i{+}1$};
\end{tikzpicture}
.
\end{gather}
\begin{gather}\label{E:Sliding}
\begin{gathered}
\begin{tikzpicture}[scale=1.2,anchorbase,smallnodes,rounded corners]
\draw[ghost](1,0)--++(0,1)node[above,yshift=-1pt]{$i$};
\draw[ghost](1.5,0)--++(0,1)node[above,yshift=-1pt]{$i$};
\draw[solid,spinach](0,0)node[below]{$i$}--++(0,1);
\draw[solid](0.5,0)node[below]{$i$}--++(0,1);
\draw[solid](1.25,0)node[below]{$i{+}1$}--++(0,1);
\end{tikzpicture}
=
-
\begin{tikzpicture}[scale=1.2,anchorbase,smallnodes,rounded corners]
\draw[ghost](1,0)--++(0.5,0.5)--++(-0.5,0.5)node[above,yshift=-1pt]{$i$};
\draw[ghost,dot](1.5,0)--++(-0.5,0.5)--++(0.5,0.5)node[above,yshift=-1pt]{$i$};
\draw[solid](0,0)node[below]{$i$}--++(0.5,0.5)--++(-0.5,0.5);
\draw[solid,dot](0.5,0)node[below]{$i$}--++(-0.5,0.5)--++(0.5,0.5);
\draw[solid](1.25,0)node[below]{$i{+}1$}--++(-0.5,0.5)--++(0.5,0.5);
\end{tikzpicture}
-
\begin{tikzpicture}[scale=1.2,anchorbase,smallnodes,rounded corners]
\draw[ghost,dot](1,0)--++(0.5,0.5)--++(-0.5,0.5)node[above,yshift=-1pt]{$i$};
\draw[ghost](1.5,0)--++(-0.5,0.5)--++(0.5,0.5)node[above,yshift=-1pt]{$i$};
\draw[solid,dot](0,0)node[below]{$i$}--++(0.5,0.5)--++(-0.5,0.5);
\draw[solid](0.5,0)node[below]{$i$}--++(-0.5,0.5)--++(0.5,0.5);
\draw[solid](1.25,0)node[below]{$i{+}1$}--++(0.5,0.5)--++(-0.5,0.5);
\end{tikzpicture},
\\
\begin{tikzpicture}[scale=1.2,anchorbase,smallnodes,rounded corners]
\draw[ghost](0.25,0)--++(0,1)node[above,yshift=-1pt]{$i{-}1$};
\draw[ghost,spinach](1,0)--++(0,1)node[above,yshift=-1pt]{$i$};
\draw[ghost](1.5,0)--++(0,1)node[above,yshift=-1pt]{$i$};
\draw[solid](0,0)node[below]{$i$}--++(0,1);
\draw[solid](0.5,0)node[below]{$i$}--++(0,1);
\end{tikzpicture}
=
+
\begin{tikzpicture}[scale=1.2,anchorbase,smallnodes,rounded corners]
\draw[ghost](1,0)--++(0.5,0.5)--++(-0.5,0.51)node[above,yshift=-1pt]{$i$};
\draw[ghost,dot](1.5,0)--++(-0.5,0.5)--++(0.5,0.5)node[above,yshift=-1pt]{$i$};
\draw[ghost](0.25,0)--++(-0.5,0.5)--++(0.5,0.5)node[above,yshift=-1pt]{$i{-}1$};
\draw[solid](0,0)node[below]{$i$}--++(0.5,0.5)--++(-0.5,0.5);
\draw[solid,dot](0.5,0)node[below]{$i$}--++(-0.5,0.5)--++(0.5,0.5);
\end{tikzpicture}
+
\begin{tikzpicture}[scale=1.2,anchorbase,smallnodes,rounded corners]
\draw[ghost,dot](1,0)--++(0.5,0.5)--++(-0.5,0.5)node[above,yshift=-1pt]{$i$};
\draw[ghost](1.5,0)--++(-0.5,0.5)--++(0.5,0.5)node[above,yshift=-1pt]{$i$};
\draw[ghost](0.25,0)--++(0.5,0.5)--++(-0.5,0.5)node[above,yshift=-1pt]{$i{-}1$};
\draw[solid,dot](0,0)node[below]{$i$}--++(0.5,0.5)--++(-0.5,0.5);
\draw[solid](0.5,0)node[below]{$i$}--++(-0.5,0.5)--++(0.5,0.5);
\end{tikzpicture}.
\end{gathered}
\end{gather}
\end{Lemma}

\begin{Remark}
We will apply this lemma to pull the green strings, or dots on green strings,
to the right. By applying \autoref{R:IIISliding} to the right-hand side of
\autoref{E:Sliding} we can move the dot on the left $i$-string to the right.
\end{Remark}

\begin{proof}
Equation \autoref{R:GhostSolidTypeA} follows immediately
from \autoref{R:GhostSolid} and our choice of $Q$-polynomial.
We only prove the first
equation in \autoref{E:Sliding} as the other can be proven {\muta}.
To this end, we use
\begin{gather*}
\begin{tikzpicture}[scale=1.2,anchorbase,smallnodes,rounded corners]
\draw[ghost](1,1)node[above,yshift=-1pt]{$\phantom{i}$}--++(1,-1)node[below]{$i$};
\draw[ghost](2,1)--++(-1,-1)node[below]{$i$};
\draw[solid,smallnodes,rounded corners](1.5,1)--++(-0.5,-0.5)--++(0.5,-0.5)node[below]{$i{+}1$};
\end{tikzpicture}
=
\begin{tikzpicture}[scale=1.2,anchorbase,smallnodes,rounded corners]
\draw[ghost](3,1)node[above,yshift=-1pt]{$\phantom{i}$}--++(1,-1)node[below]{$i$};
\draw[ghost](4,1)--++(-1,-1)node[below]{$i$};
\draw[solid,smallnodes,rounded corners](3.5,1)--++(0.5,-0.5)--++(-0.5,-0.5)node[below]{$i{+}1$};
\end{tikzpicture}
-
\begin{tikzpicture}[scale=1.2,anchorbase,smallnodes,rounded corners]
\draw[ghost](6.2,1)node[above,yshift=-1pt]{$\phantom{i}$}--++(0,-1)node[below]{$i$};
\draw[ghost](7.2,1)--++(0,-1)node[below]{$i$};
\draw[solid](6.7,1)--++(0,-1)node[below]{$i{+}1$};
\end{tikzpicture}
,\quad
\begin{tikzpicture}[scale=1.2,anchorbase,smallnodes,rounded corners]
\draw[solid](1,1)node[above,yshift=-1pt]{$\phantom{i}$}--++(1,-1)node[below]{$i$};
\draw[solid](2,1)--++(-1,-1)node[below]{$i$};
\draw[ghost,smallnodes,rounded corners](1.5,1)--++(-0.5,-0.5)--++(0.5,-0.5)node[below]{$i{+}1$};
\end{tikzpicture}
=
\begin{tikzpicture}[scale=1.2,anchorbase,smallnodes,rounded corners]
\draw[solid](3,1)node[above,yshift=-1pt]{$\phantom{i}$}--++(1,-1)node[below]{$i$};
\draw[solid](4,1)--++(-1,-1)node[below]{$i$};
\draw[ghost,smallnodes,rounded corners](3.5,1)--++(0.5,-0.5)--++(-0.5,-0.5)node[below]{$i{+}1$};
\end{tikzpicture}
+
\begin{tikzpicture}[scale=1.2,anchorbase,smallnodes,rounded corners]
\draw[solid](7.2,1)node[above,yshift=-1pt]{$\phantom{i}$}--++(0,-1)node[below]{$i$};
\draw[solid](8.2,1)--++(0,-1)node[below]{$i$};
\draw[ghost](7.7,1)--++(0,-1)node[below]{$i{+}1$};
\end{tikzpicture}
\end{gather*}
(this is \autoref{R:BraidGSG} for our choice of $Q$-polynomial),
to pull the two solid $i$-strings together:
\begin{gather*}
\begin{tikzpicture}[scale=1.2,anchorbase,smallnodes,rounded corners]
\draw[ghost](1,0)--++(0,1)node[above,yshift=-1pt]{$i$};
\draw[ghost](1.5,0)--++(0,1)node[above,yshift=-1pt]{$i$};
\draw[solid](0,0)node[below]{$i$}--++(0,1);
\draw[solid](0.5,0)node[below]{$i$}--++(0,1);
\draw[solid](1.25,0)node[below]{$i{+}1$}--++(0,1);
\end{tikzpicture}
=
-
\begin{tikzpicture}[scale=1.2,anchorbase,smallnodes,rounded corners]
\draw[ghost](1,0)--++(0.5,1)node[above,yshift=-1pt]{$i$};
\draw[ghost](1.5,0)--++(-0.5,1)node[above,yshift=-1pt]{$i$};
\draw[solid](0,0)node[below]{$i$}--++(0.5,1);
\draw[solid](0.5,0)node[below]{$i$}--++(-0.5,1);
\draw[solid](1.25,0)node[below]{$i{+}1$}--++(-0.25,0.5)--++(0.25,0.5);
\end{tikzpicture}
+
\begin{tikzpicture}[scale=1.2,anchorbase,smallnodes,rounded corners]
\draw[ghost](1,0)--++(0.5,1)node[above,yshift=-1pt]{$i$};
\draw[ghost](1.5,0)--++(-0.5,1)node[above,yshift=-1pt]{$i$};
\draw[solid](0,0)node[below]{$i$}--++(0.5,1);
\draw[solid](0.5,0)node[below]{$i$}--++(-0.5,1);
\draw[solid](1.25,0)node[below]{$i{+}1$}--++(0.25,0.5)--++(-0.25,0.5);
\end{tikzpicture}
.
\end{gather*}
The equality then follows by applying the identity
\begin{gather*}
-
\begin{tikzpicture}[scale=1.2,anchorbase,smallnodes,rounded corners]
\draw[solid,dot](0,0)node[below]{$i$}--++(0.5,0.5)--++(-0.5,0.5);
\draw[solid](0.5,0)node[below]{$i$}--++(-0.5,0.5)--++(0.5,0.5)node[above,yshift=-1pt]{\phantom{$i$}};
\end{tikzpicture}
=
\begin{tikzpicture}[scale=1.2,anchorbase,smallnodes,rounded corners]
\draw[solid](0,0)node[below]{$i$}--++(0.5,1);
\draw[solid](0.5,0)node[below]{$i$}--++(-0.5,1)node[above,yshift=-1pt]{\phantom{$i$}};
\end{tikzpicture}
=
\begin{tikzpicture}[scale=1.2,anchorbase,smallnodes,rounded corners]
\draw[solid](0,0)node[below]{$i$}--++(0.5,0.5)--++(-0.5,0.5);
\draw[solid,dot](0.5,0)node[below]{$i$}--++(-0.5,0.5)--++(0.5,0.5)node[above,yshift=-1pt]{\phantom{$i$}};
\end{tikzpicture}
,
\end{gather*}
which is a consequence of \autoref{R:DotCrossing}
and \autoref{R:SolidSolid}.
\end{proof}

Recall that $X\subset\R$. Let $\min X$ and $\max X$ be the
minimal and maximal elements of $X$, respectively.
The \emph{region defined by $X$} is $[\min X,\max X+1]\times[0,1]$.
(We need $\max X+1$, and not $\max X$, for the rightmost ghost string
since the rightmost solid strings have boundary points at $\max X$.)

\begin{Definition}
An idempotent diagram is \emph{$\varepsilon$-separated} if all of its coordinates are
within the region defined by $X$ and all of its strings are at least $\varepsilon$ apart.
\end{Definition}

\begin{Lemma}\label{L:Separated}
Suppose that $S$ is a straight line diagram.
Then $S$ factors through a unique
$\varepsilon$-separated idempotent diagram $L(S)$ that is minimal in
the dominance order.
\end{Lemma}

\begin{proof}
Up to isotopy we may assume
that $S$ is contained in the region defined by $X$.
Using isotopy, pull the leftmost
string in $S$ as far to the left as possible. Using
isotopy again, pull the second
string as far to the left as possible, so that it remains at
least $\varepsilon$ to the right of
the first string. Continuing in this way proves the lemma.
\end{proof}

The diagram $L(S)$ in \autoref{L:Separated} is the \emph{left justification} of $S$.

\begin{Remark}
Left justification gives a ``normal form'' for comparing straight line diagrams with respect to the $\ledom_{A}$-ordering. We could equally well use right justified diagrams. We use left justified diagrams for consistency with \cite{MaTu-klrw-algebras-bad}.
\end{Remark}

\begin{Notation}
As in the proof of \autoref{L:Separated}, isotopy allows us to assume
that all strings are contained in the region defined by
$X$. We will do this without further notice in all proofs
in this section.
\end{Notation}

Recall that we use affine red strings at position $\affine{\kappa}_{m}$, for
$\ell<m\leq\hell$. In a straight line diagram two strings are \emph{adjacent} if it
is possible to draw a horizontal line connecting them without crossing any other string.
An $i$-string is \emph{left adjacent} to a $j$-string, and the $j$-string is
\emph{right adjacent} to the $i$-string, if the strings are adjacent and the $i$-string is to the left of the $j$-string.

\begin{Example}
Adjacency is a local and not a multilocal condition. For example,
\begin{gather}\label{E:Close}
\begin{tikzpicture}[scale=1.2,anchorbase,smallnodes,rounded corners]
\draw[ghost](1,0)--++(0,1)node[above,yshift=-1pt]{$i$};
\draw[ghost](1.2,0)--++(0,1)node[above,yshift=-1pt]{$j$};
\draw[solid](0,0)node[below]{$i$}--++(0,1);
\draw[solid](0.2,0)node[below]{$j$}--++(0,1);
\end{tikzpicture}
\quad\text{and}
\quad
\begin{tikzpicture}[scale=1.2,anchorbase,smallnodes,rounded corners]
\draw[ghost](1,0)--++(0,1)node[above,yshift=-1pt]{$i$};
\draw[ghost](1.2,0)--++(0,1)node[above,yshift=-1pt]{$j$};
\draw[solid](0,0)node[below]{$i$}--++(0,1);
\draw[solid](0.2,0)node[below]{$j$}--++(0,1);
\draw[solid](1.1,0)node[below]{$k$}--++(0,1);
\end{tikzpicture}
\end{gather}
are both examples of adjacent solid $i$ and $j$-strings. In both diagrams the solid $i$-string is
left adjacent to the solid $j$-string. In the left-hand diagram, the ghost
$j$-string is right adjacent to the ghost $i$-string but the ghosts strings are not adjacent in the right-hand diagram.
\end{Example}

The crucial illustrations to keep in mind for the following definition are
\autoref{E:BoxConfiguration} and \autoref{E:CasesTypeA}.

\begin{Definition}\label{D:Young}
Let $S$ be a straight line diagram.
A solid $i$-string in $S$ is \emph{Young equivalent} to a solid $j$-string if
one or more of the following conditions hold:
\begin{enumerate}

\item $j=i+1$ and the ghost $i$-string is right adjacent to the solid $(i+1)$-string;

\item $j=i-1$ and the solid $i$-string is right adjacent to the ghost $(i-1)$-string.

\end{enumerate}
Let $Y$ be a Young equivalence class of solid strings and let $S_{Y}$ be the
subdiagram of $S$ that contains the solid strings in $Y$ together with their ghosts.
The \emph{suspension point} of $Y$ is the rightmost solid string in $S_{Y}$ that is
left adjacent to an (affine) red string in $S$.
\end{Definition}

\begin{Example}\label{E:YoungEquiClasses}
Maintain the notation from \autoref{Ex:BigPositioning}. In particular,
let $\blam=(3,2|\emptyset|1^{2}|\emptyset|\dots|\emptyset)$. Then the solid strings
in the diagram $\1_{\blam}$ are in one of two Young equivalence classes $Y(1)$ or $Y(3)$, depending upon whether the strings correspond
to the nodes in the first or third component of $\blam$. The subdiagrams $S_{Y(i)}$
for these two equivalences class are:
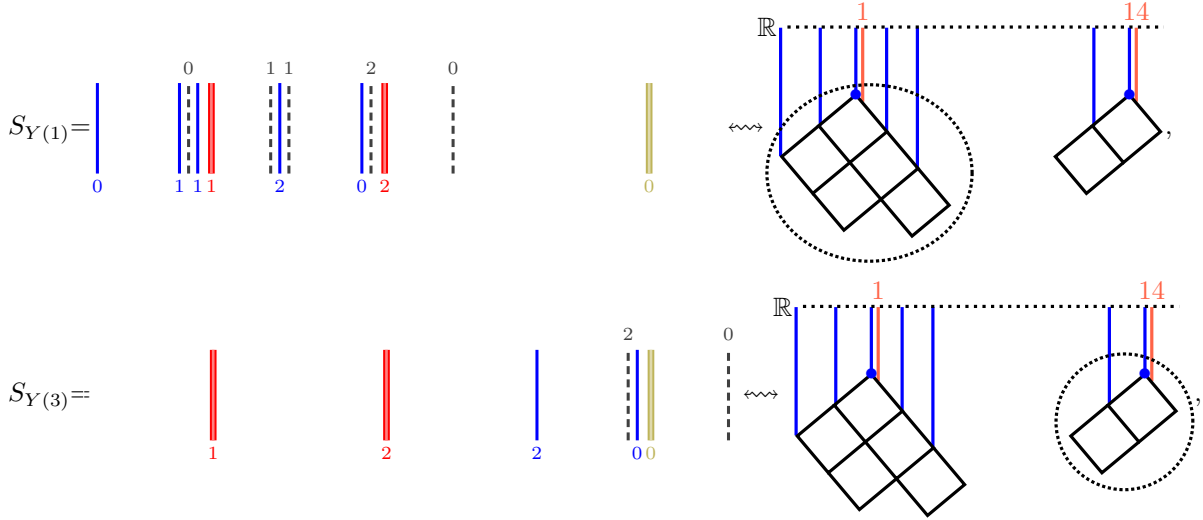
\begin{align*}
S_{Y(1)}&\hspace*{-0.1cm}=\hspace*{-0.2cm}
\begin{tikzpicture}[scale=1.2,anchorbase,smallnodes,rounded corners]
\draw[white](6.1,0)--++(0,1);
\draw[ghost](0.2,0)node[below]{$\phantom{i}$}--++(0,1)node[above,yshift=-1pt]{$0$};
\draw[ghost](1.1,0)node[below]{$\phantom{i}$}--++(0,1)node[above,yshift=-1pt]{$1$};
\draw[ghost](1.3,0)node[below]{$\phantom{i}$}--++(0,1)node[above,yshift=-1pt]{$1$};
\draw[ghost](2.2,0)node[below]{$\phantom{i}$}--++(0,1)node[above,yshift=-1pt]{$2$};
\draw[ghost](3.1,0)node[below]{$\phantom{i}$}--++(0,1)node[above,yshift=-1pt]{$0$};
\draw[solid](-0.8,0)node[below]{$0$}--++(0,1)node[above,yshift=-1pt]{$\phantom{i}$};
\draw[solid](0.1,0)node[below]{$1$}--++(0,1)node[above,yshift=-1pt]{$\phantom{i}$};
\draw[solid](0.3,0)node[below]{$1$}--++(0,1)node[above,yshift=-1pt]{$\phantom{i}$};
\draw[solid](1.2,0)node[below]{$2$}--++(0,1)node[above,yshift=-1pt]{$\phantom{i}$};
\draw[solid](2.1,0)node[below]{$0$}--++(0,1)node[above,yshift=-1pt]{$\phantom{i}$};
\draw[redstring](0.45,0)node[below]{$1$}--++(0,1)node[above,yshift=-1pt]{$\phantom{i}$};
\draw[redstring](2.35,0)node[below]{$2$}--++(0,1)node[above,yshift=-1pt]{$\phantom{i}$};
\draw[affine](5.25,0)node[below]{$0$}--++(0,1)node[above,xshift=0.25cm,yshift=-0.02]{$\phantom{i}$};
\end{tikzpicture}
\hspace*{-0.1cm}\leftrightsquigarrow\hspace*{-0.3cm}
\begin{tikzpicture}[scale=1.2,anchorbase,scale=0.75]
\draw[very thick,blue] (-0.52,-1.45) to (-0.52,-0.01);
\draw[very thick,blue] (0,-1) to (0,-0.01);
\draw[very thick,tomato] (0.1,-1.1) to (0.1,-0.01)node[above,yshift=-1pt]{$1$};
\draw[very thick,blue] (-1.1,-1.9) to (-1.1,-0.01);
\draw[very thick,blue] (0.45,-1.55) to (0.45,-0.01);
\draw[very thick,blue] (0.9,-2.1) to (0.9,-0.01);
\draw[very thick,blue] (3.48,-1.45) to (3.48,-0.01);
\draw[very thick,blue] (4,-1) to (4,-0.01);
\draw[very thick,tomato] (4.1,-1.1) to (4.1,-0.01)node[above,yshift=-1pt]{$14$};
\begin{scope}[rotate around={-5:(0,-1)}]
\draw[very thick] (-0.5,-1.5) to (-1,-2) to (-0.5,-2.5) to (0,-2);
\draw[very thick] (0,-1) to (-0.5,-1.5) to (0,-2) to (0.5,-1.5) to (0,-1);
\draw[very thick] (0.5,-1.5) to (1,-2) to (0.5,-2.5) to (0,-2);
\draw[very thick] (1,-2) to (1.5,-2.5) to (1,-3) to (0.5,-2.5);
\draw[very thick] (-0.5,-2.5) to (0,-3) to (0.5,-2.5);
\end{scope}
\begin{scope}[rotate around={-5:(4,-1)}]
\draw[very thick] (3.5,-1.5) to (3,-2) to (3.5,-2.5) to (4,-2);
\draw[very thick] (4,-1) to (3.5,-1.5) to (4,-2) to (4.5,-1.5) to (4,-1);
\end{scope}
\draw[very thick,dotted] (-1,0)node[left]{$\R$} to (4.5,0);
\node[blue] at (0,-1){$\bullet$};
\node[blue] at (4,-1){$\bullet$};
\draw[very thick,black,densely dotted] (0.2,-2.15) ellipse (1.5cm and 1.3cm);
\end{tikzpicture}
,\\
S_{Y(3)}&\hspace*{-0.1cm}=\hspace*{-0.2cm}
\begin{tikzpicture}[scale=1.2,anchorbase,smallnodes,rounded corners]
\draw[white](-0.97,0)--++(0,1);
\draw[ghost](5.0,0)node[below]{$\phantom{i}$}--++(0,1)node[above,yshift=-1pt]{$2$};
\draw[ghost](6.1,0)node[below]{$\phantom{i}$}--++(0,1)node[above,yshift=-1pt]{$0$};
\draw[solid](4.0,0)node[below]{$2$}--++(0,1)node[above,yshift=-1pt]{$\phantom{i}$};
\draw[solid](5.1,0)node[below]{$0$}--++(0,1)node[above,yshift=-1pt]{$\phantom{i}$};
\draw[redstring](0.45,0)node[below]{$1$}--++(0,1)node[above,yshift=-1pt]{$\phantom{i}$};
\draw[redstring](2.35,0)node[below]{$2$}--++(0,1)node[above,yshift=-1pt]{$\phantom{i}$};
\draw[affine](5.25,0)node[below]{$0$}--++(0,1)node[above,xshift=0.25cm,yshift=-0.02]{$\phantom{i}$};
\end{tikzpicture}
\hspace*{-0.1cm}\leftrightsquigarrow\hspace*{-0.3cm}
\begin{tikzpicture}[scale=1.2,anchorbase,scale=0.75]
\draw[very thick,blue] (-0.52,-1.45) to (-0.52,-0.01);
\draw[very thick,blue] (0,-1) to (0,-0.01);
\draw[very thick,tomato] (0.1,-1.1) to (0.1,-0.01)node[above,yshift=-1pt]{$1$};
\draw[very thick,blue] (-1.1,-1.9) to (-1.1,-0.01);
\draw[very thick,blue] (0.45,-1.55) to (0.45,-0.01);
\draw[very thick,blue] (0.9,-2.1) to (0.9,-0.01);
\draw[very thick,blue] (3.48,-1.45) to (3.48,-0.01);
\draw[very thick,blue] (4,-1) to (4,-0.01);
\draw[very thick,tomato] (4.1,-1.1) to (4.1,-0.01)node[above,yshift=-1pt]{$14$};
\begin{scope}[rotate around={-5:(0,-1)}]
\draw[very thick] (-0.5,-1.5) to (-1,-2) to (-0.5,-2.5) to (0,-2);
\draw[very thick] (0,-1) to (-0.5,-1.5) to (0,-2) to (0.5,-1.5) to (0,-1);
\draw[very thick] (0.5,-1.5) to (1,-2) to (0.5,-2.5) to (0,-2);
\draw[very thick] (1,-2) to (1.5,-2.5) to (1,-3) to (0.5,-2.5);
\draw[very thick] (-0.5,-2.5) to (0,-3) to (0.5,-2.5);
\end{scope}
\begin{scope}[rotate around={-5:(4,-1)}]
\draw[very thick] (3.5,-1.5) to (3,-2) to (3.5,-2.5) to (4,-2);
\draw[very thick] (4,-1) to (3.5,-1.5) to (4,-2) to (4.5,-1.5) to (4,-1);
\end{scope}
\draw[very thick,dotted] (-1,0)node[left]{$\R$} to (4.5,0);
\node[blue] at (0,-1){$\bullet$};
\node[blue] at (4,-1){$\bullet$};
\draw[very thick,black,densely dotted] (3.7,-1.7) circle (1cm);
\end{tikzpicture}
,
\end{align*}
where we include the red strings for comparison, as they are not included in $S_{Y(i)}$.
The Young equivalence classes are the solid strings in the left-hand diagrams, which corresponded to the Young diagrams inside the dotted circles on the right-hand side.
The suspension point of $Y(1)$ is the rightmost solid $1$-string in the illustration
above and the suspension point of $Y(3)$ is the rightmost solid $0$-string.
\end{Example}

As indicated in \autoref{E:YoungEquiClasses}, we think of Young equivalence classes as the strings associated to one Young diagram and of their suspension point as the suspension point of the associated Young diagram.

\begin{Lemma}\label{L:VerticalHell}
Let $S$ be a straight line diagram. Then $L(S)=L(\1_{\blam})$ for some
$\blam\in\hParts$ if and only if for each Young equivalence class $Y$ of solid string in $S$
we have the following.
\begin{enumerate}

\item The suspension point of $Y$ exists and is left adjacent to an (affine) red string in $S$ of the same residue, and,

\item there are no adjacent solid $i$-strings, and no adjacent ghost $i$-strings, in $S_{Y}$.

\end{enumerate}
\end{Lemma}

\begin{proof}
Let $\blam=(\lambda^{(1)}|\dots|\lambda^{(\hell)})\in\hParts$. Each integer $m$ with
$1\leq m\leq\hell$ let $Y_{m}$ be the set of strings with $\hcoord$-coordinates
$\set{\hcoord(m,r,c)|(m,r,c)\in\blam}$. Then $Y_{m}$ is a Young equivalence class of solid
strings in $\1_{\blam}$ with suspension point at position $\hcoord(m,1,1)$, which is adjacent
to an (affine) red string. Hence, if $\blam\in\hParts$, then $\1_{\blam}$ satisfies condition
(a) and condition (b) by virtue of \autoref{E:BoxConfiguration}.

To prove the converse, suppose that $S$ is an idempotent diagram satisfying conditions
(a) and (b). Let $Y$ be a Young equivalence class in $S$. By assumption, the suspension
point exists and is adjacent to an (affine) red string. If the $\hcoord$-coordinate of
the red string is $\hcoord(m,0,0)$, then associate the node $(m,1,1)$ to the suspension
point of $Y$. We now inductively associate every solid string in $Y$ with a node as follows.
Suppose that the solid $i$-string is associated to the node $(m,r,c)$. As in \autoref{D:Young}.(a),
if a solid $(i+1)$-string is adjacent to the ghost $i$-string, then the solid $(i+1)$-string
corresponds to the node $(m,r,c+1)$. Similarly, as in \autoref{D:Young}.(b), the solid $(i-1)$-string
that has its ghost left adjacent to the solid $i$-string corresponds to the node $(m,r+1,c)$.
(These two cases exactly correspond to the two diagrams in \autoref{E:CasesTypeA}.)

Finally, \autoref{E:BoxConfiguration} ensures that the diagram that is constructed in this
way has $\hell$-partition shape. More precisely, if we remove either the $(i+1)$ or the
$(i-1)$-string from the diagram in \autoref{E:BoxConfiguration}, then we have a diagram that
has adjacent solid $i$-strings, or adjacent ghost $i$-strings, respectively, which is not possible by condition (b) of the lemma.
\end{proof}

The next result is the key to proving \autoref{T:Basis}. The main idea in the proof is to use
the (honest) Reidemeister II relations to pull the solid $i$-strings as far to the right
as possible while satisfying conditions (a) and (b) from
\autoref{L:VerticalHell}. The resulting diagram then factors through an idempotent
diagram $\1_{\blam}$ that dominates $L(S)$.
For example,
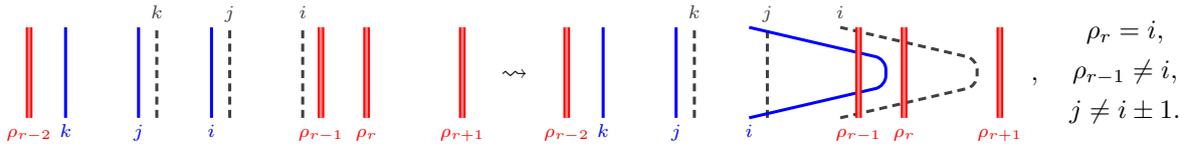
\begin{gather*}
\begin{tikzpicture}[scale=1.2,anchorbase,smallnodes]
\draw[ghost](-0.6,0)--++(0,1)node[above,yshift=-1pt]{$k$};
\draw[ghost](0.2,0)--++(0,1)node[above,yshift=-1pt]{$j$};
\draw[ghost,rounded corners](1,0)--++(0,1)node[above,yshift=-1pt]{$i$};
\draw[solid](-1.6,0)node[below]{$k$}--++(0,1);
\draw[solid](-0.8,0)node[below]{$j$}--++(0,1);
\draw[solid,rounded corners](0,0)node[below]{$i$}--++(0,1);
\draw[redstring](-2,0)node[below]{$\rho_{r-2}$}--++(0,1);
\draw[redstring](1.2,0)node[below]{$\rho_{r-1}$}--++(0,1);
\draw[redstring](1.7,0)node[below]{$\rho_{r}$}--++(0,1);
\draw[redstring](2.75,0)node[below]{$\rho_{r+1}$}--++(0,1);
\end{tikzpicture}
\rightsquigarrow
\begin{tikzpicture}[scale=1.2,anchorbase,smallnodes,rounded corners]
\draw[ghost](-0.6,0)--++(0,1)node[above,yshift=-1pt]{$k$};
\draw[ghost](0.2,0)--++(0,1)node[above,yshift=-1pt]{$j$};
\draw[ghost,rounded corners](1,0)--++(1.5,0.35)--++(0,0.3)--++(-1.5,0.35)node[above,yshift=-1pt]{$i$};
\draw[solid](-1.6,0)node[below]{$k$}--++(0,1);
\draw[solid](-0.8,0)node[below]{$j$}--++(0,1);
\draw[solid,rounded corners](0,0)node[below]{$i$}--++(1.5,0.35)--++(0,0.3)--++(-1.5,0.35);
\draw[redstring](-2,0)node[below]{$\rho_{r-2}$}--++(0,1);
\draw[redstring](1.2,0)node[below]{$\rho_{r-1}$}--++(0,1);
\draw[redstring](1.7,0)node[below]{$\rho_{r}$}--++(0,1);
\draw[redstring](2.75,0)node[below]{$\rho_{r+1}$}--++(0,1);
\end{tikzpicture}
,\quad
\begin{gathered}
\rho_{r}=i,
\\
\rho_{r-1}\neq i,
\\
j\neq i\pm 1.
\end{gathered}
\end{gather*}
As the right-hand diagram $S^{\prime}$ is obtained by pulling one string to the right we have $L(S)\ledom_{A}L(S^{\prime})$. There are some cases where we cannot naively pull strings further to the right. In these cases we apply \autoref{L:GeneralSliding} or \autoref{L:Sliding}.

\begin{Proposition}\label{P:VerticalDominance}
Suppose that $D\in\WA[n](X)$ and that $D$ factors through the idempotent diagram $S$. Then there exists $\blam\in\hParts$ such that $D$ factors through $\1_{\blam}$ and $L(S)\ledom_{A}\blam$.
\end{Proposition}

\begin{proof}
It is enough to show that $S$ factors through $\1_{\blam}$, for some $\blam\in\hParts$ with $L(S)\ledom_{A}\blam$. If $S$ satisfies the assumptions of \autoref{L:VerticalHell}, then $S=\1_{\blam}$ for some $\blam\in\hParts$, so there is nothing to prove. Let $Y$ be a Young equivalence class of solid strings in $S$. There are four cases to consider.

\noindent\textbf{Case 1.}
First assume that $Y$ does not have a suspension point.
Let $s$ be any solid string in $S_{Y}$ and assume that $s$ is an $i$-string. Then $s$ is left adjacent in $S_{Y}$ to its own ghost, to a ghost $(i-1)$-string, or
to a solid $(i+1)$-string that in turn is left adjacent to the ghost of $s$.
In the first and third situations we can use isotopy to pull $s$ to the right.
In the second situation we proceed as follows. The ghost $(i-1)$-string
has an associated solid $(i-1)$-string further to the left than the solid $i$-string. This solid $(i-1)$-string cannot be a suspension point, so we can repeat the argument for this string, again possibly finding an $(i-2)$-string further to the left. This process eventually terminates to give
a solid $j$-string that is either left adjacent in $S_{Y}$ to its own ghost or a solid $j+1$-string. Hence, in all cases we can pull this string further to the right.

\noindent\textbf{Case 2.} We assume that suspension points exist. If the suspension point of $Y$ does not satisfy condition (a) of \autoref{L:VerticalHell}, then we can pull the suspension point to the right using the Reidemeister II relations of \autoref{D:RationalCherednik}.

Hence, we can assume that every Young equivalence class of solid strings in $S$ satisfies condition (a) of \autoref{L:VerticalHell}.
Now, suppose that $S$ does not satisfy condition (b) of \autoref{L:VerticalHell} and consider the rightmost pair of solid or ghost strings that violate this condition. We consider only the case when $S$ contains two adjacent solid $i$-strings, since the case of adjacent ghost $i$-strings is similar. There are two more cases to consider.

\noindent\textbf{Case 3.} There is a solid $(i+1)$-string in between the two ghost $i$-strings. By applying \autoref{E:Sliding}, and \autoref{R:IIISliding} if necessary, we can pull the solid $i$-string to the right to give a linear combination of diagrams each of which has a dot on the rightmost $i$-string. We can now apply \autoref{R:GhostSolidTypeA} to pull the dotted $i$-string further to the right giving a more dominant diagram.

\noindent\textbf{Case 4.} A solid $(i+1)$-string does not appear in between the two ghost $i$-strings. By using isotopies and Reidemeister II, if necessary, we can assume that we are in the situation of \autoref{R:IISliding}. Hence, we can again pull the $i$-string to the right to give a linear combination of diagrams with a dot on the rightmost solid $i$-string. Repeating the argument from Case 3, we can pull the solid $i$-further to the right.

(Note that there cannot be more than one solid $(i+1)$-string in between the two ghost $i$-strings, since we took the rightmost pair of such strings to begin with.)

After finitely many iterations, all of the strings in $S$ will satisfy conditions (a) and (b) of \autoref{L:VerticalHell}, with all strings either being in the same positions or further to the right.
To see this let us analyze the situation carefully. Since we
start with the strings being in the region defined by $X$, the first
solid $i$-string that we pull rightwards in the above procedure will
stop before the rightmost affine red $i$-string. The next solid $j$-string
pulled rightwards will stop before the rightmost affine red $j$-string, and so on. Note that we never require more than $n$ affine red $i$-strings for this to work because there are only $n$ solid strings.

Hence, $S$ factors through a more dominant idempotent diagram $\1_{\blam}$.
\end{proof}

To simplify notation we write $\gdom=\gdom_{A}$ for the remainder of this section.

\begin{Lemma}\label{L:PullingDots}
Let $\blam\in\hParts$ and suppose that $1\leq m\leq n$
satisfies $\hcoord(\blam)_{m}\leq\hcoord(\ell,1,n)$. Then $y_{m}\1_{\blam}\in\WAlam*$.
\end{Lemma}

\begin{proof}
By \autoref{L:VerticalHell} all strings of $\1_{\blam}$
are in some Young equivalence class associated to an (affine) red string.
If $\hcoord(\ell,1,n)<\hcoord(\blam)_{m}$, then, by \autoref{L:RightLeft},
the Young equivalence
of the $m$th string is associated with a red string, rather than an affine red string.
By assumption, when a dot is added to the $m$th string in $\1_{\blam}$ then either relation \autoref{R:RedSolid} or \autoref{R:GhostSolidTypeA} applies.
(Note that two $i$-strings are never next to each other
in $\1_{\blam}$ by condition (b) of \autoref{L:VerticalHell}.)
If \autoref{R:RedSolid} applies, then we can pull the $m$th solid string through the red string, losing the dot and making the diagram more dominant. If \autoref{R:GhostSolidTypeA} applies, then we can pull the dotted $m$th solid string through the $(m+1)$th ghost string, or the dotted $m$th ghost string through the $(m+1)$th solid string, to give a linear combination of diagrams that are either more dominant or have the dot further to the right. Hence, the result follows by induction by either repeating this argument or by applying \autoref{P:VerticalDominance}.
\end{proof}

Recall from \autoref{D:Dw} that $D(w)$ is the permutation diagram associated to $w\in\Sym$.
For $\blam\in\hParts$ set $D_{\blam}(w)=D(w)\1_{\blam}$. Let $\Sym[\blam]$ be the corresponding Young, or parabolic, subgroup of $\Sym[n]$, where we consider $\blam$ as a composition. Given $\bmu\in\hParts$ let $\Dcal_{\blam\bmu}$ be the set of minimal length $(\Sym[\blam],\Sym[\bmu])$-double coset representatives; see, for example, \cite[Proposition 4.4]{Ma-hecke-schur}.

We have the following immediate consequence of \autoref{P:VerticalDominance}.

\begin{Lemma}\label{L:CosetReps}
Suppose that $\blam\in\hParts$ and $w\in\Sym[\blam]$. Then $D_{\blam}(w)\1_{\blam},\1_{\blam}D_{\blam}(w)\in\WAlam*$.
\end{Lemma}

\begin{proof}
It is enough to consider the case when $w=s_{r}$, where $1\leq r<n$. Let $S$ be the idempotent diagram in $D(s_{r})\1_{\blam}$ above the crossing.
Using \autoref{L:VerticalHell},
since $s_{r}\in\Sym[\blam]$, if the $r$th string is an $i$-string, then the $(r+1)$st string is an $(i+1)$-string; {\cf}
\autoref{E:CasesTypeA}. Consequently, in $S$ we can pull this $i$-string further to the right and then apply \autoref{P:VerticalDominance}. Hence $D(s_{r})1_{\blam}$, belongs to $\WAlam*$.
\end{proof}

\begin{Proposition}\label{P:KLRWSpanning}
The algebra $\WA[n](X)$ is
spanned by the diagrams in \autoref{E:TBasis}.
\end{Proposition}

\begin{proof}
Recall that $X$, the set of coordinates of endpoints, is defined by using $\blam\in\hParts$. Consequently, $\sum_{\blam}\1_{\blam}$ is the identity of $\WA[n](X)$. By \autoref{P:WABasis}, after multiplying \autoref{E:AffineBasis}
from the right by a possibly trivial permutation diagram,
$\WA[n](X)$ is spanned by the diagrams
$D(w)y^{\ba}\1_{\blam}D(v)$, for $\blam\in\hParts$, $v,w\in\Sym$ and $\ba\in\N^{n}$.

By \autoref{L:PullingDots} we can assume that $\ba\in\Affch$.
Applying \autoref{L:CosetReps} to the top and bottom of the diagram, it follows that $\WA[n](X)$ is spanned by the diagrams $D(w)y^{\ba}\1_{\blam}D(v)$,
where $\blam,\bmu,\bnu\in\hParts$,
$w\in\Dcal_{\bmu\blam}$, $v\in\Dcal_{\blam\bnu}$ and $\ba\in\Affch$.
Hence, it is enough to show that any diagram of the form $D(w)y^{\ba}\1_{\blam}D(v)$
is a linear combination of diagrams $D^{\ba}_{\bS\bT}$, where $\bS$ and $\bT$ are semistandard. Notice that we can write $D(w)=D_{\bS}$ and $D(v)=D_{\bT}$ for some, not necessary semistandard, $\blam$-tableaux, where $w=w_{\bS}$ and $v=v_{\bT}$ are the permutations defined
in \autoref{D:DST}. If the tableau $\bS$ is not semistandard, then $w_{\bS}=xs_{r}$, where $s_{r}\in\Sym[\blam]$ and $x$ and
$s_{r}$ are of shorter length than $w_{\bS}$. By the argument of \autoref{L:FiniteGeneration}, there exist scalars $a_{u}\in R$ such that
\begin{gather*}
D_{\bS}\1_{\blam}=D_{\blam}(x)D_{\blam}(s_{r})\1_{\blam}
+\sum_{u<w_{\bS}}a_{u}D_{\blam}(x)\1_{\blam}
\equiv\sum_{u<w_{\bS}}a_{u}D_{\blam}(u)\1_{\blam}
\pmod{\WAlam*},
\end{gather*}
where the last equality follows by \autoref{L:CosetReps}. Hence, by induction on $\ell(w_{\bS})$, we the diagram $D^{\ba}_{\bS\bT}$ can be written as a linear combination of diagrams of the form $D^{\ba}_{\bS^{\prime}\bT}$, for $\bS^{\prime}$ semistandard, plus a linear combination of diagrams in $\WAlam*$. By symmetry, modulo $\WAlam*$, the diagram
$D^{\ba}_{\bS\bT}$ is equal to a linear combination of diagrams of the form $D^{\ba}_{\bS^{\prime}\bT^{\prime}}$, where $\bS^{\prime}$ and $\bT^{\prime}$ are both semistandard.
\end{proof}

\begin{proof}[Proof of \autoref{T:Basis}]
We first show that \autoref{E:TBasis} is a basis of $\WA[n](X)$. By \autoref{P:KLRWSpanning} we only need to show that these diagrams are linearly independent.

Recall from \autoref{D:PolynomialAction} that $P_{\beta}(X)$ is the polynomial module of $\WA[n](X)$. To show that the diagrams in \autoref{E:TBasis} are linearly independent it suffices to show that their images in $\End(P_{\beta}(X))$ are linearly independent. By the proof of \autoref{P:WABasis}, if $\bT$ is semistandard, then the diagram $D_{\bT}$ acts non-trivially only on $\1_{\bT}P_{\beta}(X)$ where it sends a polynomial $\1_{\bT}f(y_{1},\dots,y_{n})$ to $\1_{\bT} f(y_{w_{\bT}(1)},\dots,y_{w_{\bT}(n)})$. Hence, in the action of $\gr\WA[n](W)$ on $P_{\beta}(X)$, the diagram $D^{\ba}_{\bS\bT}$ acts as
\begin{gather*}
\1_{\bT}f(y_{1},\dots,y_{n})\mapsto
\1_{\bS}y^{\ba}f(y_{w_{\bS}^{-1}w_{\bT}(1)},\dots,y_{w_{\bS}^{-1}w_{\bT}(n)}).
\end{gather*}
It follows that the diagrams in \autoref{E:TBasis} are linearly independent.

To show that the basis in \autoref{E:TBasis} is a homogeneous affine cellular basis the only axiom that does not immediately
follow by construction is (AC$_{3}$). We need to show that, if $x\in\WA(X)$, $\ba\in\Affch$ and $\bS,\bT\in\SStd(\blam)$ for $\blam\in\hParts$, then
\begin{gather*}
xD^{\ba}_{\bS\bT}\equiv\sum_{\bU\in\SStd(\blam)}
r_{\bS\bU}D^{\ba}_{\bU\bT}\pmod{\WAlam},
\qquad\text{for some }r_{\bS\bU}\in R.
\end{gather*}
If $x$ is a crossing, then this is immediate from \autoref{L:CosetReps}.
If $x$ is a dot, then we can pull it towards the diagram $\1_{\blam}$ at the equator of $D^{\ba}_{\bS\bT}$. By \autoref{R:DotCrossing} the dot slides freely through each crossing except for the $(i,i)$-crossing where an additional ``error term'' is produced that has the $(i,i)$-crossing split apart. The diagram with the error term satisfies (AC$_3$) by the argument of \autoref{P:KLRWSpanning}. By \autoref{L:PullingDots}, when a dot on the $m$th string reaches the equator of the diagram then the diagram belongs to $\WAlam*$ if $\hcoord(\blam)_{m}\leq\hcoord(\hell,1,n)$. If $\hcoord(\hell,1,n)<\hcoord(\blam)_{m}$, then the exponent $\ba$ increases by $1$ in the $m$th position so that the diagram is contained in $\WAlam*$.
\end{proof}

We can think of the final part of the argument above roughly as follows. To compute the product
$D^{\ba}_{\bS\bT}$ and $D^{\bb}_{\bS^{\prime}\bT^{\prime}}$
we can use the two statements \autoref{L:CosetReps} and \autoref{R:DotCrossing}
to simplify the middle of the picture:
\begin{gather*}
\begin{tikzpicture}[scale=1.2,anchorbase,scale=1]
\draw[line width=0.75,color=black,fill=cream] (0,1) to (0.25,0.5) to (0.75,0.5) to (1,1) to (0,1);
\node at (0.5,0.75){$T$};
\draw[line width=0.75,color=black,fill=cream] (0.25,0) to (0.25,0.5) to (0.75,0.5) to (0.75,0) to (0.25,0);
\node at (0.5,0.25){$\1_{\blam}$};
\draw[line width=0.75,color=black,fill=cream] (0,1) to (0.25,1.5) to (0.75,1.5) to (1,1) to (0,1);
\node at (0.5,1.25){$S^{\prime}$};
\draw[line width=0.75,color=black,fill=cream] (0.25,1.5) to (0.25,2) to (0.75,2) to (0.75,1.5) to (0.25,1.5);
\node at (0.5,1.75){\scalebox{0.85}{$y^{\bb}$}};
\end{tikzpicture}
\equiv
r_{TS^{\prime}}\cdot
\begin{tikzpicture}[scale=1.2,anchorbase,scale=1]
\draw[line width=0.75,color=black,fill=cream] (0.25,0) to (0.25,0.5) to (0.75,0.5) to (0.75,0) to (0.25,0);
\node at (0.5,0.25){$\1_{\blam}$};
\draw[line width=0.75,color=black,fill=cream] (0.25,0.5) to (0.25,1) to (0.75,1) to (0.75,0.5) to (0.25,0.5);
\node at (0.5,0.75){\scalebox{0.85}{$y^{\bb}$}};
\end{tikzpicture}
\equiv
r_{TS^{\prime}}\cdot
\begin{tikzpicture}[scale=1.2,anchorbase,scale=1]
\draw[line width=0.75,color=black,fill=cream] (0.25,0) to (0.25,0.5) to (0.75,0.5) to (0.75,0) to (0.25,0);
\node at (0.5,0.25){$\1_{\blam}$};
\end{tikzpicture}
\quad\pmod{\WAlam*},
\end{gather*}
which illustrates the relevant simplifications for $x=C^{\bb}_{S^{\prime}T^{\prime}}$.


\section{Homogeneous (affine) cellular bases in type \texorpdfstring{$C$}{C}}\label{S:TypeC}


This section constructs homogeneous affine cellular bases for the weighted KLRW algebras of type $C^{(1)}_{e}$. These bases descend to give homogeneous
cellular bases for the corresponding cyclotomic algebras.
The ideas are essentially the same as in \autoref{R:Strategy} for type $A$, however, the constructions look different because the position function, and hence the idempotent diagrams $\1_{\blam}$, change.


\subsection{The positioning function}


We use the same partition and tableau combinatorics as in \autoref{SS:Tableaux}.
In particular, define $\Parts$, $\hell=\ell+n(e+1)$, $\Parts(\hell)$,
$\affine{\charge}=
(\affine{\kappa}_{1},\dots,\affine{\kappa}_{\hell})\in\Z^{\hell}$
and
$\affine{\brho}=(\affine{\rho}_{1},\dots,\affine{\rho}_{\hell})\in I^{\hell}$, similar to \autoref{S:AffineCellular}. The set $I$ is still identified
with $\Z/(e+1)\Z$. We need to change the definitions of both the positioning function for partitions and the residues of tableaux.

Given an integer $k\in\Z$ write $k=2ek^{\prime}+k^{\prime\prime}$, where $k^{\prime},k^{\prime\prime}\in\Z$ are the unique integers such that $0\leq k^{\prime\prime}<2e$. Define a function $\re\colon\Z\longrightarrow I$ by
\begin{gather}\label{E:Residue}
\re(k)=
\begin{cases*}
k^{\prime\prime}+(e+1)\Z & if $0\leq k^{\prime\prime}<e$,\\
2e-k^{\prime\prime}+(e+1)\Z & if $e\leq k^{\prime\prime}<2e$.
\end{cases*}
\end{gather}
The \emph{($\affine{\brho}$-)residue} of the node $(m,r,c)$ is $\res(m,r,c)=\re(c-r+\rho_{m})+(e+1)\Z\in I$.

\begin{Remark}
A type $C$ tableaux combinatorics for the KLR algebras
appears in several sources. The above is
taken almost {\ver} from \cite[Section 1]{ArPaSp-specht-klr-typec} or \cite{EvMa-deformation-klr}.
\end{Remark}

The definition of the residue looks more demanding than it actually is:

\begin{Example}
Let $\lambda=(6^{6})$ and $\brho=(0)$.
We fill the nodes with the corresponding residues for $e=2$ and $e=3$, respectively:
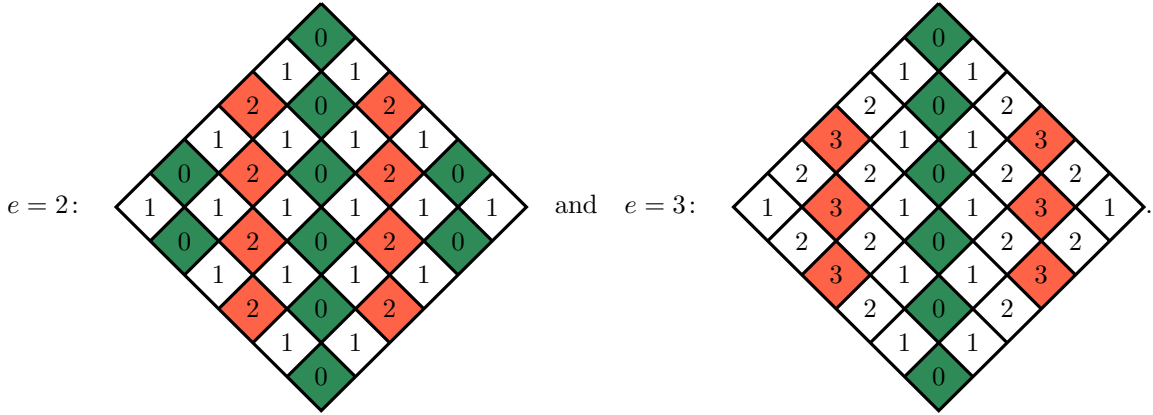
\begin{gather*}
e=2\colon\quad
\begin{tikzpicture}[scale=1.2,anchorbase,scale=0.75]
\foreach \x in {0,...,5} {
\draw[fill=spinach] (0,0-\x) to (-0.5,-0.5-\x) to (0,-1-\x) to (0.5,-0.5-\x) to (0,0-\x);
}
\foreach \x in {0,...,3} {
\draw[fill=tomato] (-1,-1-\x) to (-1.5,-1.5-\x) to (-1,-2-\x) to (-0.5,-1.5-\x) to (-1,-1-\x);
\draw[fill=tomato] (1,-1-\x) to (1.5,-1.5-\x) to (1,-2-\x) to (0.5,-1.5-\x) to (1,-1-\x);
}
\foreach \x in {0,...,1} {
\draw[fill=spinach] (-2,-2-\x) to (-2.5,-2.5-\x) to (-2,-3-\x) to (-1.5,-2.5-\x) to (-2,-2-\x);
\draw[fill=spinach] (2,-2-\x) to (2.5,-2.5-\x) to (2,-3-\x) to (1.5,-2.5-\x) to (2,-2-\x);
}
\draw[very thick] (0,0) to (-3,-3) to (0,-6) to (3,-3) to (0,0);
\draw[very thick] (-0.5,-0.5) to (2.5,-3.5);
\draw[very thick] (-1.0,-1.0) to (2.0,-4.0);
\draw[very thick] (-1.5,-1.5) to (1.5,-4.5);
\draw[very thick] (-2.0,-2.0) to (1.0,-5.0);
\draw[very thick] (-2.5,-2.5) to (0.5,-5.5);
\draw[very thick] (0.5,-0.5) to (-2.5,-3.5);
\draw[very thick] (1.0,-1.0) to (-2.0,-4.0);
\draw[very thick] (1.5,-1.5) to (-1.5,-4.5);
\draw[very thick] (2.0,-2.0) to (-1.0,-5.0);
\draw[very thick] (2.5,-2.5) to (-0.5,-5.5);
\node at (-2.5,-3.0){$1$};
\node at (-2.0,-2.5){$0$};
\node at (-2.0,-3.5){$0$};
\node at (-1.5,-2.0){$1$};
\node at (-1.5,-3.0){$1$};
\node at (-1.5,-4.0){$1$};
\node at (-1.0,-1.5){$2$};
\node at (-1.0,-2.5){$2$};
\node at (-1.0,-3.5){$2$};
\node at (-1.0,-4.5){$2$};
\node at (-0.5,-1.0){$1$};
\node at (-0.5,-2.0){$1$};
\node at (-0.5,-3.0){$1$};
\node at (-0.5,-4.0){$1$};
\node at (-0.5,-5.0){$1$};
\node at (0,-0.5){$0$};
\node at (0,-1.5){$0$};
\node at (0,-2.5){$0$};
\node at (0,-3.5){$0$};
\node at (0,-4.5){$0$};
\node at (0,-5.5){$0$};
\node at (0.5,-1.0){$1$};
\node at (0.5,-2.0){$1$};
\node at (0.5,-3.0){$1$};
\node at (0.5,-4.0){$1$};
\node at (0.5,-5.0){$1$};
\node at (1.0,-1.5){$2$};
\node at (1.0,-2.5){$2$};
\node at (1.0,-3.5){$2$};
\node at (1.0,-4.5){$2$};
\node at (1.5,-2.0){$1$};
\node at (1.5,-3.0){$1$};
\node at (1.5,-4.0){$1$};
\node at (2.0,-2.5){$0$};
\node at (2.0,-3.5){$0$};
\node at (2.5,-3.0){$1$};
\end{tikzpicture}
\quad\text{and}\quad
e=3\colon\quad
\begin{tikzpicture}[scale=1.2,anchorbase,scale=0.75]
\foreach \x in {0,...,5} {
\draw[fill=spinach] (0,0-\x) to (-0.5,-0.5-\x) to (0,-1-\x) to (0.5,-0.5-\x) to (0,0-\x);
}
\foreach \x in {0,...,2} {
\draw[fill=tomato] (-1.5,-1.5-\x) to (-2,-2-\x) to (-1.5,-2.5-\x) to (-1,-2-\x) to (-1.5,-1.5-\x);
\draw[fill=tomato] (1.5,-1.5-\x) to (2,-2-\x) to (1.5,-2.5-\x) to (1,-2-\x) to (1.5,-1.5-\x);
}
\draw[very thick] (0,0) to (-3,-3) to (0,-6) to (3,-3) to (0,0);
\draw[very thick] (-0.5,-0.5) to (2.5,-3.5);
\draw[very thick] (-1.0,-1.0) to (2.0,-4.0);
\draw[very thick] (-1.5,-1.5) to (1.5,-4.5);
\draw[very thick] (-2.0,-2.0) to (1.0,-5.0);
\draw[very thick] (-2.5,-2.5) to (0.5,-5.5);
\draw[very thick] (0.5,-0.5) to (-2.5,-3.5);
\draw[very thick] (1.0,-1.0) to (-2.0,-4.0);
\draw[very thick] (1.5,-1.5) to (-1.5,-4.5);
\draw[very thick] (2.0,-2.0) to (-1.0,-5.0);
\draw[very thick] (2.5,-2.5) to (-0.5,-5.5);
\node at (-2.5,-3.0){$1$};
\node at (-2.0,-2.5){$2$};
\node at (-2.0,-3.5){$2$};
\node at (-1.5,-2.0){$3$};
\node at (-1.5,-3.0){$3$};
\node at (-1.5,-4.0){$3$};
\node at (-1.0,-1.5){$2$};
\node at (-1.0,-2.5){$2$};
\node at (-1.0,-3.5){$2$};
\node at (-1.0,-4.5){$2$};
\node at (-0.5,-1.0){$1$};
\node at (-0.5,-2.0){$1$};
\node at (-0.5,-3.0){$1$};
\node at (-0.5,-4.0){$1$};
\node at (-0.5,-5.0){$1$};
\node at (0,-0.5){$0$};
\node at (0,-1.5){$0$};
\node at (0,-2.5){$0$};
\node at (0,-3.5){$0$};
\node at (0,-4.5){$0$};
\node at (0,-5.5){$0$};
\node at (0.5,-1.0){$1$};
\node at (0.5,-2.0){$1$};
\node at (0.5,-3.0){$1$};
\node at (0.5,-4.0){$1$};
\node at (0.5,-5.0){$1$};
\node at (1.0,-1.5){$2$};
\node at (1.0,-2.5){$2$};
\node at (1.0,-3.5){$2$};
\node at (1.0,-4.5){$2$};
\node at (1.5,-2.0){$3$};
\node at (1.5,-3.0){$3$};
\node at (1.5,-4.0){$3$};
\node at (2.0,-2.5){$2$};
\node at (2.0,-3.5){$2$};
\node at (2.5,-3.0){$1$};
\end{tikzpicture}
.
\end{gather*}
Here, we shade the $0$ and $e$ nodes to highlight the pattern.
In words, the residue increases along rows and columns until it hits $e$, and then
the residue bounces back until it hits $0$, and starts again.
The initial pattern in the first row and column, depends on the
value of the corresponding entry of $\brho$.
\end{Example}

As noted in \autoref{R:HangingPictures}, the type $A$ positioning function orders the nodes $(m,r,c)$ in a $\hell$-partition according to their height $c+r$ in the Russian diagram, by their diagonal $c-r$ and their component $m$. In type $C$ we use the row reading order $p_{\blam}(m,r,c)$, that associates to the node $(m,r,c)$ its position when reading along the rows of $\lambda^{(m)}$.
Let $\re^{\R}(k)=\re(k)-2\delta_{e\re(k)}$, where $\re$ is the function from
\autoref{E:Residue} with its output
interpreted as a real number.
As before, we need a \emph{positioning function}:

\begin{Definition}
Suppose that $0<\varepsilon<\frac{1}{4n\hell}$ and
$\blam\in\hParts$. The \emph{coordinate} of $(m,r,c)\in\blam$, is
\begin{gather*}
\hcoordc(m,r,c)=\affine{\kappa}_{m}-\frac{m}{\hell}+\re^{\R}(c-r)-\re^{\R}(\rho_{m})-p_{\blam}(m,r,c)\varepsilon,
\end{gather*}
where $p_{\blam}(m,r,c)=c+\sum_{i=1}^{r-1}\lambda^{(m)}_{i}$.
As in type $A$, $\affine{\kappa}_{m}$ is the \emph{suspension
point} of $\lambda^{(m)}$.
\end{Definition}

Note that the row reading order is not local.
{
\begin{Lemma}\label{L:RightLeftC}
We have $\hcoordc(\ell,1,n)<\hcoordc(k,r,c)$ for all $k\in\set{\ell+1,\dots,\hell}$ and all $r,c\in\set{1,\dots,n}$.
\end{Lemma}

\begin{proof}
As in type $A^{(1)}_{e}$.
\end{proof}
}

\begin{Remark}
Recall that we use the Russian convention for our Young diagrams, which
gave us a nice interpretation of the positioning function in type $A$. We are not aware of such a nice interpretation for the type $C$ positioning function. For consistency we continue to use the Russian convention when drawing tableau.
\end{Remark}

Unlike type $A$, where the strings move to the right as $n$ increases, the positioning function in type $C$ puts all of the strings into a region of length $e+1$. This is forced on us by the type $C$ quiver (for example, the solid $(e-1)$-strings do not have ghosts), and
because, as in type $A$, we want to pull strings as far to the right as possible, but the weighted KLRW relations prevent us from doing this.

\begin{Remark}\label{R:TypeCPositions}
Using the type $C$ positioning function, the
analogs of \autoref{E:BoxConfiguration} and \autoref{E:CasesTypeA} take the following form.

The important pictures are (here and below we remove ghost $(e-1)$-strings from these pictures if they appear):
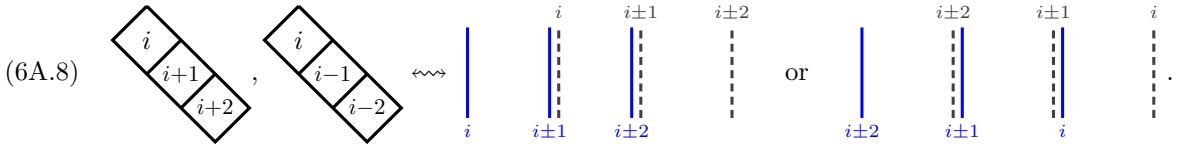
\begin{gather}\label{E:BoxConfigurationC}
\begin{tikzpicture}[scale=1.2,anchorbase,scale=0.75]
\draw[very thick] (0,0) to (-0.5,-0.5) to (1,-2.0) to (1.5,-1.5) to (0,0);
\draw[very thick] (0.5,-0.5) to (0,-1);
\draw[very thick] (1,-1) to (0.5,-1.5);
\node at (0,-0.5){$i$};
\node at (0.5,-1){\scalebox{0.85}{$i{+}1$}};
\node at (1,-1.5){\scalebox{0.85}{$i{+}2$}};
\end{tikzpicture}
,
\begin{tikzpicture}[scale=1.2,anchorbase,scale=0.75]
\draw[very thick] (0,0) to (-0.5,-0.5) to (1,-2.0) to (1.5,-1.5) to (0,0);
\draw[very thick] (0.5,-0.5) to (0,-1);
\draw[very thick] (1,-1) to (0.5,-1.5);
\node at (0,-0.5){$i$};
\node at (0.5,-1){\scalebox{0.85}{$i{-}1$}};
\node at (1,-1.5){\scalebox{0.85}{$i{-}2$}};
\end{tikzpicture}
\leftrightsquigarrow
\begin{tikzpicture}[scale=1.2,anchorbase,smallnodes,rounded corners]
\draw[ghost](1,0)node[below]{$\phantom{i}$}--++(0,1)node[above,yshift=-1pt]{$i$};
\draw[ghost](1.9,0)node[below]{$\phantom{i}$}--++(0,1)node[above,yshift=-1pt]{$i{\pm}1$};
\draw[ghost](2.9,0)node[below]{$\phantom{i}$}--++(0,1)node[above,yshift=-1pt]{$i{\pm}2$};
\draw[solid](0,0)node[below]{$i$}--++(0,1)node[above,yshift=-1pt]{$\phantom{i}$};
\draw[solid](0.9,0)node[below]{$i{\pm}1$}--++(0,1)node[above,yshift=-1pt]{$\phantom{i}$};
\draw[solid](1.8,0)node[below]{$i{\pm}2$}--++(0,1)node[above,yshift=-1pt]{$\phantom{i}$};
\end{tikzpicture}
\quad\text{or}\quad
\begin{tikzpicture}[scale=1.2,anchorbase,smallnodes,rounded corners]
\draw[ghost](-1.2,0)node[below]{$\phantom{i}$}--++(0,1)node[above,yshift=-1pt]{$i{\pm}2$};
\draw[ghost](-0.1,0)node[below]{$\phantom{i}$}--++(0,1)node[above,yshift=-1pt]{$i{\pm}1$};
\draw[ghost](1,0)node[below]{$\phantom{i}$}--++(0,1)node[above,yshift=-1pt]{$i$};
\draw[solid](-2.2,0)node[below]{$i{\pm}2$}--++(0,1)node[above,yshift=-1pt]{$\phantom{i}$};
\draw[solid](-1.1,0)node[below]{$i{\pm}1$}--++(0,1)node[above,yshift=-1pt]{$\phantom{i}$};
\draw[solid](0,0)node[below]{$i$}--++(0,1)node[above,yshift=-1pt]{$\phantom{i}$};
\end{tikzpicture}.
\end{gather}
where the residue of the middle node is not $0$, $e-1$ or $e$. (Both displayed diagrams can arise for both row configurations.) We interpret the strings as moving the right in the first diagram and to the left in the second.
When the middle node has residue $0$, $e-1$ or $e$ we have:
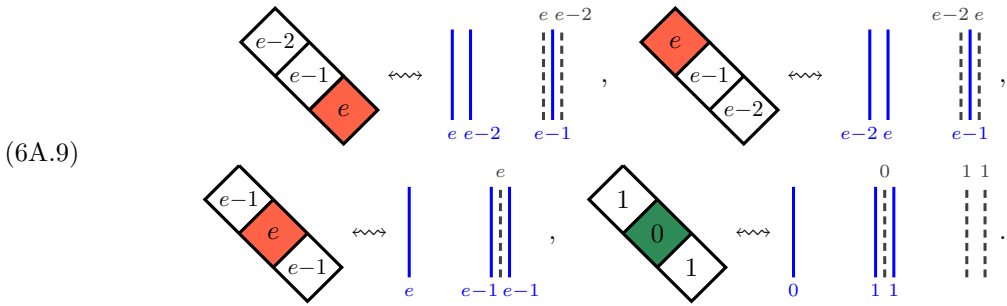
\begin{gather}\label{E:BoxConfigurationC2}
  \begin{aligned}
\begin{tikzpicture}[scale=1.2,anchorbase,scale=0.75]
\draw[fill=tomato, very thick](1,-2)--(1.5,-1.5)--(1,-1)--(0.5,-1.5)--cycle;
\draw[very thick](1,-1)--(0,0)--(-0.5,-0.5)--(0.5,-1.5);
\draw[very thick](1,-1)--(0.5,-1.5) (0.5,-0.5)--(0,-1);
\node at (0,-0.5){\scalebox{0.85}{$e{-}2$}};
\node at (0.5,-1){\scalebox{0.85}{$e{-}1$}};
\node at (1,-1.5){$e$};
\end{tikzpicture}
\leftrightsquigarrow
\begin{tikzpicture}[scale=1.2,anchorbase,smallnodes,rounded corners]
  \draw[solid](0,0)node[below]{$\phantom{aa}e{-}2$}--++(0,1)node[above,yshift=-1pt]{$\phantom{i}$};
  \draw[ghost](1,0)node[below]{$\phantom{i}$}--++(0,1)node[above,yshift=-1pt]{$\phantom{aa}e{-}2$};
  \draw[solid](0.9,0)node[below]{$e{-}1$}--++(0,1)node[above,yshift=-1pt]{$\phantom{i}$};
  \draw[ghost](0.8,0)node[below]{$\phantom{e}$}--++(0,1)node[above,yshift=-1pt]{$e$};
  \draw[solid](-0.2,0)node[below]{$e$}--++(0,1)node[above,yshift=-1pt]{$\phantom{e}$};
\end{tikzpicture}
,\quad
\begin{tikzpicture}[scale=1.2,anchorbase,scale=0.75]
\draw[fill=tomato,very thick] (0,0)--(0.5,-0.5)--(0,-1)--(-0.5,-0.5)--cycle;
\draw[very thick] (0.5,-0.5)--(1.5,-1.5)--(1,-2)--(0,-1) (1,-1)--(0.5,-1.5);
\node at (0,-0.5){$e$};
\node at (0.5,-1){\scalebox{0.85}{$e{-}1$}};
\node at (1,-1.5){\scalebox{0.85}{$e{-}2$}};
\end{tikzpicture}
\leftrightsquigarrow
\begin{tikzpicture}[scale=1.2,anchorbase,smallnodes,rounded corners]
\draw[solid](0,0)node[below]{$e$}--++(0,1)node[above,yshift=-1pt]{$\phantom{i}$};
\draw[ghost](1,0)node[below]{$\phantom{e}$}--++(0,1)node[above,yshift=-1pt]{$e$};
\draw[solid](0.9,0)node[below]{$e{-}1$}--++(0,1)node[above,yshift=-1pt]{$\phantom{i}$};
\draw[ghost](0.8,0)node[below]{$\phantom{e}$}--++(0,1)node[above,yshift=-1pt]{$e{-}2$\phantom{aa}};
\draw[solid](-0.2,0)node[below]{$e{-}2$\phantom{aa}}--++(0,1)node[above,yshift=-1pt]{$\phantom{i}$};
\end{tikzpicture}
,
\\
\begin{tikzpicture}[scale=1.2,anchorbase,scale=0.75]
\draw[fill=tomato] (0,-1) to (0.5,-0.5) to (1,-1) to (0.5,-1.5) to (0,-1);
\draw[very thick] (0,0) to (-0.5,-0.5) to (1,-2.0) to (1.5,-1.5) to (0,0);
\draw[very thick] (0.5,-0.5) to (0,-1);
\draw[very thick] (1,-1) to (0.5,-1.5);
\node at (0,-0.5){\scalebox{0.85}{$e{-}1$}};
\node at (0.5,-1){$e$};
\node at (1,-1.5){\scalebox{0.85}{$e{-}1$}};
\end{tikzpicture}
\leftrightsquigarrow
\begin{tikzpicture}[scale=1.2,anchorbase,smallnodes,rounded corners]
\draw[ghost](-0.1,0)node[below]{$\phantom{i}$}--++(0,1)node[above,yshift=-1pt]{$e$};
\draw[solid](-1.1,0)node[below]{$e$}--++(0,1)node[above,yshift=-1pt]{$\phantom{i}$};
\draw[solid](-0.2,0)node[below]{$e{-}1\phantom{aa}$}--++(0,1)node[above,yshift=-1pt]{$\phantom{i}$};
\draw[solid](0,0)node[below]{$\phantom{aa}e{-}1$}--++(0,1)node[above,yshift=-1pt]{$\phantom{i}$};
\end{tikzpicture}
,\quad
\begin{tikzpicture}[scale=1.2,anchorbase,scale=0.75]
\draw[fill=spinach] (0,-1) to (0.5,-0.5) to (1,-1) to (0.5,-1.5) to (0,-1);
\draw[very thick] (0,0) to (-0.5,-0.5) to (1,-2.0) to (1.5,-1.5) to (0,0);
\draw[very thick] (0.5,-0.5) to (0,-1);
\draw[very thick] (1,-1) to (0.5,-1.5);
\node at (0,-0.5){$1$};
\node at (0.5,-1){$0$};
\node at (1,-1.5){$1$};
\end{tikzpicture}
\leftrightsquigarrow
\begin{tikzpicture}[scale=1.2,anchorbase,smallnodes,rounded corners]
\draw[ghost](-0.1,0)node[below]{$\phantom{i}$}--++(0,1)node[above,yshift=-1pt]{$0$};
\draw[ghost](0.8,0)node[below]{$\phantom{i}$}--++(0,1)node[above,yshift=-1pt]{$1$};
\draw[ghost](1,0)node[below]{$\phantom{i}$}--++(0,1)node[above,yshift=-1pt]{$1$};
\draw[solid](-1.1,0)node[below]{$0$}--++(0,1)node[above,yshift=-1pt]{$\phantom{i}$};
\draw[solid](-0.2,0)node[below]{$1$}--++(0,1)node[above,yshift=-1pt]{$\phantom{i}$};
\draw[solid](0,0)node[below]{$1$}--++(0,1)node[above,yshift=-1pt]{$\phantom{i}$};
\end{tikzpicture}
.
\end{aligned}
\end{gather}
Again, we do not illustrate other strings that may appear in these diagrams.
As we will see, the reason why these are
minimal configurations in
type $C^{(1)}_{e}$ is
\autoref{E:SlidingC}, which is a weaker string pulling relation
than \autoref{E:Sliding}.

Another crucial configuration,
which is not local anymore as it involves two rows, is:
\begin{gather}\label{E:BoxConfigurationC3}
\begin{tikzpicture}[scale=1.2,anchorbase,scale=0.75]
\draw[very thick] (0,0) to (-0.5,-0.5) to (1,-2.0) to (1.5,-1.5) to (0,0);
\draw[very thick] (0.5,-0.5) to (0,-1);
\draw[very thick] (1,-1) to (0.5,-1.5);
\draw[very thick] (-0.5,-0.5) to (-1,-1) to (-0.5,-1.5) to (0.5,-0.5);
\node at (0,-0.5){\scalebox{0.85}{$k{\pm}1$}};
\node at (0.5,-1.05){$\dots$};
\node at (1,-1.5){$i$};
\node at (-0.5,-1){$k$};
\end{tikzpicture}
\leftrightsquigarrow
\begin{cases}
\begin{tikzpicture}[scale=1.2,anchorbase,smallnodes,rounded corners]
\draw[ghost](0.9,0)node[below]{$\phantom{i}$}--++(0,1)node[above,yshift=-1pt]{$i$};
\draw[ghost](1,0)node[below]{$\phantom{i}$}--++(0,1)node[above,yshift=-1pt]{$i$};
\draw[solid](-0.1,0)node[below]{$i$}--++(0,1)node[above,yshift=-1pt]{$\phantom{i}$};
\draw[solid](0,0)node[below]{$i$}--++(0,1)node[above,yshift=-1pt]{$\phantom{i}$};
\end{tikzpicture}
&\text{if $i=k$},
\\
\begin{tikzpicture}[scale=1.2,anchorbase,smallnodes,rounded corners]
\draw[ghost](1,0)node[below]{$\phantom{i}$}--++(0,1)node[above,yshift=-1pt]{$i$};
\draw[ghost](1.9,0)node[below]{$\phantom{i}$}--++(0,1)node[above,yshift=-1pt]{$k$};
\draw[solid](0,0)node[below]{$i$}--++(0,1)node[above,yshift=-1pt]{$\phantom{i}$};
\draw[solid](0.9,0)node[below]{$k$}--++(0,1)node[above,yshift=-1pt]{$\phantom{i}$};
\end{tikzpicture}
\quad\text{or}\quad
\begin{tikzpicture}[scale=1.2,anchorbase,smallnodes,rounded corners]
\draw[ghost](-0.1,0)node[below]{$\phantom{i}$}--++(0,1)node[above,yshift=-1pt]{$k$};
\draw[ghost](1,0)node[below]{$\phantom{i}$}--++(0,1)node[above,yshift=-1pt]{$i$};
\draw[solid](-1.1,0)node[below]{$k$}--++(0,1)node[above,yshift=-1pt]{$\phantom{i}$};
\draw[solid](0,0)node[below]{$i$}--++(0,1)node[above,yshift=-1pt]{$\phantom{i}$};
\end{tikzpicture}
&\text{if $|i-k|=1$},
\\
\begin{gathered}
\text{the solid/ghost $k$-string is close}
\\
\text{to a ghost/solid $(k\pm 1)$-string}
\end{gathered}
&\text{otherwise}
.
\end{cases}
\end{gather}
Here we have not illustrated the $(k\pm 1)$-strings, which might be far away in this picture.

There are two potential problems with \autoref{E:BoxConfigurationC3}:
\begin{enumerate}

\item Two adjacent solid or ghost strings
of the same residue can appear in these diagrams. This is
problematic because such diagrams are zero when conjugated with crossings by \autoref{R:SolidSolid}.

\item In the second case the pictures are not distinguishable
from \autoref{E:BoxConfigurationC} or \autoref{E:BoxConfigurationC2}.

\end{enumerate}
We will overcome the issues in (a) and (b) by placing a dot on the $j$ string in such cases.

The following (illegal) configurations, which do not correspond to $\hell$-partitions,
\begin{gather}\label{E:BoxConfigurationNotGood}
\begin{tikzpicture}[scale=1.2,anchorbase,scale=0.75]
\draw[very thick] (0,0) to (-1,-1) to (0,-2) to (0.5,-1.5) to (0,-1) to (0.5,-0.5) to (0,0);
\draw[very thick] (-0.5,-0.5) to (0,-1) to (-0.5,-1.5);
\node at (0,-0.5){$i$};
\node at (0,-1.5){$i$};
\node at (-0.5,-1){\scalebox{0.85}{$i{-}1$}};
\end{tikzpicture}
\leftrightsquigarrow
\begin{tikzpicture}[scale=1.2,anchorbase,smallnodes,rounded corners]
\draw[ghost](-0.1,0)node[below]{$\phantom{i}$}--++(0,1)node[above,yshift=-1pt]{$i{-}1$};
\draw[ghost](0.8,0)node[below]{$\phantom{i}$}--++(0,1)node[above,yshift=-1pt]{$i$};
\draw[ghost](1,0)node[below]{$\phantom{i}$}--++(0,1)node[above,yshift=-1pt]{$i$};
\draw[solid](-1.1,0)node[below]{$i{-}1$}--++(0,1)node[above,yshift=-1pt]{$\phantom{i}$};
\draw[solid](-0.2,0)node[below]{$i$}--++(0,1)node[above,yshift=-1pt]{$\phantom{i}$};
\draw[solid](0,0)node[below]{$i$}--++(0,1)node[above,yshift=-1pt]{$\phantom{i}$};
\end{tikzpicture}
\quad\text{or}\quad
\begin{tikzpicture}[scale=1.2,anchorbase,scale=0.75]
\draw[very thick] (0,0) to (-1,-1) to (0,-2) to (0.5,-1.5) to (0,-1) to (0.5,-0.5) to (0,0);
\draw[very thick] (-0.5,-0.5) to (0,-1) to (-0.5,-1.5);
\node at (0,-0.5){$i$};
\node at (0,-1.5){$i$};
\node at (-0.5,-1){\scalebox{0.85}{$i{+}1$}};
\end{tikzpicture}
\leftrightsquigarrow
\begin{tikzpicture}[scale=1.2,anchorbase,smallnodes,rounded corners]
\draw[ghost](-0.2,0)node[below]{$\phantom{i}$}--++(0,1)node[above,yshift=-1pt]{$i$};
\draw[ghost](0,0)node[below]{$\phantom{i}$}--++(0,1)node[above,yshift=-1pt]{$i$};
\draw[ghost](0.9,0)node[below]{$\phantom{i}$}--++(0,1)node[above,yshift=-1pt]{$i{+}1$};
\draw[solid](-1.2,0)node[below]{$i$}--++(0,1)node[above,yshift=-1pt]{$\phantom{i}$};
\draw[solid](-1,0)node[below]{$i$}--++(0,1)node[above,yshift=-1pt]{$\phantom{i}$};
\draw[solid](-0.1,0)node[below]{$i{+}1$}--++(0,1)node[above,yshift=-1pt]{$\phantom{i}$};
\end{tikzpicture}
.
\end{gather}
These configurations will be of importance later on. As before,
the crucial relation for these is \autoref{E:Sliding}.
\end{Remark}


\subsection{The (dotted) idempotent diagrams}


Let $\1_{\blam}$ be the straight line diagram associated to
the position function above. The following dot placement recipe is motivated by
\autoref{E:BoxConfigurationC3}.

\begin{Definition}\label{D:Reorder}
Suppose $\blam\in\hParts$.
If the $k$th solid string in $\1_{\blam}$ corresponds to the node $(m,r,c)$, then
\begin{gather*}
a_{k}=
\begin{cases*}
1& if $\res(m,r,c)=\res(m,r+1,1)$,\\
0& \text{otherwise}.
\end{cases*}
\end{gather*}
The \emph{dotted idempotent} associated to $\blam$
is $y_{\blam}\1_{\blam}$, where $y_{\blam}=y_{1}^{a_{1}}\dots
y_{n}^{a_{n}}\in R[y_{1},\dots,y_{n}]$.
\end{Definition}

\begin{Remark}
In the language often used in the context of (cyclotomic) KLR algebras,
such as \cite{HuMa-klr-basis}: $a_{k}\ne0$ in
\autoref{D:Reorder} if there is an addable node of the same residue in the following row of $\blam$.
\end{Remark}

As in type $A$, the idempotent diagram $\1_{\blam}$ depends mainly on the components of $\blam$ so in the examples below we only consider $1$-partitions.
Below we use the \emph{residue sequence}
$\res(\blam)$ associated to $\blam\in\hParts$, which is
the ordered tuple of the residues of the nodes of $\blam$ listed in row reading order.

\begin{Example}\label{Ex:TypeCIdempotent}
We take $\ell=1$ and $\brho=(0)$.

\begin{enumerate}

\item For the quiver of type $C^{(1)}_{3}$ and $\blam=(9)$,
reading from left to right along rows gives the residue
sequence $\res(\blam)=(0,1,2,3,2,1,0,1,2)$. Then:
\begin{gather*}
y_{\blam}\1_{\blam}=
\DottedIdempotent[2,6,8]{2}{0,1,2,3,2,1,0,1,2}
.
\end{gather*}
One can check that it is not possible to pull any string in $\1_{\blam}$ further to the right.
For type $C^{(1)}_{4}$ we have $\res(\blam)=(0,1,2,3,4,3,2,1,0)$, and we get:
\begin{gather*}
y_{\blam}\1_{\blam}=
\DottedIdempotent[2,8]{3}{0,1,2,3,4,3,2,1,0}
.
\end{gather*}

\item Let $\blam=(5,5)$ for the quiver $C^{(1)}_{2}$. Then
$\res(\blam)=(0,1,2,1,0,1,0,1,2,1)$ and
the positioning function and the dot placement give:
\begin{gather*}
y_{\blam}\1_{\blam}=
\DottedIdempotent[2,4,9]{1}{0,1,2,1,0,1,0,1,2,1}
.
\end{gather*}
The dotted idempotent in this case for the quiver $C^{(1)}_{3}$ is
(with $\res(\blam)=(0,1,2,3,2,1,0,1,2,3)$)
\begin{gather*}
y_{\blam}\1_{\blam}=
\DottedIdempotent[2,9]{2}{0,1,2,3,2,1,0,1,2,3}
.
\end{gather*}

\item If $\blam=(3,3,3)$ and we take the quiver $C^{(1)}_{2}$, then we get $\res(\blam)=(0,1,2,1,0,1,2,1,0)$ and:
\begin{gather*}
y_{\blam}\1_{\blam}=
\DottedIdempotent[2,8]{1}{0,1,2,1,0,1,2,1,0}
.
\end{gather*}
Changing the quiver to $C^{(1)}_{3}$ gives the same residue sequence, but a different $y_{\blam}\1_{\blam}$:
\begin{gather*}
y_{\blam}\1_{\blam}=
\DottedIdempotent[2]{2}{0,1,2,1,0,1,2,1,0}
.
\end{gather*}

\end{enumerate}
For readability, all diagrams are stretched in the $x$-direction because these diagrams tend to get clustered. Note that this is the typical behavior: the $(e-2)$ and the $e$-strings
always appear very close together because the $(e-1)$-strings does not have a ghost and the positioning function compensates for this by
bouncing back, which interleaves these strings and their ghosts.
\end{Example}

\begin{Remark}
The dotted idempotent $y_{\blam}\1_{\blam}$
is not as complicated as the definition
suggests. One first places the strings using the row reading order, following the rule that they are placed as far to the right as possible such that the previously placed ghost and strings stops the new strings being pulled further to the right. One then adds dots by looking at addable nodes. This is easy to feed into a machine and, in fact, the diagrams above are drawn by a TikZ macro that uses only the residue sequence. These macros are also available on GitHub \cite{MaTu-sagemath-finite-type-klrw}.
\end{Remark}

\begin{Lemma}\label{L:Distinct1lam}
Suppose that $\blam,\bmu\in\hParts$. Then $y_{\blam}\1_{\blam}=y_{\bmu}\1_{\bmu}$ if and only if $\blam=\bmu$.
\end{Lemma}

\begin{proof}
If $\res(\blam)\neq\res(\bmu)$, then this is clear since $\1_{\blam}\neq\1_{\bmu}$ in this case.
If $\blam\neq\bmu$ and $\res(\blam)=\res(\bmu)$, then the polynomials $y_{\blam}$ and $y_{\bmu}$ are different by induction, so $y_{\blam}\1_{\blam}\neq y_{\bmu}\1_{\bmu}$ by \autoref{P:WABasis}.
\end{proof}


\subsection{The homogeneous (affine) cellular basis of type $C$}


For the algebra $\WA(X)$ for type $C^{(1)}_{e}$ we will choose the $Q$-polynomials as in \autoref{E:QPoly}.

Define $\hParts$ and $X$ exactly as in \autoref{D:ThePartitions} and consider the weighted KLRW algebra $\WA(X)$ of type $C^{(1)}_{\ell}$. Repeating \autoref{D:Semistandard}, define the set $\SStd(\blam,\bmu)$ of semistandard $\blam$-tableaux of type $\bmu$ using the type $C$ positioning function above. Set $\SStd(\blam)=\bigcup_{\bmu}\SStd(\blam,\bmu)$. As for type $A$, $({}_{-})^{\star}$ is the antiinvolution from \autoref{E:StarMap}, which we will use from now on. As in \autoref{D:DST}, define diagrams $D_{\bS}$ and $D_{\bT}$ for $\bS,\bT\in\Std(\blam,\bmu)$. Given $\ba\in\Z^{n}$ define $D_{\bS\bT}^{\ba}=D_{\bS}^{\star}y^{\ba}y_{\blam}\1_{\blam}D_{\bT}$. For $\ba=(0,\dots,0)$ set $D_{\bS\bT}=D_{\bS\bT}^{\ba}$.

The definition of the diagrams $D_{\bS\bT}$ is almost exactly the same as in type $A$ except that $y_{\blam}$ appears on their equator and we are using a different choice of positioning function. Recall the definition of the set $\Affch$ from \autoref{SS:AAffineCellular}.

\begin{Definition}\label{D:CDominance}
Let $\blam,\bmu\in\hParts$. Then $\blam$ \emph{dominates} $\bmu$, written $\blam\gedom_{C}\bmu$, if there exists a bijection $d\colon\blam\to\bmu$ such that
$\hcoordc(\alpha)\geq\hcoordc\big(d(\alpha)\big)$ and the solid string in $\1_{\blam}$ at position $\hcoordc(\alpha)$ has at least as many dots as the solid string in $\1_{\bmu}$ at position $\hcoordc\big(d(\alpha)\big)$, for all $\alpha\in\blam$. Write $\blam\gdom_{C}\bmu$ if $\blam\gedom_{C}\bmu$ and $\blam\neq\bmu$.
\end{Definition}

\begin{Remark}
As in type $A$, it is not clear to the authors how the partial order $\gdom_{C}$ compares with
the dominance order on partitions, which is used in \cite{ArPaSp-specht-klr-typec} or \cite{EvMa-deformation-klr}.
\end{Remark}

Throughout this section we consider $\hParts$ as a poset ordered by $\ledom_{C}$. In particular, we use \autoref{D:CDominance} as the partial order in the following theorem. As in type $A^{(1)}_{e}$ we define $\deg\bS=\deg D_{\bS}$, $\deg\ba=\deg y^{\ba}\1_{\blam}$ and $\deg\bT=\deg D_{\bT}$.

\begin{Theorem}\label{T:BasisC}
The set $\BX=\set{D_{\bS\bT}^\ba|\blam\in\hParts,\bS,\bT\in\SStd(\blam),\ba\in\Affch}$ is a homogeneous affine cellular basis of $\WA[n](X)$ with respect to the poset $(\hParts,\ledom_{C})$.
\end{Theorem}

The proof follows the same road as the proof of \autoref{T:Basis}, see \autoref{SS:IngredientsTypeC} below.

As before we immediately obtain:

\begin{Corollary}\label{C:BasisC}
The set $\BX[{\mathscr{R}}]=\set[\big]{D_{\bS\bT}|\blam\in\Parts,\bS,\bT\in\SStd(\blam)}$
is a homogeneous cellular basis of
the cyclotomic weighted KLRW algebra $\WAc[n](X)$.
\end{Corollary}

Giving a basis for the KLR algebra of type $C^{(1)}_{e}$ requires some additional notation. Recall from \autoref{SS:AAffineCellular} that $\bom=(1^{n}|0|\dots|0)$. By definition, $\bom\ledom_{C}\blam$, for all $\blam\in\hParts$. Recall from \autoref{E:1Abom} that $\1_{A,n}$ is the sum of diagrams with positions given by~$\bom$, where the sum is over all residues. Similarly, define
\begin{gather}\label{E:1C}
\1_{C,n}=\sum_{\bi\in I^{n}}D^{\hcoordc(\bom),\bi}_{\hcoord(\bom),\bi}(1)\in\WA[n](X).
\end{gather}
That is, $\1_{C,n}$ is the sum of straight line diagrams connecting the type $A$ coordinates and type $C$ coordinates for $\bom$, where the sum is over all residues $\bi\in I^{n}$.
Then $\1_{C,n}^{\star}\WA(X)\1_{C,n}\cong\TA$ exactly as in \autoref{C:KLRCellular}, and we also have $\1_{C,n}^{\star}\WAc(X)\1_{C,n}\cong\TAc$. We identify the corresponding algebras.

As in type $A$, a \emph{standard} tableau is a semistandard tableau of type $\bom$. Let $\Std(\blam)$ be the set of standard tableaux of shape $\blam$, for $\blam\in\hParts$. For $\bs,\bt\in\Std(\blam)$ and $\ba\in\Affch$ set $E^{\ba}_{\bs\bt}=\1_{C,n}^{\star}D^{\ba}_{\bs\bt}\1_{C,n}$ and let $E_{\bs\bt}=E^{(0,\dots,0)}_{\bs\bt}$.

\begin{Proposition}\label{P:KLRCellularC}
The set $\BX[{\mathcal{W}}]=\set[\big]{E^{\ba}_{\bs\bt}|\blam\in\Parts,\bs,\bt\in\Std(\blam),\ba\in\Affch}$
is a homogeneous affine cellular basis of
$\TA$ and
$\BX[{\mathcal{R}}]=\set[\big]{E_{\bs\bt}|\blam\in\Parts,\bs,\bt\in\Std(\blam)}$
is a homogeneous cellular basis of $\TAc$.
\end{Proposition}

\begin{proof}
As noted above, the argument of \autoref{P:WebAlg}
shows that $\TA\cong\1_{C,n}^{\star}\WA(X)\1_{C,n}$. We claim that
$\TA\cong\1_{C,n}^{\star}\WA(X)\1_{C,n}$. The evident conjugation map $\TA\to\1_{C,n}^{\star}\WA(X)\1_{C,n}$ identifies the polynomial bases of \autoref{P:WABasis} on both sides, and so it is an isomorphism. By \autoref{T:BasisC}, $\set{D^{\ba}_{\bs\bt}}$ is a basis
of $\1_{C,n}^{\star}\WA(X)\1_{C,n}$ so $\set{E^{\ba}_{\bs\bt}}$ is a basis of $\TA$.
The unsteady basis elements are identified under this isomorphism, so we also get
$\TAc\cong\1_{C,n}^{\star}\WAc(X)\1_{C,n}$, which finishes the proof.
\end{proof}

As in type $A$, as a corollary we obtain dimension formulas for the cyclotomic KLRW and KLR algebras that are similar to \autoref{E:DimFormulas}. For the cyclotomic KLR algebras this recovers \cite[Theorem~2.5]{ArPaSp-specht-klr-typec}.

\begin{Remark}
As in \autoref{R:KLRComparision}, \autoref{P:KLRCellularC} constructs a basis of the KLR algebras of type $C^{(1)}_{e}$ for each choice of loading. A basis for the cyclotomic KLR algebras of type $C^{(1)}_{e}$ is constructed in \cite{EvMa-deformation-klr} but it is not clear to the authors how the basis in \cite{EvMa-deformation-klr} is related to the basis in
\autoref{P:KLRCellularC} because, at first sight, they use different partial orders. This is slightly surprising because the definition of the basis in \autoref{P:KLRCellularC} was partly motivated by \cite{EvMa-deformation-klr}.
\end{Remark}


\subsection{Some results regarding simple modules}


Classifying the simple $\WA(X)$ and $\WAc(X)$ modules in type~$C$ is more involved than in type $A$, which was discussed in \autoref{SS:SimplesA}. More precisely, most of
\autoref{SS:SimplesA} goes through without change with the crucial difference being that
\autoref{T:SimplesA} because it is not clear when pairing on
$\Delta(\blam,K)$ can be zero. We are only able to show that the cellular pairing is nonzero for the $\hell$-partitions $\hPartsgood$, that we now define.

\begin{Definition}\label{D:GoodPartitions}
Let $\hPartsgood\subset\hParts$ be the subset of all $\hell$-partitions $\blam$
such that $y_{\blam}\1_{\blam}$ has no dots or each dotted $i$-string is locally of the form:
\begin{gather*}
\begin{tikzpicture}[scale=1.2,anchorbase,smallnodes,rounded corners]
\draw[solid](-0.2,0)node[below]{$i$}--++(0,1)node[above,yshift=-1pt]{$\phantom{i}$};
\draw[solid,dot](0,0)node[below]{$i$}--++(0,1)node[above,yshift=-1pt]{$\phantom{i}$};
\end{tikzpicture}
,
\quad\text{or}\quad
\begin{tikzpicture}[scale=1.2,anchorbase,smallnodes,rounded corners]
\draw[ghost](1,0)--++(0,1)node[above,yshift=-1pt]{$i{-}1$};
\draw[ghost](1.5,0)--++(0,1)node[above,yshift=-1pt]{$i{-}1$};
\draw[solid](1.25,0)node[below]{$i$}--++(0,1);
\draw[solid,dot](1.75,0)node[below]{$i$}--++(0,1);
\end{tikzpicture}
\text{ and }i\neq 0,
\quad\text{or}\quad
\begin{tikzpicture}[scale=1.2,anchorbase,smallnodes,rounded corners]
\draw[ghost](0.25,0)--++(0,1)node[above,yshift=-1pt]{$i$};
\draw[ghost,dot](0.75,0)--++(0,1)node[above,yshift=-1pt]{$i$};
\draw[solid](0,0)node[below]{$i{+}1$}--++(0,1);
\draw[solid](0.5,0)node[below]{$i{+}1$}--++(0,1);
\end{tikzpicture}
\text{ and }i\neq e,
\end{gather*}
and the illustrated strings can be pulled arbitrarily close to one another. In particular,
there are no multilocal configurations stopping a naive isotopy.
\end{Definition}

\begin{Example}
We ignore the affine part in this example.
For type $C^{(1)}_{3}$, $n=5$, $\ell=1$ and $\brho=(0)$ we have seven $1$-partitions of five:
\begin{align*}
y_{(5)}\1_{(5)}&=
\DottedIdempotent[2]{2}{0,1,2,3,2}
,\\
y_{(4,1)}\1_{(4,1)}&=
\DottedIdempotent[2]{2}{0,1,2,3,1}
,\\
y_{(3,2)}\1_{(3,2)}&=
\DottedIdempotent[2]{2}{0,1,2,1,0}
,\\
y_{(3,1^{2})}\1_{(3,1^{2})}&=
\DottedIdempotent[2]{2}{0,1,2,1,2}
,
\\
y_{(2^{2},1)}\1_{(2^{2},1)}&=
\DottedIdempotent[2]{2}{0,1,1,0,2}
,\\
y_{(2,1^{3})}\1_{(2,1^{3})}&=
\DottedIdempotent[2]{2}{0,1,1,2,3}
,\\
y_{(1^{5})}\1_{(1^{5})}&=
\DottedIdempotent[]{2}{0,1,2,3,2}
.
\end{align*}
Hence, $\set{(3,1^{2}),(2^{2},1),(2,1^{3}),(1^{5})}=\Partsgood$
in this case.
\end{Example}

With the same notions as in type $A$, the (much weaker) analog of \autoref{T:SimplesA} reads:

\begin{Proposition}\label{P:SimplesC}
Suppose that $R$ is a field.
\begin{enumerate}

\item The set $\set{L(\blam,K)|\blam\in\hPartsgood\text{ and }K\in S\big(B(\blam)\big)}$
is a non-redundant set of simple $\WA[n](X)$-modules.

\item The set $\set{q^{s}L(\blam)|\blam\in\hPartsgood\text{ and }s\in\Z}$ is
a non-redundant set of graded simple $\WA[n](X)$-modules.

\item The set $\set{L(\blam)|\blam\in\Partsgood}$ is a non-redundant set of simple $\WAc[n](X)$-modules.

\item The set $\set{q^{s}L(\blam)|\blam\in\Partsgood\text{ and }s\in\Z}$ is a non-redundant set of graded simple $\WAc[n](X)$-modules.

\end{enumerate}
\end{Proposition}

\begin{proof}
In type $A$, the proof of \autoref{T:SimplesA} shows that because every $\blam\in\hParts$
has an associated idempotent the
cellular pairing is nonzero. We argue in the same way for $\blam\in\hPartsgood$. If $y_{\blam}\1_{\blam}$ has no dots,
then $y_{\blam}\1_{\blam}=\1_{\blam}$ is an idempotent and as before $L(\blam)\neq 0$. Otherwise, if only the first
condition in \autoref{D:GoodPartitions} appears for $y_{\blam}\1_{\blam}$, then
our construction of the basis elements $D^{\ba}_{\bs\bt}$ implies that we get idempotents that are locally of the form
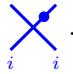
\begin{gather*}
\begin{tikzpicture}[scale=1.2,anchorbase,smallnodes,rounded corners]
\draw[solid,dot=0.25](0.5,0.5)node[above,yshift=-1pt]{$\phantom{i}$}--(0,0) node[below]{$i$};
\draw[solid](0,0.5)--(0.5,0) node[below]{$i$};
\end{tikzpicture}
.
\end{gather*}
That these diagrams are indeed idempotents follows from \autoref{R:DotCrossing}.
In all other situations we first use \autoref{E:Sliding}. Note that
in both of the right-hand diagrams in \autoref{E:Sliding} we have locally
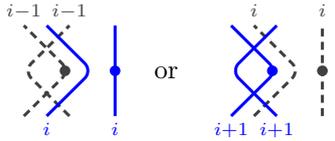
\begin{gather*}
\begin{tikzpicture}[scale=1.2,anchorbase,smallnodes,rounded corners]
\draw[ghost,dot](1,0)--++(0.5,0.5)--++(-0.5,0.5)node[above,yshift=-1pt]{$i{-}1$};
\draw[ghost](1.5,0)--++(-0.5,0.5)--++(0.5,0.5)node[above,yshift=-1pt]{$i{-}1$};
\draw[solid](1.25,0)node[below]{$i$}--++(0.5,0.5)--++(-0.5,0.5);
\draw[solid,dot](2,0)node[below]{$i$}--++(0,1);
\end{tikzpicture}
\quad\text{or}\quad
\begin{tikzpicture}[scale=1.2,anchorbase,smallnodes,rounded corners]
\draw[ghost](0.25,0)--++(0.5,0.5)--++(-0.5,0.5)node[above,yshift=-1pt]{$i$};
\draw[ghost,dot](1,0)--++(0,1)node[above,yshift=-1pt]{$i$};
\draw[solid,dot](0,0)node[below]{$i{+}1$}--++(0.5,0.5)--++(-0.5,0.5);
\draw[solid](0.5,0)node[below]{$i{+}1$}--++(-0.5,0.5)--++(0.5,0.5);
\end{tikzpicture}
,
\end{gather*}
so these summands secretly contain idempotents by \autoref{R:DotCrossing}.
This implies that the cellular pairing is nonzero in these cases, so the associated simple modules are nonzero. The general theory of cellular algebras then implies that they are pairwise non-isomorphic, so this completes the proof.
\end{proof}

\begin{Remark}
In small examples, if $\blam\in\Partsgood$, then the corresponding simple module of the KLR algebra $\TA$ is nonzero, although not all simple modules arise in this way. This is consistent with \cite[Corollary 2.16]{We-weighted-klr}, which claims that $\WAc(X)$ is Morita equivalent to $\TA$.
\end{Remark}

The analog of \autoref{C:qhereditary} fails in type $C^{(1)}_{e}$:

\begin{Proposition}
Suppose that $R$ is a field.
The algebra $\WA(X)$ is not a graded affine quasi-hereditary algebra and $\WAc(X)$ is not a graded quasi-hereditary algebra.
\end{Proposition}

\begin{proof}
The cell modules for $\hell$-partitions with a component of the form $(n)$
are one dimensional and concentrated in strictly positive degree, so they cannot carry a simple module. Hence, these algebras cannot be quasi-hereditary by the general theory of cellular algebras; {\cf} \cite[Theorem 4.1]{KoXi-affine-cellular}.
\end{proof}


\subsection{Proof of cellularity in type $C$}\label{SS:IngredientsTypeC}


Most of the arguments that we use to prove \autoref{T:BasisC} are the same as those in type $A$, for example
\autoref{L:GeneralSliding} holds without change,
and we focus on the differences.
Recall that the $Q$-polynomials for type $C$ that we use are
the ones in \autoref{E:QPoly}.
We have the following analog of \autoref{L:Sliding}.
The proof is the same as the one for \autoref{L:Sliding}
but taking the changed $Q$-polynomials into account.

\begin{Lemma}\label{L:SlidingC}
The relations in \autoref{R:GhostSolidTypeA} hold
unless (note that there are no ghost $(e-1)$-strings):
\begin{gather*}
\begin{tikzpicture}[scale=1.2,anchorbase,smallnodes,rounded corners]
\draw[ghost,dot,spinach](0,1)node[above,yshift=-1pt]{$\phantom{i}$}--++(0,-1)node[below]{$i$};
\draw[solid](0.5,1)--++(0,-1)node[below]{$j$};
\end{tikzpicture}
=
\begin{tikzpicture}[scale=1.2,anchorbase,smallnodes,rounded corners]
\draw[ghost](0,1)node[above,yshift=-1pt]{$\phantom{i}$}--++(0,-1)node[below]{$i$};
\draw[solid,dot=0.4,dot=0.6](0.5,1)--++(0,-1)node[below]{$j$};
\end{tikzpicture}
+
\begin{tikzpicture}[scale=1.2,anchorbase,smallnodes,rounded corners]
\draw[ghost](0,1)node[above,yshift=-1pt]{$\phantom{i}$}--++(0.5,-0.5)--++(-0.5,-0.5) node[below]{$i$};
\draw[solid](0.5,1)--++(-0.5,-0.5)--++(0.5,-0.5) node[below]{$j$};
\end{tikzpicture}
,\quad
\begin{tikzpicture}[scale=1.2,anchorbase,smallnodes,rounded corners]
\draw[ghost](0.5,1)node[above,yshift=-1pt]{$\phantom{i}$}--++(0,-1)node[below]{$i$};
\draw[solid,dot=0.4,dot=0.6,spinach](0,1)--++(0,-1)node[below]{$j$};
\end{tikzpicture}
=
\begin{tikzpicture}[scale=1.2,anchorbase,smallnodes,rounded corners]
\draw[ghost,dot](0.5,1)node[above,yshift=-1pt]{$\phantom{i}$}--++(0,-1)node[below]{$i$};
\draw[solid](0,1)--++(0,-1)node[below]{$j$};
\end{tikzpicture}
-
\begin{tikzpicture}[scale=1.2,anchorbase,smallnodes,rounded corners]
\draw[ghost](0.5,1)node[above,yshift=-1pt]{$\phantom{i}$}--++(-0.5,-0.5)--++(0.5,-0.5) node[below]{$i$};
\draw[solid](0,1)--++(0.5,-0.5)--++(-0.5,-0.5) node[below]{$j$};
\end{tikzpicture}
\text{ for either of }
\begin{cases}
i=0,j=1,
\\
i=e,j=e-1,
\end{cases}
\end{gather*}
where the ghost $i$-string should be ignored in these diagrams when $i=e-1$. Moreover, \autoref{E:Sliding} holds (deleting ghost $(e-1)$-strings if necessary) unless:
\begin{gather}\label{E:SlidingC}
\begin{gathered}
\begin{tikzpicture}[scale=1.2,anchorbase,smallnodes,rounded corners]
\draw[ghost,dot,spinach](1,0)--++(0,1)node[above,yshift=-1pt]{$1$};
\draw[ghost](1.5,0)--++(0,1)node[above,yshift=-1pt]{$1$};
\draw[solid,dot,spinach](0,0)node[below]{$1$}--++(0,1);
\draw[solid](0.5,0)node[below]{$1$}--++(0,1);
\draw[solid](1.25,0)node[below]{$0$}--++(0,1);
\end{tikzpicture}
=
\begin{tikzpicture}[scale=1.2,anchorbase,smallnodes,rounded corners]
\draw[ghost](1,0)--++(0,1)node[above,yshift=-1pt]{$1$};
\draw[ghost,dot](1.5,0)--++(0,1)node[above,yshift=-1pt]{$1$};
\draw[solid](0,0)node[below]{$1$}--++(0,1);
\draw[solid,dot](0.5,0)node[below]{$1$}--++(0,1);
\draw[solid](1.25,0)node[below]{$0$}--++(0,1);
\end{tikzpicture}
-
\begin{tikzpicture}[scale=1.2,anchorbase,smallnodes,rounded corners]
\draw[ghost](1,0)--++(0.5,0.5)--++(-0.5,0.5)node[above,yshift=-1pt]{$1$};
\draw[ghost,dot](1.5,0)--++(-0.5,0.5)--++(0.5,0.5)node[above,yshift=-1pt]{$1$};
\draw[solid](0,0)node[below]{$1$}--++(0.5,0.5)--++(-0.5,0.5);
\draw[solid,dot](0.5,0)node[below]{$1$}--++(-0.5,0.5)--++(0.5,0.5);
\draw[solid](1.25,0)node[below]{$0$}--++(-0.5,0.5)--++(0.5,0.5);
\end{tikzpicture}
-
\begin{tikzpicture}[scale=1.2,anchorbase,smallnodes,rounded corners]
\draw[ghost,dot](1,0)--++(0.5,0.5)--++(-0.5,0.5)node[above,yshift=-1pt]{$1$};
\draw[ghost](1.5,0)--++(-0.5,0.5)--++(0.5,0.5)node[above,yshift=-1pt]{$1$};
\draw[solid,dot](0,0)node[below]{$1$}--++(0.5,0.5)--++(-0.5,0.5);
\draw[solid](0.5,0)node[below]{$1$}--++(-0.5,0.5)--++(0.5,0.5);
\draw[solid](1.25,0)node[below]{$0$}--++(0.5,0.5)--++(-0.5,0.5);
\end{tikzpicture},
\\
\begin{tikzpicture}[scale=1.2,anchorbase,smallnodes,rounded corners]
\draw[ghost](0.25,0)--++(0,1)node[above,yshift=-1pt]{$0$};
\draw[ghost,dot,spinach](1,0)--++(0,1)node[above,yshift=-1pt]{$1$};
\draw[ghost](1.5,0)--++(0,1)node[above,yshift=-1pt]{$1$};
\draw[solid,dot,spinach](0,0)node[below]{$1$}--++(0,1);
\draw[solid](0.5,0)node[below]{$1$}--++(0,1);
\end{tikzpicture}
=
\begin{tikzpicture}[scale=1.2,anchorbase,smallnodes,rounded corners]
\draw[ghost](0.25,0)--++(0,1)node[above,yshift=-1pt]{$0$};
\draw[ghost](1,0)--++(0,1)node[above,yshift=-1pt]{$1$};
\draw[ghost,dot](1.5,0)--++(0,1)node[above,yshift=-1pt]{$1$};
\draw[solid](0,0)node[below]{$1$}--++(0,1);
\draw[solid,dot](0.5,0)node[below]{$1$}--++(0,1);
\end{tikzpicture}
+
\begin{tikzpicture}[scale=1.2,anchorbase,smallnodes,rounded corners]
\draw[ghost](1,0)--++(0.5,0.5)--++(-0.5,0.51)node[above,yshift=-1pt]{$1$};
\draw[ghost,dot](1.5,0)--++(-0.5,0.5)--++(0.5,0.5)node[above,yshift=-1pt]{$1$};
\draw[ghost](0.25,0)--++(-0.5,0.5)--++(0.5,0.5)node[above,yshift=-1pt]{$0$};
\draw[solid](0,0)node[below]{$1$}--++(0.5,0.5)--++(-0.5,0.5);
\draw[solid,dot](0.5,0)node[below]{$1$}--++(-0.5,0.5)--++(0.5,0.5);
\end{tikzpicture}
+
\begin{tikzpicture}[scale=1.2,anchorbase,smallnodes,rounded corners]
\draw[ghost,dot](1,0)--++(0.5,0.5)--++(-0.5,0.5)node[above,yshift=-1pt]{$1$};
\draw[ghost](1.5,0)--++(-0.5,0.5)--++(0.5,0.5)node[above,yshift=-1pt]{$1$};
\draw[ghost](0.25,0)--++(0.5,0.5)--++(-0.5,0.5)node[above,yshift=-1pt]{$0$};
\draw[solid,dot](0,0)node[below]{$1$}--++(0.5,0.5)--++(-0.5,0.5);
\draw[solid](0.5,0)node[below]{$1$}--++(-0.5,0.5)--++(0.5,0.5);
\end{tikzpicture},
\\
\begin{tikzpicture}[scale=1.2,anchorbase,smallnodes,rounded corners]
\draw[ghost](0.25,0)--++(0,1)node[above,yshift=-1pt]{$e$};
\draw[solid,dot,spinach](0,0)node[below]{$e{-}1$}--++(0,1);
\draw[solid](0.5,0)node[below]{$e{-}1$}--++(0,1);
\end{tikzpicture}
=
\begin{tikzpicture}[scale=1.2,anchorbase,smallnodes,rounded corners]
\draw[ghost](0.25,0)--++(0,1)node[above,yshift=-1pt]{$e$};
\draw[solid](0,0)node[below]{$e{-}1$}--++(0,1);
\draw[solid,dot](0.5,0)node[below]{$e{-}1$}--++(0,1);
\end{tikzpicture}
+
\begin{tikzpicture}[scale=1.2,anchorbase,smallnodes,rounded corners]
\draw[ghost](0.25,0)--++(-0.5,0.5)--++(0.5,0.5)node[above,yshift=-1pt]{$e$};
\draw[solid](0,0)node[below]{$e{-}1$}--++(0.5,0.5)--++(-0.5,0.5);
\draw[solid,dot](0.5,0)node[below]{$e{-}1$}--++(-0.5,0.5)--++(0.5,0.5);
\end{tikzpicture}
+
\begin{tikzpicture}[scale=1.2,anchorbase,smallnodes,rounded corners]
\draw[ghost](0.25,0)--++(0.5,0.5)--++(-0.5,0.5)node[above,yshift=-1pt]{$e$};
\draw[solid,dot](0,0)node[below]{$e{-}1$}--++(0.5,0.5)--++(-0.5,0.5);
\draw[solid](0.5,0)node[below]{$e{-}1$}--++(-0.5,0.5)--++(0.5,0.5);
\end{tikzpicture}
.
\end{gathered}
\end{gather}
\end{Lemma}

As usual, we pull the green strings, or the dot on the green strings, to the right.

\begin{proof}
These follow using similar arguments to those in \autoref{L:Sliding}.
\end{proof}

\begin{Remark}
In contrast to type $A^{(1)}_{e}$, relation \autoref{E:SlidingC}
does not allow any strings to be pulled further to the right for the configurations
\begin{gather*}
\begin{tikzpicture}[scale=1.2,anchorbase,smallnodes,rounded corners]
\draw[ghost](1,0)--++(0,1)node[above,yshift=-1pt]{$1$};
\draw[ghost](1.5,0)--++(0,1)node[above,yshift=-1pt]{$1$};
\draw[solid](0,0)node[below]{$1$}--++(0,1);
\draw[solid](0.5,0)node[below]{$1$}--++(0,1);
\draw[solid](1.25,0)node[below]{$0$}--++(0,1);
\end{tikzpicture}
,\quad
\begin{tikzpicture}[scale=1.2,anchorbase,smallnodes,rounded corners]
\draw[ghost](0.25,0)--++(0,1)node[above,yshift=-1pt]{$0$};
\draw[ghost](1,0)--++(0,1)node[above,yshift=-1pt]{$1$};
\draw[ghost](1.5,0)--++(0,1)node[above,yshift=-1pt]{$1$};
\draw[solid](0,0)node[below]{$1$}--++(0,1);
\draw[solid](0.5,0)node[below]{$1$}--++(0,1);
\end{tikzpicture}
,\quad
\begin{tikzpicture}[scale=1.2,anchorbase,smallnodes,rounded corners]
\draw[ghost](0.25,0)--++(0,1)node[above,yshift=-1pt]{$e$};
\draw[solid](0,0)node[below]{$e{-}1$}--++(0,1);
\draw[solid](0.5,0)node[below]{$e{-}1$}--++(0,1);
\end{tikzpicture}
.
\end{gather*}
This is the reason why \autoref{E:BoxConfigurationC2} are
minimal configurations.
\end{Remark}

A \emph{dotted straight line diagram} is a diagram of the form $y_{1}^{a_{1}}\dots y_{n}^{a_{n}}S$, where $S$ is a straight line diagram and $a_{k}\in\set{0,1}$ for $1\leq k\leq n$. (Note that there is at most one dot per string.)

We need a type $C$ version of \autoref{L:VerticalHell}, for which we need analogs of Young equivalence classes.

In type $A$, Young equivalence classes of strings were defined using adjacency. In type $C$, we define Young equivalence classes using ``closeness'' of strings.
In a dotted straight line diagram two strings are \emph{close} if
we can pull them
arbitrarily close to one another in a
neighborhood that does not contain any other strings using only (multilocal) isotopies. Note that two solid strings can be adjacent without being close because their ghost strings can prevent them from being pulled arbitrarily close together. We use \emph{close and to the left/right} in the evident way.

\begin{Example}
In contrast to being adjacent,
the property of being close is multilocal. For example, in \autoref{E:Close}
only the left diagram has a solid $i$-string
that is close to a solid $j$-string.
\end{Example}

For $i,j\in I$, a \emph{close triple of the form $(i,j,i)$} is a collection of three strings of the form
\begin{gather*}
\begin{tikzpicture}[scale=1.2,anchorbase,smallnodes,rounded corners]
\draw[ghost](1.1,0)--++(0,1)node[above,yshift=-1pt]{$j$};
\draw[solid](1,0)node[below]{$i$}--++(0,1);
\draw[solid](1.2,0)node[below]{$i$}--++(0,1);
\end{tikzpicture}
\quad\text{or}\quad
\begin{tikzpicture}[scale=1.2,anchorbase,smallnodes,rounded corners]
\draw[ghost](1,0)--++(0,1)node[above,yshift=-1pt]{$i$};
\draw[ghost](1.2,0)--++(0,1)node[above,yshift=-1pt]{$i$};
\draw[solid](1.1,0)node[below]{$j$}--++(0,1);
\end{tikzpicture}
,
\end{gather*}
that are all close to their direct neighbors.
Next, motivated by \autoref{E:BoxConfigurationC3}
we define:

\begin{Definition}\label{D:RowEquivalence}
Let $S$ be a dotted straight line diagram.
Fix $k\in I$. Solid $i$ and $j$-strings of $S$ are \emph{$k$-row equivalent} if the $i$-string carries a dot if and only if $i=k$, and either:
\begin{enumerate}

\item $0\leq i<e-1$, $j=i+1$ and the ghost $i$-string is close and to the right of the solid $(i+1)$-string;

\item $i=e$, $j=e-1$ and the ghost $e$-string is close and to the right of the solid $(e-1)$-string;

\item $0<i<e$, $j=i-1$ and the solid $i$-string is close and to the right of the ghost $(i-1)$-string;

\item $i=e-1$, $j=e$ and the solid $(e-1)$-string is close and to the right of the ghost $e$-string.

\end{enumerate}
A $k$-row equivalence class is \emph{admissible} if all close triples in the equivalence class are of the form $(1,0,1)$ or $(e-1,e,e-1)$.
\end{Definition}

When we apply \autoref{D:RowEquivalence}, if the $i$-string corresponds to a node $(m,r,c)$, then $k=\res(m,r+1,1)$ is the residue of the first now in the next row, as in \autoref{E:BoxConfigurationC3}.

By definition, the strings in a $k$-row equivalence class are ordered by closeness,
starting with a string that is not close and to the right of any other in the same class. Hence, we can refer to the \emph{first} and \emph{last} strings in a row equivalence class. Recall the function $\re\map{\Z}{I}$ from \autoref{E:Residue}.

\begin{Definition}\label{D:CYoung}
A \emph{Young equivalence class} $Y$ of solid strings in $S$ is a disjoint union of admissible $k$-row equivalence classes, $R_{1}\cup\dots\cup R_{z}$ where $R_{a}$ is a $k_{a}$-row equivalence class, such that:
\begin{enumerate}

\item The first string in $R_{1}$ has residue $i=\re(k_{1}+1)$ and is close to an (affine) red $i$-string;

\item $|R_{1}|\geq|R_{2}|\geq\dots\geq|R_{z}|$;

\item $k_{a+1}=\re(k_{a}-1)$ and the first string in $R_{a+1}$ is a $\re(k_{a})$-string and close to a dotted string of the same residue in $R_{a}$, or there is a $j$-string in $R_{a}$ that satisfies one of conditions (a)--(d) in \autoref{D:RowEquivalence} with respect to this string.

\end{enumerate}
\end{Definition}

\begin{Remark}
As in type $A$, \autoref{D:CYoung} is motivated by the box configurations in \autoref{R:TypeCPositions}
\end{Remark}

As in \autoref{L:Separated}, let $L(S)$ be the left justification of a (dotted) straight line diagram, defined with respect to the dominance order $\gdom_{C}$.
Note that not all solid strings in a dotted straight line diagram necessarily belong to a Young equivalence class, as implied, for example, by \autoref{D:CYoung}(a). Instead, we have:

\begin{Lemma}\label{L:CVerticalHell}
Let $S$ be a dotted straight line diagram.
Then $L(S)=L(y_{\blam}\1_{\blam})$, for some
$\blam\in\hParts$, if and only if every solid string in $S$ belongs to a Young equivalence class.
\end{Lemma}

\begin{proof}
The dotted straight line diagram $y_{\blam}\1_{\blam}$ is constructed by reading along the rows of the components of $\blam$. By \autoref{D:RowEquivalence},
\autoref{R:TypeCPositions} and \autoref{D:Reorder}, the solid strings corresponding to nodes in the same row of $\blam$ are $k$-row equivalent
(where $k$ is the residue of the first node in the row below the one that defines the equivalence) and the strings corresponding to the component $\lambda^{(m)}$ are Young equivalent. Hence, each solid string in $y_{\blam}\1_{\blam}$ belongs to a Young equivalence class.

Conversely, suppose that $S$ is a dotted straight line diagram and that the solid strings in $S$ are a disjoint union of Young equivalence classes.
Let $Y$ be a Young equivalence class of solid strings in $S$. As in \autoref{L:VerticalHell} we inductively associate the solid strings in $Y$ with a $\hell$-partition $\blam$. By \autoref{D:CYoung}.(a), the first string of $Y$ is left adjacent to an (affine) red string of the same residue. If this is the $m$th red string, then identify the solid string with the node $(m,1,1)$. By induction, suppose that the solid $i$-string corresponds to the node $(m,r,c)\in\blam$.

First suppose that this $i$-string is not the last solid string in its row equivalence class then,
in view of \autoref{R:TypeCPositions}, conditions (a)--(d) in \autoref{D:RowEquivalence} ensure that the next solid $j$-string corresponds to the node $(m,r,c+1)$. Furthermore, $\res(m,r,c+1)=j$ since the row equivalence classes are admissible. By \autoref{D:CYoung}.(b), $\blam\cup\set{(m,r,c+1)}$ is still of $\hell$-partition shape.

Finally, suppose that the $i$-string is the last string in its row equivalence class $R_{a}$ consider the $j$-string that is the first string in $R_{a+1}$. The two residue conditions in \autoref{D:CYoung}.(c) ensure that $(m,r+1,1)$ is an addable $j$-node of $\blam$. It remains to observe that the closeness requirements in \autoref{D:CYoung}.(c) correspond to the string configurations of \autoref{E:BoxConfigurationC3}.
\end{proof}

\begin{Proposition}\label{P:CVerticalDominance}
Suppose that $D\in\WA[n](X)$ and that $D$ factors through the dotted idempotent diagram $S$. Then there exists $\blam\in\hParts$ such that $D$ factors through $y_{\blam}\1_{\blam}$ and $L(S)\ledom_{C}\blam$.
\end{Proposition}

The proof of \autoref{P:CVerticalDominance} is similar to the proof of
\autoref{P:VerticalDominance}.

\begin{proof}
It is enough to show that $S$ factors through $y_{\blam}\1_{\blam}$ for some $\blam\in\hParts$. If $S$ satisfies the conditions of \autoref{L:CVerticalHell}, then $S=y_{\blam}\1_{\blam}$ for some $\blam\in\hParts$, so the proposition follows.
Hence, we can assume that the solid strings in $S$ do not form a disjoint union of Young equivalence classes.

Take the rightmost solid string not in any Young equivalence class, call it $s$.
We proceed with a case-by-case check. There are now six cases to consider.

\noindent\textbf{Case 1.} First, if there is no string close and
to the right of $s$, or its ghost, then we could pull $s$ arbitrarily far to the right. So, we can assume that this string is left adjacent to an affine red string of the same residue.

Hence, we are reduced to considering the case where $s$ is close and to the left of a solid or ghost string~$t$.
We claim that we can pull $s$ rightwards or link it to the
Young equivalence class of $t$, if applicable.
Below we assume
that $s$ is close to $t$. The case where
the ghost of $s$ is close to $t$ can be treated {\muta}.

\noindent\textbf{Case 2.} First assume that $s$ is not
in the Young equivalence class of $t$ because the admissibility
condition on close triples is violated. In this situation
we can either apply a Reidemeister II relation to pull
$s$ further rightwards, or apply \autoref{E:Sliding} together with
\autoref{R:IIISliding} to jump the dot.

\noindent\textbf{Case 3.} We now assume that $s$ is not
in the Young equivalence class of $t$ because \autoref{D:CYoung}.(b)
is violated. If the row equivalence class of $t$ does not contain an earlier string of the same residue as $s$, then we are in the situation
illustrated in \autoref{E:BoxConfigurationNotGood}. In this case
we can apply either \autoref{E:Sliding} or \autoref{E:SlidingC} (note that if  $i\pm 1\in\set{0,e}$ in leftmost string so \autoref{E:SlidingC} applies), to pull $s$ further to the right.
Note that this is the only situation where this could happen.
To see this consider
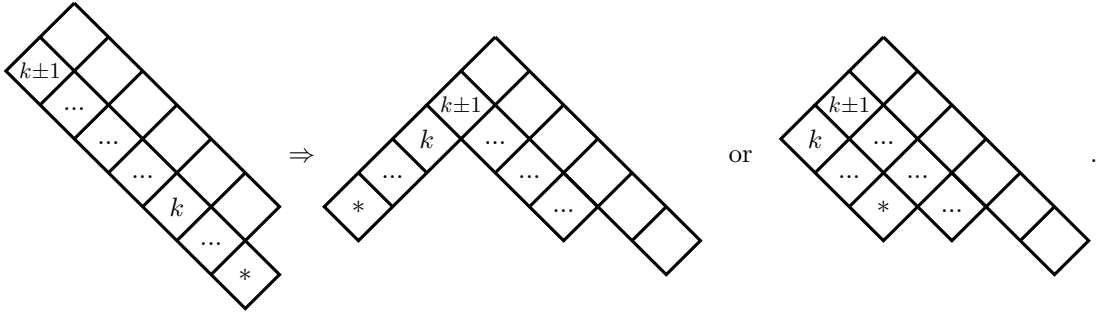
\begin{gather*}
\begin{tikzpicture}[scale=1.2,anchorbase,scale=0.75]
\draw[very thick] (0,0) to (-0.5,-0.5) to (2.5,-3.5) to (3,-3) to (0,0);
\draw[very thick] (0.5,-0.5) to (0,-1);
\draw[very thick] (1,-1) to (0.5,-1.5);
\draw[very thick] (1.5,-1.5) to (1,-2);
\draw[very thick] (2,-2) to (1.5,-2.5);
\draw[very thick] (2.5,-2.5) to (2,-3);
\draw[very thick] (-0.5,-0.5) to (-1,-1) to (2.5,-4.5) to (3,-4) to (2.5,-3.5);
\draw[very thick] (0,-1) to (-0.5,-1.5);
\draw[very thick] (0.5,-1.5) to (0,-2);
\draw[very thick] (1,-2) to (0.5,-2.5);
\draw[very thick] (1.5,-2.5) to (1,-3);
\draw[very thick] (2,-3) to (1.5,-3.5);
\draw[very thick] (2.5,-3.5) to (2,-4);
\node at (-0.5,-1){\scalebox{0.85}{$k{\pm}1$}};
\node at (0,-1.55){$\dots$};
\node at (0.5,-2.05){$\dots$};
\node at (1,-2.55){$\dots$};
\node at (1.5,-3){$k$};
\node at (2,-3.55){$\dots$};
\node at (2.5,-4){$\ast$};
\end{tikzpicture}
\Rightarrow
\begin{tikzpicture}[scale=1.2,anchorbase,scale=0.75]
\draw[very thick] (0,0) to (-0.5,-0.5) to (2.5,-3.5) to (3,-3) to (0,0);
\draw[very thick] (0.5,-0.5) to (0,-1);
\draw[very thick] (1,-1) to (0.5,-1.5);
\draw[very thick] (1.5,-1.5) to (1,-2);
\draw[very thick] (2,-2) to (1.5,-2.5);
\draw[very thick] (2.5,-2.5) to (2,-3);
\draw[very thick] (-0.5,-0.5) to (-1,-1) to (1,-3) to (1.5,-2.5) to (1,-2);
\draw[very thick] (0,-1) to (-0.5,-1.5);
\draw[very thick] (0.5,-1.5) to (0,-2);
\draw[very thick] (1,-2) to (0.5,-2.5);
\draw[very thick] (-1,-1) to (-2.5,-2.5) to (-2,-3) to (-0.5,-1.5);
\draw[very thick] (-2,-2) to (-1.5,-2.5);
\draw[very thick] (-1.5,-1.5) to (-1,-2);
\node at (-0.5,-1){\scalebox{0.85}{$k{\pm}1$}};
\node at (0,-1.55){$\dots$};
\node at (0.5,-2.05){$\dots$};
\node at (1,-2.55){$\dots$};
\node at (-1,-1.5){$k$};
\node at (-1.5,-2.05){$\dots$};
\node at (-2,-2.5){$\ast$};
\end{tikzpicture}
\quad\text{or}\quad
\begin{tikzpicture}[scale=1.2,anchorbase,scale=0.75]
\draw[very thick] (0,0) to (-0.5,-0.5) to (2.5,-3.5) to (3,-3) to (0,0);
\draw[very thick] (0.5,-0.5) to (0,-1);
\draw[very thick] (1,-1) to (0.5,-1.5);
\draw[very thick] (1.5,-1.5) to (1,-2);
\draw[very thick] (2,-2) to (1.5,-2.5);
\draw[very thick] (2.5,-2.5) to (2,-3);
\draw[very thick] (-0.5,-0.5) to (-1,-1) to (1,-3) to (1.5,-2.5) to (1,-2);
\draw[very thick] (0,-1) to (-0.5,-1.5);
\draw[very thick] (0.5,-1.5) to (0,-2);
\draw[very thick] (1,-2) to (0.5,-2.5);
\draw[very thick] (-1,-1) to (-1.5,-1.5) to (0,-3) to (0.5,-2.5);
\draw[very thick] (-0.5,-1.5) to (-1,-2);
\draw[very thick] (0,-2) to (-0.5,-2.5);
\node at (-0.5,-1){\scalebox{0.85}{$k{\pm}1$}};
\node at (0,-1.55){$\dots$};
\node at (0.5,-2.05){$\dots$};
\node at (1,-2.55){$\dots$};
\node at (-1,-1.5){$k$};
\node at (-0.5,-2.05){$\dots$};
\node at (0,-2.5){$\ast$};
\end{tikzpicture}
.
\end{gather*}
Note that there are four subcases here since the residues could be either increasing or decreasing as we read down the first column or along this row. In the case of the diagram on the right, if necessary we iterate this process to move the nodes as low as possible in the component so that the corresponding dotted idempotent does not require additional dots. Hence, in all cases, $S$ factors through a more dominant dotted idempotent diagram.

We now assume that we are not in Cases 1, 2 and 3.

\noindent\textbf{Case 4.} If $t$ is an (affine) red string, then
it has to have a different residue to $s$ by \autoref{D:CYoung}.(a).
In this case, a standard Reidemeister II relation applies and we can move
$s$ rightwards.

\noindent\textbf{Case 5.} If $t$ is a solid string, then
the only situations preventing us from pulling $s$ rightwards
via a Reidemeister II relation are the ones displayed in
\autoref{R:IISliding} and \autoref{R:IIISliding}.
In both cases we can pull $s$ rightwards except for the diagram
to the right of the equals sign in \autoref{R:IIISliding}. However, for this diagram we observe that either $s$ and $t$ are in the same Young equivalence class by \autoref{D:CYoung}.(c), or $t$ is now the rightmost solid string that is not in any Young equivalence class. Note that it can happen that $t$ is no longer contained in a Young equivalence class because the dot has moved onto $t$. As $s$ and $t$ are row equivalence by repeating the argument applied to $t$ sufficiently many times shows that $s$ (and $t$), belong to a Young equivalence class.

\noindent\textbf{Case 6.} Assume that $t$ is a ghost string.
Note that we are not in the situations of (c) and (d) of
\autoref{D:RowEquivalence} because, by assumption, $s$ is
rightmost with respect to being not in any Young equivalence class. Hence, a Reidemeister II relation applies.

Note that after pulling $s$ rightwards it is either contained in a Young equivalence class (cases 1, 2, 3 or 5) or we can continue pulling it two the right (cases 4 and 6). Repeating the argument sufficiently many times shows that $s$ is contained in a Young equivalence class because we cannot pull $s$ arbitrarily far to the right as it will eventually stop next to an (affine)
red string of the same residue. Hence, the result follows by induction.
\end{proof}

We again simplify notation and write $\gdom=\gdom_{C}$.

\begin{Lemma}\label{L:PullingDotsC}
Suppose that $\blam\in\Parts$ and $1\leq m\leq n$
satisfies $\hcoordc(\blam)_{m}\leq\hcoordc(\ell,1,n)$. Then $y_{m}y_{\blam}\1_{\blam}\in\WAlam*(C)$.
\end{Lemma}

\begin{proof}
The proof is essentially the same as \autoref{L:PullingDots}. Suppose that the $m$th string is an $i$-string. By \autoref{L:CVerticalHell} and \autoref{D:Reorder}, $y_{m}y_{\blam}\1_{\blam}$ is a dotted idempotent diagram if and only if $y_{m}y_{\blam}\1_{\blam}=y_{\bmu}\1_{\bmu}$ for some $\hell$-partition $\bmu$ with $\bmu\gdom_{C}\blam$. If $y_{m}y_{\blam}\1_{\blam}$ is not a dotted straight line diagram, then, by \autoref{P:CVerticalDominance}, it factors through a dotted straight line diagram $y_{\bmu}\1_{\bmu}$, for some $\hell$-partition $\bmu$ with $\bmu\gdom_{C}\blam$. So $y_{m}y_{\blam}\1_{\blam}\in\WAlam*(C)$, as required.
\end{proof}

For $w\in\Sym$ and $\blam\in\hParts$ define the diagram $D_{\blam}(w)=D(w)\1_{\blam}$ as in \autoref{SS:TypeA}.

\begin{Lemma}\label{L:CosetRepsC}
Suppose that $\blam\in\hParts$ and $w\in\Sym[\blam]$. Then $D_{\blam}(w)\1_{\blam},\1_{\blam}D_{\blam}(w)\in\WAlam*(C)$.
\end{Lemma}

\begin{proof}
This follows by repeating the argument of \autoref{L:CosetReps}.
\end{proof}

\begin{proof}[Proof of \autoref{T:BasisC}]
Using \autoref{P:CVerticalDominance} and the last two lemmas, the (completely formal) argument of \autoref{P:KLRWSpanning} shows that $\WA(X)$ is spanned by the elements in the statement of \autoref{T:BasisC}. Similarly, linear independence and cellularity follow exactly as in the proof of \autoref{T:Basis}. Note that the argument for linear independence relies on the faithful module in \autoref{P:WABasis}, which is valid in any type, whilst the proof of cellularity is a purely formal argument using \autoref{P:CVerticalDominance}, \autoref{L:PullingDotsC} and \autoref{L:CosetRepsC}.
\end{proof}

\section{Table of notation and central concepts}

The following is the list of the most important concepts in this paper.

\begin{xltabular}{\textwidth}{p{3.5cm}p{3.5cm}p{7cm}}\toprule
Name & Symbol & Description
\\
\midrule\endhead
\bottomrule\endfoot
The quiver
&
$\Gamma=(I,E)$, $i\to j$, $i\Rightarrow j$, $i\Rrightarrow j$ and $i\rightsquigarrow j$ (multiplicity is unimportant)
&
The quiver and its various edges, see
\autoref{SS:Quiver}
\\
\hline\mystrut
$\ell$=level, $\charge$=charge
&
$n$, $\ell$, $e$, $\charge$, $\brho$, $X$, $\bsig$
&
The number of solid strings, the number of red strings, the number of vertices in the quiver, the positions of the red strings, the labels of the red strings, the positions and labels of the solid strings, the ghost shift; all of these are fixed from the start, see \autoref{S:WebsterAlgebras}
\\
\hline\mystrut
Solid, ghost, red and affine strings
&
$\begin{tikzpicture}[scale=1.2,anchorbase,smallnodes,rounded corners]
\draw[solid] (0,0)node[below]{$i$} to (0,0.5)node[above,yshift=-1pt]{$\phantom{i}$};
\end{tikzpicture}$
,
$\begin{tikzpicture}[scale=1.2,anchorbase,smallnodes,rounded corners]
\draw[ghost] (0,0)node[below]{$\phantom{i}$} to (0,0.5)node[above,yshift=-1pt]{$i$};
\end{tikzpicture}$
,
$\begin{tikzpicture}[scale=1.2,anchorbase,smallnodes,rounded corners]
\draw[redstring] (0,0)node[below]{$i$} to (0,0.5)node[above,yshift=-1pt]{$\phantom{i}$};
\end{tikzpicture}$
,
$\begin{tikzpicture}[scale=1.2,anchorbase,smallnodes,rounded corners]
\draw[affine] (0,0)node[below]{$i$} to (0,0.5)node[above,yshift=-1pt]{$\phantom{i}$};
\end{tikzpicture}$
&
Strings in diagrams, see
\autoref{SS:Diagrams}
\\
\hline\mystrut
The infinitesimal and the asymptotic case
&
No specific symbol
&
The infinitesimal case, where $\bsig=(\varepsilon,\dots,\varepsilon)$, and the asymptotic case, where $\bsig=(1/\varepsilon,\dots,1/\varepsilon)$, both for $0<\varepsilon\ll 1$, see \autoref{SS:Diagrams}
\\
\hline\mystrut
$Q$-polynomials and algebraic dots
&
$Q_{ij}(u,v)$, $Q_{ijk}(u,v,w)$, $y_{i}$
&
The $Q$-polynomials, see \autoref{E:QPoly}, and variables that act as dots
\\
\hline\mystrut
Affine notation
&
$\hell$, $\affine{\charge}$, $\affine{\brho}$
&
We set $\hell=\ell+n(e+1)$. These ensure that we have enough affine red string to catch every solid string, see \autoref{SS:Tableaux}
\\
\hline\mystrut
Weighted KLRW algebra
&
$\WA[n](X)$
&
A diagram algebra with $n$ solid strings (sometimes we use this with a $\beta\in Q^{+}$) with top and bottom $x$-coordinates in $X$, see \autoref{SS:WebsterAlgebras}
\\
\hline\mystrut
The cyclotomic version
&
$\WAc[n](X)$
&
A quotient of $\WA[n](X)$, {\cf}
\autoref{SS:CyclotomicQuotients}
\\
\hline\mystrut
Antiinvolution
&
$({}_{-}){}^{\star}$
&
The antiinvolution mirrors
diagrams in a horizontal axis, see \autoref{SS:FirstProperties}
\\
\hline\mystrut
Idempotent diagram for $(\bx,\bi)$ or $\blam$
&
$\1_{\bx,\bi}$, $\1_{\blam}$
&
The idempotent diagram created by inductively placing strings
while reading along the rows of
the $\affine{\brho}$-partition $\blam$, see \autoref{SS:Tableaux}
\\
\hline\mystrut
Truncation idempotents
&
$\1_{A,n}$, $\1_{C,n}$
&
The idempotent diagrams for truncation to the KLR(W) algebras, see \autoref{E:1Abom} and \autoref{E:1C}
\\
\hline\mystrut
Permutation diagram for $w\in\Sym$ or $\bS$
&
$D_{w}$, $D_{\bS}$
&
A permutation diagram associated to, for example,
the $\affine{\brho}$-tableaux $\bS$, see
\autoref{SS:Tableaux}
\\
\hline\mystrut
KLRW and its cyclotomic quotient
&
$\TA[n]$, $\TAc[n]$
&
The `classical' KLRW algebras, see \autoref{SS:KLRW}
\\
\hline\mystrut
Notation for subdivision
&
$\overline{{}_{-}}$
&
Everything related to the subdivision is overlined in \autoref{S:Subdivision}
\\
\hline\mystrut
Partitions
&
$\hParts$, $\Parts$
&
Indexing sets for the middle of the cellular bases of $\WA[n](X)$ and $\WAc[n](X)$, see \autoref{SS:AAffineCellular}
\\
\hline\mystrut
Nodes
&
$(m,r,c)$
&
Nodes in (shifted) partitions: $m$ is the component index, $r$ the row index and $c$ the column index, see \autoref{SS:AAffineCellular}
\\
\hline\mystrut
Positioning function
&
$\hcoord(m,r,c)$, $\hcoordc(m,r,c)$
&
The functions (in types $A$ and $C$) that give the $x$-coordinates of solid strings,
see \autoref{SS:Tableaux}
\\
\hline\mystrut
Cellular basis elements
&
$D_{\bS\bT}^{\ba}$, $D_{\bS\bT}$
&
The homogeneous cellular basis elements, see \autoref{SS:AAffineCellular}, for $\WA[n](X)$ and $\WAc[n](X)$
\\
\hline\mystrut
Cellular basis sets
&
$\BX[{\mathscr{W}}]$, $\BX[{\mathscr{R}}]$
&
The homogeneous cellular bases sets, see \autoref{SS:AAffineCellular}, for $\WA[n](X)$ and $\WAc[n](X)$
\\
\hline\mystrut
Cellular order
&
$\ledom_{A}$, $\ledom_{C}$
&
The cellular orders (type $A$ and $C$) measuring how far strings
are to the right, see \autoref{SS:AAffineCellular}
\\
\hline\mystrut
(Affine) cell and simple modules
&
$\Delta(\blam,K)$, $\Delta(\blam)$, and $L(\blam,K)$, $L(\blam)$
&
The cell modules coming for the cellular structure and their associated simple modules, see \autoref{SS:SimplesA}
\end{xltabular}
\label{table:notation}

\end{document}